Dissertation zur Erlangung des Doktorgrades der Fakultät für Mathematik und Physik der Albert-Ludwigs-Universität Freiburg im Breisgau

PhD Thesis

# Comparison of the Categories of Motives defined by Voevodsky and Nori

Daniel Harrer

**April 2016**

Betreuer: Prof. Dr. Annette Huber-Klawitter

# Contents

# Chapter 0

# Introduction

**Abstract**

In this thesis we compare V. Voevodsky's geometric motives to the derived category of M. Nori's abelian category of mixed motives by constructing a triangulated tensor functor between them. It will be compatible with the Betti realizations on both sides. We allow an arbitrary noetherian ring of coefficients, but require it to be a field or a Dedekind domain for the tensor structure on derived Nori motives to exist.

There are three key ingredients: we present a theory of Nisnevich covers on finite acyclic diagrams of finite correspondences, explain, following D. Rydh, how to interpret finite correspondences as multivalued morphisms and elaborate on M. Nori's cohomological cell structures. For the first two, we will be working over an arbitrary regular scheme, but the last one will require that we restrict ourselves to a subfield of the complex numbers.

On the way we also show that smooth commutative group schemes over a normal base automatically admit transfers, generalizing a result by M. Spiess and T. Szamuely.

## 0.1 Main results

The central goal of this thesis is to prove the following theorem originally proposed and briefly sketched in the affine case by M. Nori:

**Theorem 7.3.1, 7.4.17.** *Assume that $k \subseteq \mathbb{C}$ is a field and let $\Lambda$ be a noetherian ring.*

*There exists a contravariant triangulated functor*

$$\mathrm{C}^{\mathrm{eff}}\colon \mathrm{DM}^{\mathrm{eff}}_{\mathrm{gm}}(k,\Lambda) \to D^b(\mathcal{MM}^{\mathrm{eff}}_{\mathrm{Nori}}(k,\Lambda))$$

*between Voevodsky's geometric motives and derived Nori motives which calculates singular cohomology $H_{\mathrm{sing}}$ in the following sense:*

*If* $\omega_{\mathrm{sing}}\colon \mathcal{MM}^{\mathrm{eff}}_{\mathrm{Nori}}(k,\Lambda) \to \Lambda\text{-}\mathrm{Mod}$ *is the forgetful fibre functor of Nori motives, then there is a natural isomorphism*

$$\omega_{\mathrm{sing}}\Big(H^n\big(\,\mathrm{C}^{\mathrm{eff}}(X[0])\big)\Big) \cong H^n_{\mathrm{sing}}(X^{\mathrm{an}},\Lambda)$$

*for all smooth varieties* $X$ *over* $k$*.*

*If furthermore* $\Lambda$ *is a Dedekind domain or a field, then* $\mathrm{C}^{\mathrm{eff}}$ *is a tensor functor and descends to a contravariant triangulated tensor functor*

$$\mathrm{C}\colon\ \mathrm{DM}_{\mathrm{gm}}(k,\Lambda) \to D^b(\mathcal{MM}_{\mathrm{Nori}}(k,\Lambda))$$

*between the non-effective versions.*

We will see that the conditions on $\Lambda$ are only imposed to guarantee the existence of Nori motives and a tensor structure on their derived category. While showing the main theorem we also prove several other interesting results. For example, we generalize the interpretation of finite correspondences as multivalued morphisms originally put forward by A. Suslin and V. Voevodsky:

**Theorem 3.7.7.** *Let* $S$ *be a noetherian normal scheme.*

*Then there exists an additive, strict monoidal, graded and faithful functor*

$$\mathrm{sym}_{\mathbb{Z}}\colon\ \mathrm{SmCor}_S \to \mathrm{Multi}_S$$

*between the preadditive, symmetric monoidal and* $\mathbb{Z}$*-graded categories of smooth finite correspondences and multivalued morphisms.*

*Additionally, this functor is compatible with the embeddings of* $\mathrm{Sm}_S$ *into both sides.*

This allows us to generalize a theorem by M. Spiess and T. Szamuely:

**Theorem 3.8.1.** *Let* $S$ *be a noetherian normal scheme. Let* $\mathrm{Sm}_S$ *be the category of schemes that are smooth, separated and of finite type over* $S$ *and let* $\mathrm{AlgSp}^\flat_S$ *be the category of algebraic spaces that are flat and separated over* $S$*. Let* $\mathcal{G}$ *be an abelian group object in* $\mathrm{AlgSp}^\flat_S$*.*

*Then the corresponding presheaf*

$$\mathcal{G}(-) = \mathrm{AlgSp}^\flat_S(-,\mathcal{G})\colon\ \mathrm{Sm}_S \hookrightarrow \mathrm{AlgSp}^\flat_S \to \mathrm{Ab}$$

*admits transfers, i.e. it extends to a presheaf*

$$\widetilde{\mathcal{G}}\colon\ \mathrm{SmCor}(S,\mathbb{Z}) \to \mathrm{Ab}\,.$$

## 0.2 Historical overview

There exists a plethora of cohomology theories in algebraic geometry: de Rham, Betti/singular, étale, $\ell$-adic cohomology, and many more. Their existence is not always universal and may among others depend on the characteristic of the base or on the smoothness of the variety. Nonetheless they are often found to satisfy very similar axioms, for example that of a Weil cohomology theory on smooth projective varieties, which includes Poincaré duality and the Künneth formula.

This led A. Grothendieck to the question if every variety $X$ has an associated *motive* $M(X)$, a single object encoding all the information found within all possible reasonable notions of cohomology, which already satisfies these properties. This and similar expectations ultimately culminated in the search for what is now called the category $\mathcal{MM}$ of *mixed motives*, as independently suggested by A. Beilinson and P. Deligne. All motives should be part of this abelian tensor category $\mathcal{MM}$ and every suitable cohomology functor $\mathrm{Var} \to \mathcal{A}$ into an abelian category should factor over a *realization* functor $r\colon \mathcal{MM} \to \mathcal{A}$. One might also hope that it is rigid, i.e. contains duals, and Tannakian, which would give rise to a motivic Galois group.

A first result was given by A. Grothendieck in the mid of the twentieth century by constructing the category of *pure motives*, those that correspond to smooth projective schemes over a field $k$. The situation becomes more complicated for general varieties, though, as witnessed by P. Deligne's theory of mixed Hodge structures of [Del71a], [Del71b] and [Del74]. It refines singular cohomology via factorization over the abelian category $\mathrm{MHS}(\mathbb{Z})$ of integral mixed Hodge structures, whose objects come with a natural weight filtration. Hence we might expect motives to also carry this finer information.

As of now, the existence of the category $\mathcal{MM}$ is still mostly conjectural. In an attempt to construct it, or at least mimic its most important properties, another concept was born, that of *motivic cohomology*, calculated by what is called the category DM of *triangulated motives.* The latter should simply be the derived category of mixed motives, but the idea behind this attempt is to go in the opposite direction: DM should be a tensor triangulated category carrying a $t$-structure, allowing us to define $\mathcal{MM}$ as its heart.

Several candidates for both $\mathcal{MM}$ and DM have been brought forward. On the triangulated side we have in particular the categories defined by M. Hanamura [Han99], M. Levine [Lev98] and V. Voevodsky [VSF00, chapter 5]. The latter was then generalized by F. Ivorra in [Ivo05] and by D.-C. Cisinski and F. Déglise in [CD12] to more general bases, and refined into another triangulated category DA of motives by J. Ayoub ([Ayo07a], [Ayo07b], [Ayo14a]).

On the abelian side we have categories described by P. Deligne and U. Jannsen (cf. [Del90]) and M. Nori (unpublished, but see [HM16]). Furthermore, we have U. Jannsen's mixed realizations of [Jan90]. Lastly, we also have a category of perverse Nori motives with rational coefficients defined by F. Ivorra in [Ivo14].

None of the above are known to satisfy all the desired axioms. There is, however, good evidence that we are on the right track: several of the suggested triangulated categories of motives are known to agree, at least with rational coefficients. Additionally, there is a result [CG15] by U. Choudhury and M. Gallauer Alves de Souza that the Galois group of M. Nori's Tannakian category agrees with that resulting from J. Ayoub's weak Tannakian formalism applied to his motives.

## 0.3 Summary of the thesis

The proposed constructions for the various categories of mixed or triangulated motives are often quite different in both approach and execution. It is thus natural to search for connections between them. In particular, proving a specific property may be rather hard in one, but becomes easy in another. Therefore a way to transfer some of their qualities may prove to be useful.

On one side of Theorem 7.4.17 we find V. Voevodsky's geometric motives $\mathrm{DM}_{\mathrm{gm}}(k)$. They are very geometric and deeply connected to cycles and intersection theory. On the other side is M. Nori's abelian tensor category $\mathcal{MM}_{\mathrm{Nori}}(k)$. It is of a more combinatorial nature, being induced by representations of quivers via singular cohomology, and is a direct candidate for the aforementioned category $\mathcal{MM}(k)$.

Hence results such as Theorem 7.4.17 that compare categories of motives, in this case prospective candidates for triangulated motives, are desirable to clarify the conjectural picture. We have the following overarching strategy of proof:

**Strategy 0.3.1.**

1. Introduce a theory of Nisnevich covers on finite acyclic diagrams of finite correspondences between smooth varieties (Definition 6.5.1 and Theorem 6.5.3) and a way to associate Čech complexes to them (Definition 5.4.4). Together with a finiteness result (Proposition 6.6.3), this allows us to replace the objects of V. Voevodsky's triangulated category $\mathrm{DM}^{\mathrm{eff}}_{\mathrm{gm}}$ by affine ones.

2. Define the notion of a multivalued morphism, representing the concept of a map that associates to every point not just one, but several points (Definition 3.4.10). Functorially translate finite correspondences into multivalued morphisms (Theorem 3.7.5).

3. Construct cohomological CW-structures on affine varieties (Definition 4.1.1) and multivalued morphisms (Definition 4.5.3) between them. Show that they can be functorially assigned to any finite acyclic diagram of multivalued morphisms (Theorem 4.8.1).

4. Combine the above to assign to every complex of finite correspondences between smooth varieties the pro-ind-system consisting of cohomological CW-structures on Čech complexes of affine Nisnevich covers.

Calculating M. Nori's cohomology of pairs then defines the functor $\mathrm{C}^{\mathrm{eff}}$. One then has to verify that this is independent of the choices and check compatibility with the tensor structure.

This approach extends M. Nori's original suggestion on how to construct such a functor $\mathrm{DM}^{\mathrm{eff}}_{\mathrm{gm}}(k, \mathbb{Z}) \to D^b(\mathcal{MM}^{\mathrm{eff}}_{\mathrm{Nori}}(k, \mathbb{Z}))$, but he only gave a brief sketch on how to deal with a single finite correspondence between smooth affine varieties. Our proof is based on the elaboration of his work presented by A. Huber and S. Müller-Stach in [HM14] and [HM16], where the following is shown:

**Theorem 0.3.2.** *Let $k \subseteq \mathbb{C}$ be a field. There exists a functor*

$$\widetilde{\mathrm{C}} \colon C^b(\mathbb{Z}[\mathrm{Var}_k]) \to D^b\left(\mathcal{MM}^{\mathrm{eff}}_{\mathrm{Nori}}(k)\right)$$

*such that*

$$\omega_{\mathrm{sing}}\left(H^n\left(\widetilde{\mathrm{C}}(X[0])\right)\right) \cong H^n_{\mathrm{sing}}(X^{\mathrm{an}}, \Lambda)$$

*for all varieties $X$ over $k$.*

They, however, did not have to deal with finite correspondences, which complicate almost every step of the argument.

Only Step 3. of Strategy 0.3.1 will actually require us to work over a field $k \subseteq \mathbb{C}$, while the remainder requires at most a regular base. We hope to extend our results to fields of positive characteristic, and possibly more general bases, in the future. Hence we opted to work in greater generality whenever possible.

Let us explain in more detail how the above steps will be achieved and what other results will be shown along the way:

**Step 1:** We introduce the new notion of *bridges* in Definition 5.2.1. They allow us to treat finite correspondences between covers akin to morphisms. To be precise, we show that they offer suitable notions of composition and fibre products. Furthermore, they behave well under refinements of the covers. Their morphism-like behaviour will allow us to make Čech complexes functorial under finite correspondences.

Following the ideas of [HM16] we *rigidify* Nisnevich covers $\mathcal{U} \twoheadrightarrow X$ by choosing point-wise splittings $X \to \mathcal{U}$ which satisfy a constructibility condition. This makes morphisms between the covers unique when requiring them to be compatible with the splittings. Due to their nature, this extends to bridges, rendering finite correspondences between rigidified Nisnevich covers over normal schemes unique. A pointed version of the usual yoga of henselian local schemes will conversely show that such a bridge always exists, assuming the cover of the domain is chosen sufficiently fine.

We then use rigidifications and bridges to introduce a notion of Nisnevich covers on diagrams of finite correspondences. In Theorem 6.5.3 we show that there are plenty such covers, even when adding further restrictions, on finite acyclic diagrams. This will follow quite formally from the aforementioned constructions. We can then interpret bounded complexes and morphisms thereof as finite acyclic diagrams, to which we assign the pro-system of all unifibrant Nisnevich covers on them.

We also show a finiteness result (Proposition 5.5.7), roughly stating that there are only finitely many isomorphism classes of connected components in the Čech complex of an étale cover. Furthermore, there exist only finitely many possible finite correspondences between any two such components that come from a bridge (cf. Proposition 6.6.3). This will prove useful when applying the results of Step 3.

We may remark upon the circumstance that Definition 5.3.1 introduces *unifibrant* Nisnevich covers. They form a pretopology between the usual Nisnevich pretopology and the one induced by Nisnevich cd-squares. Standard results then imply that the resulting sheaves will be the same. Nonetheless it satisfies properties not shared with general Nisnevich covers. Most prominently, a smooth variety, understood as an abelian sheaf via the Yoneda embedding, becomes isomorphic to the repetition-free variant of the Čech complex of a Nisnevich cover if and only if the cover is unifibrant, as seen in Theorem A.0.10. This has an analogue at the level of topological spaces (cf. Proposition B.0.11). A similar result for the full Čech complex is true for all Nisnevich covers as witnessed by Proposition 3.1.3 of [VSF00, chapter 5] or Theorem 10.3.3 of [CD12].

**Step 2:** We define the category $\mathrm{Multi}(S)$ of *multivalued morphisms* over a scheme $S$ in Definition 3.4.10. This extends on Appendix A of [Ayo14b] and requires certain functors $\mathfrak{S}^n$, where $n \in \mathbb{N}_0$. Such functors will be given by the *symmetric products* $\mathrm{S}^n(Y)$ which can be understood as unordered $n$-tuples of elements of $Y$. To properly work with them we elaborate on standard results regarding their existence: it is unconditional as an algebraic space by a result of P. Deligne, but requires the morphisms to be affin finie if a scheme is desired.

Multivalued morphisms then correspond by definition to maps $X \to \mathrm{S}^n(Y)$, i.e. functions associating to every point of $X$ an unordered $n$-tuple in $Y$. Composition is easily defined as soon as we get a natural transformation $\mathrm{S}^m(\mathrm{S}^n(Y)) \to \mathrm{S}^{mn}(Y)$, i.e. if we can interpret an $m$-tuple of $n$-tuples as an $mn$-tuple.

This leads us to the functorialities of the functors $\mathrm{S}^n(-)$, often requiring us to assume flatness. Here we meet the broad theory of *divided powers* put forward by N. Roby for rings (cf. [Rob63] and [Rob80]) and D. Rydh for schemes (cf. [Ryd08a], [Ryd08b], [Ryd08c] and [Ryd08d]). They satisfy the very same functorialities without any assumptions of flatness, giving us a second candidate $\Gamma^n(-)$ for the endofunctors $\mathfrak{S}^n$. Furthermore, divided powers simply agree with the symmetric products on flat schemes, thus should be understood as the better-suited notion. We, however, avoid them for our main results due to the greater technicality of their construction and because we will not require this broader generality.

Instead, we give a detailed description of the affine versions of symmetric products which will allow us to transfer ideas from divided powers to them. Among other things we present in Theorem 2.1.18 an explicit set of generators as an algebra over the base, which we call the *elementary symmetric tensors* due to their relations to elementary symmetric polynomials. In fact, the latter are generalized by the former as witnessed by the main result of [Vac05].

Divided powers, on the other hand, satisfy a universal property (cf. Proposition 2.6.9), linking them to so-called polynomial laws. The Galois-theoretic norm, or more generally the determinant, defines such a law, which in return defines *Grothendieck-Deligne norm maps*. We describe, following A. Suslin and V. Voevodsky, the analogous construction for symmetric powers. Additionally, we check their behaviour on the elementary symmetric tensors, which get mapped to coefficients of characteristic polynomials.

On schemes, this associates to every finite surjective morphism $f\colon Y \to X$ into a connected normal scheme a morphism $X \to \mathrm{S}^n(Y)$, where $n$ is the generic degree of $f$. We check the compatibility of this construction with the functorialities emerging from A. Suslin's and V. Voevodsky's theory of relative cycles, such as pushforward, pullback and the correspondence morphism. Ultimately this means that we get for every noetherian normal scheme $S$ a natural additive faithful tensor functor

$$\mathrm{SmCor}^{\mathrm{eff}}(S) \to \mathrm{Multi}^{\mathrm{eff}}(S)$$

which interprets finite correspondences as multivalued morphisms. The right hand side involves the symmetric powers $\mathrm{S}^n(X)$, which are in general no longer smooth. We therefore replaced finite correspondences with morphisms at the cost of regularity. This result also follows from the vast work of D. Rydh on cycles, in particular [Ryd08b] and [Ryd08d]. Over a field of characteristic 0 this was, except for its compatibility with the tensor structures, already shown by J. Ayoub (cf. Appendix A of [Ayo14b]).

**Step 3:** Using the previous steps, we can now translate finite correspondences, or even finite acyclic diagrams of them, into finite acyclic diagrams of multivalued morphisms between affine varieties. Elaborating on M. Nori's original suggestion we use his Basic Lemma 4.3.1 to associate to every affine variety $X$ over a field $k \subseteq \mathbb{C}$ a directed system of *cellular filtrations*. They are constructed to mimic the singular cohomological behaviour of cell structures in topology. In particular, their $k$-skeleta $X_k$ induce a complex with objects $H^k_{\mathrm{sing}}(X_k, X_{k-1})$ which calculates the singular cohomology $H_{\mathrm{sing}}(X)$.

The defining property of Nori motives, being the *diagram category* associated to the singular cohomology representation to the quiver of *effective pairs*, then gives us complexes in $C^b(\mathcal{MM}_{\mathrm{Nori}}(k,\Lambda))$. They, however, depend on the chosen cellular filtration, but this dependence vanishes in $D^b(\mathcal{MM}_{\mathrm{Nori}}(k,\Lambda))$.

This was enough for the results of [HM16]. Dealing with multivalued morphisms requires further arguments. We need to extend the notion of cellular filtrations to multivalued morphisms themselves. This adds significant complexity:

We formulate an equivariant version of M. Nori's Basic Lemma, suggested in his original sketch. It allows us to show the existence of cellular filtrations on many diagrams, including those necessary to check the existence, additivity, functoriality and multiplicativity of assigning cellular filtrations and their cellular cohomology. A recurring theme is that many properties are only eventually true in the resulting pro-systems.

**Step 4:** Everything left to do is putting the above three steps together and checking that the resulting functor is well-defined. This uses standard arguments involving pro-ind-objects and singular cohomology, but needs a careful treatment at times due to the finiteness assumptions that are involved.

On the way we will also prove several results that are interesting on their own. One of them is Theorem 3.8.1 on group schemes already mentioned above, which is a natural side result of the extensive theory of multivalued morphisms. Furthermore, we show some consequences of the Main Theorem 7.4.17 such as the existence of weight structures on $D^b(\mathcal{MM}_{\mathrm{Nori}}(k,\Lambda))$ and the existence of a chain

$$\begin{aligned}\mathrm{DM}_{\mathrm{gm}}(k,\mathbb{Q}) \to D^b(\mathcal{MM}_{\mathrm{Nori}}(k,\mathbb{Q})) \to D^b(\mathcal{MM}_{\mathrm{AH}}) \to\\ \to \mathcal{D}_{\mathcal{MR}} \to D^b(\mathrm{MHS}(\mathbb{Q})) \to D^b(\mathbb{Q}\text{-}\mathrm{Mod})\end{aligned}$$

of realizations. We also offer a new proof for the rigidity of Nori motives.

## 0.4 Structure of the thesis

### Fundamentals on relative cycles and motives (Chapter 1)

We give a recapitulation of the basics of Nori motives very closely based on that of A. Huber and S. Müller-Stach in [HM16]. Then we fully reintroduce A. Suslin's and V. Voevodsky's calculus of relative cycles and finite correspondences, using a variation developed by S. Kelly in [Kel13].

We continue with a reminder on V. Voevodsky's motives. We repeat the definition of the category of geometric motives $\mathrm{DM}^{\mathrm{eff}}_{\mathrm{gm}}(S,\Lambda)$ and the category of unbounded motivic complexes $\mathrm{DM}^{\mathrm{eff}}(S,\Lambda)$.

### Symmetric products, divided powers and multivalued morphisms (Chapters 2 and 3)

We give an overview over the concept of a multivalued morphisms and the underlying functors $\mathrm{S}^n(-)$ and $\Gamma^n(-)$.

In preparation to the general definitions, we will start with the affine versions, which will occupy all of Chapter 2. There we will give a very explicit description of the involved norm maps and their relation to the theory of invariants.

This then readily extends to normal schemes in Chapter 3, reproving and slightly extending results of D. Rydh.

### Cellular filtrations on finite correspondences (Chapter 4)

We give a reminder on the techniques behind the 'Yoga of good pairs' and our closely related global choice of cellular filtrations on varieties and multivalued morphisms.

Afterwards we prove an equivariant version of M. Nori's crucial Basic Lemma and show several theorems based on it.

The primary goal of this chapter is to provide us with well-behaved, i.e. additive (Proposition 4.7.2), functorial (Theorem 4.8.4) and monoidal (Theorem 4.9.6), directed systems of cellular filtrations.

### Bridges, functorial Čech complexes and Nisnevich covers of diagrams (Chapters 5 and 6)

We introduce bridges and explain how they relate to Čech complexes and rigidifications. We use them to define Nisnevich covers on finite acyclic diagrams. We then show a finiteness result, allowing us to apply the results of Chapter 4.

We also define the new pretopology of unifibrant Nisnevich covers. It becomes linked to properties of Čech complexes in Appendix A and has a topological analogue explored in Appendix B.

### Collection of main theorems (Chapter 7)

In this final part we combine the previous results to build the functor

$$\mathrm{C}^{\mathrm{eff}}\colon \mathrm{DM}^{\mathrm{eff}}_{\mathrm{gm}}(k) \to D^b(\mathcal{MM}^{\mathrm{eff}}_{\mathrm{Nori}}(k)).$$

Afterwards we verify that the functor descends to the tensor-localizations at the Tate and Lefschetz motives and thereby finish the proof of the main theorem. We also give a sketch for an alternative proof suggested by J. Ayoub.

Lastly, we state several interesting consequences of our Main Theorem 7.4.17. A few are already known, but we offer new proofs.

## 0.5 Acknowledgements

First and foremost I would like to thank my advisor, Annette Huber, for her personal support and plentiful advice in matters concerning this thesis and mathematics in general.

Further thanks go to Joseph Ayoub, who provided feedback on an early draft of this document, and his research group at the University of Zürich, in particular his student Martin Gallauer Alves de Souza. I am indebted for the hospitality extended to me during my regrettably short stay with their group.

I am very grateful to Darij Grinberg for proofreading Chapter 2. Similarly, I am much obliged to Oliver Bräunling, Paul-Jonas Hamacher, Shane Kelly, Matvey Soloviev and Konrad Völkel who provided me with many helpful comments on preliminary versions of this thesis.

Furthermore, I wish to thank Sebastian Goette for answering my numerous questions regarding algebraic topology. I also want to express my gratitude towards Viktor Kleen, Helene Sigloch, Konrad Völkel and Matthias Wendt for their help in understanding motives, homotopy theory and subtleties surrounding triangulated categories. Lastly, I benefited greatly from several fruitful discussions with Jens Eberhardt, Fritz Hörmann and Wolfgang Soergel.

Once again, I would like to express my gratitude to all of the aforementioned. Without their mathematical and moral support this thesis would likely never have come to completion.

# Chapter 1

# Review of Motives

The aim of this chapter is to remind the reader of the definitions and constructions encountered around the motives defined by Voevodsky and Nori.

In particular, we repeat Nori's definition of his motives as a diagram category. His construction associates to every quiver representation an abelian category which satisfies a universal property with respect to the given representation. Singular cohomology, or any other adequate cohomological functor, gives such a representation on pairs of varieties, defining the desired category. Our exposition is closely based on [HM16], where Nori's work was elaborated in full.

Afterwards, we aim to introduce Voevodsky's triangulated category of geometric motives. It is fundamentally based on a theory of cycles by him and Suslin. As we will use these cycles very often in this thesis, we give a full introduction to them. We will base this on the existing literature, most prominently the original article, found as [VSF00, chapter 2]. Our notations and arguments are, however, often taken from the newer work [Kel13] by Kelly, which offers quicker access to the definition of a relative cycle. Two more texts which we rely on are [Ivo07] and [CD12]. The latter is often significantly more general than what we need, though.

Ultimately, we want to work with coefficients in a ring $\Lambda$. We avoid re-proving the results of [VSF00] and [Kel13], who only considered integral coefficients. Instead, we start with integral coefficients and then make a change of coefficients. This leads in general to a different category of motives, but only over non-regular bases, which will be of no concern. Hence we decided to use this shortcut.

## 1.1 Notation and conventions

We now introduce notation that we use throughout this thesis. Most of them are standard or only slight variations of established ones.

**Categories:**

- Set: sets.
- Top: topological spaces with continuous maps.
- Ab: abelian groups.
- $\Lambda$-$\widetilde{\mathrm{Mod}}$: modules over a ring $\Lambda$.
- $\Lambda$- Mod: finitely generated modules over a ring $\Lambda$.
- $\Lambda$- Alg: algebras over a ring $\Lambda$.
- $\widetilde{Sch}_S$: schemes separated over a base scheme $S$.
- $\mathrm{Sch}_S$: schemes separated and of finite type over a base scheme $S$.
- $\mathrm{Var}_k$: varieties over a field $k$, i.e. reduced and separated schemes of finite type over $\mathrm{Spec}(k)$.
- $\mathrm{Sm}_S$: schemes which are smooth, separated and of finite type over a base scheme $S$.
- $\mathrm{AlgSp}_S$: algebraic spaces separated over an algebraic space $S$.
- $\mathrm{PreShv}(-)$: presheaves.
- $\mathrm{Shv}_\tau(-)$: sheaves with respect to a pretopology $\tau$.
- $\mathrm{Hom}_{\mathcal{C}}(A, B)$, sometimes simply $\mathcal{C}(A, B)$: hom-set of morphisms $A \to B$ in a category $\mathcal{C}$.
- $C^b(\mathcal{A})$, $C^+(\mathcal{A})$, $C^-(\mathcal{A})$, $C(\mathcal{A})$: (cohomological, i.e. increasing) chain complexes over an additive category $\mathcal{A}$ which are bounded, bounded to the left, bounded to the right or unbounded, respectively.
- $K^b(\mathcal{A})$, $K^+(\mathcal{A})$, $K^-(\mathcal{A})$, $K(\mathcal{A})$: homotopy category of chain complexes over an additive category $\mathcal{A}$ which are bounded, bounded to the left, bounded to the right or unbounded, respectively.
- $D^b(\mathcal{A})$, $D^+(\mathcal{A})$, $D^-(\mathcal{A})$, $D(\mathcal{A})$: derived category of chain complexes over an abelian category $\mathcal{A}$ which are bounded, bounded to the left, bounded to the right or unbounded, respectively.
- $\Lambda[\mathcal{C}]$: $\Lambda$-linearisation of a category $\mathcal{C}$, with $\Lambda$ a ring, i.e. the category which has the same objects as $\mathcal{C}$, but has as hom-sets the free $\Lambda$-modules over the hom-sets of $\mathcal{C}$.

**Superscripts:**

- $-^{\mathrm{eff}}$ effective versions of the respective objects and categories.
- $-^{\mathrm{aff}}$ affine schemes (in a category of schemes).
- $-^{\mathrm{sep}}$ objects separated over the base (in a category of algebraic spaces or schemes).

**Other notation:**

- For any integer $n \in \mathbb{N}_0$ we let $[n] = \{1, 2, \ldots, n\}$, both as an ordered or ordinary set.
- In a category $\mathcal{C}$ admitting fibre products we use $(f, g)_S$ to denote the morphism $A \to B \times_S C$ resulting from morphisms $f\colon A \to B$ and $g\colon A \to C$ over an object $S$ in $\mathcal{C}$.
- If applicable, and there is no danger of confusion, we will abbreviate the identity $\mathrm{id}_A\colon A \to A$ of an object $A$ of a category simply as $A$.
- If the base $S$ is clear from the context we use the shortcut $X_1 X_2 \ldots X_r$ for $X_1 \times_S X_2 \times_S \times \cdots \times_S X_r$. Under this circumstances we write $\mathrm{pr}^{XYZ}_{XZ}$ to denote the projection $XYZ \to XZ$.
- A *generic point* of a scheme $X$ is one of an irreducible component of $X$.

**Conventions:**

**Convention 1.1.1.** We assume the following unless explicitly stated otherwise:

- All rings are commutative and with unity.
- All schemes and morphisms of schemes are separated.

**Convention 1.1.2.** Let non-negative integers $m$, $n$, $r$ and $n_1, \ldots, n_r$ be given. We interpret $S_m$ as the group of bijections $[m] \to [m]$.

We fix the following conventions to interpret certain permutation groups as subgroups of symmetric groups:

- $S_m \times S_n$ will be identified with the subgroup of $S_{m+n}$ that leaves the subsets $[m]$ and $[m+n] \backslash [m]$ of $[m+n]$ invariant.

  This inductively extends to an embedding

$$S_{n_1} \times \ldots \times S_{n_r} \hookrightarrow S_{n_1 + \ldots + n_r}.$$

- $[m] \times [n]$ will be identified with $[mn]$ via the map $(x, y) \mapsto x + m \cdot (y - 1)$, thereby arranging the elements of $[mn]$ from left to right, top to bottom in a rectangular shape. If $m$ and $n$ are non-zero, this allows us to identify $S_m \times S_n$ with the subgroup of $S_{mn} \cong S_{[m]\times[n]}$ that individually preserves the two relations of lying in the same row and of lying in the same column.

  This inductively extends to a morphism

  $$S_{n_1} \times \ldots \times S_{n_r} \to S_{n_1 \cdot \ldots \cdot n_r}$$

  which is injective if the $n_i$ are non-zero.

- If $m$ is non-zero, the semi-direct product $S_m \ltimes S_n^m$, where $S_m$ acts by permuting the factors, will be interpreted as the subgroup

  $$\{\sigma \in S_{mn} \mid a \equiv b \mod m \implies \sigma(a) \equiv \sigma(b) \mod m\} \subseteq S_{mn}.$$

  In the interpretation of the previous point this is the subgroup preserving the relation of lying in the same column. Hence we got a chain

  $$S_m \times S_n \hookrightarrow S_m \ltimes S_n^m \to S_{mn},$$

  where the left inclusion corresponds to the diagonal morphism $S_n \to S_n^m$. Note that these morphisms exist even if $m$ is 0.

**Convention 1.1.3.** We need a global sign convention for total complexes:

For us, double complexes have commuting squares, in contrast to anti-commuting ones. For a double complex $A^{\bullet,\bullet}$ with differentials

$$\delta_h^{\bullet,\bullet} \colon A^{\bullet,\bullet} \to A^{\bullet+1,\bullet},$$
$$\delta_v^{\bullet,\bullet} \colon A^{\bullet,\bullet} \to A^{\bullet,\bullet+1}$$

we set

$$\mathrm{Tot}^\bullet(A^{\bullet,\bullet}) = \mathrm{Tot}^\bullet_{i,j}(A^{i,j}) := \left( \bigoplus_{i+j=\bullet} A^{i,j}, \delta^\bullet \right),$$

where $\delta^{i+j} = \delta_v^{i,j} + (-1)^j \delta_h^{i,j}$ on the summand $A^{i,j}$.

For complexes $(A^\bullet, \delta_A^\bullet), (B^\bullet, \delta_B^\bullet) \in C^b(\mathcal{C})$ over a preadditive tensor category $\mathcal{C}$ we therefore set

$$(A^\bullet, \delta_A^\bullet) \otimes (B^\bullet, \delta_B^\bullet) := \left( \bigoplus_{i+j=\bullet} A^i \otimes B^j, \delta_{A\otimes B}^\bullet \right)$$

with

$$\delta_{A\otimes B}^{i+j}(a_i \otimes b_j) = a_i \otimes \delta_B b_j + (-1)^j \delta_A a_i \otimes b_j$$

for $a_i \in A^i, b_j \in B^j$.

This is the same sign convention as in Section 1.3.3 of [HM16], but opposite to the one used for the Künneth formula in [Hat02] (see in particular its Lemma 3.6).

**Convention 1.1.4.** It is readily checked that forming total complexes is associative. Hence we can without ambiguity define total complexes $\mathrm{Tot}^{\bullet}_{i,j,k}(A^{i,j,k})$ if triple complexes, or higher dimensional ones.

## 1.2 Diagrams

Let us quickly recapitulate a few terminologies regarding categories.

**Definition 1.2.1.** Let $\mathcal{C}$ be a category and let $\rightarrow$ be the category consisting of two objects 'source' $s$ and 'target' $t$, together with three morphisms $\mathrm{id}_s$, $\mathrm{id}_t$ and $s \rightarrow t$. We define the *arrow category* over $\mathcal{C}$ as the functor category $\mathcal{C}^{\rightarrow}$ which thus has as objects the morphisms $Y \rightarrow X$ in $\mathcal{C}$, simply denoted by $(Y|X)$, and commutative squares

$$\begin{array}{ccc} Y & \longrightarrow & X \\ \downarrow{\scriptstyle g} & & \downarrow{\scriptstyle f} \\ Y' & \longrightarrow & X' \end{array}$$

as morphisms $(Y|X) \rightarrow (Y'|X')$, simply denoted by $(g|f)$.

**Definition 1.2.2** (Quiver)**.** As usual, a *quiver* $Q$ is a generalization of a directed graph where multiple arrows between the same vertices as well as loops are allowed. Formally, a quiver is a functor $Q\colon \hat{Q} \rightarrow \mathrm{Set}$, where $\hat{Q}$ is the free quiver, which is the category with two objects $\mathcal{V}, \mathcal{E}$ and four morphisms $\mathrm{id}_{\mathcal{V}}$, $\mathrm{id}_{\mathcal{E}}$ and $s, t\colon \mathcal{E} \rightarrow \mathcal{V}$.

The sets $V(Q) := Q(\mathcal{V})$ and $E(Q) := Q(\mathcal{E})$ will be called the *vertices* and *edges* of $Q$, respectively.

A morphism of quivers is defined to be a natural transformation of functors.

**Remark 1.2.3.** Equivalently, a quiver consists of sets $V$ and $E$ as well as two functions $s, t\colon E \rightarrow V$ called *source* and *target.* The interpretation of a quiver as a graph corresponds to interpreting $e \in E$ as an edge from $s(e)$ to $t(e)$.

**Remark 1.2.4.** In Nori's original context and several others, a quiver is often called a 'diagram'. We wish to avoid this unnecessary confusion and the collision with the more common Definition 1.2.9 below. See also Remark 1.2.6.

**Definition 1.2.5** (Path category)**.** Let $Q$ be a quiver. Its *path category* $\mathcal{P}(Q)$ (also called the *free category* on $Q$) is the category with objects $V(Q)$. The morphisms $\mathrm{Hom}_{\mathcal{P}(Q)}(v, w)$ are finite formal chains

$$v = v_0 \xrightarrow{e_1} v_1 \xrightarrow{e_2} \dots \xrightarrow{e_n} v_n = w$$

of edges $e_i \in E(Q)$ such that $s(e_i) = v_{i-1}$ and $t(e_i) = v_i$. Composition is defined by concatenating chains and the identities are the trivial chains of length 0. There are by definition no further relations between the chains.

**Remark 1.2.6.** From the definition it is immediate that $\mathcal{P}$ is a left adjoint to the forgetful functor from small strict categories, i.e. both objects and morphisms forming sets, to the category of quivers.

**Definition 1.2.7** (Quiver representations)**.** A *representation* $R$ of a quiver $Q = (V, E, s, t)$ in a category $\mathcal{C}$ is defined by choosing for each $v \in V(Q)$ an object $R(v) \in \mathcal{C}$ and for each $e \in E(Q)$ a morphism $R(e)\colon R(s(e)) \to R(t(e))$.

**Remark 1.2.8.** In the context of Remark 1.2.6, a representation of the quiver $Q$ in the category $\mathcal{C}$ is a functor $\mathcal{P}(Q) \to \mathcal{C}$.

**Definition 1.2.9** (Diagrams)**.**

- A *finite category* is a category where the collection of objects and every hom-set is a finite set. Equivalently, the collection of all morphisms forms a finite set.

- An *acyclic category* is a category that, when interpreted as a graph, does not contain any directed loop except the identities. Formally, it is a small category $\mathcal{C}$ such that for all objects $A, B \in \mathcal{C}$ the non-emptiness of both hom-sets $\mathrm{Mor}_{\mathcal{C}}(A, B)$ and $\mathrm{Mor}_{\mathcal{C}}(B, A)$ implies that $A = B$ and $\mathrm{Mor}_{\mathcal{C}}(A, A) = \{\mathrm{id}_A\}$.

- An *almost acyclic category* is a small category $\mathcal{C}$ such that for all objects $A, B \in \mathcal{C}$ the non-emptiness of both hom-sets $\mathrm{Mor}_{\mathcal{C}}(A, B)$ and $\mathrm{Mor}_{\mathcal{C}}(B, A)$ implies that $A = B$ and that $\mathrm{Mor}_{\mathcal{C}}(A, A)$ is a finite group of automorphisms of $A$.

- An *initial segment* of an acyclic category $\mathcal{C}$ is a full subcategory $I \subseteq \mathcal{C}$ such that if $A \to B$ is a morphism of $\mathcal{C}$ with $B \in I$, then $A$ is in $I$. Dually, a *terminal segment* of an acyclic category $\mathcal{C}$ is a full subcategory $I \subseteq \mathcal{C}$ such that if $A \to B$ is a morphism of $\mathcal{C}$ with $A \in I$, then $B$ is in $I$.

- A *diagram in a category* $\mathcal{D}$ is a covariant functor $D : \mathcal{C} \to \mathcal{D}$, where $\mathcal{C}$ is a small category. A *subdiagram* is a diagram induced by restricting the functor to a subcategory of $\mathcal{C}$. By abuse of notation we will sometimes call quiver representations a diagram, by which we mean the induced diagram on its path category (cf. Remark 1.2.6).

- We call a diagram $D\colon \mathcal{C} \to \mathcal{D}$ *finite* (respectively *(almost) acyclic*) if the underlying category $\mathcal{C}$ is. We extend the notions of *initial* and *terminal segments* to diagrams $D$, them being the restriction of the functor $D$ to an initial or terminal segment of $\mathcal{C}$, respectively.

**Remark 1.2.10.** A quiver $Q$ is acyclic as a graph if and only if its path category $\mathcal{P}(Q)$ is an acyclic category.

**Definition 1.2.11.** Two diagrams $D_1: C_1 \to A$ and $D_2: C_2 \to A$ over the same category and a bifunctor $\otimes: X \times X \to X$ give by composition rise to a *product diagram* $D_1 \otimes D_2: C_1 \times C_2 \to X$.

## 1.3 Nori's diagram categories

**Convention 1.3.1.** We fix a base field $k \subseteq \mathbb{C}$. Recall that a *variety* over $k$ is a reduced and separated scheme of finite type over $k$. We denote the corresponding category by $\mathrm{Var}_k$. We also fix a noetherian ring $\Lambda$ of coefficients.

For a variety $X$ and a Zariski-closed $Y \hookrightarrow X$ we use the abbreviation $H^n_{\mathrm{sing}}(X, Y)$ for the singular cohomology $H^n_{\mathrm{sing}}(X^{\mathrm{an}}, Y^{\mathrm{an}}, \Lambda)$ of their complex analytifications.

The main step in the construction of Nori motives is the following Tannakian construction due to Nori that associates an abelian category to a quiver representation:

**Theorem 1.3.2.** *Let $R: Q \to \Lambda$-$\mathrm{Mod}$ be a quiver representation in the category of finitely generated $\Lambda$-modules, where $\Lambda$ is a noetherian ring.*

*Then there exists*

- *a $\Lambda$-linear abelian category $\mathcal{C}(Q, R)$,*
- *a $\Lambda$-linear faithful exact* forgetful functor $\omega_R: \mathcal{C}(Q, R) \to \Lambda$-$\mathrm{Mod}$,
- *a quiver representation $\widetilde{R}: Q \to \mathcal{C}(Q, R)$*

*which factors $R = \omega_R \circ \widetilde{R}$ and such that the given data are universal with respect to these properties:*

*Given*

- *another $\Lambda$-linear abelian category $\mathcal{A}$,*
- *a $\Lambda$-linear faithful exact functor $\omega_{\mathcal{A}}: \mathcal{A} \to \Lambda$-$\mathrm{Mod}$,*
- *a quiver representation $R': Q \to \mathcal{A}$*

*with $R = \omega_{\mathcal{A}} \circ R'$, then there exists a $\Lambda$-linear functor $a: \mathcal{C}(Q, R) \to \mathcal{A}$, which*

*is unique up to unique isomorphism, making the resulting diagram*

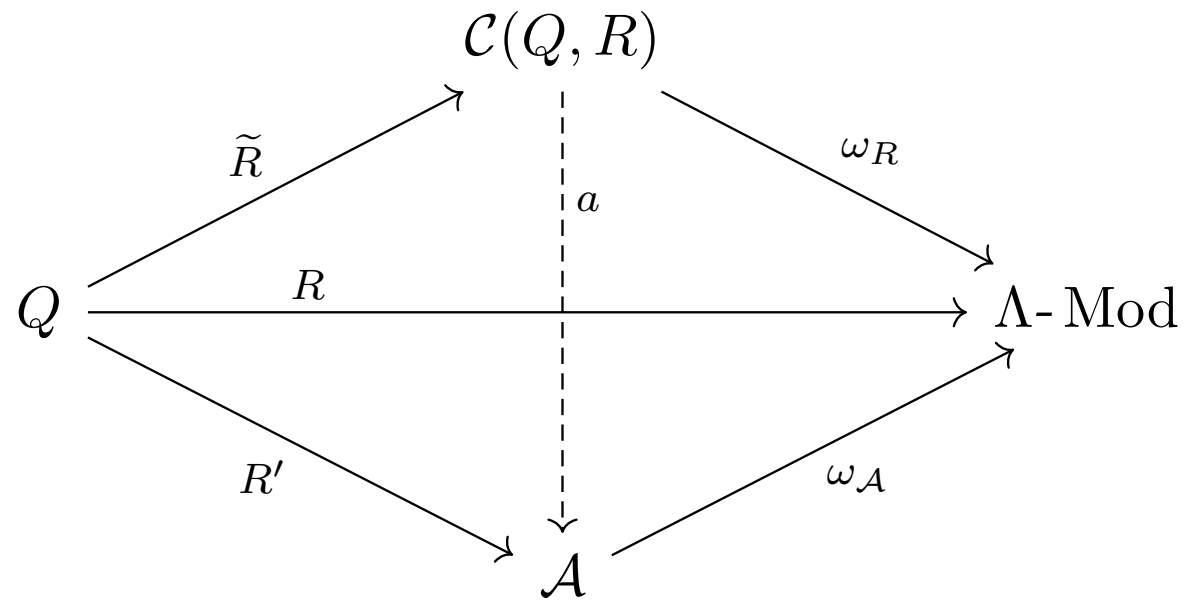

*commutative.*

*Proof.* This is Theorem 7.1.13 of [HM16]. □

**Definition 1.3.3.** The category $\mathcal{C}(Q,R)$ of Theorem 1.3.2 will be called *Nori's diagram category* associated to the quiver representation $R : Q \to \Lambda\text{-}\mathrm{Mod}$.

Nori's definition of motives now arises from a specific quiver representation, that of singular cohomology:

**Definition 1.3.4.** The quiver of *effective pairs* $\mathrm{Pairs}^{\mathrm{eff}}$ over a field $k \subseteq \mathbb{C}$ has as vertices the triples $(X, Y, i)$, where $X$ is a variety over $k$, $Y$ is a closed subvariety of $X$ and $i \in \mathbb{Z}$ is an integer. It has two types of edges for all integers $i$:

(a) for each morphism $f\colon X \to X'$ of varieties over $k$ and closed subvarieties $Y \hookrightarrow X$, $Y' \hookrightarrow X'$ with $f(Y) \subseteq Y'$ a *pullback edge*

$$f^*\colon (X', Y', i) \to (X, Y, i),$$

(b) for each triple $Z \hookrightarrow Y \hookrightarrow X$ of closed subvarieties of $X \in \mathrm{Var}_k$ a *boundary edge*

$$\delta\colon (Y, Z, i) \to (X, Y, i+1).$$

We then have a *singular cohomology representation*

$$\begin{aligned} H_{\mathrm{sing}}\colon \mathrm{Pairs}^{\mathrm{eff}} &\to \Lambda\text{-}\mathrm{Mod} \\ (X, Y, i) &\mapsto H^i_{\mathrm{sing}}(X, Y, \Lambda), \end{aligned}$$

interpreting the pullback edges as the contravariant functoriality and the boundary edges as the boundary maps of singular cohomology.

**Definition 1.3.5.** The category $\mathcal{MM}^{\mathrm{eff}}_{\mathrm{Nori}}(k,\Lambda)$ of *effective Nori motives* is the diagram category $\mathcal{C}(\mathrm{Pairs}^{\mathrm{eff}}, H_{\mathrm{sing}})$.

We denote the resulting quiver representation $\mathrm{Pairs}^{\mathrm{eff}} \to \mathcal{MM}^{\mathrm{eff}}_{\mathrm{Nori}}(k,\Lambda)$ by $H_{\mathrm{Nori}}$ and the image of $(X,Y,i)$ by $H^i_{\mathrm{Nori}}(X,Y) = H^i_{\mathrm{Nori}}(X,Y,\Lambda)$. The resulting faithful exact forgetful functor will be denoted by

$$\omega = \omega_{\mathrm{sing}} \colon \mathcal{MM}^{\mathrm{eff}}_{\mathrm{Nori}}(k,\Lambda) \to \Lambda\text{-}\mathrm{Mod}\,.$$

**Remark 1.3.6.** It is important to consider $k$ as a subfield of $\mathbb{C}$, i.e. equipped with a fixed embedding $k \hookrightarrow \mathbb{C}$. This happens because singular cohomology depends on this embedding, not only on $k$. When we talk about a field $k \subseteq \mathbb{C}$ we therefore always imply to use this embedding.

**Definition 1.3.7.** Let $(X,Y,i) \in \mathrm{Pairs}^{\mathrm{eff}}$.

- We call $(X,Y,i)$ *good* if
  - $H^j_{\mathrm{sing}}(X,Y,\Lambda) = 0$ for all $j \neq i$,
  - $H^i_{\mathrm{sing}}(X,Y,\Lambda)$ is a finitely generated projective $\Lambda$-module.
- We call $(X,Y,i)$ *very good* if
  - $(X,Y,i)$ is good,
  - $X$ is affine
  - $X\backslash Y$ is smooth over the base field $k$,
  - either $X = Y$ as well as $i = \dim(X)$ or $i = \dim(X) = \dim(Y)$.

This induces two smaller quivers $\mathrm{Good}^{\mathrm{eff}}$ and $\mathrm{VGood}^{\mathrm{eff}}$ of *effective good pairs* and *effective very good pairs*, respectively. The edges are given as before. We restrict the singular cohomology representation to those quivers, and denote the resulting quiver representations again by $H_{\mathrm{sing}}$.

**Theorem 1.3.8.** *Let $k \subseteq \mathbb{C}$ be a field and let $\Lambda$ be a noetherian ring.*

*Then the inclusions $\mathrm{VGood}^{\mathrm{eff}} \subseteq \mathrm{Good}^{\mathrm{eff}} \subseteq \mathrm{Pairs}^{\mathrm{eff}}$ induce equivalences*

$$\mathcal{C}(\mathrm{VGood}^{\mathrm{eff}}, H_{\mathrm{sing}}) \cong \mathcal{C}(\mathrm{Good}^{\mathrm{eff}}, H_{\mathrm{sing}}) \cong \mathcal{C}(\mathrm{Pairs}^{\mathrm{eff}}, H_{\mathrm{sing}})$$

*of diagram categories.*

*Proof.* For $\Lambda = \mathbb{Z}$ this is Theorem 9.2.21 of [HM16]. The general version follows by the same arguments, but using $\Lambda$-coefficients everywhere. See also the following Remark 1.3.9. □

**Remark 1.3.9.** The functor $\mathrm{C}^{\mathrm{eff}} \colon \mathrm{DM}^{\mathrm{eff}}_{\mathrm{gm}}(k,\Lambda) \to D^b(\mathcal{MM}^{\mathrm{eff}}_{\mathrm{Nori}}(k,\Lambda))$ we construct utilizes very good pairs (cf. Remark 4.3.4).

With some slight modifications, our arguments re-prove Theorem 1.3.8. The most important change would be to work with non-smooth varieties in Chapter 4, in particular Definition 4.5.3. This is possible due to Remark 4.7.4. It, however, overcomplicates things, as one does not need to consider finite correspondences, nor multivalued morphisms, at all to prove Remark 1.3.9. Indeed, the proof of Theorem 9.2.21 in [HM16] only considers morphisms. Nonetheless, in Remark 4.5.4 both arguments are found to be very close, which comes at no surprise because our Main Strategy 0.3.1 is based on [HM16].

**Proposition 1.3.10.** *Let $k \subseteq \mathbb{C}$ be a field and let $\Lambda$ be a noetherian ring.*

*(a) Every object of $\mathcal{MM}^{\mathrm{eff}}_{\mathrm{Nori}}(k,\Lambda)$ is a subquotient of one of the form*

$$\bigoplus_{i=1}^{r} H^{n_i}_{\mathrm{Nori}}(X_i, Y_i),$$

*where the $(X_i, Y_i, n_i)$ are very good pairs.*

*(b) The elements $H^{n}_{\mathrm{Nori}}(X, Y)$ corresponding to very good pairs $(X, Y, n)$ generate $\mathcal{MM}^{\mathrm{eff}}_{\mathrm{Nori}}(k,\Lambda)$ as an abelian category.*

*Proof.* Combine Proposition 7.1.16 of [HM16] with Theorem 1.3.8 above. □

**Theorem 1.3.11.** *Let $k \subseteq \mathbb{C}$ be a field and let $\Lambda$ be a noetherian ring.*

*Then the category $\mathcal{MM}^{\mathrm{eff}}_{\mathrm{Nori}}(k,\Lambda)$ carries a commutative tensor structure for which the forgetful functor*

$$\omega_{\mathrm{sing}} \colon \mathcal{MM}^{\mathrm{eff}}_{\mathrm{Nori}}(k,\Lambda) \to \Lambda\text{-}\mathrm{Mod}$$

*is a tensor functor.*

**Remark 1.3.12.** The tensor structure is induced by a larger, *graded* diagram that carries what [HM16] calls a *commutative product structure* which is reflected in a *(graded) multiplicative representation*:

The central idea is to add distinguished arrows mimicking the properties of tensor categories. We refer to op. cit., sections 8.1 and 9.3 for the details. A noteworthy subtlety is that one must sometimes extend to the path categories of the quivers involved, see Remark 8.1.6 of [HM16].

*Proof of Theorem 1.3.11.* For $\Lambda = \mathbb{Z}$ this is Theorem 9.1.5 of [HM16]. The given proof works for any noetherian coefficients as it relies solely on Theorem 8.1.5 of op. cit., which was shown for any such ring. See also Remark 9.1.8 of op. cit. □

## 1.4 Basics of relative cycles

This section explains parts of [VSF00, chapter 2], where the notions of relative cycles, pushforward and pullback are introduced. Unlike more classical approaches such as [Ful98], the pushforward and pullback are always defined, not only for proper and flat morphisms, respectively. This adds some interesting additional structures. It forces us to restrict to *universally integral relative cycles (of relative dimension* 0*)* in the terminology of [VSF00].

We try to keep our notations close to more recent treatments such as [Kel13] or [CD12]. In particular we will introduce the notion $f^{\circledast}$ to denote the pullback along a morphism $f$. It goes back to [Ivo07].

The following is taken from [CD12]:

**Definition 1.4.1.** A morphism $X \to Y$ is called *pseudo-dominant* if each generic point of $X$ is mapped to one of $Y$.

**Lemma 1.4.2.** *Let*

$$\begin{array}{ccc} X' & \xrightarrow{f \times_S S'} & S' \\ \downarrow & & \downarrow p \\ X & \xrightarrow{f} & S \end{array}$$

*be a cartesian square of noetherian schemes. If $f$ is pseudo-dominant and $p$ is universally open, then $f \times_S S'$ is pseudo-dominant.*

*Proof.* This is Lemma 3.3.7 of [VSF00, chapter 2]. □

**Definition 1.4.3** (Naive cycles)**.** Let $X \to S$ be a morphism of finite type between noetherian schemes. A *basic cycle* on $X|S$ is a point $\xi \in X$ such that the closure $\overline{\xi}$, equipped with the reduced induced structure, is finite and surjective onto an irreducible component of $S$, i.e. $\overline{\xi} \to S$ is finite and pseudo-dominant.

The abelian monoid $c^{\mathrm{nai,eff}}(X|S)$ of *effective naive cycles* is the free abelian monoid generated by the basic cycles on $X|S$. The abelian group $c^{\mathrm{nai}}(X|S)$ of *naive cycles* is the free abelian group generated by the basic cycles on $X|S$, i.e. the group completion of $c^{\mathrm{nai,eff}}(X|S)$.

**Definition 1.4.4** (Support)**.** Let $X \to S$ be a morphism of finite type between noetherian schemes.

The *support* $\mathrm{supp}(\alpha)$ of a cycle $\alpha \in c^{\mathrm{nai}}(X|S)$ is the finite union of closures $\overline{\xi}$ of all points $\xi$ that appear with a non-zero coefficient in $\alpha$. We equip it with the reduced induced structure, turning it into a closed subscheme of $X$.

**Remark 1.4.5.** The definition in [VSF00, chapter 2] is slightly different: they use proper subschemes of relative dimension 0 in their Definition 3.1.3. The following Lemma 1.4.6 states that this is equivalent.

**Lemma 1.4.6.** *Let $g\colon \Gamma \to X$ be a morphism of finite type between noetherian schemes.*

*Then $g$ is finite and pseudo-dominant if and only if it is proper and equidimensional of relative dimension* 0 *(cf. Definition 2.1.2 of [VSF00, chapter 2]).*

*Proof.* Assume that $g$ is proper and equidimensional of relative dimension 0. Thus $g$ is pseudo-dominant and closed, therefore its image is a finite union of irreducible components of $X$. Replacing $X$ by this finite union we may therefore assume that $g$ is surjective. By [EGA4-3], Théorème 8.11.1 we have to show that $g$ is quasi-finite. But the fibre over a point is a noetherian scheme of dimension 0, thus a finite discrete union of points.

If conversely $g$ is finite, then it is proper and of constant relative dimension 0. Adding pseudo-dominance makes it equidimensional of relative dimension 0. □

**Definition 1.4.7.** Let $X \to S$ be a morphism of finite type between noetherian schemes and let $Z \hookrightarrow X$ be a closed immersion. Also assume that the resulting morphism $f\colon Z \to S$ is finite.

We can *decompose* $Z$ into a cycle

$$\mathrm{cycl}_{X|S}(Z) = \sum_{i=1}^{r} m_i \zeta_i \in c^{\mathrm{nai}}(X|S),$$

where the $\zeta_i$ are the different generic points of $Z$ that lie over a generic point of $S$. The *multiplicity* $m_i \in \mathbb{N}$ is defined as the length of the ring $\mathcal{O}_{f(\zeta_i)\times_S Z,\zeta_i}$ as a module over itself.

**Remark 1.4.8.** The original definition of $\mathrm{cycl}_{X|S}(Z)$ in [VSF00, chapter 2] used the length of $\mathcal{O}_{Z,\zeta_i}$ instead. This is problematic as it is not functorial as elaborated in [Kel13] or [CD12]. Over a reduced base $S$ this change has no effect, especially not for our applications over a normal scheme in Chapter 3.

**Definition 1.4.9** (Pushforward)**.** Let $f\colon X \to Y$ be a morphism of schemes of finite type over a noetherian scheme $S$.

We define the *pushforward*

$$f_*\colon c^{\mathrm{nai}}(X|S) \to c^{\mathrm{nai}}(Y|S)$$

by sending a basic cycle $\xi \in X$ to $m \cdot f(\xi)$, where $m = [\kappa(\xi)\colon \kappa(f(\xi))]$ is the (automatically finite) degree of the induced extension of residue fields, and then extending linearly.

**Remark 1.4.10.** By [MVW06], Lemma 1.4, the pushforward indeed preserves the defining property of naive cycles.

**Definition 1.4.11** (Naive pullback)**.** Let $f\colon S' \to S$ be a morphism of noetherian schemes and let $X \to S$ be a morphism of finite type.

We define the *naive pullback*

$$f^*\colon c^{\mathrm{nai}}(X|S) \to c^{\mathrm{nai}}(X \times_S S'|S')$$

by sending a basic cycle $\xi \in c^{\mathrm{nai}}(X|S)$ to

$$\mathrm{cycl}_{X\times_S S'|S'}(\mathrm{supp}(\xi) \times_S S')$$

and extending linearly.

**Remark 1.4.12.** This is not yet the correct version of the pullback. We will see that the naive and genuine pullbacks agree in many cases (cf. Propositions 1.5.4 and 1.5.6), but not always.

We now briefly explain some of the central ideas behind the definition of relative cycles and their pullback. This is mostly based on Kelly's elaboration [Kel13] of the original work [VSF00, chapter 2], but we add another shortcut on top of his.

The naive pullback of Definition 1.4.11 is not functorial, in particular it induces the wrong multiplicities even in some moderately well-behaved situations. This is reflected in the classical theory by Serre's Tor formula, which exposes the deficit between the naive and the correct multiplicity.

On the other hand, there are some cases where we expect the naive approach to give the correct results. For example, pullbacks along flat morphisms or of flat cycles should be the naive ones. Another example are birational morphisms between bases, which have no effect on the generic points and hence should leave cycles unchanged.

This leads to the following approach:

**Definition 1.4.13.** Let $X \to S$ be a morphism of finite type between noetherian schemes, let $k$ be a field and let $\kappa\colon \mathrm{Spec}(k) \to S$ be a $k$-point in $S$. Also let

$$\alpha = \sum_{i=1}^{r} a_i \xi_i \in c^{\mathrm{nai}}(X|S)$$

be a naive cycle.

A *good factorization* of $\kappa$ with respect to $\alpha$ is a factorization

$$\kappa = p \circ i\colon \mathrm{Spec}(k) \xrightarrow{i} S' \xrightarrow{p} S$$

such that

- $p$ is proper and birational, i.e. an abstract blow-up,
- the proper transforms $\overline{\zeta_i \times_S S'}$ of the supports $\mathrm{supp}(\zeta_i)$ are flat over $S'$.

By our prior reasoning we expect the correct pullback $\kappa^{\circledast}\alpha$ to be given by naive pullbacks as $i^*p^*\alpha$ (cf. Definition 1.4.17 below). Somewhat problematic is that good factorizations exist in general only after a finite extension:

**Proposition 1.4.14** (Platification theorem)**.** *Let $X \to S$ be a morphism of finite type between noetherian schemes. Let $\kappa\colon \operatorname{Spec}(k) \to S$ be a $k$-point, where $k$ is a field, and let $\alpha \in c^{\mathrm{nai}}(X|S)$ be a naive relative cycle.*

*Then there exists a finite field extension $L|k$ such that the induced $L$-point $\operatorname{Spec}(L) \to S$ has a good factorization.*

*Proof.* This is Theorem 2.2.17 of [Kel13], which goes back to Théorème 5.2.2 of [RG71] □

Worse, we do not know if the pullback $i^*p^*\alpha$ is independent of the choice of a good factorization. This justifies the following technical definition, which is an adaptation of the original Definition 3.1.3 of [VSF00, chapter 2] and is based on Definition 2.4.1 of [Kel13].

**Definition 1.4.15.** Let $X \to S$ be a morphism of finite type between noetherian schemes. Then the set $c(X|S)$ of *relative cycles* consists of those naive relative cycles $\alpha \in c^{\mathrm{nai}}(X|S)$ which satisfy the following property:

For every point $s \in S$ there exists a naive cycle

$$\alpha_{|s} \in c^{\mathrm{nai}}(X \times_S s|s)$$

which satisfies $q^*\alpha_{|s} = i^*p^*\alpha$ for every commutative square

$$\begin{array}{ccc} \operatorname{Spec}(k) & \xrightarrow{\;i\;} & S' \\ \downarrow{\scriptstyle q} & & \downarrow{\scriptstyle p} \\ s & \hookrightarrow & S \end{array}$$

involving a field $k$ and a good factorization $p \circ i$ with respect to $\alpha$.

An *effective relative cycle* is an effective naive cycle satisfying the above property. We denote the corresponding set by $c^{\mathrm{eff}}(X|S)$.

**Remark 1.4.16.** Given $s$, a square as in Definition 1.4.15 exists by Proposition 1.4.14. Thus the naive cycle $\alpha_{|s}$ is unique because the naive pullback $q^*$ is injective.

By the nature of this definition we have already found a candidate $\alpha_{|s}$ for the pullback along the inclusion of a point. Thus, according to our self-set guidelines, the following Definition 1.4.17 should be the correct one:

**Definition 1.4.17.** Let $f\colon T \to S$ be a morphism of noetherian schemes and let $X \to S$ be a morphism of finite type.

Then the *pullback* $f^{\circledast}\colon c(X|S) \to c(X \times_S T|T)$ is defined in three steps as follows:

1) If $f\colon T = \{s\} \hookrightarrow S$ is the inclusion of a point $s \in S$, we let $f^{\circledast}(\alpha) := \alpha_{|s}$ be the naive cycle of Definition 1.4.15.

2) If $T = \{t\}$ is a single reduced point, we let $i\colon s := f(t) \hookrightarrow S$ be the inclusion and denote the resulting morphism $t \to s$ also by $f$. Using the naive pullback of Definition 1.4.11 we then set $f^{\circledast}(\alpha) = f^* i^{\circledast}\alpha$.

3) In general, we $z_k\colon \zeta_k \hookrightarrow T$, $k \in \{1, 2, \ldots, r\}$, be the generic points of $T$. Then we define

$$f^{\circledast}\alpha := \sum_{k=1}^{r} (f \circ z_k)^{\circledast}\alpha.$$

Note that $f^{\circledast}\alpha$ as just defined is only a formal linear combination of points in $X \times_S T$. One thus has to check that it is indeed an element of $c(X \times_S T|T)$. Furthermore, we have to show that this gives rise to a reasonable theory of cycles, in particular that relative cycles are preserved by pushforwards and that the set $c(X|S)$ is closed under addition. For all of this and more we point to Section 2.4 of [Kel13] or to sections 3.3 to 3.6 of [VSF00, chapter 2].

Under certain assumptions on the base we have simpler descriptions for the set of relative cycles:

**Proposition 1.4.18.** *Let $X \to S$ be a morphism of finite type between noetherian schemes.*

*(a) If $S$ is normal, then $c(X|S) = c^{\mathrm{eff}}(X|S) \otimes_{\mathbb{N}} \mathbb{Z}$, i.e. the abelian group $c(X|S)$ is generated by its effective elements.*

*(b) If $S$ is regular, then $c(X|S) = c^{\mathrm{nai}}(X|S)$.*

*Proof.* The second part is Corollary 3.4.6 of [VSF00, chapter 2].

Let now $S$ be normal, or more generally geometrically unibranch, and let $\alpha \in c(X|S)$. In the notations of loc. cit. we have by definition $c(X|S) = c_{equi}(X/S, 0) = z(X/S, 0) \cap PropCycl_{equi}(X/S, 0) \subseteq Cycl(X/S, 0)$. Thus combining Corollary 3.4.4 of [VSF00, chapter 2] with its Proposition 3.3.14 shows that there exist a positive integer $N$ and effective $\alpha_+, \alpha_- \in c^{\mathrm{eff}}(X|S)$ with $N\alpha = \alpha_+ - \alpha_-$. Therefore $\alpha + \alpha_-$ is effective, proving the first part. □

**Definition 1.4.19.** Let $f\colon X \to S$ be a finite and pseudo-dominant morphism of noetherian schemes and let $Y \to X$ be a morphism of finite type. Then every naive cycle $\alpha \in c^{\mathrm{nai}}(Y|X)$ is trivially, as a formal linear combination of points, also a naive cycle over $S$. We denote this linear map as $f_{\#}\colon c^{\mathrm{nai}}(Y|X) \to c^{\mathrm{nai}}(Y|S)$ and call it *forgetting*.

**Definition 1.4.20.** Let $X \to S$ be a morphism of finite type between noetherian schemes and let $i\colon Y \hookrightarrow X$ be an immersion of a subscheme. Also let $\alpha \in c^{\mathrm{nai}}(X|S)$ be a naive cycle with support $\mathrm{supp}(\alpha) \subseteq Y$.

Then we can also interpret $\alpha$ as a naive cycle in $c^{\mathrm{nai}}(Y|S)$ by using the same linear combination of points. We denote this *restriction* by $\mathrm{res}_Y(\alpha)$.

**Lemma 1.4.21.** *If in Definition 1.4.20 we have a relative cycle $\alpha \in c(X|S)$, then $\mathrm{res}_Y(\alpha) \in c(Y|S)$, i.e. its restriction to $Y$ is a relative cycle as well.*

*Proof.* It suffices to check that the properties of Definition 1.4.15 are preserved, which follows easily from $Y \hookrightarrow X$ being an immersion. $\square$

Lastly, there is one more essential construction in Section 3.7 of [VSF00, chapter 2]:

**Definition 1.4.22.** Let $Z \to Y$ and $f\colon Y \to X$ be morphisms of finite type between noetherian schemes. Let $\alpha = \sum_{k=1}^{r} a_k \alpha_k \in c^{\mathrm{nai}}(Y|X)$ be a naive cycle and let $\beta \in c(Z|Y)$ be a relative cycle. Denote the inclusions $\mathrm{supp}(\alpha_k) \hookrightarrow Y$ by $i_k$.

The *correspondence map*

$$\mathrm{Cor}\colon c(Z|Y) \otimes c^{\mathrm{nai}}(Y|X) \to c^{\mathrm{nai}}(Z|X)$$

is defined by

$$\mathrm{Cor}(\beta, \alpha) := \sum_{k=1}^{r} a_k \, (i_k \times_Y Z)_* \, (f \circ i_k)_\# i_k^{\circledast}(\beta) \in c^{\mathrm{nai}}(Z|X).$$

**Remark 1.4.23.** Let us illustrate this with a commutative diagram:

$$\begin{array}{ccccc}
\mathrm{supp}(\alpha_k) \times_Y Z & \longrightarrow & \mathrm{supp}(\alpha_k) & & \\
\Big\downarrow {\scriptstyle \mathrm{pr}_Z = i_k \times_Y Z} & & \Big\downarrow {\scriptstyle i_k} & \searrow & \\
Z & \longrightarrow & Y & \xrightarrow{f} & X.
\end{array}$$

The correspondence morphism gives by construction a naive cycle. That it preserves relative cycles is much more difficult:

**Proposition 1.4.24.** *Let $Z \to Y \to X$ be morphisms of finite type between noetherian schemes. Let $\alpha \in c(Y|X)$ and $\beta \in c(Z|Y)$ be relative cycles.*

*Then $\mathrm{Cor}(\beta, \alpha)$ lies in $c(Z|X)$, i.e. is a relative cycle.*

*Proof.* This is [VSF00, chapter 2], Theorem 3.7.3. $\square$

**Remark 1.4.25.** It is immediate from the definitions that pushforward, pullback, forgetting, restriction and the correspondence map preserve being effective.

## 1.5 Recollection of functorialities of cycles

We offer a collection of the most important relations between relative cycles, almost all of which are based on [VSF00, chapter 2] and chapter 2 of [Kel13].

**Proposition 1.5.1** (Functoriality of pushforward)**.** *Let $S$ be a noetherian scheme and let*

$$X \xrightarrow{f} Y \xrightarrow{g} Z \longrightarrow S$$

*be morphisms of finite type.*

*Then*

$$(g \circ f)_* = g_* \circ f_*$$

*as maps $c^{\mathrm{nai}}(X|S) \to c^{\mathrm{nai}}(Z|S)$.*

*Proof.* This is obvious from the definition. □

**Proposition 1.5.2** (Functoriality of pullback)**.** *Let $X \to S$ be a morphism of finite type and let*

$$S'' \xrightarrow{q} S' \xrightarrow{p} S$$

*be morphisms of noetherian schemes.*

*Then*

$$(p \circ q)^{\circledast} = q^{\circledast} \circ p^{\circledast}$$

*as maps $c(X|S) \to c(X \times_S S''|S'')$.*

*Proof.* This is Lemma 2.4.7 of [Kel13]. □

**Proposition 1.5.3** (Compatibility between pushforward and pullback)**.** *Let*

$$\begin{array}{ccccc} Y' & \xrightarrow{p \times_S S'} & X' & \longrightarrow & S' \\ \downarrow & & \downarrow & & \downarrow f \\ Y & \xrightarrow{p} & X & \longrightarrow & S \end{array}$$

*be two cartesian squares of noetherian schemes and assume that the horizontal morphisms are of finite type. Then we have the compatibility*

$$f^{\circledast} \circ p_* = (p \times_S S')_* \circ f^{\circledast}$$

*as maps $c(Y|S) \to c(X'|S')$.*

*Proof.* This is part 2 of Proposition 3.6.2 of [VSF00, chapter 2]. □

**Proposition 1.5.4** (Pseudo-dominant pullback)**.** *Let*

$$\begin{array}{ccc} X' & \longrightarrow & S' \\ \big\downarrow & & \big\downarrow{\scriptstyle f} \\ X & \longrightarrow & S \end{array}$$

*be a cartesian square of noetherian schemes where the horizontal morphisms are of finite type. Assume that $f$ is pseudo-dominant. Let $\alpha \in c(X|S)$ be a relative cycle, written as $\alpha = \sum_{i=1}^{r} a_i \alpha_i$ for $a_i \in \mathbb{Z}$ and $\alpha_i \in X$.*

*Then the pullback $f^{\circledast}(\alpha)$ is equal to the naive pullback*

$$f^*(\alpha) = \sum_{i=1}^{r} a_i \, \mathrm{cycl}_{X'|S'}(\mathrm{supp}(\alpha_i) \times_S S').$$

*Proof.* This is Proposition 2.4.6 of [Kel13]. See also Proposition 3.3.12 of [VSF00] for a special case. □

**Proposition 1.5.5.** *Let*

$$\begin{array}{ccc} X' & \longrightarrow & S' \\ \big\downarrow & & \big\downarrow{\scriptstyle f} \\ X & \longrightarrow & S \end{array}$$

*be a cartesian square of noetherian schemes where the horizontal morphisms are of finite type. Also assume that $f$ is pseudo-dominant and dominant.*

*Then the pullback $f^{\circledast}$ is injective as a map $c(X|S) \to c(X'|S')$.*

*Proof.* This follows by combining Proposition 1.5.4 above with Lemma 2.2.12 of [Kel13]. □

**Proposition 1.5.6** (Flat decomposition)**.** *Let*

$$\begin{array}{ccc} X' & \longrightarrow & S' \\ \big\downarrow & & \big\downarrow{\scriptstyle f} \\ X & \longrightarrow & S \end{array}$$

*be a cartesian square of noetherian schemes and assume that the horizontal morphisms are of finite type. Let $Z \hookrightarrow X$ be a closed subscheme that is flat and finite over $S$.*

*Then the decomposition $\mathrm{cycl}_{X|S}(Z)$ is a relative cycle and its pullback satisfies*

$$f^{\circledast}(\mathrm{cycl}_{X|S}(Z)) = f^*(\mathrm{cycl}_{X|S}(Z)) = \mathrm{cycl}_{X'|S'}(Z \times_S S').$$

*Proof.* This is Proposition 2.5.1 of [Kel13]. □

**Proposition 1.5.7** (Compatibility between Cor and pushforward on the right)**.** *Let*

$$\begin{array}{ccccc} Y' & \longrightarrow & X' & & \\ \big\downarrow{\scriptstyle p\times_X Y} & & \big\downarrow{\scriptstyle p} & & \\ Y & \longrightarrow & X & \longrightarrow & S \end{array}$$

*be a cartesian square of schemes of finite type over a noetherian scheme $S$. Then for all $\alpha \in c^{\mathrm{nai}}(X'|S)$ and $\beta \in c(Y|X)$ we have the equality*

$$\mathrm{Cor}(\beta, p_*(\alpha)) = (p \times_X Y)_* \, \mathrm{Cor}(p^{\circledast}\beta, \alpha)$$

*in $c^{\mathrm{nai}}(Y|S)$.*

*Proof.* This is Lemma 3.7.1 of [VSF00, chapter 2]. □

**Proposition 1.5.8** (Compatibility between Cor and pushforward on the left)**.** *Let*

$$Z \xrightarrow{q} Y \longrightarrow X \longrightarrow S$$

*be morphisms of schemes of finite type over a noetherian scheme $S$. Then for all $\alpha \in c^{\mathrm{nai}}(X|S)$ and $\beta \in c(Z|X)$ we have the equality*

$$\mathrm{Cor}(q_*\beta, \alpha) = q_* \, \mathrm{Cor}(\beta, \alpha)$$

*in $c(Y|S)$.*

*Proof.* We may by linearity assume that $\alpha$ is basic. Let $f : X \to S$ denote the given morphism and let $i\colon \mathrm{supp}(\alpha) \hookrightarrow X$ be the inclusion. Then, using the functoriality of pushforward (cf. Proposition 1.5.1) and Proposition 1.5.3, we get

$$\begin{aligned} q_* \, \mathrm{Cor}(\beta, \alpha) =& q_* \, (i \times_X Z)_* \, (f \circ i)_\# i^{\circledast}\beta = \\ =& \, (q \circ (i \times_X Z))_* \, (f \circ i)_\# i^{\circledast}\beta = \\ =& \, (i \times_X Y)_* \, (\mathrm{supp}(\alpha) \times_X q)_* \, (f \circ i)_\# i^{\circledast}\beta = \\ =& \, (i \times_X Y)_* \, (f \circ i)_\# \, (\mathrm{supp}(\alpha) \times_X q)_* \, i^{\circledast}\beta = \\ =& \, (i \times_X Y)_* \, (f \circ i)_\# i^{\circledast} q_*\beta = \mathrm{Cor}(q_*\beta, \alpha) \end{aligned}$$

as required. □

**Proposition 1.5.9** (Compatibility between Cor and pullback)**.** *Let*

$$\begin{array}{ccccc} Y' & \longrightarrow & X' & \longrightarrow & S' \\ \big\downarrow & & \big\downarrow{\scriptstyle f\times_S X} & & \big\downarrow{\scriptstyle f} \\ Y & \longrightarrow & X & \longrightarrow & S \end{array}$$

*be two cartesian squares of schemes of finite type over a noetherian scheme* $S$*. Then for all* $\alpha \in c(X|S)$ *and* $\beta \in c(Y|X)$ *we have the equality*

$$f^{\circledast} \operatorname{Cor}(\beta, \alpha) = \operatorname{Cor}((f \times_S X)^{\circledast}\beta, f^{\circledast}\alpha)$$

*in* $c(Y'|S')$*.*

*Proof.* This is Proposition 3.7.3 of [VSF00]. □

**Proposition 1.5.10** (Associativity of Cor)**.** *Let*

$$Z \longrightarrow Y \longrightarrow X \longrightarrow S$$

*be morphisms of schemes of finite type over a noetherian scheme* $S$*. Then for all* $\alpha \in c^{\mathrm{nai}}(X|S)$*,* $\beta \in c(Y|X)$ *and* $\gamma \in c(Z|Y)$ *we have the equality*

$$\operatorname{Cor}(\operatorname{Cor}(\gamma, \beta), \alpha) = \operatorname{Cor}(\gamma, \operatorname{Cor}(\beta, \alpha))$$

*in* $c^{\mathrm{nai}}(Z|S)$*.*

*Proof.* This is Proposition 3.7.7 of [VSF00]. □

## 1.6 Finite correspondences

**Definition 1.6.1** (Finite correspondences)**.** Let $\Lambda$ be a ring, let $S$ be a noetherian scheme and let $X, Y \in \mathrm{Sch}_S$. A *finite* $\Lambda$*-correspondence* over $S$ from $X$ to $Y$ is an element of $c(X \times_S Y|X) \otimes_{\mathbb{Z}} \Lambda$. We will denote them as zigzagged arrows $X \rightsquigarrow Y$.

If $\Lambda = \mathbb{Z}$ we simply call it a *finite correspondence* over $S$. An *effective finite correspondence* over $S$ is a finite correspondence with non-negative coefficients, i.e. an element of $c^{\mathrm{eff}}(X \times_S Y|X)$. We also call them finite $\mathbb{N}$-correspondences.

**Remark 1.6.2.** Definition 1.6.1 gives an ad-hoc definition of finite correspondences with coefficients in a ring $\Lambda$. Because it extends everything $\Lambda$-linearly we will meet no problems when using the results of Sections 1.4 and 1.5.

This definition is, however, not the most general one: one could use $\Lambda$-coefficients from the beginning. This can, depending on the base and the ring of coefficients, lead to larger modules $c_\Lambda(X|S) \supsetneq c(Y|X) \otimes_{\mathbb{Z}} \Lambda$ of cycles. We point the interested reader to [CD12], in particular to Remark 9.1.3, which contains a full treatise of this generality.

Most of our results will concern regular schemes $X$, where by Proposition 1.4.18 every basic cycle is already a relative cycle. Hence the definitions agree in this case, identifying $c_\Lambda(X \times_S Y|X)$ and $c(X \times_S Y|X) \otimes_{\mathbb{Z}} \Lambda$ with the free $\Lambda$-module $c^{\mathrm{nai}}(X \times_S Y|X) \otimes_{\mathbb{Z}} \Lambda$ generated by the basic cycles.

**Remark 1.6.3.** Let $f : Y \to X$ be a morphism in $\mathrm{Sch}_S$, where $S$ is a noetherian scheme. Let $\alpha \in c(Y|X)$. Then $\alpha$ can be interpreted as a finite correspondence $X \rightsquigarrow Y$ via the pushforward $((f, \mathrm{id}_Y)_S)_* \, \alpha \in c(X \times_S Y|X)$ along the graph $(f, \mathrm{id}_Y)_S \colon Y \hookrightarrow X \times_S Y$.

**Proposition 1.6.4.** *Let $\Lambda$ be a ring and let $S$ be a noetherian scheme. Let $I$ and $J$ be finite sets. Also let $X_i$, $i \in I$, and $Y_j$, $j \in J$, in $\mathrm{Sch}_S$ be given.*

*Then there is a natural additive decomposition*

$$\mathrm{SchCor}_{S,\Lambda}\left(\coprod_{i\in I} X_i, \coprod_{j\in J} Y_j\right) \cong \bigoplus_{\substack{i\in I\\ j\in J}} \mathrm{SchCor}_{S,\Lambda}(X_i, Y_j)$$

*induced by the decomposition*

$$\left(\coprod_{i\in I} X_i\right) \times_S \left(\coprod_{j\in J} Y_j\right) = \coprod_{\substack{i\in I\\ j\in J}} (X_i \times_S Y_j).$$

*Proof.* The decomposition into the individual $X_i$ follows from taking the pullback and Proposition 1.5.4. The decomposition along the $Y_j$ follows easily from Lemma 3.6.4 of [VSF00, chapter 2]. □

**Definition 1.6.5** (Composition of finite correspondences)**.** Let $\Lambda$ be a ring, let $S$ be a noetherian scheme and let $X, Y, Z \in \mathrm{Sch}_S$.

If $\alpha \colon X \rightsquigarrow Y$ and $\beta \colon Y \rightsquigarrow Z$ are finite $\mathbb{Z}$-correspondences, we define, following p. 59 of [VSF00, chapter 2], their composition as

$$\beta \circ \alpha := \left(\mathrm{pr}^{XYZ}_{XZ}\right)_* \mathrm{Cor}\left(\left(\mathrm{pr}^{XY}_{Y}\right)^{\circledast} \beta, \alpha\right),$$

which is easily seen to be a finite $\mathbb{Z}$-correspondence $X \rightsquigarrow Z$. This operation is then extended $\Lambda$-linearly to one of finite $\Lambda$-correspondences.

**Remark 1.6.6.** When working with smooth schemes over a regular base, the composition reduces to Serre's Tor formula and hence to classical intersection theory, as demonstrated by [VSF00, chapter 2], Lemma 3.5.9.

**Theorem 1.6.7** (Composition of finite correspondences is associative:)**.** *Let $S$ be a noetherian scheme and let $W, X, Y, Z$ be schemes of finite type over $S$. Let $\alpha \colon W \rightsquigarrow X$, $\beta \colon X \rightsquigarrow Y$ and $\gamma \colon Y \rightsquigarrow Z$ be finite correspondences. Then*

$$(\gamma \circ \beta) \circ \alpha = \gamma \circ (\beta \circ \alpha).$$

*Proof.* This can be found as Lemme 2.1.2 in [Ivo05] and as Proposition 9.1.7 in [CD12]. Let us elaborate the argument, which is a calculation using the results of Section 1.5:

$$
\begin{aligned}
&\gamma \circ (\beta \circ \alpha) = \\
&\overset{a)}{=} (\mathrm{pr}^{WYZ}_{WZ})_* \operatorname{Cor}\left((\mathrm{pr}^{WY}_{Y})^{\circledast}\gamma, (\mathrm{pr}^{WXY}_{WY})_* \operatorname{Cor}\left((\mathrm{pr}^{WX}_{X})^{\circledast}\beta, \alpha\right)\right) = \\
&\overset{b)}{=} (\mathrm{pr}^{WYZ}_{WZ})_*(\mathrm{pr}^{WXYZ}_{WYZ})_* \operatorname{Cor}\left((\mathrm{pr}^{WXY}_{WY})^{\circledast}(\mathrm{pr}^{WY}_{Y})^{\circledast}\gamma, \operatorname{Cor}\left((\mathrm{pr}^{WX}_{X})^{\circledast}\beta, \alpha\right)\right) = \\
&\overset{c)}{=} (\mathrm{pr}^{WXYZ}_{WZ})_* \operatorname{Cor}\left((\mathrm{pr}^{WXY}_{Y})^{\circledast}\gamma, \operatorname{Cor}\left((\mathrm{pr}^{WX}_{X})^{\circledast}\beta, \alpha\right)\right) = \\
&\overset{d)}{=} (\mathrm{pr}^{WXYZ}_{WZ})_* \operatorname{Cor}\left(\operatorname{Cor}\left((\mathrm{pr}^{WXY}_{Y})^{\circledast}\gamma, (\mathrm{pr}^{WX}_{X})^{\circledast}\beta\right), \alpha\right) = \\
&\overset{e)}{=} (\mathrm{pr}^{WXZ}_{WZ})_*(\mathrm{pr}^{WXYZ}_{WXZ})_* \operatorname{Cor}\left(\operatorname{Cor}\left((\mathrm{pr}^{WXY}_{XY})^{\circledast}(\mathrm{pr}^{XY}_{Y})^{\circledast}\gamma, (\mathrm{pr}^{WX}_{X})^{\circledast}\beta\right), \alpha\right) = \\
&\overset{f)}{=} (\mathrm{pr}^{WXZ}_{WZ})_* \operatorname{Cor}\left((\mathrm{pr}^{WXYZ}_{WXZ})_*(\mathrm{pr}^{WX}_{X})^{\circledast} \operatorname{Cor}\left((\mathrm{pr}^{XY}_{Y})^{\circledast}\gamma, \beta\right), \alpha\right) = \\
&\overset{g)}{=} (\mathrm{pr}^{WXZ}_{WZ})_* \operatorname{Cor}\left((\mathrm{pr}^{WX}_{X})^{\circledast}(\mathrm{pr}^{XYZ}_{XZ})_* \operatorname{Cor}\left((\mathrm{pr}^{XY}_{Y})^{\circledast}\gamma, \beta\right), \alpha\right) = \\
&\overset{h)}{=} (\gamma \circ \beta) \circ \alpha.
\end{aligned}
$$

Here we suppressed forgetting (cf. Definition 1.4.19) and used the following:

a) Definition 1.4.22 of Cor,

b) compatibility of Cor$(-,-)$ with pushforward on the right (Proposition 1.5.7),

c) functoriality of pushforward and pullback (Proposition 1.5.1 and Proposition 1.5.2),

d) associativity of Cor$(-.-)$ (Proposition 1.5.10),

e) functoriality of pushforward and pullback (Propositions 1.5.1 and 1.5.2),

f) compatibility of Cor$(-,-)$ with pushforward on the left (Proposition 1.5.8) and with pullback (Proposition 1.5.9),

g) compatibility between pushforward and pullback (Proposition 1.5.3),

h) Definition 1.4.22 of Cor.

□

**Definition 1.6.8.** Let $S$ be a noetherian scheme and let $\Lambda$ be a ring.

The $\Lambda$-linear category $\mathrm{SchCor}_{S,\Lambda} = \mathrm{SchCor}(S, \Lambda)$ *of finite correspondences* over $S$ with coefficients in $\Lambda$ has the same objects as $\mathrm{Sch}_S$, i.e. schemes that are separated and of finite type over $S$.

Its morphisms are the finite correspondences

$$\mathrm{SchCor}_{S,\Lambda}(X, Y) = c(X \times_S Y | X) \otimes_{\mathbb{Z}} \Lambda$$

of Definition 1.6.1, and we continue to use zigzagged arrows $X \rightsquigarrow Y$ to denote them. Their composition is given by Definition 1.6.5.

The identities $\mathrm{id}_X = \mathrm{cycl}_{X\times_S X|X}(\Delta_X)$ are obtained by decomposing the diagonals, as is easily verified.

By Remark 1.4.25 we also have the category $\mathrm{SchCor}_{S,\mathbb{N}} = \mathrm{SchCor}^{\mathrm{eff}}(S)$ of effective finite correspondences.

**Definition 1.6.9.** Let $S$ be a noetherian scheme. We define a natural embedding $[-]\colon \mathrm{Sch}_S \to \mathrm{SchCor}_S$:

On objects it is the identity. A morphism $f\colon X \to Y$ over $S$ is sent to its graph $\Gamma_f$, understood as the finite correspondence

$$\mathrm{cycl}_{X\times_S Y|X}(\Gamma_f) \in c(X \times_S Y|X) = \mathrm{SchCor}_S(X, Y)$$

via Proposition 1.5.6. Functoriality is part of Proposition 1.6.10 below.

By abuse of notation we will often omit $[-]$ and hence use $f$ to denote both the morphism as well as the associated finite correspondence.

**Proposition 1.6.10.** *Let $S$ be a noetherian scheme.*

*If one or both finite correspondences come from a morphism of schemes, then composition simplifies:*

(a) *If $\alpha\colon X \rightsquigarrow Y$ is a finite correspondence in* $\mathrm{SchCor}_S$ *and $g\colon Y \to Z$ is a morphism in* $\mathrm{Sch}_S$*, then*

$$[g] \circ \alpha = (X \times_S g)_* \, \alpha.$$

(b) *If $f\colon X \to Y$ is a morphism in* $\mathrm{Sch}_S$ *and $\beta\colon Y \rightsquigarrow Z$ is a finite correspondence in* $\mathrm{SchCor}_S$*, then*

$$\beta \circ [f] = f^{\circledast}\beta.$$

(c) *If $f\colon X \to Y$ and $g\colon Y \to Z$ are morphisms in* $\mathrm{Sch}_S$*, then*

$$[g] \circ [f] = [g \circ f].$$

*Proof.* The first two are part of Proposition 2.5.7 of [Kel13]. The last one is immediate from the first part and Proposition 1.5.6. □

**Proposition 1.6.11.** *Let $S$ be a noetherian scheme and let $f\colon X \to Y$ be a morphism between schemes of finite type over $S$.*

(a) *If $f$ is an immersion, then $[f]$ is a monomorphism in* $\mathrm{SchCor}_S$.

(b) *If $f$ is dominant and pseudo-dominant (cf. Definition 1.4.1), then $[f]$ is an epimorphism in* $\mathrm{SchCor}_S$.

*Proof.* The pushforward of cycles along an immersion is injective. Hence part (a) follows from Proposition 1.6.10 (a).

We get part (b) from Proposition 1.5.5 and Proposition 1.6.10 (b). □

## 1.7 The exterior product

Following the end of Section 3.7 of [VSF00, chapter 2] we have:

**Definition 1.7.1.** Let $S$ be a noetherian scheme. Let $z_1 \colon Z_1 \to S$ and $Z_2 \to X$ be schemes of finite type over $S$.

We define the *exterior product* of cycles over $S$ as the linear map

$$c(Z_1|S) \otimes c(Z_2|S) \to c(Z_1 \times_S Z_2|S)$$

given by

$$\alpha \otimes_S \beta := \mathrm{Cor}(z_1^{\circledast}\beta, \alpha).$$

**Lemma 1.7.2.** *The exterior product of cycles over $S$ is associative and commutative.*

**Remark 1.7.3.** The commutativity is somewhat tricky because the definition of $\mathrm{Cor}(-,-)$ is highly asymmetric. We loosely follow the argument given for Lemme 1.4.16 of [Ivo05].

*Proof of Lemma 1.7.2.* Let $x \colon X \to S$, $y \colon Y \to S$ and $Z \to S$ be schemes of finite type over $S$ and let $\alpha \in c(X|S)$, $\beta \in c(Y|S)$ and $\gamma \in c(Z|S)$. Then we get associativity by

$$\begin{aligned}
(\alpha \otimes_S \beta) \otimes_S \gamma &\overset{\text{a)}}{=} \mathrm{Cor}((x \times_S y)^{\circledast}\gamma, \mathrm{Cor}(x^{\circledast}\beta, \alpha)) = \\
&\overset{\text{b)}}{=} \mathrm{Cor}(\mathrm{Cor}((x \times_S y)^{\circledast}\gamma, x^{\circledast}\beta), \alpha)) = \\
&\overset{\text{c)}}{=} \mathrm{Cor}(\mathrm{Cor}((x \times_S Y)^{\circledast}y^{\circledast}\gamma, x^{\circledast}\beta), \alpha)) = \\
&\overset{\text{d)}}{=} \mathrm{Cor}(x^{\circledast}\,\mathrm{Cor}(y^{\circledast}\gamma, \beta), \alpha) = \\
&\overset{\text{e)}}{=} \alpha \otimes_S (\beta \otimes_S \gamma),
\end{aligned}$$

where we used:

a) Definition 1.7.1 of the exterior product,

b) associativity of Cor (Proposition 1.5.10),

c) functoriality of pullback (Proposition 1.5.2),

d) compatibility of pullback and Cor (Proposition 1.5.9),

e) Definition 1.7.1 of the exterior product.

Due to the injectivity of Proposition 1.5.5 it suffices, using Proposition 1.5.9, to check the commutativity after a pullback to each of the generic points of $S$. Hence it suffices to check it for $S = \mathrm{Spec}(k)$ the spectrum of a field. By Proposition 1.4.18 and the linearity of the exterior product we can then assume that $\alpha$ and $\beta$ are both basic, given by closed points $a \in X$ and $b \in Y$.

Using Proposition 1.5.6 and the Definitions we see that the exterior product is given by $a \otimes_S b = \mathrm{cycl}_{X \times_S Y|S}(a \times_S b)$. The right hand side is symmetric, finishing the proof. □

**Definition 1.7.4.** Let $S$ be a noetherian scheme. Let $\alpha_1 \colon X_1 \rightsquigarrow Y_1$ and $\alpha_2 \colon X_2 \rightsquigarrow Y_2$ be finite correspondences between schemes of finite type over $S$.

We define the *tensor product* of finite correspondences

$$- \otimes - \colon \ \mathrm{SchCor}_S(X_1, Y_1) \otimes \mathrm{SchCor}_S(X_2, Y_2) \to \mathrm{SchCor}_S(X_1 \times_S X_2, Y_1 \times_S Y_2)$$

by

$$\alpha_1 \otimes \alpha_2 := \left( \left(\mathrm{pr}^{X_1 X_2}_{X_1}\right)^{\circledast} \alpha_1 \right) \otimes_{X_1 X_2} \left( \left(\mathrm{pr}^{X_1 X_2}_{X_2}\right)^{\circledast} \alpha_2 \right).$$

This linearly extends to finite correspondences with coefficients in a ring $\Lambda$, inducing a *tensor product* on $\mathrm{SchCor}_{S,\Lambda}$ again denoted by $- \otimes -$.

**Remark 1.7.5.** Explicitly, we thus have

$$\alpha_1 \otimes \alpha_2 = \mathrm{Cor} \left( \left(\mathrm{pr}^{X_1 X_2 Y_1}_{X_2}\right)^{\circledast} \alpha_2, \left(\mathrm{pr}^{X_1 X_2}_{X_1}\right)^{\circledast} \alpha_1 \right).$$

**Proposition 1.7.6.** *Let $S$ be a noetherian scheme and let $\Lambda$ be a ring.*

*The tensor product of finite $\Lambda$-correspondences turns* $\mathrm{SchCor}_{S,\Lambda}$ *into a $\Lambda$-linear symmetric tensor category.*

*Proof.* The $\Lambda$-linearity is trivial. The symmetry and associativity of $- \otimes -$ follow immediately from Lemma 1.7.2. The functoriality of $- \otimes -$ is a lengthy but straightforward calculation. We wish to point towards Section 9.2 of [CD12] for a detailed proof. □

**Remark 1.7.7.** The embedding $\mathrm{Sch}_S \hookrightarrow \mathrm{SchCor}_S$ of Definition 1.6.9 is compatible with the product structures: if $f_1 \colon X_1 \to Y_1$ and $f_2 \colon X_2 \to Y_2$ are morphisms over $S$, then $[f_1 \times_S f_2] = [f_1] \otimes [f_2]$, as is easily checked.

## 1.8 Degree

**Definition 1.8.1** (Degree)**.** Let $S$ be an irreducible noetherian scheme with generic point $\eta$ and let $f\colon X \to S$ be a morphism of finite type. Then clearly $c(S|S) = c^{\mathrm{nai}}(S|S) = \mathbb{Z} \cdot \eta \cong \mathbb{Z}$.

Taking the pushforward $f_*$ hence defines a *degree map*

$$\deg = \deg_{X|S}\colon c^{\mathrm{nai}}(X|S) \to \mathbb{Z}.$$

**Remark 1.8.2.** Assume that $S$ is noetherian and integral. Due to its Definition 1.4.9 the degree of a basic cycle $\alpha$ is the generic degree of the finite surjective morphism $\mathrm{supp}(\alpha) \to S$.

Conversely one can define the degree of a morphism $X \to S$ of finite type, where $S$ is noetherian and irreducible, as $\deg(\mathrm{cycl}_{X|S}(X))$. It gets identified with the usual definition by taking the pullback $S^{\mathrm{red}} \to S$ along the reduction.

**Definition 1.8.3** (Constant degree)**.** Let $S$ be a noetherian scheme with irreducible components $s_i\colon S_i \hookrightarrow S$, $i \in \{1, \ldots, r\}$. Let $f\colon X \to S$ be a scheme of finite type over $S$ and let $\alpha \in c^{\mathrm{nai}}(X|S)$.

We say that $\alpha$ has *constant degree* $n$ if each of the naive pullbacks $s_i^*\alpha \in c^{\mathrm{nai}}(X \times_S S_i|S_i)$ has degree $n$. Under these circumstances we write $\deg(\alpha) = n$.

**Remark 1.8.4.** If $\alpha \in c(X|S)$ is a not necessarily naive relative cycle, then $s_i^*\alpha = s_i^{\circledast}\alpha$ by Proposition 1.5.4. If furthermore $\sigma_i$ is the inclusion of the generic point of $S_i$ into $S$, we thus observe that

$$\deg(s_i^*\alpha) = \deg(s_i^{\circledast}\alpha) = \deg(\sigma_i^{\circledast}\alpha).$$

**Lemma 1.8.5.** *Let $S$ be a noetherian scheme and let $g\colon Y \to X$ be a morphism between schemes of finite type over $S$. Let $\alpha \in c^{\mathrm{nai}}(Y|S)$ be of constant degree.*

*Then $g_*\alpha \in c^{\mathrm{nai}}(X|S)$ has constant degree*

$$\deg(g_*\alpha) = \deg(\alpha).$$

*Proof.* The naive pullback and the degree only depend on the behaviour over the generic points of $S$, hence we may assume that $S$ is irreducible. Now the statement is simply the functoriality of pushforward, i.e. Proposition 1.5.1. □

**Lemma 1.8.6.** *Let $p\colon S' \to S$ be a morphism between irreducible noetherian schemes and let $X \to S$ be a morphism of finite type. Let $\alpha \in c(X|S)$.*

*Then*

$$\deg(p^{\circledast}\alpha) = \deg(\alpha).$$

*Proof.* We use $f$ to denote the given morphism $X \to S$.

Let $\eta$ and $\eta'$ be the generic points of $S$ and $S'$, respectively. We also use them to denote the respective generators of $c(S|S)$ and $c(S'|S')$. As

$$p^{\circledast}\eta = p^{\circledast}\operatorname{cycl}_{S|S}(S) = \operatorname{cycl}_{S'|S'}(S') = \eta'$$

by Proposition 1.5.6, we get from Proposition 1.5.3 that indeed

$$\deg(p^{\circledast}\alpha)\cdot\eta' = (f\times_S S')_* p^{\circledast}\alpha = p^{\circledast} f_*\alpha = p^{\circledast}(\deg(\alpha)\cdot\eta) = \deg(\alpha)\cdot\eta'.$$

□

**Lemma 1.8.7.** *Let $p\colon S' \to S$ be a morphism of noetherian schemes and assume that every connected component of $S$ is irreducible. Let $X \to S$ be a morphism of finite type and let $\alpha \in c(X|S)$ be of constant degree.*

*Then $p^{\circledast}\alpha$ has constant degree*

$$\deg(p^{\circledast}\alpha) = \deg(\alpha).$$

*Proof.* Decomposing $S$ into connected components we may assume that $S$ is connected, hence by assumption irreducible. Let $s_i\colon S_i' \hookrightarrow S'$ be the inclusion of an irreducible component. Then Remark 1.8.4, Proposition 1.5.2 and Lemma 1.8.6 imply that

$$\deg(s_i^* p^{\circledast}\alpha) = \deg(s_i^{\circledast} p^{\circledast}\alpha) = \deg((p\circ s_i)^{\circledast}\alpha) = \deg(\alpha)$$

is independent of the irreducible component $S_i$ of $S'$, showing the result. □

The following is trivial from Remark 1.8.2:

**Lemma 1.8.8.** *Let $f\colon X \to S$ be a finite and surjective morphism between irreducible noetherian schemes and let $Y \to X$ be a morphism of finite type. Let $\alpha \in c^{\mathrm{nai}}(Y|X)$ and denote the generic degree of $f$ by $m$.*

*Then $f_{\#}\alpha \in c^{\mathrm{nai}}(Y|S)$ (cf. Definition 1.4.19) has degree $m\cdot\deg(\alpha)$.*

**Lemma 1.8.9.** *Let $Z \to Y \to X$ be morphisms of finite type between noetherian schemes and assume that the connected components of $X$ and $Y$ are irreducible. Let $\alpha \in c(Y|X)$ and $\beta \in c(Z|Y)$ be both of constant degree.*

*Then $\operatorname{Cor}(\beta,\alpha)$ is of constant degree*

$$\deg(\operatorname{Cor}(\beta,\alpha)) = \deg(\beta)\deg(\alpha).$$

*Proof.* By a pullback along the inclusion of a connected, i.e. irreducible, component of $X$ we may by Proposition 1.5.9 assume that $X$ is irreducible.

Note that both sides of the claimed equality are defined even if $\alpha$ is only a naive cycle. Hence by linearity we may assume that $\alpha$ is basic. Let $i\colon \operatorname{supp}(\alpha) \hookrightarrow Y$ be the inclusion and let $f\colon Y \to X$ be the given morphism. Hence by Definition 1.4.22

$$\operatorname{Cor}(\beta,\alpha) = (i \times_Y Z)_* (f \circ i)_\# i^{\circledast}(\beta).$$

The degree remains unchanged under $(i \times_Y Z)_*$ by Lemma 1.8.5 and under $i^{\circledast}$ by Lemma 1.8.7. The degree of $f \circ i\colon \operatorname{supp}(\alpha) \to X$ is by Remark 1.8.2 the degree of $\alpha$. Hence the lemma now follows from Lemma 1.8.8. □

Due to Definition 1.7.1 we get from Lemmas 1.8.7 and 1.8.9:

**Lemma 1.8.10.** *Let $Z_1$, $Z_2$ be schemes of finite type over a noetherian scheme $S$. Assume that the connected components of $S$ and $Z_1$ are irreducible. Let $\alpha_1 \in c(Z_1|S)$ and $\alpha_2 \in c(Z_2|S)$ be relative cycles of constant degree.*

*Then their exterior product $\alpha_1 \otimes_S \alpha_2 \in c(Z_1 \times_S Z_2|S)$ has constant degree*

$$\deg(\alpha_1 \otimes_S \alpha_2) = \deg(\alpha_1)\deg(\alpha_2).$$

**Definition 1.8.11.** Let $S$ be a noetherian scheme and let $\alpha\colon X \rightsquigarrow Y$ be a finite correspondence. We say that $\alpha$ is of *constant degree* if the underlying relative cycle $\alpha \in c(X \times_S Y|X)$ is. In this case we define its *degree*

$$\deg(\alpha) := \deg_{X \times_S Y|X}(\alpha)$$

as the degree of the underlying relative cycle.

**Remark 1.8.12.** We give an alternative description inspired by [BV08]:

Let $x\colon X \to S$ and $y\colon Y \to S$ be the structure morphisms into the base and assume that $\alpha\colon X \rightsquigarrow Y$ is of constant degree. Hence Proposition 1.6.10 shows that

$$\deg(\alpha)\cdot[x] = \left(\operatorname{pr}_X^{XY}\right)_* \alpha = [y] \circ \alpha$$

as finite correspondences $X \rightsquigarrow S$, which can therefore be used to calculate or define the degree.

**Lemma 1.8.13.** *Let $S$ be a noetherian scheme and let $X$, $Y$ and $Z$ be schemes of finite type over $S$. Assume that the connected components of $Y$ are irreducible. Let $\alpha\colon X \rightsquigarrow Y$ and $\beta\colon Y \rightsquigarrow Z$ be finite correspondences of constant degree over $S$.*

*Then $\beta \circ \alpha$ is of constant degree*

$$\deg(\beta \circ \alpha) = \deg(\beta)\deg(\alpha).$$

*Proof.* Let $\iota_i\colon X_i \hookrightarrow X$, $i \in \{1, 2, \dots, r\}$, be the inclusions of the irreducible components of $X$. By Remark 1.8.4 and Proposition 1.6.10 it suffices to check the lemma for the individual $\alpha \circ \iota_i\colon X_i \rightsquigarrow Y$ instead of $\alpha$. Hence we may assume that $X$ is irreducible.

If $\alpha$ is basic, then unravelling the definitions shows that

$$\beta \circ \alpha = (j \times_S Z)_* j_\# i^\circledast(\beta),$$

where $i\colon \operatorname{supp}(\alpha) \to Y$ is the projection to $Y$ and $j\colon \operatorname{supp}(\alpha) \to X$ is the projection to $X$. Thus the lemma can be proven by the same arguments as Lemma 1.8.9. □

**Remark 1.8.14.** Remark 1.8.12 offers a different approach, at least for irreducible $X$ and $Y$:

Denote the three structure morphism from $X, Y, Z$ to $S$ as $x, y, z$, respectively. Then

$$\deg(\beta \circ \alpha) \cdot [x] = [z] \circ \beta \circ \alpha = \deg(\beta) \cdot [y] \circ \alpha = \deg(\beta)\deg(\alpha) \cdot [x].$$

**Lemma 1.8.15.** *Let $S$ be a noetherian scheme and let $\alpha_1\colon X_1 \rightsquigarrow Y_1$ and $\alpha_2\colon X_2 \rightsquigarrow Y_2$ be finite correspondences over $S$ of constant degree. Assume that the connected components of $X_1$, $X_2$ and $X_1 \times_S X_2$ are irreducible.*

*Then the tensor product $\alpha_1 \otimes \alpha_2\colon X_1 \times_S X_2 \rightsquigarrow Y_1 \times_S Y_2$ is of constant degree*

$$\deg(\alpha_1 \otimes \alpha_2) = \deg(\alpha_1)\deg(\alpha_2).$$

*Proof.* This is a simple consequence of Lemmas 1.8.7 and 1.8.10. □

## 1.9 Voevodsky motives

Let us elaborate and extend the definitions of Voevodsky's geometric motives as found in [VSF00, chapter 5].

**Definition 1.9.1.** Let $S$ be a noetherian scheme. The category $\mathrm{SmCor}(S, \Lambda)$ is defined as the full subcategory of $\mathrm{SchCor}(S, \Lambda)$ consisting of schemes which are smooth over $S$.

**Definition 1.9.2.** Let $S$ be a noetherian scheme. We define the category $\mathrm{SmCor}^{\mathrm{aff}}(S, \Lambda)$ as the full subcategory of $\mathrm{SmCor}(S, \Lambda)$ consisting of those schemes which are smooth and affine over $S$.

**Definition 1.9.3.** Let $S$ be a noetherian scheme. Consider the homotopy category $K^b(\mathrm{SmCor}(S, \Lambda))$, which for every smooth $X \to S$ contains the following complexes:

- Homotopy Invariance, HI: the complex

$$\mathbb{A}^1_S \times_S X = \mathbb{A}^1_X \xrightarrow{\mathrm{pr}_X} X,$$

- Mayer-Vietoris-Nisnevich, MVN: for every Nisnevich cd-square (cf. Definition 5.3.1)

$$\begin{array}{ccc} U \times_X V & \hookrightarrow & V \\ \downarrow & & \downarrow g \\ U & \overset{f}{\hookrightarrow} & X \end{array}$$

the complex

$$U \times_X V \xrightarrow{\mathrm{pr}_U - \mathrm{pr}_V} U \sqcup V \xrightarrow{f \sqcup g} X.$$

The localization

$$\underline{\mathrm{DM}}^{\mathrm{eff}}_{\mathrm{gm}}(S, \Lambda) := K^b(\mathrm{SmCor}(S, \Lambda))/\langle \mathrm{HI}, \mathrm{MV} \rangle$$

at the thick subcategory generated by the above two types of complexes is called *Voevodsky's effective geometric pre-motives.*

Its pseudo-abelian envelope

$$\mathrm{DM}^{\mathrm{eff}}_{\mathrm{gm}}(S, \Lambda)$$

is called *Voevodsky's effective geometric motives.*

**Remark 1.9.4.** In the original work [VSF00, chapter 5], $S$ is the spectrum of a perfect field and only open covers $X = U \cup V$ were allowed for the complexes of type MVN. As is widely known, Theorem 3.2.6 of [VSF00, chapter 5] and Theorem 1.9.10 below imply that both variants produce equivalent categories.

**Definition 1.9.5.** Following [CD12] in generalizing [VSF00, chapter 5], we consider the $\Lambda$-linear abelian category

$$\mathrm{PreShv}(\mathrm{SmCor}(S, \Lambda))$$

of *presheaves with* $\Lambda$*-transfers*, i.e. contravariant $\Lambda$-linear functors

$$\mathrm{SmCor}(S, \Lambda) \to \Lambda\text{-}\mathrm{Mod}$$

It contains the full subcategory $\mathrm{Shv}_{\mathrm{Nis}}(\mathrm{SmCor}(S, \Lambda))$ of *Nisnevich sheaves with* $\Lambda$*-transfers*: those presheaves which are sheaves in the Nisnevich topology (see Definition 5.3.1) when restricted to the subcategory $\mathrm{Sm}_S \subseteq \mathrm{SmCor}(S, \Lambda)$ via the embedding of Definition 1.6.9.

We have the covariant Yoneda embedding

$$L = L_\Lambda \colon \mathrm{SmCor}(S, \Lambda) \to \mathrm{PreShv}(\mathrm{SmCor}(S, \Lambda)),$$

$X \mapsto \mathrm{SmCor}_{S,\Lambda}(-, X)$. The following is then easily checked:

**Lemma 1.9.6.** *The Yoneda embedding*

$$L\colon \mathrm{SmCor}(S,\Lambda) \to \mathrm{PreShv}(\mathrm{SmCor}(S,\Lambda))$$

*lands inside* $\mathrm{Shv}_{\mathrm{Nis}}(\mathrm{SmCor}(S,\Lambda))$.

Hence the Yoneda embedding extends to a functor

$$KL\colon K^b(\mathrm{SmCor}(S,\Lambda)) \to K^b(\mathrm{Shv}_{\mathrm{Nis}}(\mathrm{SmCor}(S,\Lambda))$$

and by localization gives a functor

$$DL\colon K^b(\mathrm{SmCor}(S,\Lambda)) \to D(\mathrm{Shv}(\mathrm{SmCor}(S,\Lambda))).$$

**Definition 1.9.7.** *Voevodsky's category of (unbounded) effective motivic complexes* is the Bousfield localization

$$\mathrm{DM}^{\mathrm{eff}}(S,\Lambda) := D(\mathrm{Shv}_{\mathrm{Nis}}(\mathrm{SmCor}(S,\Lambda)))/\langle DL(\mathrm{HI})\rangle^{\oplus}$$

at the localizing subcategory generated by the images of the complexes $\mathbb{A}^1_X \to X$ of Definition 1.9.3 under the Yoneda embedding.

By the standard properties of Bousfield localizations we can identify $\mathrm{DM}^{\mathrm{eff}}(S,\Lambda)$ with the $\mathbb{A}^1$-invariant Nisnevich sheaves with transfers. We also have the following description:

**Proposition 1.9.8.** *Let $S$ be regular and of finite dimension.*
*Then there is a natural equivalence*

$$\begin{aligned}\mathrm{DM}^{\mathrm{eff}}(S,\Lambda) =& D(\mathrm{Shv}_{\mathrm{Nis}}(\mathrm{SmCor}(S,\Lambda)))/\langle DL(\mathrm{HI})\rangle^{\oplus} \cong \\ \cong& D(\mathrm{PreShv}(\mathrm{SmCor}(S,\Lambda)))/\langle DL(\mathrm{HI}), DL(\mathrm{MVN})\rangle^{\oplus}\end{aligned}$$

*of Bousfield localizations of derived categories of (pre)sheaves.*

*Proof.* This follows from the expositions in [Voe10a] and [Voe10b], as elaborated in Lemme 3.3.5 of [Ivo05] for $\Lambda = \mathbb{Z}$. It can also be found in [CD12]. □

**Remark 1.9.9.** The same arguments show that the equivalence of Proposition 1.9.8 holds even without taking localizations at $DL(\mathrm{HI})$, i.e. we have an equivalence

$$D(\mathrm{Shv}_{\mathrm{Nis}}(\mathrm{SmCor}(S,\Lambda))) \cong D(\mathrm{PreShv}(\mathrm{SmCor}(S,\Lambda)))/\langle DL(\mathrm{MVN})\rangle^{\oplus}.$$

**Theorem 1.9.10.** *Let $S$ be a scheme of finite dimension and let $\Lambda$ be a ring. Assume that $S$ is regular or that $\Lambda$ is flat over $\mathbb{Z}$.*
*Then there exists a fully faithful embedding*

$$i_\Lambda\colon \mathrm{DM}^{\mathrm{eff}}_{\mathrm{gm}}(S,\Lambda) \to \mathrm{DM}^{\mathrm{eff}}(S,\Lambda)$$

*such that the square*

$$\begin{array}{ccc} K^b(\mathrm{SmCor}(S,\Lambda)) & \xrightarrow{DL} & D(\mathrm{Shv}(\mathrm{SmCor}(S,\Lambda))) \\ \downarrow & & \downarrow \\ \mathrm{DM}^{\mathrm{eff}}_{\mathrm{gm}}(S,\Lambda) & \xrightarrow{i_\Lambda} & \mathrm{DM}^{\mathrm{eff}}(S,\Lambda) \end{array}$$

*commutes. It identifies* $\mathrm{DM}^{\mathrm{eff}}_{\mathrm{gm}}(S,\Lambda)$ *with the compact objects of* $\mathrm{DM}^{\mathrm{eff}}(S,\Lambda)$.

*Proof.* For $\Lambda = \mathbb{Z}$, this follows quite formally from Proposition 1.9.8 and standard results on Bousfield localizations, as seen found in Proposition 4.1.23 of [Ivo05]. The case where $\Lambda$ is flat over $\mathbb{Z}$ then follows by applying the exact functor $- \otimes_{\mathbb{Z}} \Lambda$ to all hom-sets.

The version where $S$ is regular can, recalling Remark 1.6.2, be found as Theorem 11.1.13 in [CD12]. □

Note that if $S = \mathrm{Spec}(k)$ is the spectrum of a field and $\Lambda = \mathbb{Z}$, then most of Theorem 1.9.10 can already be found as Theorem 3.2.6 of [VSF00, chapter 5].

**Remark 1.9.11.** In Theorem 1.9.10, the assumption that $S$ is regular or that $\Lambda$ is flat over $\mathbb{Z}$ was only due to Remark 1.6.2 and can be removed when using the correct definitions as hinted there. We will, however, never need Theorem 1.9.10 in such greater generality.

**Remark 1.9.12.** One can repeat the definition of $\mathrm{DM}^{\mathrm{eff}}_{\mathrm{gm}}(S,\Lambda)$, but restrict both the objects and relations to affine schemes. This defines an affine version $\mathrm{DM}^{\mathrm{eff,aff}}_{\mathrm{gm}}(S,\Lambda)$ of Voevodsky's geometric motives.

We expect the obvious functor $\mathrm{DM}^{\mathrm{eff,aff}}_{\mathrm{gm}}(S,\Lambda) \to \mathrm{DM}^{\mathrm{eff}}_{\mathrm{gm}}(S,\Lambda)$ to be an equivalence of categories. This, however, is not as obvious as one might expect. The difficulty lies in the functor being fully faithful, or equivalently in the informal question whether the affine Nisnevich cd-squares induce all the relations coming from arbitrary ones.

We hope to tackle this question in the future. Recent results by Asok, Hoyois and Wendt (see [AHW15]) already go in this direction. A promising strategy of proof would be to show an affine analogue of Theorem 1.9.10, as it is quite easy to compare $\mathrm{Shv}_{\mathrm{Nis}}(\mathrm{SmCor}^{\mathrm{aff}}(S,\Lambda))$ with $\mathrm{Shv}_{\mathrm{Nis}}(\mathrm{SmCor}(S,\Lambda))$.

Such a result would significantly simplify the proofs of Chapter 7.

## 1.10 Tensor-localizations

**Definition 1.10.1.** Let $\mathcal{A}$ be a tensor category and let $L \in \mathcal{A}$. The *tensor-localization* $\mathcal{A}[L(-1)]$ of $\mathcal{A}$ at $L$ is the 2-colimit

$$\mathcal{A}[L^{\otimes -1}] := 2\text{-}\varinjlim_{n \in \mathbb{N}} \mathcal{A}(n),$$

where $\mathcal{A}(n) = \mathcal{A}$ and the transition functors for $m \geq n$ are

$$- \otimes L^{\otimes(m-n)} \colon \mathcal{A}(n) \to \mathcal{A}(m).$$

**Remark 1.10.2.** Less abstractly, the objects of $\mathcal{A}[L^{\otimes -1}]$ are pairs

$$A(n) := (A, n) \in \mathcal{A} \times \mathbb{Z}.$$

One could informally interpret them as $A \otimes L^{\otimes n}$, even for negative $n$, in the sense of duals. The morphisms in $\mathcal{A}[L^{\otimes -1}]$ are described by

$$\mathrm{Hom}_{\mathcal{A}[L^{\otimes -1}]}(A(m), B(n)) = \varinjlim_{k \geq -m, -n} \mathrm{Hom}_{\mathcal{A}}(A(m+k), B(n+k)).$$

**Definition 1.10.3.** Let $\mathcal{A}$ be a tensor category and let $L \in \mathcal{A}$.

We call $L$ *symmetric* if the morphism $L \otimes L \to L \otimes L$ swapping the two tensor factors equals the identity.

If $L$ is symmetric we extend the tensor product of $\mathcal{A}$ to $\mathcal{A}[L^{\otimes -1}]$ by defining $A(m) \otimes B(n) := (A \otimes B)(m+n)$ for all $A, B \in \mathcal{A}$ and $m, n \in \mathbb{Z}$.

**Remark 1.10.4.** The condition on $L$ may at first glance seem unnecessary. It comes from the suggestive notion $A(n)$, which hides the required compatibilities. We point the interested reader to Appendix A of [MVW06], or to [Tho80] for a similar problem in K-theory.

**Lemma 1.10.5.** *Let $\mathcal{A}$ be an abelian tensor category and let $L \in \mathcal{A}$. Assume that $L$ is flat (cf. Definition 7.4.1).*

*Then there exists a natural equivalence*

$$D^b(\mathcal{A}[L^{\otimes -1}]) \cong D^b(\mathcal{A})[(L[0])^{\otimes -1}]$$

*of triangulated categories.*

*Proof.* The lemma holds at the level of homotopy categories:

$$K^b(\mathcal{A}[L^{\otimes -1}]) \cong K^b(\mathcal{A})[(L[0])^{\otimes -1}].$$

It descends to the localizations at quasi-isomorphisms by the flatness of $T$. □

**Definition 1.10.6** (Tate/Lefschetz motives)**.** We have the following special objects in our categories of motives:

- Let $k$ be a field and let $\Lambda$ be a ring. We call

$$\Lambda_{\mathrm{DM}}(1) := \{\mathbb{P}^1_k \to \mathrm{Spec}(k)\}[-2] \in \mathrm{DM}^{\mathrm{eff}}_{\mathrm{gm}}(k, \Lambda),$$

with $\mathbb{P}^1$ in degree 2, the *Tate motive* of $\mathrm{DM}^{\mathrm{eff}}_{\mathrm{gm}}(k, \Lambda)$.

- Let $k \subseteq \mathbb{C}$ be a field and let $\Lambda$ be a noetherian ring. We call

$$\Lambda_{\mathrm{Nori}}(-1) := H^1_{\mathrm{Nori}}(\mathbb{G}_m) \in \mathcal{MM}^{\mathrm{eff}}_{\mathrm{Nori}}(k, \Lambda)$$

the *Lefschetz motive* of $\mathcal{MM}^{\mathrm{eff}}_{\mathrm{Nori}}(k, \Lambda)$.

**Remark 1.10.7.** The naming conventions in Definition 1.10.6 follow those found in [VSF00] and [HM16]. Usually, the Tate motive is understood to be dual to the Lefschetz motive. In our case, singular cohomology and hence Definition 1.3.5 are contravariant with respect to morphisms of varieties, while $\mathrm{DM}^{\mathrm{eff}}_{\mathrm{gm}}(S, \Lambda)$ is clearly covariant, justifying this choice.

**Proposition 1.10.8.** *Let $k$ be a field and let $\Lambda$ be a ring. Then the Tate motive $\Lambda_{\mathrm{DM}}(1)$ is symmetric.*

*Proof.* For $\Lambda = \mathbb{Z}$ this is Corollary 2.1.5 of [VSF00, chapter 5]. The general case follows by a change of coefficients. □

**Remark 1.10.9.** The restriction to the base $S = \mathrm{Spec}(k)$ is not necessary for this to hold.

**Definition 1.10.10.** Let $k$ be a field and let $\Lambda$ be a ring. The triangulated tensor category $\mathrm{DM}_{\mathrm{gm}}(k, \Lambda)$ of *Voevodsky's geometric motives* is the tensor-localization of $\mathrm{DM}^{\mathrm{eff}}_{\mathrm{gm}}(k, \Lambda)$ at the Tate motive $\Lambda_{\mathrm{DM}}(1)$.

**Proposition 1.10.11.** *Let $k \subseteq \mathbb{C}$ be a field and let $\Lambda$ be a noetherian ring. Then the Lefschetz motive $\Lambda_{\mathrm{Nori}}(-1)$ is symmetric.*

*Proof.* This can be checked after applying the faithful forgetful tensor functor $\omega_{\mathrm{sing}} \colon \mathcal{MM}^{\mathrm{eff}}_{\mathrm{Nori}}(k, \Lambda) \to \Lambda\text{-}\mathrm{Mod}$. As $\omega_{\mathrm{sing}}(\Lambda_{\mathrm{Nori}}(-1)) \cong \Lambda$, the statement hence reduces to the symmetry of $\Lambda \in \Lambda\text{-}\mathrm{Mod}$, which is clearly satisfied. □

**Lemma 1.10.12.** *The Tate motive $\Lambda_{\mathrm{DM}}(1) \in \mathrm{DM}^{\mathrm{eff}}_{\mathrm{gm}}(k, \Lambda)$ is isomorphic to $(\mathbb{G}_m \to \mathrm{Spec}(k))[-1]$, where $\mathbb{G}_m$ is placed in degree 1.*

*Proof.* We can by Theorem 1.9.10 check this in $\mathrm{DM}(k, \Lambda)$. Then it reduces to Lemme 4.1.5 of [Ivo05], whose proof works without any further change for arbitrary coefficients. □

**Definition 1.10.13.** Let $k \subseteq \mathbb{C}$ be a field and let $\Lambda$ be a noetherian ring. The abelian tensor category of *Nori motives* $\mathcal{MM}_{\mathrm{Nori}}(k, \Lambda)$ is the tensor-localization of $\mathcal{MM}^{\mathrm{eff}}_{\mathrm{Nori}}(k, \Lambda)$ at the Lefschetz motive $\Lambda_{\mathrm{Nori}}(-1)$.

**Remark 1.10.14.** Because $\Lambda$ has a tensor-inverse as a module over itself, we find that

$$(\Lambda\text{-}\mathrm{Mod})[\omega_{\mathrm{sing}}(\Lambda_{\mathrm{Nori}}(-1))^{\otimes -1}] \cong (\Lambda\text{-}\mathrm{Mod})[\Lambda^{\otimes -1}] \cong \Lambda\text{-}\mathrm{Mod}\,.$$

Therefore $\omega_{\mathrm{sing}}$ descends to a faithful and exact tensor functor

$$\omega_{\mathrm{sing}} : \mathcal{MM}_{\mathrm{Nori}}(k, \Lambda) \to \Lambda\text{-}\mathrm{Mod}\,.$$

**Remark 1.10.15.** As witnessed by Proposition 8.2.5 of [HM16], one can identify $\mathcal{MM}_{\mathrm{Nori}}(k,\Lambda)$ with a diagram category of an enlargement of the quiver $\mathrm{Pairs}^{\mathrm{eff}}$ of effective pairs.

The following is a direct consequence of Proposition 1.3.10:

**Proposition 1.10.16.** *Let $k \subseteq \mathbb{C}$ be a field and let $\Lambda$ be a noetherian ring.*

*(a) Every object of $\mathcal{MM}_{\mathrm{Nori}}(k,\Lambda)$ is a subquotient of one of the form*

$$\bigoplus_{i=1}^{r} H^{n_i}_{\mathrm{Nori}}(X_i, Y_i)(m_i),$$

*where the $(X_i, Y_i, n_i)$ are very good pairs and the $m_i$ are integers.*

*(b) The elements $H^n_{\mathrm{Nori}}(X,Y)(m)$ corresponding to integers $m$ and very good pairs $(X,Y,n)$ generate $\mathcal{MM}_{\mathrm{Nori}}(k,\Lambda)$ as an abelian category.*

We will also need the following technical result in the proof of Theorem 7.6.10:

**Lemma 1.10.17.** *Let $k \subseteq \mathbb{C}$ be a field and let $\Lambda$ be a noetherian ring. Let $\mathcal{V}$ be the class of objects in $\mathcal{MM}_{\mathrm{Nori}}(k,\Lambda)$ of the form $H^n_{\mathrm{Nori}}(X,Y)(m)$, where $m$ is an integer and $(X,Y,n)$ is a very good pair. Let $\langle\mathcal{V}\rangle^{\mathrm{psab}}$ be the pseudo-abelian subcategory of $\mathcal{MM}_{\mathrm{Nori}}(k,\Lambda)$ generated by $\mathcal{V}$.*

*Then the embedding $\langle\mathcal{V}\rangle^{\mathrm{psab}} \hookrightarrow \mathcal{MM}_{\mathrm{Nori}}(k,\Lambda)$ induces an equivalence*

$$\mathcal{C}\left(\langle\mathcal{V}\rangle^{\mathrm{psab}}, \omega_{\mathrm{sing}}\right) \cong \mathcal{MM}_{\mathrm{Nori}}(k,\Lambda)$$

*between the resulting diagram category and Nori motives.*

*Proof.* Proposition 8.2.5 of [HM16] shows that $\mathcal{MM}_{\mathrm{Nori}}(k,\Lambda)$ is the diagram category associated to a certain quiver with vertices the quadruples of the form $(X,Y,n)(m)$, where again $m$ is an integer and $(X,Y,n)$ is a very good pair. The result then follows from Lemma 8.1.11 of [HM16]. We, however, offer a proof independent from the description of $\mathcal{MM}_{\mathrm{Nori}}(k,\Lambda)$ as a diagram category, solely relying on our definition as a tensor-localisation:

We denote by $\mathcal{V}^{\mathrm{eff}}$ the class of objects in $\mathcal{MM}^{\mathrm{eff}}_{\mathrm{Nori}}(k,\Lambda)$ of the form $H^n_{\mathrm{Nori}}(X,Y)$ for very good pairs $(X,Y,n)$. Then the inclusion of $\langle\mathcal{V}^{\mathrm{eff}}\rangle^{\mathrm{psab}}$ into $\mathcal{MM}^{\mathrm{eff}}_{\mathrm{Nori}}(k,\Lambda)$ induces an equivalence

$$\mathcal{C}\left(\langle\mathcal{V}^{\mathrm{eff}}\rangle^{\mathrm{psab}}, \omega_{\mathrm{sing}}\right) \cong \mathcal{MM}^{\mathrm{eff}}_{\mathrm{Nori}}(k,\Lambda)$$

due to [HM16], Lemma 8.1.11. Thus the functoriality of diagram categories (cf. Lemma 7.2.6 of [HM16]) induces a functor

$$\mathcal{F}\colon \mathcal{MM}^{\mathrm{eff}}_{\mathrm{Nori}}(k,\Lambda) \cong \mathcal{C}\left(\langle\mathcal{V}^{\mathrm{eff}}\rangle^{\mathrm{psab}}, \omega_{\mathrm{sing}}\right) \to \mathcal{C}\left(\langle\mathcal{V}\rangle^{\mathrm{psab}}, \omega_{\mathrm{sing}}\right)$$

when applied to the map $\mathcal{V}^{\text{eff}} \to \mathcal{V}$, $H^n_{\text{Nori}}(X,Y) \mapsto H^n_{\text{Nori}}(X,Y)(0)$.

For $m \in \mathbb{Z}$ we let $U(-m)$ be the image of

$$\Lambda_{\text{Nori}}(-1)^{\otimes m} \in \mathcal{V} \subseteq \langle \mathcal{V} \rangle^{\text{psab}}$$

in $\mathcal{C}(\langle \mathcal{V} \rangle^{\text{psab}}, \omega_{\text{sing}})$. We then define a functor

$$\widetilde{\mathcal{F}} \colon \mathcal{MM}_{\text{Nori}}(k,\Lambda) \to \mathcal{C}(\langle \mathcal{V} \rangle^{\text{psab}}, \omega_{\text{sing}})$$

as follows:

It sends an object $A(-m) \cong A \otimes \Lambda_{\text{Nori}}(-1)^{\otimes m}$, where $A \in \mathcal{MM}^{\text{eff}}_{\text{Nori}}(k,\Lambda)$ and $m \in \mathbb{Z}$, to $\mathcal{F}(A) \otimes U(-m)$. On a morphism $f \colon A(-m) \to B(-n)$ in

$$\mathcal{MM}_{\text{Nori}}(k,\Lambda) = \mathcal{MM}^{\text{eff}}_{\text{Nori}}(k,\Lambda)[\Lambda_{\text{Nori}}(-1)^{\otimes -1}],$$

represented by a morphism $f' \colon A(-m-k) \to B(-n-k)$ in $\mathcal{MM}^{\text{eff}}_{\text{Nori}}(k,\Lambda)$ for some $k \in \mathbb{N}_0$, it is defined as $\widetilde{\mathcal{F}}(f) := \mathcal{F}(f') \otimes U(k)$.

Conversely, we get a functor

$$\mathcal{G} \colon \mathcal{C}\left(\langle \mathcal{V} \rangle^{\text{psab}}, \omega_{\text{sing}}\right) \to \mathcal{MM}_{\text{Nori}}(k,\Lambda)$$

from the universal property of the diagram category on the left hand side applied to $\omega_{\text{sing}}$. It is now readily checked that $\mathcal{G} \circ \mathcal{F}$ is by construction the identity on $\mathcal{MM}_{\text{Nori}}(k,\Lambda)$. Hence $\mathcal{G}$ is full and essentially surjective. Furthermore, $\omega_{\text{sing}}$ and thus $\mathcal{G}$ are faithful, proving the lemma. $\square$

# Chapter 2

# Symmetric Tensors and Divided Power Algebras

This chapter is to some extent a preparation for the next one as we introduce the affine versions of what is to come in Chapter 3. This is, however, not only to glue them later: we can often restrict to generic points as explained in Remark 3.5.6, hence reducing to the results of this chapter.

The main goal is to introduce symmetric tensors and Grothendieck-Deligne norm maps. At the end we link everything to Roby's theory of divided powers, which offers a somewhat different approach. It simplifies some proofs, but involves certain technicalities. Both views agree under an assumption of flatness as witnessed by Theorem 2.6.13 and Corollary 2.6.14. Because flat algebras are the focus of this chapter we will thus rarely see a difference.

We explain methods and theorems that simplify the usage of symmetric tensors. To achieve this we give, following Vaccarino, explicit sets of generators for our algebras and describe how they are mapped by norm maps. This requires us to delve into linear algebra and the elementary theory of invariants. Ultimately, this approach allows us slick proofs, rivalling the analogous proofs for divided powers.

From an informal point of view one could say that our presentation lives between both worlds. Symmetric tensors are easier to define, especially in the scheme-theoretic setting, but lack the versatility found within divided powers. This versatility comes from their relation to so-called polynomial laws, which make them inherently better suited to deal with norm-like maps. We have thus chosen to present an approach avoiding some technicalities, instead using the best of both worlds, while aiming to point out all the important relations between both sides.

Note that many statements in this chapter are already known or folklore.

## 2.1 Symmetric tensors

**Definition 2.1.1.** Let $n$ be a non-negative integer and let $M$ be a module over a ring $A$. The *$n$-fold tensor product* of $M$ over $A$ will be denoted as

$$(M|A)^{\otimes n} := M^{\otimes_A n} := \underbrace{M \otimes_A \ldots \otimes_A M}_{n \text{ times}}.$$

The $n$-fold tensor product has a natural action of the symmetric group $S_n$ by permuting the $n$ tensor factors. The corresponding module of invariants

$$\mathrm{S}_n(M|A) := \left((M|A)^{\otimes n}\right)^{S_n}$$

is called the *$n$-th symmetric tensors* of $M$ over $A$.

**Remark 2.1.2.** In the existing literature, the symbol $\mathrm{TS}^{\mathrm{n}}_{\mathrm{A}}(\mathrm{M})$ is often used instead of $\mathrm{S}_n(M|A)$. We prefer our notation as it seems to be the more general, categorical one, and because it is dual to that of symmetric products (cf. Definition 3.2.1).

**Remark 2.1.3.** If $B$ is an $A$-algebra, then so are $(B|A)^{\otimes n}$ and $\mathrm{S}_n(B|A)$.

The universal property of the tensor product and the equivariance of the $S_n$-action shows that $(-|A)^{\otimes n}$ and $\mathrm{S}_n(-|A)$ are endofunctors on each of the categories $A\text{-}\widetilde{\mathrm{Mod}}$ and $A\text{-}\mathrm{Alg}$. Even more, the same arguments readily imply that they are target-preserving endofunctors $(-|-)^{\otimes n}$ and $\mathrm{S}_n(-|-)$ on the arrow category $(A\text{-}\mathrm{Alg})^{\to}$ of Definition 1.2.1.

We will often use this functoriality in both source and target of a morphism. In accordance to our conventions regarding arrow categories we denote the morphism $\mathrm{S}_n(B|A) \to \mathrm{S}_n(B'|A')$ induced by a commutative square

$$\begin{array}{ccc} A & \longrightarrow & B \\ \downarrow{\scriptstyle f} & & \downarrow{\scriptstyle g} \\ A' & \longrightarrow & B' \end{array}$$

by $\mathrm{S}_n(g|f)$. Also recall our conventions from Section 1.1 regarding identities by which we simply write $\mathrm{S}_n(g|A)$ if $f = \mathrm{id}_A$ and $A' = A$, and similarly use $\mathrm{S}_n(B|f)$.

This has a generalization to modules: a morphism $f\colon A \to A'$ of rings, a morphism $g\colon M \to M'$ of $A$-modules and an $A'$-module structure on $M'$ such that $f(a)g(m) = g(am)$ for all $a \in A$, $m \in M$ induce a natural morphism $\mathrm{S}_n(M|A) \to \mathrm{S}_n(M'|A')$ of $A$-modules.

To gain a better understanding we aim to give an explicit yet small set of generators for the symmetric tensors, both as a module and as an algebra. Similar notions as below, but mostly restricted to polynomials, can be found in [Vac05] where the name *monomial multisymmetric functions* is attributed to Dalbec. We instead use a different naming scheme to highlight both the relations to the $n$-fold tensor products and the elementary symmetric functions.

**Definition 2.1.4.** Let $B$ be an algebra over a ring $A$ and let $n$ be a non-negative integer.

- For every $k \in \{1, 2, \ldots, n\}$ the *$k$-th formal $n$-tensor conjugate of* $b \in B$ *over* $A$ is the pure tensor

$$\iota_k(b) := \iota_k^n(b) := 1 \otimes \ldots \otimes 1 \otimes b \otimes 1 \otimes \ldots \otimes 1 \in (B|A)^{\otimes n}$$

  with $b$ at position $k$ and 1 at all other places.

- For every $k \in \{0, 1, \ldots, n\}$ the *$k$-th elementary symmetric $n$-tensor* $\rho_k(b) := \rho_k^n(b)$ of $b \in B$ over $A$ is the coefficient of $x^{n-k}$ in the polynomial

$$(x + b)^{\otimes n} = \prod_{i=1}^{n}(x + \iota_i(b)) \in (B[x]|A[x])^{\otimes n} \cong (B|A)^{\otimes n}[x].$$

  It is clearly an element of $\mathrm{S}_n(B|A)$.

- An element of $(B|A)^{\otimes n}$ is called an *elementary symmetric $n$-tensor* if it is the $k$-th elementary symmetric $n$-tensor of some $b \in B$ for some $k$. If it is clear from the context, we will omit $n$.

Thus if $\sigma_k$ is the $k$-th elementary symmetric polynomial we have

$$\rho_k(b) = \sigma_k\left(\iota_1(b), \ldots, \iota_n(b)\right).$$

**Example 2.1.5.** In particular we have $\rho_n(b) = b^{\otimes n}$ and

$$\rho_1(b) = b \otimes 1 \otimes \cdots \otimes 1 + 1 \otimes b \otimes \cdots \otimes 1 + \ldots + 1 \otimes 1 \otimes \cdots \otimes b.$$

**Example 2.1.6.** The different types of elementary symmetric 3-tensors are

$$\begin{aligned}
\rho_1(x) &= x \otimes 1 \otimes 1 + 1 \otimes x \otimes 1 + 1 \otimes 1 \otimes x,\\
\rho_2(x) &= x \otimes x \otimes 1 + x \otimes 1 \otimes x + 1 \otimes x \otimes x,\\
\rho_3(x) &= x \otimes x \otimes x.
\end{aligned}$$

**Definition 2.1.7.** Let $n$, $r$ and $w$ be non-negative integers.

- An *$n$-type of length $r$ and total weight $w$* is a tuple $a = (a_1, a_2, \ldots, a_r)$ of non-negative integers $a_1, a_2 \ldots, a_r$ where $w = \sum_{i=1}^{r} a_i \leq n$. If $a$ is an $n$-type we write $l(a)$ to denote its length and $|a|$ to denote its weight. We refer to the set of all $n$-types as $T_n$.

- Let $B$ be an algebra over a ring $A$ and let $a$ be an $n$-type of length $r$. Let $b_1, \ldots, b_r \in B$ and set $b_0 = 1 \in B$. We define

$$\rho_a(b_1, \ldots, b_r) := \sum_{\substack{f\colon \{1,2,\ldots,n\} \to \{0,1,\ldots,r\} \\ \#(f^{-1}(i)) = a_i \text{ for all} \\ i \in \{1,2,\ldots,r\}}} b_{f(1)} \otimes b_{f(2)} \otimes \cdots \otimes b_{f(n)}.$$

  The action of the group $S_n$ only permutes the summands, hence we have $\rho_a(b_1, \ldots, b_r) \in \mathrm{S}_n(B|A)$. We call those elements the *symmetric tensors of type $a$*.

  Let more generally $M$ be an $A$-module and let $m_1, \ldots, m_r \in M$. If the type $a$ has the maximal possible weight $|a| = n$, then $b_0$ does not appear in the defining sum for $\rho_a$ and hence $\rho_a(m_1, \ldots, m_r)$ is a well-defined element of $\mathrm{S}_n(M|A)$.

**Remark 2.1.8.** Let $x_1, x_2, \ldots, x_r$ be distinct indeterminates, take a type $a = (a_1, a_2, \ldots, a_r) \in T_n$ and define

$$x^{\otimes a} := 1^{\otimes(n-|a|)} \otimes x_1^{\otimes a_1} \otimes x_2^{\otimes a_2} \otimes \ldots \otimes x_r^{\otimes a_r} \in (\mathbb{Z}[x_1, \ldots, x_r]|\mathbb{Z})^{\otimes n}.$$

Then the symmetric tensor $\rho_a(x_1, \ldots, x_r)$ of type $a$ is the sum

$$\sum_{\widetilde{x} \in S_n \cdot x^{\otimes a}} \widetilde{x} \in \mathrm{S}_n\left(\mathbb{Z}[x_1, x_2, \ldots, x_r]|\mathbb{Z}\right)$$

over the $S_n$-orbit of $x^{\otimes a}$. We recover the general version by plugging in $x_i = b_i$ while also applying the morphism $\mathrm{S}_n(B|\mathbb{Z}) \to \mathrm{S}_n(B|A)$.

**Example 2.1.9.** Symmetric tensors of a given type generalize elementary symmetric tensors because $\rho_{(a_1)}(x) = \rho_{a_1}(x)$, i.e. elementary symmetric tensors are the symmetric tensors of a type of length 1.

The archetypical non-elementary examples for $n = 3$ are

$$\begin{aligned}
\rho_{(1,1)}(x, y) &= x \otimes y \otimes 1 + 1 \otimes x \otimes y + y \otimes 1 \otimes x + \\
&+ y \otimes x \otimes 1 + 1 \otimes y \otimes x + x \otimes 1 \otimes y, \\
\rho_{(1,1,1)}(x, y, z) &= x \otimes y \otimes z + x \otimes z \otimes y + y \otimes x \otimes z + \\
&+ y \otimes z \otimes x + z \otimes x \otimes y + z \otimes y \otimes x, \\
\rho_{(2,1)}(x, y) &= x \otimes x \otimes y + x \otimes y \otimes x + y \otimes x \otimes x.
\end{aligned}$$

**Remark 2.1.10.** Symmetric tensors of type $a$ satisfy an $a$-weighted binomial-like theorem in the sense that

$$\rho_a(x_1, \ldots, \mu x_i + \nu \widetilde{x}_i, \ldots, x_r) = \sum_{m+n=a_i} \mu^m \nu^n \rho_{a:[m,n]}(x_1, \ldots, x_i, \widetilde{x}_i, \ldots, x_r).$$

Here we took $a : [m, n]$ to be the type $a$, but with the $i$-th entry $a_i$ being replaced by the two entries $m$ and $n$.

**Example 2.1.11.** There are non-trivial relations between symmetric tensors. For example set $n = 3$, where we observe that

$$\begin{aligned}\rho_1(x)\rho_1(y) =&(x \otimes 1 \otimes 1 + 1 \otimes x \otimes 1 + 1 \otimes 1 \otimes x)\cdot\\ &\cdot(y \otimes 1 \otimes 1 + 1 \otimes y \otimes 1 + 1 \otimes 1 \otimes y) =\\ =&(xy \otimes 1 \otimes 1 + 1 \otimes xy \otimes 1 + 1 \otimes 1 \otimes xy)+\\ +&(x \otimes y \otimes 1 + 1 \otimes x \otimes y + y \otimes 1 \otimes x+\\ &+ y \otimes x \otimes 1 + 1 \otimes y \otimes x + x \otimes 1 \otimes y) =\\ =&\rho_1(xy) + \rho_{(1,1)}(x, y)\end{aligned}$$

and

$$\begin{aligned}\rho_1(x)\rho_2(y) =&(x \otimes 1 \otimes 1 + 1 \otimes x \otimes 1 + 1 \otimes 1 \otimes x)\cdot\\ &\cdot(y \otimes y \otimes 1 + y \otimes 1 \otimes y + 1 \otimes y \otimes y) =\\ =&(x \otimes y \otimes y + y \otimes x \otimes y + y \otimes y \otimes x)+\\ +&(xy \otimes y \otimes 1 + 1 \otimes xy \otimes y + y \otimes 1 \otimes xy+\\ &+ y \otimes xy \otimes 1 + 1 \otimes y \otimes xy + xy \otimes 1 \otimes y) =\\ =&\rho_{(1,2)}(x, y) + \rho_{(1,1)}(xy, y).\end{aligned}$$

Additionally there are the trivial ones highlighted in Remark 2.1.15.

**Remark 2.1.12.** Tensor conjugates and (elementary) symmetric tensors are functorial in the sense that if $f : B \to B'$ is a morphism of $A$-algebras, $b, b_1, \ldots, b_r \in A$, $k \in \{1, 2, \ldots, n\}$ and $a$ is an $n$-type of length $r$, then

$$\begin{aligned}f^{\otimes_A n}(\iota_k(b)) &= \iota_k(f(b))\\ f^{\otimes_A n}(\rho_k(b)) &= \rho_k(f(b))\\ f^{\otimes_A n}(\rho_a(b_1, \ldots, b_r)) &= \rho_a(f(b_1), \ldots, f(b_r)).\end{aligned}$$

Recall:

**Theorem 2.1.13** (Govorov-Lazard Theorem)**.** *Let $A$ be a ring. An $A$-module is flat if and only if it can be written as a direct limit of free $A$-modules, which may be chosen to be finitely generated.*

*Proof.* This is Théorème 1.2 of [Laz69]. See also [Gov65]. A more recent version can e.g. be found as Theorem (4.34) in [Lam99]. □

We now give explicit sets of generators, at least for flat modules and algebras.

The following Theorem 2.1.14 is well-known (see e.g. [Ive70]) and the flatness of part (b) is sketched as Proposition 1.1 of op. cit.

**Theorem 2.1.14.** *Let $A$ be a ring and let $n$ be a non-negative integer.*

(a) *Let $M$ be a free $A$-module with basis $E$, which may possibly be infinite. Then the $A$-module $\mathrm{S}_n(M|A)$ of invariants is free. An explicit basis is*

$$\mathrm{S}_n(E|A) := \left\{ \rho_a(e_1, \ldots, e_{l(a)}) \mid {}^{a \in T_n, |a| = n}_{e_1, \ldots, e_{l(a)} \in E \text{ distinct}} \right\}.$$

(b) *Let $M$ be a flat $A$-module with an arbitrary set of generators $E$. Then the $A$-module $\mathrm{S}_n(M|A)$ is flat and generated by*

$$\mathrm{S}_n(E|A) := \left\{ \rho_a(e_1, \ldots, e_{l(a)}) \mid {}^{a \in T_n, |a| = n}_{e_1, \ldots, e_{l(a)} \in E \text{ distinct}} \right\}.$$

**Remark 2.1.15.** In part (a) of Theorem 2.1.14 it is important to understand $\mathrm{S}_n(E|A)$ as a set as some $\rho_a(e_1, \ldots, e_{l(a)})$ may be equal. For example, one has $\rho_{(2,1)}(e_1, e_2) = \rho_{(1,2)}(e_2, e_1)$ and $\rho_{(2,2)}(e_1, e_2) = \rho_{(2,2)}(e_2, e_1)$. These two equalities already demonstrate the worst that can happen:

Let $a$ be an $n$-type of maximal weight $n$ and let $e \in E^{l(a)}$. Then $\mathrm{S}_{l(a)}$ acts on both $a$ and $e$ by permuting the entries of the tuples. Assume now that the entries of $e$ are linearly independent. Then the proof of Theorem 2.1.14 (a) shows that an equality $\rho_a(e) = \rho_{a'}(e')$ with a second symmetric tensor happens if and only if $l(a) = l(a')$ and there is a permutation $\sigma \in \mathrm{S}_{l(a)}$ such that $\sigma(a) = a'$ and $\sigma(e) = e'$.

*Proof of Theorem 2.1.14.*

(a) Set $[n] = \{1, 2, \ldots, n\}$. If $f \colon [n] \to E$ is any map, we define

$$\widetilde{e}_f := f(1) \otimes f(2) \otimes \ldots \otimes f(n) \in (M|A)^{\otimes n}.$$

It is standard that the $\widetilde{e}_f$ for all such maps $f$ form an $A$-basis of $(M|A)^{\otimes n}$. Thus an element $m \in (M|A)^{\otimes n}$ can be uniquely written as

$$m = \sum_{f \colon [n] \to E} m_f \widetilde{e}_f$$

for certain $m_f \in A$, all but finitely many of them being 0. The $S_n$-action is then given by

$$\sigma m = \sum_{f \colon [n] \to E} m_f \widetilde{e}_{f \circ \sigma^{-1}} = \sum_{f \colon [n] \to E} m_{f \circ \sigma} \widetilde{e}_f.$$

Therefore such an element $m$ is $S_n$-invariant if and only if $m_f = m_{f \circ \sigma}$ for all $\sigma \in S_n$ and all maps $f \colon [n] \to E$. Hence the different orbit sums $\sum_{g \in f \circ S_n} \widetilde{e}_g$ form an $A$-basis of $\mathrm{S}_n(M|A)$. It is now sufficient to check that those orbit sums are exactly the $\rho_a(e_1, \ldots, e_r) \in \mathrm{S}_n(E|A)$.

Indeed, if $f\colon [n] \to E$ is arbitrary, we let $e_1, \dots, e_r$ be the distinct elements in its image and take $a_i$ to be the cardinality of the preimage $f^{-1}(e_i)$. Thus $a = (a_1, \dots, a_r) \in T_n$ is an $n$-type of weight $n$. By the definitions and Remark 2.1.8 we conclude that

$$\rho_a(e_1, \dots, e_r) = \sum_{g \in f \circ S_n} \widetilde{e}_g.$$

Conversely, it is easy to see that every such type $a$ and elements $e_1, \dots, e_r$ come from a function $f\colon [n] \to E$, finishing the proof of part (a).

(b) By the Govorov-Lazard Theorem 2.1.13 we may write

$$M = \varinjlim_{k \in D} M_k$$

as a direct limit of finitely generated free $A$-modules $M_k$.

As tensor products commute with direct limits and because $D$ is a final subsystem of $D^n$ via the diagonal embedding, we get

$$(M|A)^{\otimes n} \cong \varinjlim_{(k_1,\dots,k_n)\in D^n} M_{k_1} \otimes_A \dots \otimes_A M_{k_n} \cong \varinjlim_{k \in D} (M_k|A)^{\otimes n}. \quad (2.1)$$

If $G$ is a group acting on an $A$-module $N$, then the invariants $N^G$ are the kernel of the morphism

$$N \to \prod_{g \in G} N$$
$$x \mapsto (x - \sigma x)_{\sigma \in G}.$$

Direct limits are left exact, therefore we conclude that they commute with taking $G$-invariants. Setting $G = S_n$ and combining this observation with isomorphism (2.1) we have shown that

$$\mathrm{S}_n(M|A) \cong \varinjlim_{k \in D} \mathrm{S}_n(M_k|A).$$

By part (a) we find that $\varinjlim_{k \in D} \mathrm{S}_n(M_k|A)$ is a direct limit of finitely generated free $A$-modules, thus by the Govorov-Lazard Theorem 2.1.13 a flat $A$-module. Hence $\mathrm{S}_n(M|A)$ is flat.

Now let $E_k = \{\widetilde{e}_1, \dots, \widetilde{e}_s\}$ be an $A$-basis of $M_k$. By part (a) the set $\mathrm{S}_n(E_k|A)$ generates the $A$-module $\mathrm{S}_n(M_k|A)$. Thus the images of all such sets $\mathrm{S}_n(E_k|A)$ in $\mathrm{S}_n(M|A) \cong \varinjlim_{k \in D} \mathrm{S}_n(M_k|A)$ generate the latter. Writing the image of each $\widetilde{e}_i$ in $\mathrm{S}_n(M|A)$ as a finite linear combination of elements of $E$ and using Remark 2.1.10 shows the claim on the generators. □

**Remark 2.1.16.** Theorem 2.1.14 and its proof remain true for the $G$-invariants $(M|A)^{\otimes n}$ of a subgroup $G \subseteq S_n$ if the $\rho_a$ are replaced by appropriate orbit sums as in Remark 2.1.8.

**Remark 2.1.17.** Some assumption regarding flatness is necessary as seen in the following example from [Lun08]: the $\mathbb{Z}[x, y]$-module $M$ generated by two elements $s, t$ and satisfying $xs = yt$ is not flat, and it is shown in op. cit. that the resulting $\mathrm{S}_3(M|\mathbb{Z}[s, t])$ is not generated by the elements of Theorem 2.1.14.

**Theorem 2.1.18.** *Let $A$ be a ring and let $n$ be a non-negative integer.*

*(a) Let $B$ be a flat $A$-algebra which is, as an $A$-module, generated by a subset $E \subseteq B$.*

*Then the ring $\mathrm{S}_n(B|A)$ of invariants is generated as an $A$-algebra by the elementary symmetric $n$-tensors (cf. Definition 2.1.4) of the $e \in E$.*

*(b) If, in addition to (a), the integer $n!$ is invertible in $A$, then it suffices to take the first elementary symmetric tensors $\rho_1(e)$, $e \in E$, to generate $\mathrm{S}_n(B|A)$ as an $A$-algebra.*

**Remark 2.1.19.** This theorem and its proof are closely related to [Vac05] and [Ryd07]. There it is amongst others shown that the elementary symmetric tensors of monomials generate $\mathrm{S}_n(A[x_1, \ldots, x_m]|A)$. Theorem 2.1.18 is then deduced from this in [Vac06]. Both sources furthermore give explicit bounds on the degree of the monomials needed and [Vac05] also has an analogue of part (b) for the ring of polynomials, which implies the more general version of Theorem 2.1.18.

*Proof of Theorem 2.1.18.* This is Proposition 3 of [Vac06]. Let us give a direct proof:

(a) Let $C$ be the $A$-subalgebra of $\mathrm{S}_n(B|A)$ generated by the elementary symmetric $n$-tensors $\rho_k(e) = \rho_{(k)}(e)$ with $e \in E$. By Theorem 2.1.14 (b) we have to check that $\mathrm{S}_n(E|A) \subseteq C$. For this purpose we show by induction on the weight $w = |a|$ that $C$ contains all symmetric tensors $\rho_a(e_1, \ldots, e_{l(a)})$, where $a \in T_n$ is an $n$-type and $e_1, \ldots, e_{l(a)} \in E$:

All types of weight 1 are of the form $a = (0, \ldots, 0, 1, 0, \ldots, 0)$ containing a single 1 and being 0 everywhere else. Their corresponding symmetric tensors are the $\rho_a(e) = \rho_1(e)$, $e \in E$, which lie in $C$ by definition.

Let now $a = (a_1, \ldots, a_r) \in T_n$ be any $n$-type of weight $|a|$. We may assume that we have already shown that $C$ contains all $\rho_{\widetilde{a}}(\widetilde{e}_1, \ldots, \widetilde{e}_s)$, with $\widetilde{e}_1, \ldots, \widetilde{e}_s \in E$, where $\widetilde{a} \in T_n$ is an $n$-type of weight $|\widetilde{a}| < |a|$ and arbitrary length $s$.

A straightforward formal expansion of the left hand side as in Example 2.1.11 reveals that

$$\rho_{a_1}(e_1) \cdot \ldots \cdot \rho_{a_r}(e_r) = \rho_{(a_1,\ldots,a_r)}(e_1, \ldots, e_r) + \sum_{(\widetilde{a},\widetilde{e})} \lambda_{(\widetilde{a},\widetilde{e})} \rho_{\widetilde{a}}(\widetilde{e}). \quad (2.2)$$

Here the $\lambda_{(\widetilde{a},\widetilde{e})}$ are non-negative integers and the sum is over finitely many pairs $(\widetilde{a}, \widetilde{e})$, where $\widetilde{a}$ is a type of weight $|\widetilde{a}| < |a|$ and $\widetilde{e}$ is an $l(\widetilde{a})$-tuple whose entries are monomials in the $e_i$.

We can express every monomial entry of $\widetilde{e}$ as a finite $A$-linear combination of some elements of $E$, not necessarily only the chosen $e_1, \ldots, e_r$. Using Remark 2.1.10 we can then expand every $\rho_{\widetilde{a}}(\widetilde{e})$ as a finite $A$-linear combination of finitely many $\rho_{\widetilde{b}}(\widetilde{e}')$, where $\widetilde{b}$ is an $n$-type of the same weight as $\widetilde{a}$ and $\widetilde{e}' \in E^{l(\widetilde{b})}$. Therefore we know by induction that each $\rho_{\widetilde{a}}(\widetilde{e})$ belongs to $C$. Hence there is an element $\beta \in C$ such that

$$\rho_{(a_1,\ldots,a_r)}(e_1, \ldots, e_r) = \rho_{a_1}(e_1) \cdot \ldots \cdot \rho_{a_r}(e_r) - \beta.$$

Thus $\rho_{(a_1,\ldots,a_r)}(e_1, \ldots, e_r) \in C$ as claimed.

(b) Let $k \in \{1, 2, \ldots, n\}$ and let $e \in E$. Another formal expansion shows

$$\rho_1(e)^k = k! \cdot \rho_k(e) + \sum_{(\widetilde{k},\widetilde{e})} \lambda_{(\widetilde{k},\widetilde{e})} \rho_{\widetilde{k}}(\widetilde{e}) \quad (2.3)$$

for some non-negative integers $\lambda_{(\widetilde{k},\widetilde{e})}$, the sum running over finitely many pairs $(\widetilde{k}, \widetilde{e})$, where $\widetilde{k}$ is an $n$-type of weight $|\widetilde{k}| < k$ and $\widetilde{e} \in B^{l(\widetilde{k})}$ consists of powers of $e$.

Let now $a = (a_1, \ldots, a_r)$ be an $n$-type of length $r$. As before we want to show by induction on its weight $|a|$ that all $\rho_a(e)$, $e \in E^r$, are in the $A$-subalgebra of $\mathrm{S}_n(B|A)$ generated by the $\rho_1(e')$, $e' \in E$. We actually show something slightly stronger:

By equation (2.2) of part (a) we can, setting $c = \frac{|a|!}{a_1! \cdots a_r!} \in \mathbb{N}$ to be a multinomial coefficient, write

$$|a|!\rho_a(e_1, \ldots, e_r) = c \cdot \prod_{i=1}^{r} a_i!\rho_{a_i}(e_i) - \sum_{(\widetilde{a},\widetilde{e})} \lambda'_{(\widetilde{a},\widetilde{e})} |\widetilde{a}|! \rho_{\widetilde{a}}(\widetilde{e}).$$

Here we put $\lambda'_{(\widetilde{a},\widetilde{e})} := \lambda_{(\widetilde{a},\widetilde{e})} \cdot \frac{|a|!}{|\widetilde{a}|!}$, which is still an integer as $|\widetilde{a}| < |a|$.

Hence, using the expansion of equation (2.3), we can by induction on $|a|$ and arguments similar to part (a) express $|a|!\rho_a(e_1, \ldots, e_r)$ using only some $\rho_1(e')$ with $e' \in E$.

The assumed invertibility of $n!$ and thus of $|a|!$ finishes the proof. □

**Remark 2.1.20.** The example $A = \mathbb{F}_2$, $B = \mathbb{F}_4$, $n = 3$ shows that the $\rho_1(b)$, $b \in B$, do not always generate the $A$-algebra $\mathrm{S}_n(B|A)$. The same example also shows that neither the $\rho_2(b)$ nor the $\rho_3(b)$ do suffice on their own.

## 2.2 Further functorialities between symmetric tensors

We list some important properties and naturalities concerning the functor $\mathrm{S}_n(-|-)$:

**Lemma 2.2.1.** *Let $A$ be a ring and let $M, N$ be $A$-modules. Let $G$ be a finite group acting on $M$ and trivially on $N$. Assume that $N$ is flat.*

*Then the induced $G$-action on $M \otimes_A N$ induces a natural isomorphism*

$$M^G \otimes_A N \cong (M \otimes_A N)^G$$

*of $G$-invariant $A$-modules, induced by taking the tensor product of the inclusion $M^G \hookrightarrow M$ with $N$.*

*Proof.* This is easy when $N$ is finitely generated and free. The general case follows easily from a direct limit argument using the Govorov-Lazard Theorem 2.1.13 as in the proof of Theorem 2.1.14 (b). □

**Lemma 2.2.2.** *Let $A'$ be an algebra over a ring $A$ which is flat as an $A$-module and let $M$ be an $A$-module.*

*Then the functoriality of Remark 2.1.3 induces a natural isomorphism*

$$\mathrm{S}_n(M|A) \otimes_A A' \cong \mathrm{S}_n(M \otimes_A A'|A')$$

*of $A'$-modules.*

*Proof.* We have an isomorphism

$$(M|A)^{\otimes n} \otimes_A A' \cong (M \otimes_A A'|A')^{\otimes n}.$$

Taking $S_n$-invariants and using Lemma 2.2.1 shows the result. □

There is also a closely related version which shifts the assumption of flatness towards the module:

**Lemma 2.2.3.** *Let $A'$ be an algebra over a ring $A$ and let $M$ be a flat $A$-module.*

*Then the functoriality of Remark 2.1.3 induces a natural isomorphism*

$$\mathrm{S}_n(M|A) \otimes_A A' \cong \mathrm{S}_n(M \otimes_A A'|A')$$

*of $A'$-modules.*

*Proof.* By the Govorov-Lazard Theorem 2.1.13 and a short direct limit argument it suffices to show this under the stronger assumption that $M$ is free.

Let $E$ be an $A$-basis of $M$. Then $E \otimes_A A' := \{e \otimes 1 \mid e \in E\}$ is an $A'$-basis of $M \otimes_A A'$.

The morphism $\mathrm{S}_n(M|A) \to \mathrm{S}_n(M \otimes_A A'|A')$ of Remark 2.1.3 sends $\mathrm{S}_n(E|A)$ to $\mathrm{S}_n(E \otimes_A A'|A')$. This map is bijective by Remark 2.1.15. Hence the lemma follows from Theorem 2.1.14 (a). $\square$

**Lemma 2.2.4.** *Let $M$ be a flat module over a ring $A$ and let $m, n$ be non-negative integers.*

*Then the inclusion $S_m \times S_n \subseteq S_{m+n}$ of Convention 1.1.2 induces a natural isomorphism*

$$\mathrm{S}_m(M|A) \otimes_A \mathrm{S}_n(M|A) \cong \left((M|A)^{\otimes(m+n)}\right)^{S_m \times S_n}$$

*of flat $A$-modules.*

*Proof.* This is just a two-fold application of Lemma 2.2.1, noting that everything is flat by Theorem 2.1.14 (b):

$$\begin{aligned}
\mathrm{S}_m(M|A) \otimes_A \mathrm{S}_n(M|A) &= \mathrm{S}_m(M|A) \otimes_A \left(M^{\otimes_A n}\right)^{S_n} \cong \\
&\cong \left(\mathrm{S}_m(M|A) \otimes_A M^{\otimes_A n}\right)^{S_n} = \\
&= \left(\left(M^{\otimes_A m}\right)^{S_m} \otimes_A M^{\otimes_A n}\right)^{S_n} \cong \\
&\cong \left(\left(M^{\otimes_A m} \otimes_A M^{\otimes_A n}\right)^{S_m}\right)^{S_n} \cong \\
&\cong \left(M^{\otimes_A (m+n)}\right)^{S_m \times S_n}.
\end{aligned}$$ $\square$

**Corollary 2.2.5.** *Let $M$ be a flat module over a ring $A$ and let $r$ and $n_1, n_2, \ldots, n_r$ be non-negative integers.*

*Then the inclusion $S_{n_1} \times \cdots \times S_{n_r} \subseteq S_{\sum_{i=1}^r n_i}$ of Convention 1.1.2 induces a natural isomorphism*

$$\mathrm{S}_{n_1}(M|A) \otimes_A \cdots \otimes_A \mathrm{S}_{n_r}(M|A) \cong \left((M|A)^{\otimes \sum_{i=1}^r n_i}\right)^{S_{n_1} \times \cdots \times S_{n_r}}$$

*of flat $A$-modules.*

*Proof.* This follows inductively by the same argument as Lemma 2.2.4. $\square$

**Lemma 2.2.6.** *Let $M$ be a flat module over a ring $A$ and let $m, n$ be non-negative integers.*

*Then the morphism $S_m \ltimes S_n^m \to S_{mn}$ of Convention 1.1.2 induces a natural isomorphism*

$$\mathrm{S}_m\left(\mathrm{S}_n(M|A)|A\right) \cong \left(M^{\otimes_A mn}\right)^{S_m \ltimes S_n^m}$$

*of $A$-modules.*

*Proof.* By Corollary 2.2.5 there is a natural isomorphism

$$\left(\left(M^{\otimes_A n}\right)^{S_n}\right)^{\otimes_A m} \cong \left(\left(M^{\otimes_A n}\right)^{\otimes_A m}\right)^{S_n^m}.$$

Hence by Convention 1.1.2 we get

$$\begin{aligned}\mathrm{S}_m\left(\mathrm{S}_n(M|A)|A\right) &= \left(\left(\left(M^{\otimes_A n}\right)^{S_n}\right)^{\otimes_A m}\right)^{S_m} \cong \\ &\cong \left(\left(\left(M^{\otimes_A n}\right)^{\otimes_A m}\right)^{S_n^m}\right)^{S_m} \cong \\ &\cong \left(M^{\otimes_A mn}\right)^{S_m \ltimes S_n^m}.\end{aligned}$$

□

**Definition 2.2.7.** Convention 1.1.2 induces due to Corollary 2.2.5 and Lemma 2.2.6 natural transformations of endofunctors of flat $A$-modules, which restrict to flat $A$-algebras:

(a) For any non-negative integers $n_1, \ldots, n_r$ a natural transformation

$$\sigma_{n_1,\ldots,n_r} \colon \mathrm{S}_{\sum_{i=1}^r n_i}(-|A) \implies \bigotimes_{i=1}^r \mathrm{S}_{n_i}(-|A).$$

(b) For all non-negative integers $m$ and $n$ a natural transformation

$$\tau_{m,n} \colon \mathrm{S}_{mn}(-|A) \implies \mathrm{S}_m(-|A) \circ \mathrm{S}_n(-|A).$$

**Remark 2.2.8.** Let $A$ be a ring and $M$ be a flat $A$-module. The morphism

$$\sigma_{m,n} \colon \mathrm{S}_{m+n}(M|A) \to \mathrm{S}_m(M|A) \otimes_A \mathrm{S}_n(M|A)$$

is easily described at the level of elementary symmetric tensors:

Lemma 2.2.4 identifies $\mathrm{S}_m(M|A) \otimes_A \mathrm{S}_n(M|A)$ with an $A$-submodule of $(M|A)^{\otimes(m+n)}$, within which the domain and the codomain of $\sigma_{m,n}$ reside. At this level $\sigma_{m,n}$ is nothing else than a simple inclusion of submodules, in particular it is injective. Similarly we see that $\tau_{m,n}$ is injective.

Understanding both sides as elements of $(M|A)^{\otimes(m+n)}$ it is not difficult to check that

$$\sigma_{m,n}\left(\rho_k^{m+n}(b)\right) = \sum_{j=0}^{k} \rho_j^m(b) \otimes \rho_{k-j}^n(b)$$

for all $b \in M$.

**Remark 2.2.9.** There are many commutative diagrams arising from these two natural transformations. We will explore them more thoroughly in Section 3.3 for schemes, where we list many properties of their dual versions.

The two natural transformations $\sigma_{m,n}$ and $\tau_{m,n}$ are already defined at the level of flat $A$-modules. Many of the respective diagrams in Section 3.3 generalize to this setting, but we will not need this.

## 2.3 Norm maps

**Definition 2.3.1.** Let $A$ be a ring and $n$ a positive integer. An *$n$-ic ring extension* of $A$ is an $A$-algebra $B$ which as an $A$-module is free of rank $n$.

**Definition 2.3.2.** We define a *good triple* $(M|B|A)$ *of rank* $n$ to consist of a ring $A$, an $A$-algebra $B$ and a $B$-module $M$ which is a free $A$-module of finite rank $n$.

An $n$-ic ring extension $B|A$ is therefore the same as a good triple $(M|B|A)$ of rank $n$ with $M = B$.

If $A$ is a ring, $M$ is a free $A$-module of rank $n$ and $\alpha\colon M \to M$ is an $A$-endomorphism of $M$, usual linear algebra induces the notions of determinant $\det(\alpha) \in A$ and characteristic polynomial $\chi_\alpha \in A[x]$. They are connected by the defining property $\chi_\alpha(x) = \det(x - \alpha \otimes_A A[x])$. For our purposes the alternating characteristic polynomial

$$\widetilde{\chi}_\alpha(x) = (-1)^n \chi_\alpha(-x) = \det(x + \alpha \otimes A[x])$$

will be more useful.

**Definition 2.3.3.** Let $(M|B|A)$ be a good triple of rank $n$ and let $b \in B$. Then the *$k$-th characteristic coefficient* $\chi_k(b|M)$ of $b$ is the coefficient of $x^{n-k}$ in the alternating characteristic polynomial $\widetilde{\chi}_b(x)$ of the multiplication by $b$ map $M \to M$. As a special case we recover the determinant $\det(b|M) = \chi_n(b|M)$.

If $B = M$ is an $n$-ic ring extension of $A$ we will simply write $\chi_k(b)$ for $\chi_k(b|M)$.

As found in [SV96], [Fer98], [Ryd08b] and others, a good triple induces a norm map:

**Definition 2.3.4.** Let $(M|B|A)$ be a good triple of rank $n$ as given by Definition 2.3.2 and let $\Lambda^n(M|A)$ denote the $n$-th exterior product of $M$ over $A$. By Proposition 1.3 of [LT07], or from a short argument using Theorem 2.1.18 (a), we get a well-defined $A$-linear map

$$\psi_{M|B|A}\colon \mathrm{S}_n(B|A) \otimes_A \Lambda^n(M|A) \to \Lambda^n(M|A)$$

given by

$$\sum_{i=1}^{s} (\alpha_{1,i} \otimes \cdots \otimes \alpha_{n,i}) \otimes (m_1 \wedge \cdots \wedge m_n) \mapsto \sum_{i=1}^{s} (\alpha_{1,i} m_1 \wedge \cdots \wedge \alpha_{n,i} m_n).$$

Due to the tensor-hom-adjunction it therefore induces an $A$-linear map

$$\vartheta_{M|B|A}\colon \mathrm{S}_n(B|A) \to \mathrm{End}_A(\Lambda^n(M|A)) \cong A$$

which is directly seen to be a morphism of $A$-algebras. We call it the *symmetrization morphism* of $(M|B|A)$.

If $B = M$ is an $n$-ic ring extension of $A$ we simply write $\vartheta_{B|A}$.

**Remark 2.3.5.** Some freeness condition on $M$ or $B$ is essential for this morphism to be well-defined as seen in [Lun08], Example 7.3.

**Lemma 2.3.6.** *Let $(M|B|A)$ be a good triple of rank $n$.*

*Then $\vartheta_{M|B|A}$ commutes with base change, i.e. if $T$ is an $A$-algebra then $(M \otimes_A T|B \otimes_A T|T)$ is a good triple of rank $n$ and the diagram*

$$\begin{array}{ccc} \mathrm{S}_n(B|A) & \xrightarrow{\vartheta_{M|B|A}} & A \\ \downarrow & & \downarrow \\ \mathrm{S}_n(B \otimes_A T|T) & \xrightarrow{\vartheta_{M\otimes_A T|B\otimes_A T|T}} & T \end{array}$$

*of $A$-algebras commutes.*

*Proof.* Freeness and rank $n$ are trivial. The morphism $\psi_{M|B|A}$ is induced by the natural module structure of $M^{\otimes_A n}$ as a $B^{\otimes_A n}$-module by taking the $S_n$-invariants in the latter and then descending to the quotient $\Lambda^n(M|A)$ on the former. Thus $\psi_{M|B|A}$ commutes with base change and hence the same is true for $\vartheta_{M|B|A}$ by the adjunction. □

**Proposition 2.3.7.** *Let $(M|B|A)$ be a good triple of rank $n$, let $k \in \{1, 2, \ldots, n\}$ and let $b \in B$.*

*Then the symmetrization of the $k$-th elementary symmetric $n$-tensor of $b$ is the $k$-th characteristic coefficient of $b$, i.e.*

$$\vartheta_{M|B|A}(\rho_k(b)) = \chi_k(b|M).$$

*Proof.* If $k = n$ we have to show that $\vartheta_{M|B|A}(b^{\otimes n}) = \det(b|M)$. But this is simply the usual way to calculate, or even define, the determinant of an $A$-endomorphism of $M$ via the exterior product.

If $k$ is arbitrary, we apply the base-change $A \hookrightarrow A[x]$: we consider the good triple $(M[x]|B[x]|A[x])$, where we wrote $M[x]$ for $M \otimes_A A[x]$, and use the already considered case $k = n$ to get

$$\begin{aligned} \sum_{k=0}^{n} \vartheta_{M|B|A}(\rho_k(b))x^{n-k} &= \vartheta_{M[x]|B[x]|A[x]}\left(\sum_{k=0}^{n} \rho_k(b)x^{n-k}\right) = \\ &= \vartheta_{M[x]|B[x]|A[x]}\left((x+b)^{\otimes n}\right) = \\ &= \det(x+b|M[x]) = \widetilde{\chi}_b(x) = \sum_{k=0}^{n} \chi_k(b|M)x^{n-k}, \end{aligned}$$

the first equality by Lemma 2.3.6. Comparison of coefficients yields the result. □

**Remark 2.3.8.** The proof of Proposition 2.3.7 demonstrates how properties regarding characteristic coefficients or elementary symmetric tensors can, using base change, often be reduced to checking them for the determinant by the defining formulas

$$\det(x+b) = \sum_{k=0}^{n} \chi_k(b|M) x^{n-k}$$
$$(x+b)^{\otimes n} = \sum_{k=0}^{n} \rho_k(b) x^{n-k}.$$

In combination with Proposition 2.3.7 and Theorem 2.1.18 this allows functorialities to be checked at elements of type $(x+b)^{\otimes n}$, which are rather easy to work with. In the upcoming Sections 2.4 and 2.5 we will see this effect in action.

## 2.4 Functorial behaviour of the norm maps

**Lemma 2.4.1.** *Let $A$ be a ring and let $B$ be an $A$-algebra. Also let*

$$0 \to L \to M \to N \to 0$$

*be an exact sequence of $B$-modules $L$, $M$ and $N$ which are free $A$-modules of finite ranks $l$, $m$ and $n$, respectively.*

*Then $m = l + n$ and the diagram*

$$\begin{array}{ccc} \mathrm{S}_m(B|A) & \xrightarrow{\quad \vartheta_{M|B|A} \quad} & A \\ \big\downarrow{\scriptstyle \sigma_{l,n}} & & \Big\| \\ \mathrm{S}_l(B|A) \otimes_A \mathrm{S}_n(B|A) & \xrightarrow{\vartheta_{L|B|A} \otimes_A \vartheta_{N|B|A}} & A \otimes_A A \end{array}$$

*commutes.*

*Proof.* The equality of ranks is standard.

By Theorem 2.1.18 it suffices to check the commutativity for elementary symmetric $m$-tensors $\rho_k(b)$, $b \in B$. Using the base change $A \hookrightarrow A[x]$, we can by Remark 2.3.8 restrict to the case $k = m$.

We recall that $\rho_m(b) = b^{\otimes m}$ and observe, for example from Remark 2.2.8, that $\sigma_{l,n}(b^{\otimes m}) = b^{\otimes l} \otimes b^{\otimes n}$. Due to Proposition 2.3.7 we hence have to check that $\det(b|M) = \det(b|L) \cdot \det(b|N)$. But this is just the well-known fact that the determinant is multiplicative on short exact sequences. □

**Remark 2.4.2.** Recall that the length of an $A$-module $M$ is by definition the maximal size of a strictly increasing chain of $A$-submodules of $M$. It is thus automatically the maximal possible length $m$ of a composition series of $M$, i.e. a strictly increasing chain

$$\{0\} = M_m \subsetneq M_{m-1} \subsetneq \ldots \subsetneq M_1 \subsetneq M_0 = M$$

of $A$-submodules with simple quotients $M_i/M_{i+1}$.

The Jordan-Hölder Theorem states that all composition series have the same size, hence any such chain already calculates the length.

**Proposition 2.4.3.** *Let $k$ be field. Let $A$ be a finite-dimensional local $k$-algebra of degree $d$ over $k$ and length $m$ as a module over itself. Let $K$ denote the residue field of $A$, let $\pi\colon A \to K$ be the reduction and let $n$ be degree of $K$ over $k$.*

*Then $d = mn$ and the diagram*

$$\begin{array}{ccc}
\mathrm{S}_{mn}(A|k) & \xrightarrow{\quad \vartheta_{A|k} \quad} & k \\
\Big\downarrow {\scriptstyle \sigma_{\underbrace{n,n,\ldots,n}_{m \text{ times}}}} & & \Big\| \\
(\mathrm{S}_n(A|k)|k)^{\otimes_k m} & & \Big\| \\
\Big\downarrow {\scriptstyle (\mathrm{S}_n(\pi|k))^{\otimes_k m}} & & \Big\| \\
(\mathrm{S}_n(K|k)|k)^{\otimes_k m} & \xrightarrow{\quad \vartheta_{K|k}^{\otimes m} \quad} & (k|k)^{\otimes m}
\end{array}$$

*of $k$-algebras commutes.*

*Proof.* Let $\mathfrak{m}$ be the maximal ideal of $A$. Then $\mathfrak{m}^N = (0)$ for some positive integer $N$ by Krull's intersection theorem. Now we have a chain

$$(0) = \mathfrak{m}^N \subseteq \mathfrak{m}^{N-1} \subseteq \ldots \subseteq \mathfrak{m} \subseteq A$$

where each subquotient $\mathfrak{m}^i/\mathfrak{m}^{i+1}$ is a $K = A/\mathfrak{m}$-vector space. By refining it we arrive at a chain

$$(0) = A_{m'} \subsetneq A_{m'-1} \subsetneq \ldots \subsetneq A_1 = \mathfrak{m} \subsetneq A_0 = A$$

such that each subquotient $Q_i := A_i/A_{i+1}$ is a 1-dimensional $K$-vector space. In particular, the $A_i$ are ideals of $A$ and the quotients $Q_i$ are simple $A$-modules. Hence this forms a composition series, so we have $m' = m$ by Remark 2.4.2. Looking at the dimensions over $k$ we see that $d = m'n = mn$.

Inductively applying Lemma 2.4.1 now shows that

$$\begin{aligned}\vartheta_{A|k} = \vartheta_{A_0|A|k} &= \left(\vartheta_{Q_0|A|k} \otimes_k \vartheta_{A_1|A|k}\right) \circ \sigma_{n,n(m-1)} = \ldots = \\ &= \left( \left( \bigotimes_{i=0}^{j-1} \vartheta_{Q_i|A|k} \right) \otimes_k \vartheta_{A_j|A|k} \right) \circ \sigma_{\underbrace{n,n,\ldots,n}_{j \text{ times}},n(m-j)} = \ldots = \\ &= \left( \bigotimes_{i=0}^{m-1} \vartheta_{Q_i|A|k} \right) \circ \sigma_{\underbrace{n,n,\ldots,n}_{m \text{ times}}}.\end{aligned}$$

Here we used the insight of Remark 2.2.8 that the morphisms $\sigma_{-}$ can be interpreted as simple inclusions.

As $a \in A$ and $\pi(a) \in K$ induce the same $k$-endomorphisms on $Q_i \cong K$ we get

$$\begin{aligned}\vartheta_{A|k} = \left( \bigotimes_{i=0}^{m-1} \vartheta_{Q_i|A|k} \right) \circ \sigma_{\underbrace{n,n,\ldots,n}_{m \text{ times}}} &= \left( \bigotimes_{i=0}^{m-1} \left(\vartheta_{Q_i|K|k} \circ \mathrm{S}_n(\pi|k)\right) \right) \circ \sigma_{\underbrace{n,n,\ldots,n}_{m \text{ times}}} = \\ = \left(\vartheta_{K|k} \circ \mathrm{S}_n(\pi|k)\right)^{\otimes_k m} \circ \sigma_{\underbrace{n,n,\ldots,n}_{m \text{ times}}} &= \vartheta_{K|k}^{\otimes_k m} \circ \left(\mathrm{S}_n(\pi|k)\right)^{\otimes_k m} \circ \sigma_{\underbrace{n,n,\ldots,n}_{m \text{ times}}}.\end{aligned}$$

□

**Proposition 2.4.4.** *Let $B|A$ be an $m$-ic ring extension and let $C|B$ be an $n$-ic ring extension.*

*Then $C|A$ is an $mn$-ic ring extension and the diagram*

$$\begin{array}{ccc} \mathrm{S}_{mn}(C|A) & \xrightarrow{\tau_{m,n}} & \mathrm{S}_m\left(\mathrm{S}_n(C|A)|A\right) \\ & & \downarrow{\scriptstyle \mathrm{S}_m(\mathrm{S}_n(C|f)|A)} \\ \downarrow{\scriptstyle \vartheta_{C|A}} & & \mathrm{S}_m\left(\mathrm{S}_n(C|B)|A\right) \\ & & \downarrow{\scriptstyle \mathrm{S}_m(\vartheta_{C|B}|A)} \\ A & \xleftarrow[\vartheta_{B|A}]{} & \mathrm{S}_m(B|A) \end{array}$$

*of $A$-algebras commutes.*

*Proof.* Theorem 2.1.18 allows us to check this at the elementary symmetric tensors. By Remark 2.3.8 and a base change $A \hookrightarrow A[x]$ it suffices to consider pure tensor powers $\rho_{mn}(c) = c^{\otimes mn}$, $c \in C$.

Remark 2.2.8 implies that $\tau_{m,n}(c^{\otimes mn}) = (c^{\otimes n})^{\otimes m}$. Furthermore, we get from Proposition 2.3.7 that $\vartheta_{C|A}(c^{\otimes mn}) = \det_{C|A}(c)$ and

$$\vartheta_{B|A}\left(\vartheta_{C|B}^{\otimes m}(c^{\otimes mn})\right) = \vartheta_{B|A}\left(\det_{C|B}(c)^{\otimes m}\right) = \det_{B|A}\left(\det_{C|B}(c)\right)$$

for all $c \in C$. Here we wrote $\det_{T|S}(t)$ to denote the determinant of the $S$-linear multiplication-by-$t$ map $T \to T$, where $t \in T$. Therefore we have to show that $\det_{C|A} = \det_{B|A} \circ \det_{C|B}$ as maps $C \to A$.

We choose a basis $\alpha_1, \ldots, \alpha_m$ of the $A$-module $B$ and a basis $\beta_1, \ldots, \beta_n$ of the $B$-module $C$. Then the $\alpha_i\beta_j$ with $i \in \{1, 2, \ldots, m\}$ and $j \in \{1, 2, \ldots, n\}$ constitute a basis of the $A$-module $C$. This choice of bases identifies $B$ with a subring of the non-commutative matrix ring $A^{m\times m}$, $C$ with a subring of $B^{n\times n}$ and thus also $C$ with a subring of $(A^{m\times m})^{n\times n} \cong A^{mn\times mn}$, where the last isomorphism corresponds to our choice for the $A$-basis of $C$. Hence we have reduced the claimed commutativity to Theorem 1 of [KSW99]. □

**Remark 2.4.5.** Note that this generalizes the well-known

$$\mathrm{Nm}_{M|K} = \mathrm{Nm}_{L|K} \circ \mathrm{Nm}_{M|L},$$

where $M|L|K$ are finite field extensions and Nm denotes the respective norms.

## 2.5 Multiplying cycles

We also consider multiplicative, i.e. tensor, structures in Chapter 3. Our approach there is more direct then the one presented in this section. Instead one should understand this section as an alternative approach to an affine version of Theorem 3.7.4. It is still worthwhile because it offers some interesting insights not immediately apparent in the scheme-theoretic versions found in Chapter 3.

We also demonstrate how to use the formal tensor conjugates $\iota_k^n(b)$ in the proof of Theorem 2.5.2. One should think of them as conjugates in a Galois-theoretic sense. The reader interested in the formalities behind this should look at [BS14].

**Lemma 2.5.1.** *Let $m$ and $n$ be non-negative integers. Then there exist $mn+1$ polynomials $w_0 = 1, w_1, \ldots, w_{mn} \in \mathbb{Z}[u_1, \ldots, u_m, v_1, \ldots, v_n]$ in $m+n$ variables with the following properties:*

*(a) Let $\overline{a} = (a_1, \ldots., a_m)$ and $\overline{b} = (b_1, \ldots, b_n)$ be tuples of variables. Then the formal identity*

$$\prod_{\substack{i=1,\ldots,m\\ j=1,\ldots,n}} (t + a_i b_j) =$$
$$= \sum_{k=0}^{mn} w_k \left(\sigma_1(\overline{a}), \sigma_2(\overline{a}), \ldots, \sigma_m(\overline{a}), \sigma_1(\overline{b}), \sigma_2(\overline{b}), \ldots, \sigma_n(\overline{b})\right) t^{mn-k}$$

*holds. Here $\sigma_k$ is the $k$-th elementary symmetric polynomial.*

*(b) Let $X$ and $Y$ be $m\times m$ and $n\times n$ matrices of indeterminates, respectively, with alternating characteristic polynomials*

$$\det(t+X) = \sum_{i=0}^{m} \chi_i(X) t^{m-i}$$
$$\det(t+Y) = \sum_{j=0}^{n} \chi_j(Y) t^{n-j}.$$

*Then we have the formal relation*

$$\det(t+X\otimes Y) = \\ = \sum_{k=0}^{mn} w_k\left(\chi_1(X),\chi_2(X),\ldots,\chi_m(X),\chi_1(Y),\chi_2(Y),\ldots,\chi_n(Y)\right)t^{mn-k}$$

*between characteristic polynomials.*

*Furthermore, each of the properties a) and b) makes the $w_k$ unique.*

*Proof.* Using diagonal matrices, i.e. setting all off-diagonal entries to 0, shows that part (a) follows from part (b).

The coefficient of $t^{mn-k}$ in

$$\prod_{\substack{i=1,\ldots,m\\ j=1,\ldots,n}} (t+a_i b_j)$$

is invariant under the actions of $S_m$ permuting the $a_i$ and $S_n$ permuting the $b_j$. Thus a twofold application of the main theorem on elementary symmetric polynomials shows that the $w_k$ as in part (a) exist and are unique. We are thus left to show that they also satisfy the property of part (b):

Let $K$ be any algebraically closed field of characteristic 0 containing all the indeterminates occuring in $X$ and $Y$. As $\mathbb{Z}$ embeds into $K$ it suffices to check the identities over $K$. There we can upper triangularize $X$ and $Y$ into $X'$ and $Y'$, and doing so also triangularizes $X\otimes Y$ into $X'\otimes Y'$.

If we let $a_1,\ldots,a_m$ and $b_1,\ldots,b_n$ be the diagonal entries of $X'$ and $Y'$, respectively, then the diagonal entries of $X'\otimes Y'$ are the pairwise products $a_i b_j$. Thus the result for $X'$ and $Y'$ follows from part (a), and by reversing the triangularizations we recover it for $X$ and $Y$. □

**Theorem 2.5.2.** *Let $A$ and $\widetilde{A}$ be flat algebras over the same ring $R$. Let $B|A$ and $\widetilde{B}|\widetilde{A}$ be $n$-ic and $\widetilde{n}$-ic algebras, respectively.*

*Let furthermore*

$$\begin{aligned} \delta\colon \quad (B|A)^{\otimes \widetilde{n}} &\to B \\ b_1 \otimes \ldots \otimes b_{\widetilde{n}} &\mapsto \prod_{i=1}^{\widetilde{n}} b_i, \\ \widetilde{\delta}\colon \quad (\widetilde{B}|\widetilde{A})^{\otimes n} &\to \widetilde{B} \\ \widetilde{b}_1 \otimes \ldots \otimes \widetilde{b}_n &\mapsto \prod_{i=1}^{n} \widetilde{b}_i \end{aligned}$$

*be the diagonal morphisms.*

*Then $\delta^{\otimes n} \otimes \widetilde{\delta}^{\otimes \widetilde{n}}$ induces a natural morphism*

$$\rho_{n,\widetilde{n}}\colon \mathrm{S}_{n\widetilde{n}}\left(B \otimes_R \widetilde{B}|A \otimes_R \widetilde{A}\right) \to \mathrm{S}_n(B|A) \otimes_R \mathrm{S}_{\widetilde{n}}(\widetilde{B}|\widetilde{A}).$$

*Furthermore, $B \otimes_R \widetilde{B}$ is an $n\widetilde{n}$-ic algebra over $A \otimes_R \widetilde{A}$ and the diagram*

$$\begin{array}{ccc} \mathrm{S}_{n\widetilde{n}}\left(B \otimes_R \widetilde{B}|A \otimes_R \widetilde{A}\right) & \xrightarrow{\vartheta_{B\otimes_R \widetilde{B}|A\otimes_R \widetilde{A}}} & A \otimes_R \widetilde{A} \\ & {\scriptstyle \rho_{n,\widetilde{n}}} \searrow \qquad \nearrow {\scriptstyle \vartheta_{B|A} \otimes_R \vartheta_{\widetilde{B}|\widetilde{A}}} & \\ & \mathrm{S}_n(B|A) \otimes_R \mathrm{S}_{\widetilde{n}}(\widetilde{B}|\widetilde{A}) & \end{array}$$

*commutes.*

*Proof.* It is standard that $B \otimes_R \widetilde{B}$ is a free $A \otimes_R \widetilde{A}$-module of rank $n\widetilde{n}$. We use, without further mention, that if $D|C$ and $C|R$ are flat algebras, then $D|R$ is flat as well.

For the remainder of the proof it is advisable to understand $n\widetilde{n}$-fold tensor powers as 'two-dimensional' tensor products, corresponding by Convention 1.1.2 to identifications

$$\left(C^{\otimes n}\right)^{\otimes \widetilde{n}} \cong C^{\otimes n\widetilde{n}} \cong \left(C^{\otimes \widetilde{n}}\right)^{\otimes n}.$$

Each of $\delta^{\otimes n}$ and $\widetilde{\delta}^{\otimes \widetilde{n}}$ is equivariant for both $S_n$ and $S_{\widetilde{n}}$, thus the same is true for $\delta^{\otimes n} \otimes_R \widetilde{\delta}^{\otimes \widetilde{n}}$. Therefore $\delta^{\otimes n} \otimes_R \widetilde{\delta}^{\otimes \widetilde{n}}$ restricts by Remark 2.2.8 to give the morphism $\rho_{m,n}$.

If $e_1, \ldots, e_n$ are an $A$-basis of $B$ and $\widetilde{e}_1, \ldots, \widetilde{e}_{\widetilde{n}}$ are an $\widetilde{A}$-basis of $\widetilde{B}$, then the $e_i \otimes \widetilde{e}_j$ are an $A \otimes_R \widetilde{A}$-basis of $B \otimes_R \widetilde{B}$. Thus by Theorem 2.1.18 (a) and Remark 2.3.8 we only need to check that, after base change to $R[x]$, the diagram commutes for the elements of type $(x + b \otimes \widetilde{b})^{\otimes n\widetilde{n}}$ with $b \in B$ and $\widetilde{b} \in \widetilde{B}$.

We define

$$\iota^{n\widetilde{n}}_{i,j}(b\otimes\widetilde{b}) := \iota^{\widetilde{n}}_{j}\big(\iota^{n}_{i}(b\otimes\widetilde{b})\big) = \iota^{n}_{i}\big(\iota^{\widetilde{n}}_{j}(b\otimes\widetilde{b})\big) \in (B\otimes_R \widetilde{B}|A\otimes_R \widetilde{A})^{\otimes n\widetilde{n}}$$

which is the pure tensor with ones everywhere, except a $b\otimes\widetilde{b}$ at position $(i,j)$. We then have

$$(x+b\otimes\widetilde{b})^{\otimes n\widetilde{n}} = \prod_{\substack{i=1,\dots,n\\ j=1,\dots,\widetilde{n}}} \left(x+\iota^{n\widetilde{n}}_{i,j}(b\otimes\widetilde{b})\right),$$

which is mapped by $\delta^{\otimes n}\otimes\widetilde{\delta}^{\otimes\widetilde{n}}$ to

$$\prod_{\substack{i=1,\dots,n\\ j=1,\dots,\widetilde{n}}} \left(x+\iota^{n}_{i}(b)\otimes\iota^{\widetilde{n}}_{j}(\widetilde{b})\right). \tag{2.4}$$

Set

$$\begin{aligned}\bar{\iota}_n(b) &= \big(\iota^n_1(b),\dots,\iota^n_n(b)\big),\\ \bar{\iota}_{\widetilde{n}}(\widetilde{b}) &= \big(\iota^{\widetilde{n}}_1(\widetilde{b}),\dots,\iota^{\widetilde{n}}_{\widetilde{n}}(\widetilde{b})\big)\end{aligned}$$

and recall the polynomials $w_k$ from Lemma 2.5.1. Then the product (2.4) can by the defining property (a) be rewritten as

$$\sum_{k=0}^{n\widetilde{n}} w_k\left(\sigma_1\big(\bar{\iota}_n(b)\big),\dots,\sigma_n\big(\bar{\iota}_n(b)\big),\sigma_1\big(\bar{\iota}_n(\widetilde{b})\big),\dots,\sigma_{\widetilde{n}}\big(\bar{\iota}_{\widetilde{n}}(\widetilde{b})\big)\right)x^{n\widetilde{n}-k}.$$

Application of $\vartheta_{B[x]|A[x]}\otimes_{R[x]}\vartheta_{\widetilde{B}[x]|\widetilde{A}[x]}$ turns this by Proposition 2.3.7 into

$$\begin{aligned}&\sum_{k=0}^{n\widetilde{n}} w_k\left(\chi_1(b),\dots,\chi_n(b),\chi_1(\widetilde{b}),\dots,\chi_{\widetilde{n}}(\widetilde{b})\right)x^{n\widetilde{n}-k} =\\ &= \det\left(x+b\otimes\widetilde{b}\right) = \vartheta_{B[x]\otimes_{R[x]}\widetilde{B}[x]|A[x]\otimes_{R[x]}\widetilde{A}[x]}\left(x+b\otimes\widetilde{b}\right),\end{aligned}$$

the equalities by property (b) of Lemma 2.5.1 and by Proposition 2.3.7. This shows the claimed commutativity. □

## 2.6 Divided powers

There exists a close analogue to symmetric tensors, the modules and algebras of divided powers introduced by Roby in [Rob63]. They satisfy a universal property (cf. Proposition 2.6.8) better tailored to the existence and behaviour of norm maps, but at the disadvantage of having a more technical construction. In many cases these two objects are naturally isomorphic, cf. Theorem 2.6.13.

We will only briefly explain the general setting and refer to [Rob63], [Rob80], [Fer98] and [Ryd08a] for a full treatise.

Note that we will never rely on the results of this section, with the exception of Section 3.9 which explains the scheme-theoretic analogues of this section.

**Definition 2.6.1.** Let $A$ be a ring and let $M$ be any $A$-module.

Following [Rob63], the graded $A$-algebra

$$\Gamma(M|A) = \bigoplus_{d=0}^{\infty} \Gamma_d(M|A)$$

of *divided powers* with multiplication denoted by $*$ is defined as follows:

It is the universal $A$-algebra generated by symbols $\gamma_d(m) \in \Gamma_d(M|A)$, where $d \in \mathbb{N}_0$ and $m \in M$, satisfying the following relations for all $d, e \in \mathbb{N}_0$, all $a \in A$ and all $x, y \in M$:

- $\gamma_0(x) = 1$,
- $\gamma_d(am) = a^d \gamma_d(m)$,
- $\gamma_d(x + y) = \sum_{i+j=d} \gamma_i(x) * \gamma_j(y)$,
- $\gamma_d(x) * \gamma_e(x) = \binom{d+e}{d} \gamma_{d+e}(x)$.

We let $\gamma_d \colon M \to \Gamma_d(M|A)$ be the obvious maps.

If $M$ is furthermore an $A$-algebra, then there exists by [Rob80] a second multiplicative structure, induced by $\gamma_d(x)\gamma_d(y) = \gamma_d(xy)$, which turns each graded summand $\Gamma_d(M|A)$ into an $A$-algebra.

**Remark 2.6.2.** Assume that $M$ is an $A$-module. As done in [Rob63], one can check that

- there is an isomorphism $\Gamma_0(M|A) \cong A$ of $A$-algebras,
- $\gamma_1 \colon M \to \Gamma_1(M|A)$ is an isomorphism of $A$-modules, which is an isomorphism of $A$-algebras if $M$ is one.

**Remark 2.6.3.** Informally, the element $\gamma_d(x)$ satisfies the same properties as $\frac{x^d}{d!}$, assuming the latter notion to make sense. But note that this is not entirely akin to the divided power structures associated to a ring $A$ with fixed ideal $I$ because the missing relation

$$\gamma_e(\gamma_d(x)) = \frac{(de)!}{d!^e e!} \gamma_{de}(x)$$

does not even make sense.

**Remark 2.6.4.** Analogous to functorialities in Remark 2.1.3 we find $\Gamma_d(-|A)$ to give endofunctors on the categories $A\text{-}\widetilde{\mathrm{Mod}}$ and $A\text{-}\mathrm{Alg}$. It again extends to a target-preserving endofunctor on the arrow category $(A\text{-}\mathrm{Alg})^{\rightarrow}$.

This once more generalizes to modules: a morphism $f\colon A \to A'$ of rings, a morphism $g\colon M \to M'$ of $A$-modules and an $A'$-module structure on $M'$ such that $f(a)g(m) = g(am)$ for all $a \in A$, $m \in M$ induce a natural morphism $\Gamma_d(M|A) \to \Gamma_d(M'|A')$ of $A$-modules.

**Lemma 2.6.5.** *Let $A$ be a ring, let $B$ be an $A$-algebra and let $M$ be an $A$-module.*

*Then there exists a natural isomorphism*

$$\Gamma_d(M|A) \otimes_A B \cong \Gamma_d(M \otimes_A B|B)$$

*of $B$-modules.*

*Proof.* This is Théorème III.3 of [Rob63]. □

**Remark 2.6.6.** Lemma 2.6.5 should be seen as an analogue of Lemma 2.2.2 and Lemma 2.2.3. Note that there is no assumption of flatness necessary for the divided powers, which is the main reason why they work in a more general context than symmetric tensors.

The following definitions are taken from [Ryd08a] and go back to [Rob63]:

**Definition 2.6.7** (Polynomial laws)**.** Let $A$ be a ring. For any $A$-module $M$ we denote the functor $A\text{-}\mathrm{Alg} \to \mathrm{Set}$, $A' \mapsto M \otimes_A A'$, by $\mathcal{F}_M$. Note that we are therefore forgetting the module structure on $M \otimes_A A'$.

- Let $M$ and $N$ be $A$-modules. Then a *polynomial law* from $M$ to $N$ is a natural transformation $F\colon \mathcal{F}_M \implies \mathcal{F}_N$. We call it *homogeneous of degree $d$* if $F_{A'}(am) = a^d F_{A'}(m)$ for all $a \in A' \in A\text{-}\mathrm{Alg}$ and all $m \in M \otimes_A A'$.

  We denote the set of all polynomial laws from $M$ to $N$ by $\mathrm{Poly}_A(M, N)$. Its subset of homogeneous polynomial laws of degree $d$ is $\mathrm{Poly}^d_A(M, N)$.

- Let $B$ and $C$ be $A$-algebras. A polynomial law $F\colon \mathcal{F}_B \to \mathcal{F}_C$ is called *multiplicative* if $F_{A'}(1) = 1$ and $F_{A'}(xy) = F_{A'}(x)F_{A'}(y)$ for all $A' \in A\text{-}\mathrm{Alg}$ and all $x, y \in B \otimes_A A'$.

  We denote the set of all multiplicative polynomial laws from $B$ to $C$ by $\text{mult-}\mathrm{Poly}_A(B, C)$. Its subset of homogeneous multiplicative polynomial laws of degree $d$ is $\text{mult-}\mathrm{Poly}^d_A(B, C)$.

**Proposition 2.6.8.** *Let $M$ and $N$ be modules over a ring $A$ and let $d$ be a non-negative integer.*

*There exists a natural bijection*

$$\mathrm{Hom}_{\widetilde{A\text{-Mod}}}(\Gamma_d(M|A), N) \cong \mathrm{Poly}^d_A(M, N)$$

*induced by sending a morphism* $f\colon \Gamma_d(M|A) \to N$ *of* $A$*-modules to the polynomial law*

$$\begin{aligned} M \otimes_A B &\to N \otimes_A B \\ m \otimes b &\mapsto f(\gamma_d(m)) \otimes b. \end{aligned}$$

*In other words: the* $A$*-module* $\Gamma_d(M|A)$ *represents homogeneous polynomial laws* $\mathrm{Poly}^d_A(M, -)$ *of degree* $d$*.*

*Proof.* This is Théorème IV.1 of [Rob63]. □

**Proposition 2.6.9.** *Let* $B$ *and* $C$ *be algebras over a ring* $A$*.*

*Then the bijection of Proposition 2.6.8 restricts to a natural bijection*

$$\mathrm{Hom}_{A\text{-}\mathrm{Alg}}(\Gamma_d(B|A), C) \cong \text{mult-}\mathrm{Poly}^d_A(B, C)$$

*In other words: the* $A$*-algebra* $\Gamma_d(B|A)$ *represents homogeneous multiplicative polynomial laws* $\text{mult-}\mathrm{Poly}^d_A(B, -)$ *of degree* $d$*.*

*Proof.* This can be found in [Rob80] or as Proposition 2.5.1 of [Fer98]. □

Unlike symmetric tensors, they satisfy, in extension to Lemma 2.6.5, several natural isomorphisms without any assumptions of flatness:

**Proposition 2.6.10.** *Let* $A$ *be a ring and let* $M$*,* $N$ *be* $A$*-modules.*

*Then we have for all non-negative integers* $d$ *a natural isomorphism of* $A$*-modules*

*(a)*

$$\Gamma_d(M \times N|A) \cong \bigoplus_{i+j=d} \Gamma_i(M|A) \otimes_A \Gamma_j(N|A)$$

*as well as for all non-negative integers* $m, n$ *natural morphisms of* $A$*-modules*

*(b)*

$$\sigma_{m,n}\colon \Gamma_{m+n}(M|A) \to \Gamma_m(M|A) \otimes_A \Gamma_n(M|A).$$

*(c)*

$$\tau_{m,n}\colon \Gamma_{mn}(M|A) \to \Gamma_m(\Gamma_n(M|A)|A).$$

*If furthermore* $M$ *and* $N$ *are* $A$*-algebras, then the above morphisms are such of* $A$*-algebras.*

*Proof.* These are, in order, Théorème III.4 of [Rob63], (1.2.14) of [Ryd08a] and (7.1) of [Ryd08b]. All three follow from the universal property of Proposition 2.6.8. For example, the last one is induced by the polynomial law $x \mapsto \gamma_m(\gamma_n(x))$.

The versions for algebras follow analogously from Proposition 2.6.9. □

**Remark 2.6.11.** Proposition 2.6.10 gives us analogues of the natural transformations $\sigma_{m,n}$ and $\tau_{m,n}$ of Definition 2.2.7. They can also be constructed analogously to Section 2.2 by using Lemma 2.6.5 instead of Lemma 2.2.2. By construction, they correspond via Proposition 2.6.8 to the polynomial laws $x \mapsto \gamma_m(x) \otimes \gamma_n(x)$ and $x \mapsto \gamma_m\bigl(\gamma_n(x)\bigr)$.

**Definition 2.6.12.** Let $(M|B|A)$ be a good triple of rank $n$ (cf. Definition 2.3.2).

The determinant of $A$-endomorphisms of $M$ induces a multiplicative universal polynomial law $B \to A$, $b \mapsto \det(b|M)$, of degree $n$ which by Proposition 2.6.9 gives rise to a morphism

$$\theta_{M|B|A}\colon \Gamma_n(B|A) \to A$$

of $A$-algebras.

**Theorem 2.6.13.** *Let $A$ be a ring, let $M$ be an $A$-module and let $n$ be a non-negative integer.*

*There exists a natural morphism*

$$\gamma_{M|A,n}\colon \Gamma_n(M|A) \to \mathrm{S}_n(M|A)$$

*of $A$-modules which is an isomorphism if $M$ is a flat $A$-module or $\mathbb{Q} \subseteq A$. If $M$ is an $A$-algebra, then this morphism is one of $A$-algebras.*

*Proof.* To get such a morphism $\gamma_{M|A,d}$ it is by Proposition 2.6.8 enough to give a multiplicative universal polynomial law $M \to \mathrm{S}_d(M|A)$ of degree $d$. Such a law is clearly given by $b \mapsto b^{\otimes n}$. The claimed isomorphism is Corollary (4.2.5) of [Ryd08a].

If $M$ is an $A$-algebra, then $b \mapsto b^{\otimes n}$ is clearly multiplicative and hence $\gamma_{M|A,d}$ corresponds by Proposition 2.6.9 to a morphism of $A$-algebras. □

**Corollary 2.6.14.** *If $(M|B|A)$ is a good triple of rank $n$, then we have the equality*

$$\theta_{M|B|A} = \vartheta_{M|B|A} \circ \gamma_{B|A,n}.$$

*Proof.* Let $b \in B$ be arbitrary. By Proposition 2.6.9 and the definition of $\theta_{M|B|A}$ we only need to check that $\vartheta_{M|B|A}$ sends $b^{\otimes n}$ to $\det(b)$, which is Proposition 2.3.7. □

**Remark 2.6.15.** As announced at the beginning of this chapter, we focus on symmetric tensors, mostly due to them being better known and less technical to define. As we have seen, there are no relevant differences as long as everything is flat. Let us brielfy look back to highlight the differences:

The advantage of symmetric tensors is being defined via group actions, which adds a functoriality amiss in the divided powers, namely that of changing the group. In many prominent cases, including those of Convention 1.1.2, this can be mimicked by Proposition 2.6.10.

The advantage of divided powers lies, as already mentioned, in their ties to norm maps. The results of Sections 2.4 and 2.5 are for example slightly easier when using them.

In Chapter 3 we will again meet both settings, dealing with the symmetric world first and then, as in this chapter, explaining the few differences to the divided one. The main difference will again be that flatness is often crucial in the symmetric case as several functorialities rely on it, but we require no such assumption when working with divided powers.

# Chapter 3

# Translating Cycles into Multivalued Morphisms

This chapter deals with the notion of multivalued morphisms and their relation to finite correspondences. Such a comparison is to be expected as both sides are formalizations of the concept of a map that sends a single point to several. This relation was originally described over a field of characteristic 0 in Section 6 of [SV96]. They did, however, not talk about composition or other structures. We aim to go further than op. cit. and wish to fully compare the natural structures on both sides: additivity, decomposition, pushforward, pullback, correspondence map, composition, exterior product and tensor structure. This extends the results of [Ayo14b], Appendix A, where additivity and functoriality over a field of characteristic 0 were considered.

Let us describe now describe some central ideas:

An unordered $n$-tuple of elements of a set $X$ can be defined as an element of the quotient $X^n/S_n$, i.e. orbits of $n$-tuples from $X^n$ under the permuting action of the symmetric group $S_n$. Note that this allows the same element to occur multiple times and hence encompasses multiplicities.

On varieties, or more generally, schemes, this is suggestive: if this quotient exists, then the *symmetric product* $\mathrm{S}^n(X) := X^n/S_n$ (cf. Definition 3.2.1), where the symmetric group $S_n$ permutes the factors, is a natural candidate for what we want. There are some pitfalls, though: this quotient might not exist, or, more precisely, only as an algebraic space. Furthermore, it might not always satisfy the properties we want.

In particular, we will see that for schemes $X, Y \in \mathrm{Sch}_S$ a natural isomorphism $(X \times_S Y^n)/S_n \cong X \times (Y^n/S_n)$ is desirable, which in general only exists if $X$ is flat over $S$. Additionally, flatness is not always preserved by group quotients. A solution comes from Theorem 2.1.14, which implies that flatness is nonetheless preserved by symmetric products, despite this being false for arbitrary group quotients.

We then carry on to the main objective of this chapter: interpreting finite correspondences as multivalued morphisms and how this effects the aforementioned additional structures on both sides. As an immediate consequence we show that smooth commutative group schemes over a normal base admit transfers, generalizing a result of Spiess and Szamuely.

**Remark 3.0.1.** Our results, but not the proofs, can also be extracted from those of Rydh in [Ryd08b] and [Ryd08d]. Indeed, he works in a direction opposite to ours, mapping multivalued morphisms to cycles, and ultimately reaches similar conclusions. He also shows most of the desired compatibilities, but in a slightly different language: that of divided powers. Similar to Chapter 2, we will therefore link everything to Rydh's work on divided powers in Section 3.9 and briefly explain how his results interact with ours. In particular, we explain where his statements generalize those of this chapter.

## 3.1 Quotients by actions of finite groups

We briefly replicate some standard facts about quotients by finite group actions in preparation of the next Section 3.2. We have no use for the more general notions and constructions regarding algebraic group schemes, hence we stay close to [SGA1], Exposé V, and omit them. We will, however, deal with algebraic spaces.

**Definition 3.1.1.** Let $\mathcal{C}$ be a category and let $X$ be an object of $\mathcal{C}$. A *(discrete) action* of a group $G$ on $X$ is a group homomorphism $G \to \mathrm{Aut}_{\mathcal{C}}(X)$ into the automorphisms of $X$. We call it an *action on* $X|S$, where $S$ is some object of $\mathcal{C}$ under $X$, if this group homomorphism maps to the subgroup $\mathrm{Aut}_{\mathcal{C}}(X|S)$ of $S$-automorphisms, i.e. those preserving the given structure morphism $X \to S$.

By abuse of notation we will often identify $\sigma \in G$ with its image in $\mathrm{Aut}_{\mathcal{C}}(X)$.

Let now $f\colon X \to Y$ be a morphism in $\mathcal{C}$ and let the group $G$ act on $X$ and $Y$. We say that $f$ is *$G$-equivariant* if $\sigma \circ f = f \circ \sigma$ for all $\sigma \in G$.

**Definition 3.1.2.** Let $\mathcal{C}$ be a category and let $X$ be an object of $\mathcal{C}$. Let $G$ be a group acting on $X$.

A *(categorical) quotient* with respect to this action is an object $X/G$ of $\mathcal{C}$ and a morphism $\pi = \pi_{X,G}\colon X \to X/G$ invariant under the action of $G$, satisfying the following universal property: every morphism $f\colon X \to Y$ in $\mathcal{C}$ invariant under the action of $G$ on $X$ factors uniquely through $\pi$.

**Remark 3.1.3.** By standard abstract nonsense, the universal property immediately implies uniqueness of the quotient up to unique isomorphism. Existence, on the other side, is not guaranteed. Our cases of interest lie within schemes $X$, where such quotients might indeed not exist as e.g. demonstrated by Hironaka's example. We therefore extend our scope to algebraic spaces:

**Proposition 3.1.4.** *Let $X \to S$ be a separated morphism of algebraic spaces and let $G$ be a finite group acting on $X|S$.*

*Then the categorical quotient $\pi\colon X \to X/G$ exists as an algebraic space. It is a topological and geometric quotient, i.e.:*

*(a) The fibres of $\pi$ correspond to the $G$-orbits on $X$.*

*(b) The topology of $X/G$ is the quotient topology induced by $\pi$.*

*(c) $\pi$ is a finite surjective morphism.*

*(d) $\pi$ induces an isomorphism $\mathcal{O}_{X/G} \cong (\pi_* \mathcal{O}_X)^G$.*

*(e) If $V \hookrightarrow X/G$ is an open immersion, then $\pi$ restricts to an isomorphism $\pi^{-1}(V)/G \cong V$.*

*Proof.* This is Theorem 5.4 of [Ryd13], which goes back to Deligne. □

**Remark 3.1.5.** In particular, if $X = \operatorname{Spec}(A)$ is affine, then the quotient by a finite group $G$ exists as a scheme. Indeed, it is then given as $X/G = \operatorname{Spec}(A^G)$ by Proposition 3.1.4 (d), as could also be verified directly.

**Definition 3.1.6.** The action of a finite group $G$ on a scheme $X$ is called *admissible* if the $G$-orbit of every point of $X$ is contained in an open affine subset of $X$.

In [SGA1], Exposé V, Définition 1.7, admissibility is defined differently, but is equivalent to ours thanks to loc. cit. Proposition 1.8. This shows that quotients by admissible group actions behave well:

**Proposition 3.1.7.** *Let $X$ be a scheme with an admissible action by a finite group $G$. Then the quotient $\pi\colon X \to X/G$ exists as a scheme.*

**Remark 3.1.8.** Remark 4.5 of [Ryd13] shows a converse of Proposition 3.1.7: if $X \to S$ is a separated morphism of schemes and $G$ is a finite group acting on $X|S$ such that the quotient $X/G$ exists as a scheme, then the action is admissible.

**Proposition 3.1.9.** *Let $X$ and $Y$ be algebraic spaces separated over an algebraic space $S$. Assume that $Y$ is flat over $S$. Let $G$ be a finite group acting on $X|S$.*

*Then the induced action of $G$ on $X \times_S Y|Y$ induces a natural isomorphism*

$$(X \times_S Y)/G \cong (X/G) \times_S Y$$

*of algebraic spaces. If the algebraic spaces are schemes and the action of $G$ on $X$ is admissible, then so is its induced action on $X \times_S Y$ and the isomorphism is one of schemes.*

*Proof.* The existence of a morphism $(X \times_S Y)/G \to (X/G) \times_S Y$ is immediate from the universal properties. The scheme-theoretic case is [SGA1], V, Proposition 1.9. The general case follows from this by étale descent. We omit the details and point instead to [Ryd13], Proposition 2.10 in particular. □

The following observation follows directly from the definitions and the universal property of group quotients:

**Lemma 3.1.10.** *Let $H \to G$ be a morphism of finite groups. Let $X \to Y$ be a $G$-equivariant morphism of algebraic spaces separated over an algebraic space $S$.*

*Then the induced action of $H$ on $X$ and $Y$ induces a natural commutative diagram*

$$\begin{array}{ccc} X/H & \longrightarrow & X/G \\ \downarrow & & \downarrow \\ Y/H & \longrightarrow & Y/G. \end{array}$$

Let us offer a property that ensures admissibility while also having good permanence properties.

**Definition 3.1.11.** Let $f\colon X \to S$ be a morphism of schemes.

We call $X$ *(weakly) AF over $S$*, or say that the morphism $f$ is *(weakly) AF*, if every finite set of points $x_1, x_2, \ldots, x_r \in X$ over the same point $s \in S$ is contained in an open subscheme $U \hookrightarrow X$ which is affine over $S$.

We will generally omit 'weakly', but refer to Remark 3.1.13 for its relevance and a *strong* version.

AF morphisms satisfy several useful properties as pointed out in Remarks 3.1.3 and 3.1.4 of [Ryd08c]. Most importantly we have:

**Lemma 3.1.12.** *We have the following permanence properties of AF morphisms:*

*(a) Every AF morphism is separated.*

*(b) If $X \to S$ is an AF morphism and $\Gamma \subseteq X$ is a closed subscheme, then $\Gamma \to S$ is AF.*

*(c) Let $X \to S$ be a separated morphism of schemes and let $Y$, $Y'$ be schemes over $X$ that are AF over $S$. Then $Y \times_X Y'$ is AF over $S$.*

*(d) Let $f\colon Y \to X$ and $g\colon X \to S$ be morphisms of schemes. If $g \circ f$ is AF and $g$ is separated, then $f$ is AF.*

*(e) Let $f\colon X \to S$ be an AF morphism of schemes. If $G$ is a finite group acting on $X|S$, then this action on $X$ is admissible. The quotient $X/G$, which therefore exists as a scheme by Proposition 3.1.7, is furthermore AF over $S$.*

*Proof.*

(a) Let $X \to S$ be AF, let $x \in X \times_S X$ and let $x_1, x_2 \in X$ be its projections to the factors. We have to show that the diagonal $\Delta_X$ in $X \times_S X$ is closed. By assumption we find an open subscheme $U \hookrightarrow X$ which is affine over $S$ and contains $x_1, x_2$. Thus $x$ lies in the open subscheme $U \times_S U$ whose diagonal $\Delta_U = (U \times_S U) \cap \Delta X$ is closed. The result now follows by varying $x$ over $X \times_S X$.

(b) This follows directly from closed immersions being affine.

(c) The natural morphism $Y \times_X Y' \to Y \times_S Y'$ is a closed immersion because it is the pullback of the closed immersion $\Delta : X \to X \times_S X$ along $Y \times_S Y' \to X \times_S X$. Thus by the previous part we may assume that $X = S$.

If the points $z_1, z_2, \dots, z_r \in Y \times_S Y'$ lie over the same point $s \in S$ we let $y_i$ and $y'_i$ be their projections to $Y$ and $Y'$, respectively. Then we find open subschemes $U \hookrightarrow Y$ and $U' \hookrightarrow Y'$ which are affine over $S$ and contain the respective points $y_i$ and $y'_i$. Then $U \times_S U'$ is an open neighbourhood of each $z_i$. It is affine over $S$ as fibre product of affine schemes over $S$.

(d) Let $y_1, y_2, \dots, y_r$ be over the same point $x \in X$, thus especially over the same point $g(x) \in S$. Then by $g \circ f$ being AF we find an open $V \subseteq Y$ affine over $S$ and containing each $y_i$. Hence $V$ is affine over $X$ by Proposition 12.3 (3) of [GW10]. Alternatively, this follows easily from the previous part by considering the closed immersion $Y \to Y \times_X Y$.

(e) Every orbit of the action of $G$ is finite and lies over a single point $s \in S$. Because $X$ is AF over $S$ we find an open subscheme $U \hookrightarrow X$ affine over $S$ and containing the orbit. Choosing an affine open neighbourhood of $s$ and taking its preimage in $U$ shows that every orbit is contained in an affine open subscheme, i.e. the $G$-action is admissible. Thus the quotient $X/G$ exists by Proposition 3.1.7.

Now let $Q$ be any finite set of points in $X/G$ over the same point $s \in S$. Then, because $\pi \colon X \to X/G$ is finite by Proposition 3.1.4 (c), its preimage $P = \pi^{-1}(Q)$ in $X$ is finite. Trivially, all of $P$ lies above $s$, thus by $X$ being AF over $S$ we find an open $U \subseteq X$ affine over $S$ and containing $P$. We may replace $U$ by $\bigcap_{g \in G} gU$. Indeed, it is open, contains $P$ and is affine over $S$, where the last part can be shown similarly to Lemma 7.1.1. Thus we may assume $U$ to be $G$-invariant.

We claim that the image $\pi(U)$ is the desired open subscheme of $X/G$. It is open by Proposition 3.1.4 (c) and, as a scheme, isomorphic to $U/G$ by Proposition 3.1.4 (e). It also clearly contains $Q$. If $S' \hookrightarrow S$ is an open affine subscheme, we let $V = f^{-1}(S') \cap U$ be its preimage in $U$, which by assumption is affine, say $V = \operatorname{Spec}(A)$. Then by Proposition 3.1.4 (d) the corresponding open subscheme $\pi(V) = V/G = \operatorname{Spec}(A^G)$ of $X/G$ is affine as well, as required. □

**Remark 3.1.13.** There also exists a more widespread version of AF, clearly implying the weak one: a morphism $X \to S$ of schemes is *(strongly) AF* if every finite set of points $x_1, \dots, x_r \in X$, not necessarily over the same point of $S$, is contained in an open subscheme affine over $S$.

After minor changes the proof of Lemma 3.1.12 works for strongly AF morphisms. Furthermore it can be shown that if $f\colon X \to Y$ and $g\colon Y \to Z$ are morphisms such that $g$ is strongly AF, then $g \circ f$ is weakly (respectively strongly) AF if and only if $g$ is weakly (respectively strongly) AF.

Finally, quasi-projective morphisms $X \to S$ are strongly and hence weakly AF, as can be seen using an easy hyperplane argument. We omit the details as we will not need this property. If the base $S$ is noetherian, then the quotient $X/G$ by an action of a finite group $G$ on $X/S$ is quasi-projective again, see [Ryd08c], Remark 2.3.3.

## 3.2 Symmetric products

This and the next section are primarily a scheme-theoretic version of Chapter 2: we introduce symmetric products and their basic properties.

**Definition 3.2.1.** Let $Y \to X$ be a morphism of algebraic spaces and let $n$ be a non-negative integer. We have a natural action of the symmetric group $S_n$ on the $n$-fold fibre product $(Y|X)^{\times n} := Y \times_X Y \times_X \cdots \times_X Y$ by permuting the $n$ factors.

We then define the *symmetric product*

$$\mathrm{S}^n(Y|X) := (Y|X)^{\times n}/S_n,$$

which by Proposition 3.1.4 exists as an algebraic space,

We will simply write $\mathrm{S}^n(Y)$ if the base $X$ is clear from the context.

Note that this definition makes sense if $n = 0$, as then $S_0$ is the trivial group and $X = (Y|X)^0 = \mathrm{S}^0(Y|X)$. If the schemes are affine, or if more generally the morphism $Y \to X$ is affine, then the symmetric products are according to Proposition 3.1.4 (d) given by

$$\mathrm{S}^n(Y|X) \cong \operatorname{Spec}\left(\left((\mathcal{O}_Y|\mathcal{O}_X)^{\otimes n}\right)^{S_n}\right) = \operatorname{Spec}\left(\mathrm{S}_n\left(\mathcal{O}_Y|\mathcal{O}_X\right)\right).$$

**Example 3.2.2.** The fundamental theorem of symmetric polynomials states that $\mathrm{S}^n(\mathbb{A}^1_S|S) \cong \mathbb{A}^n_S$. Consequently we also find $\mathrm{S}^n(\mathbb{P}^1_S|S) \cong \mathbb{P}^n_S$.

**Lemma 3.2.3.** *Let $n$ be a non-negative integer. Let $X \to S$ be a flat and separated morphism of algebraic spaces.*

*Then $\mathrm{S}^n(X|S)$ is flat over $S$.*

*Proof.* This can be checked locally, where it amounts to Theorem 2.1.14 (b). □

**Lemma 3.2.4.** *Let $n$ be a non-negative integer. Let $X \to S$ be a separated morphism of schemes and let $Y$ be a scheme over $X$ which is AF over $S$.*

*Then the symmetric product $\mathrm{S}^n(Y|X)$ exists as a scheme and is furthermore AF over $S$.*

*Proof.* By Lemma 3.1.12 (c) we find $(Y|X)^n$ to be AF over $S$. Therefore by Lemma 3.1.12 (e) the quotient $\mathrm{S}^n(Y|X)$ exists and is AF over $S$. □

As an exemplary consequence we see that $\mathrm{S}^m(\mathrm{S}^n(X|S')|S)$ exists as a scheme for any chain of morphisms $X \to S' \to S$ such that $S'$ is separated over $S$ and $X$ is AF over $S$.

**Remark 3.2.5.** Let $S$ be a scheme. Then the $n$-fold fibre product $(Y|X) \mapsto (Y|X)^n$, the latter seen as a scheme over $X$, defines base-preserving endofunctors $(-|-)^n$ on $\mathrm{Sch}_S^{\to}$, $\widetilde{\mathrm{Sch}}_S^{\to}$ and $\mathrm{AlgSp}_S^{\to,\mathrm{sep}}$. The endofunctor on the latter also exists when $S$ is only an algebraic space.

The universal properties of group quotients as well as Lemma 3.2.4 then readily imply that $\mathrm{S}^n(-|-)$ gives base-preserving endofunctors on $(\mathrm{Sch}_S^{\mathrm{aff}})^{\to}$, $\mathrm{Sch}_S^{\to,\mathrm{AF}}$ and $\mathrm{AlgSp}^{\to,\mathrm{sep}}$. They restrict by Lemma 3.2.3 to the subcategories of flat morphisms.

To avoid confusion we will often make the base explicit, for example as $\mathrm{S}^n(Y|X)|X$.

## 3.3 Functorialities of symmetric products

This section is dedicated to three natural transformations $\tau^{m,n}$, $\sigma^{m,n}$ and $\rho^{m,n}$ as well as their properties and compatibilities. We will see in Sections 3.5 to 3.7 that they correspond on cycles to composition, addition and the tensor product.

**Remark 3.3.1.** Over a field, most of the contents of this section, with the exception of the existence and the properties of the morphisms $\rho^{m,n}$ corresponding to the tensor structure, were stated without explicit proofs in Appendix A of [Ayo14b]. There, a different setting was used: that of a category closed under finite products and coproducts satisfying a formal analogue of the decomposition into connected components and admitting all quotients by actions of finite groups. Furthermore, the quotients were assumed to satisfy $(X \times Y)/G \cong X \times (Y/G)$ for all finite groups $G$ acting on an object $Y$. The latter, however, is not automatically true for arbitrary schemes, as it requires flatness, which in return is not preserved by finite group quotients over general bases. Our proofs will work with minor changes in Ayoub's setting as well. See also Remark 3.3.10 below for additional details.

**Lemma 3.3.2.** *Let non-negative integers $m$, $n$, $r$ and $n_1, \ldots n_r$ be given and let $N = \sum_{i=1}^r n_i$. Recall that Convention 1.1.2 identified $[m] \times [n]$ with $[mn]$ and thereby induced a morphism from the semi-direct product $S_m \ltimes S_n^m$ to $S_{mn}$. It also identified the product $\prod_{i=1}^r S_{n_i}$ with a subgroup of $S_{n_1+\ldots+n_r}$.*

*Furthermore let $S$ be an algebraic space and let $X$ and $X_1, \ldots, X_r$ be algebraic spaces which are flat and separated over $S$. Then the identifications of Convention 1.1.2 induce natural isomorphisms*

*(a)*

$$\prod_{i=1}^r \mathrm{S}^{n_i}(X_i|S) \cong \left(\prod_{i=1}^r (X_i|S)^{n_i}\right) / \left(\prod_{i=1}^r S_{n_i}\right),$$

*(b)*

$$\mathrm{S}^m(\mathrm{S}^n(X|S)|S) \cong (X|S)^{mn} / \left(S_m \ltimes S_n^m\right).$$

*Proof.* We will refrain from repeatedly writing down the base $S$ over which all symmetric products and fibre products are to be understood.

Let $G$ be a subgroup of $S_k$ acting on $X/S$. Each $\mathrm{S}^{n_i}(X_i|S)$ is flat over $S$ by Lemma 3.2.3. Thus, using Proposition 3.1.9 and the universal properties, we get natural isomorphisms

$$\begin{aligned}(X/G) \times \mathrm{S}^{n_r}(X_r) &\cong \left(X \times (X_r^{n_r}/S_{n_r})\right)/G \cong \\ &\cong \left((X \times X_r^{n_r})/S_{n_r}\right)/G \cong (X \times X_r^{n_r})/(G \times S_{n_r}).\end{aligned}$$

Therefore (a) follows by induction on $r$.

Now (b) follows from (a) via

$$\begin{aligned}\mathrm{S}^m(\mathrm{S}^n(X)) &= \left((X^n)/S_n\right)^m/S_m \cong \\ &\cong \left(((X^n)^m)/S_n^m\right)/S_m = X^{mn}/\left(S_m \ltimes S_n^m\right),\end{aligned}$$

where the last equality is the identification of Convention 1.1.2.

The functoriality is immediate from that of $\mathrm{S}^n(-|-)$ and fibre products. □

**Proposition 3.3.3.** *Let $X \to S' \to S$ be separated morphisms of algebraic spaces. Assume that $X \to S$ is flat. Furthermore let non-negative integers $m$ and $n$ be given. We identify $[m] \times [n]$ with $[mn]$ as in Convention 1.1.2.*

*Then there exists a unique morphism*

$$\tau^{m,n}\colon\ \mathrm{S}^m(\mathrm{S}^n(X|S')|S) \to \mathrm{S}^{mn}(X|S)$$

*of algebraic spaces over $S$ such that the diagram*

$$\begin{array}{ccc}
((X|S')^n|S)^m & \longrightarrow & ((X|S)^n|S)^m \\
\big\downarrow\!\!\!\downarrow & & \cong\big\downarrow \text{Conv. 1.1.2} \\
\mathrm{S}^m\left((X|S')^n|S\right) & & (X|S)^{mn} \\
\big\downarrow\!\!\!\downarrow & & \big\downarrow\!\!\!\downarrow \\
\mathrm{S}^m(\mathrm{S}^n(X|S')|S) & \xrightarrow{\tau^{m,n}} & \mathrm{S}^{mn}(X|S)
\end{array}$$

*commutes.*

*It induces a natural transformation, i.e. if*

$$\begin{array}{ccccc}
X & \longrightarrow & S' & \longrightarrow & S \\
\big\downarrow & & \big\downarrow & & \big\downarrow \\
Y & \longrightarrow & T' & \longrightarrow & T
\end{array}$$

*is a commutative diagram of algebraic spaces with the above conditions on each row, then the diagram*

$$\begin{array}{ccc}
\mathrm{S}^m(\mathrm{S}^n(X|S')|S) & \xrightarrow{\tau^{m,n}} & \mathrm{S}^{mn}(X|S) \\
\big\downarrow & & \big\downarrow \\
\mathrm{S}^m(\mathrm{S}^n(Y|T')|T) & \xrightarrow{\tau^{m,n}} & \mathrm{S}^{mn}(Y|T)
\end{array}$$

*commutes.*

*Proof.* The functorialities of $\mathrm{S}^m(-|-)$ and $\mathrm{S}^n(-|-)$ give a natural morphism $\mathrm{S}^m(\mathrm{S}^n(X|S')|S) \to \mathrm{S}^m(\mathrm{S}^n(X|S)|S)$, allowing us to assume $S = S'$.

We identify $\mathrm{S}_m \ltimes \mathrm{S}_n^m$ with a subgroup of $\mathrm{S}_{mn}$ as in Convention 1.1.2. Then the existence follows directly from Lemma 3.3.2 (b) and Lemma 3.1.10. The uniqueness, and therefore the naturality, follow from the universal property of quotients by a group action. □

**Proposition 3.3.4.** *Let non-negative integers $r$ and $n_1, n_2, \ldots, n_r$ be given and set $N = \sum_{i=1}^r n_i$. Let $X_1, X_2, \ldots, X_r$ be algebraic spaces which are flat and separated over an algebraic space $S$. Denote the inclusions $X_i \hookrightarrow \coprod_{j=1}^r X_j$ by $\iota_i$.*

*Then there exists a unique morphism*

$$\sigma^{(n_i)}\colon\ \prod_{i=1}^r \mathrm{S}^{n_i}(X_i|S) \to \mathrm{S}^N\left(\coprod_{j=1}^r X_j|S\right)$$

*of algebraic spaces over $S$ such that the diagram*

$$\begin{array}{ccccc}
\prod_{i=1}^{r}(X_i|S)^{n_i} & \xrightarrow{\prod_{i=1}^{r}(\iota_i|S)^{n_i}} & \prod_{i=1}^{r}\left(\left(\coprod_{j=1}^{r}X_j|S\right)^{n_i}\right) & \xrightarrow[\cong]{\text{Conv. 1.1.2}} & \left(\coprod_{j=1}^{r}X_j|S\right)^{N} \\
\downarrow & & & & \downarrow \\
\prod_{i=1}^{r}\mathrm{S}^{n_i}(X_i|S) & & \xrightarrow{\sigma^{(n_i)}} & & \mathrm{S}^{N}\left(\coprod_{j=1}^{r}X_j|S\right)
\end{array}$$

*commutes.*

*It induces a natural transformation, i.e. if we are given commutative diagrams*

$$\begin{array}{ccc}
X_i & \longrightarrow & S \\
\downarrow & & \downarrow \\
Y_i & \longrightarrow & T
\end{array}$$

*of algebraic spaces for all $i \in \{1, 2, \ldots, r\}$, satisfying the above conditions on each row, then the diagram*

$$\begin{array}{ccc}
\prod_{i=1}^{r}\mathrm{S}^{n_i}(X_i|S) & \xrightarrow{\sigma^{(n_i)}} & \mathrm{S}^{N}\left(\coprod_{j=1}^{r}X_j|S\right) \\
\downarrow & & \downarrow \\
\prod_{i=1}^{r}\mathrm{S}^{n_i}(Y_i|T) & \xrightarrow{\sigma^{(n_i)}} & \mathrm{S}^{N}\left(\coprod_{j=1}^{r}Y_j|T\right)
\end{array}$$

*commutes.*

*Proof.* We identify $\prod_{i=1}^{r}\mathrm{S}_{n_i}$ with a subgroup of $\mathrm{S}_N$ as in Convention 1.1.2. Then the existence follows directly from Lemma 3.3.2 (a) and Lemma 3.1.10. The uniqueness, and therefore the naturality, follow from the universal property of quotients by a group action. □

**Proposition 3.3.5.** *Let $S$ be an algebraic space. Let $X$, $Y$ and $X_i$, $i \in I$, be algebraic spaces which are separated and flat over $S$. Let $k$, $m$ and $n$ be non-negative integers. Denote the inclusions of $X$ and $Y$ into their disjoint union $X \sqcup Y$ by $x$ and $y$, respectively.*

*Then the composition*

$$\sigma^{m,n} \circ (\mathrm{S}^m(x|S), \mathrm{S}^n(y|S))_S :$$
$$\mathrm{S}^m(X|S) \times_S \mathrm{S}^n(Y|S) \to \mathrm{S}^m(X \sqcup Y|S) \times_S \mathrm{S}^n(X \sqcup Y|S) \to \mathrm{S}^{m+n}(X \sqcup Y|S)$$

*is an open and closed immersion.*

*Furthermore, it induces a natural decomposition*

$$\mathrm{S}^k\left(\coprod_{i\in I} X_i|S\right) \cong \coprod_{\substack{(d_i)\in(\mathbb{N}_0)^I\\ \sum_{i\in I} d_i=k}} \prod_{\substack{i\in I\\ d_i\neq 0}} \mathrm{S}^{d_i}(X_i|S).$$

*Proof.* For $j=(j_1,\ldots,j_k)\in I^k$ we set $X_j := X_{j_1}\times_S\cdots\times_S X_{j_k}$.

We hence get a decomposition

$$\left(\coprod_{j\in I} X_i|S\right)^k \cong \coprod_{j\in I^k} X_j \cong \coprod_{s\in I^k/S_k}\coprod_{j\in s} X_i$$

into $S_k$-orbits $Z_s := \coprod_{j\in s} X_j$, $s\in I^k/S_k$. Then the universal property of group quotients implies

$$\mathrm{S}^k\left(\coprod_{j\in I} X_j|S\right) \cong \coprod_{s\in I^k/S_k} (Z_s/S_k).$$

We now fix an $s\in I^k/S_k$. For $i\in I$ and an arbitrary $j\in s$ we let $d_i := \#\{m\in[k] \mid j_m = i\}$ be the number of occurrences of $i$ in $j$. Note that it does not depend on the choice of $j$, but only on $s$ and $i$. Also set $I_s := \{i\in I \mid d_i\neq 0\}$.

Thus it suffices to check that $\sigma^{(d_i)}$ induces an isomorphism

$$\prod_{i\in I_s} \mathrm{S}^{d_i}(X_i|S) \cong Z_s/S_k.$$

To see this, we let $H=\prod_{i\in I_s} S_{d_i}$, which we identify with the stabilizer of $j$. Hence we have to show that

$$\prod_{i\in I_s} \mathrm{S}^{d_i}(X_i|S) \cong \left(\prod_{i\in I_s} X_i^{d_i}\right)/H,$$

which follows inductively from Lemma 3.3.2 (a). □

**Remark 3.3.6.** Note that Proposition 3.3.5 can alternatively be used to define $\sigma^{m,n}$ as the inclusion of the component $\mathrm{S}^m(X|S)\times_S \mathrm{S}^n(X|S)$ into $\mathrm{S}^{m+n}(X\sqcup X|S)$, followed by the morphism $\mathrm{S}^{m+n}(X\sqcup X|S)\to \mathrm{S}^{m+n}(X|S)$.

**Proposition 3.3.7.** *Let non-negative integers $r$ and $n_1,n_2,\ldots,n_r$ be given. Set $N=\prod_{i=1}^r n_i$ and let $m_i$ be the integers such that $m_i n_i = N$. Let $X_1,X_2,\ldots,X_r$ be algebraic spaces which are flat and separated over an algebraic space $S$. Let*

$$\Delta_i\colon (X_i|S)^{n_i} \to ((X_i|S)^{n_i}|S)^{m_i} \cong (X_i|S)^N$$

*be the $m_i$-fold diagonal morphism, where the isomorphism comes from Convention 1.1.2.*

*Then there exists a unique morphism*

$$\rho^{(n_i)} \colon \prod_{i=1}^{r} \mathrm{S}^{n_i}(X_i|S) \longrightarrow \mathrm{S}^N\left(\prod_{i=1}^{r} X_i|S\right)$$

*of schemes over $S$ such that the diagram*

$$\begin{array}{ccccc} \prod_{i=1}^{r}(X_i|S)^{n_i} & \xrightarrow{\prod_{i=1}^{r}\Delta_i} & \prod_{i=1}^{r}(X_i|S)^N & \xrightarrow{\cong} & \left(\prod_{i=1}^{r}(X_i|S)\right)^N \\ \downarrow & & & & \downarrow \\ \prod_{i=1}^{r}\mathrm{S}^{n_i}(X_i|S) & & \xrightarrow{\rho^{(n_i)}} & & \mathrm{S}^N\left(\prod_{i=1}^{r} X_i|S\right) \end{array}$$

*commutes.*

*It induces a natural transformation, i.e. if we are given commutative diagrams*

$$\begin{array}{ccc} X_i & \longrightarrow & S \\ \downarrow & & \downarrow \\ Y_i & \longrightarrow & T \end{array}$$

*for all $i \in \{1, 2, \ldots, r\}$, with the above conditions on each row, then the diagram*

$$\begin{array}{ccc} \prod_{i=1}^{r}\mathrm{S}^{n_i}(X_i|S) & \xrightarrow{\rho^{(n_i)}} & \mathrm{S}^N\left(\prod_{i=1}^{r} X_i|S\right) \\ \downarrow & & \downarrow \\ \prod_{i=1}^{r}\mathrm{S}^{n_i}(Y_i|T) & \xrightarrow{\rho^{(n_i)}} & \mathrm{S}^N\left(\prod_{i=1}^{r} Y_i|T\right) \end{array}$$

*commutes.*

*Proof.* Convention 1.1.2 fixes a morphism $G := \prod_{i=1}^{r} S_{n_i} \to S_N$. We let $G$ act on $(X_i|S)^{n_i}$ via its projection to $S_{n_i}$ and note that the diagonal morphisms $\Delta_i$ are $G$-equivariant. This induces a $G$-equivariant morphism

$$\prod_{i=1}^{r} X_i^{n_i} \to \left(\prod_{i=1}^{r} X_i|S\right)^N.$$

Combining it with Lemma 3.3.2 (a) and Lemma 3.1.10 we thus get a natural morphism

$$\prod_{i=1}^{r} \mathrm{S}^{n_i}(X_i|S) \cong \left(\prod_{i=1}^{r}(X_i|S)^{n_i}\right)/G \to \left(\prod_{i=1}^{r} X_i|S\right)^{N}/G \to \mathrm{S}^N\left(\prod_{i=1}^{r} X_i|S\right)$$

as desired. Uniqueness and naturality follow from the universal property of quotients by a group action. □

**Lemma 3.3.8.** *The $S$ be a scheme and let $X, Y$ be flat over $S$. Then for all non-negative integers $m, n$ we have a commutative diagram*

$$\begin{array}{ccc}
\mathrm{S}^m\big(X \times_S \mathrm{S}^n(Y|S)|S\big) & \xleftarrow{\;\cong\;} & \mathrm{S}^m\big(\mathrm{S}^n(X \times_S Y|X)|S\big) \\
\uparrow & & \\
\mathrm{S}^m\big(X \times_S \mathrm{S}^n(Y|S)|\,\mathrm{S}^n(Y|S)\big) & & \Big\downarrow \tau^{m,n} \\
\Big\downarrow \cong & & \\
\mathrm{S}^m(X|S) \times_S \mathrm{S}^n(Y|S) & \xrightarrow{\;\rho^{m,n}\;} & \mathrm{S}^{mn}(X \times_S Y|S)
\end{array}$$

*induced by the functorialities of symmetric products, where the isomorphisms are those of Lemma 3.3.2.*

*Proof.* All occurring objects are flat by Lemma 3.2.3. In particular, all the (iso)morphisms are defined.

Due to the defining properties of $\tau^{m,n}$ and $\rho^{m,n}$, the commutativity of the lemma reduces to that of the diagram

$$\begin{array}{ccc}
\big(X \times_S (Y|S)^n|S\big)^m & \xleftarrow{\;\cong\;} & \big((X \times_S Y|X)^n|S\big)^m \\
\uparrow & & \\
\big(X \times_S (Y|S)^n|(Y|S)^n\big)^m & & \Big\downarrow \text{Conv. 1.1.2} \\
\Big\downarrow \cong & & \\
(X|S)^m \times_S (Y|S)^n & \xrightarrow[\cong]{\text{Conv. 1.1.2}} & (X \times_S Y|S)^{mn},
\end{array}$$

which is trivial. □

**Proposition 3.3.9** (Compatibilities of $\tau^{m,n}$, $\sigma^{m,n}$ and $\rho^{m,n}$)**.** *Let $S$ be an algebraic space and let $k, l, m, n$ be non-negative integers. Let $X$, $Y$ and $Z$ be algebraic spaces which are flat and separated over $S$. We now omit the base $S$ in the sense that we use $\times$ for $\times_S$ and $\mathrm{S}^n(X)$ for $\mathrm{S}^n(X|S)$.*

*The morphisms $\sigma^{m,n}$, $\tau^{m,n}$ and $\rho^{m,n}$ defined by Propositions 3.3.3, 3.3.4 and 3.3.7 satisfy, in addition to the naturalities stated in the aforementioned propositions, also the following compatibilities, whose names will become clear by Theorem 3.4.11:*

- *Commutativity of addition:*

$$\begin{array}{ccc}
\mathrm{S}^m(X)\times\mathrm{S}^n(X) & \xrightarrow{\sigma^{m,n}} & \mathrm{S}^{m+n}(X) \\
\downarrow{\scriptstyle \mathit{swap}} & & \Vert \\
\mathrm{S}^n(X)\times\mathrm{S}^m(X) & \xrightarrow{\sigma^{n,m}} & \mathrm{S}^{m+n}(X)
\end{array}$$

- *Associativity of addition:*

$$\begin{array}{ccc}
\mathrm{S}^k(X)\times\mathrm{S}^m(X)\times\mathrm{S}^n(X) & \xrightarrow{\mathrm{S}^k(X)\times\sigma^{m,n}} & \mathrm{S}^k(X)\times\mathrm{S}^{m+n}(X) \\
\downarrow{\scriptstyle \sigma^{k,m}\times\mathrm{S}^n(X)} & & \downarrow{\scriptstyle \sigma^{k,m+n}} \\
\mathrm{S}^{k+m}(X)\times\mathrm{S}^n(X) & \xrightarrow{\sigma^{k+m,n}} & \mathrm{S}^{k+m+n}(X)
\end{array}$$

- *Associativity of composition*

$$\begin{array}{ccc}
\mathrm{S}^k(\mathrm{S}^m(\mathrm{S}^n(X))) & \xrightarrow{\mathrm{S}^k(\tau^{m,n})} & \mathrm{S}^k(\mathrm{S}^{mn}(X)) \\
\downarrow{\scriptstyle \tau^{k,m}} & & \downarrow{\scriptstyle \tau^{k,mn}} \\
\mathrm{S}^{km}(\mathrm{S}^n(X)) & \xrightarrow{\tau^{km,n}} & \mathrm{S}^{kmn}(X)
\end{array}$$

- *Left distributivity of addition and composition:*

$$\begin{array}{ccc}
\mathrm{S}^m(\mathrm{S}^k(X))\times\mathrm{S}^n(\mathrm{S}^k(X)) & \xrightarrow{\tau^{m,k}\times\tau^{n,k}} & \mathrm{S}^{km}(X)\times\mathrm{S}^{kn}(X) \\
\downarrow{\scriptstyle \sigma^{m,n}} & & \downarrow{\scriptstyle \sigma^{km,kn}} \\
\mathrm{S}^{m+n}(\mathrm{S}^k(X)) & \xrightarrow{\tau^{m+n,k}} & \mathrm{S}^{k(m+n)}(X)
\end{array}$$

- *Right distributivity of addition and composition:*

$$\begin{array}{ccc}
\mathrm{S}^k(\mathrm{S}^m(X)\times\mathrm{S}^n(X)) & \xrightarrow[\times(\tau^{k,n}\circ\mathrm{S}^k(\mathrm{pr}^2))]{(\tau^{k,m}\circ\mathrm{S}^k(\mathrm{pr}^1))\times} & \mathrm{S}^{km}(X)\times\mathrm{S}^{kn}(X) \\
\downarrow{\scriptstyle \sigma^{m,n}} & & \downarrow{\scriptstyle \mathrm{S}^k(\sigma^{km,kn})} \\
\mathrm{S}^k(\mathrm{S}^{m+n}(X)) & \xrightarrow{\tau^{k,m+n}} & \mathrm{S}^{k(m+n)}(X)
\end{array}$$

- *Symmetry of multiplication:*

$$\begin{array}{ccc}
\mathrm{S}^m(X)\times\mathrm{S}^n(Y) & \xrightarrow{\rho^{m,n}} & \mathrm{S}^{mn}(X\times Y) \\
\downarrow{\scriptstyle \mathit{swap}} & & \downarrow{\scriptstyle \mathrm{S}^{mn}(\mathit{swap})} \\
\mathrm{S}^n(Y)\times\mathrm{S}^m(X) & \xrightarrow{\rho^{n,m}} & \mathrm{S}^{mn}(Y\times X)
\end{array}$$

- *Associativity of multiplication:*

$$\begin{array}{ccc}
\mathrm{S}^k(X) \times \mathrm{S}^m(Y) \times \mathrm{S}^n(Z) & \xrightarrow{\mathrm{S}^k(X)\times\rho^{m,n}} & \mathrm{S}^k(X) \times \mathrm{S}^{mn}(Y \times Z) \\
\Big\downarrow{\scriptstyle \rho^{k,m}\times\mathrm{S}^n(Z)} & & \Big\downarrow{\scriptstyle \rho^{k,mn}} \\
\mathrm{S}^{km}(X \times Y) \times \mathrm{S}^n(Z) & \xrightarrow{\rho^{km,n}} & \mathrm{S}^{kmn}(X \times Y \times Z)
\end{array}$$

- *Distributivity of addition and multiplication:*

$$\begin{array}{ccc}
\mathrm{S}^k(X) \times (\mathrm{S}^m(Y) \times \mathrm{S}^n(Y)) & \xrightarrow{\Delta(\mathrm{S}^k(X))} & \begin{array}{c}(\mathrm{S}^k(X) \times \mathrm{S}^m(Y))\times \\ \times(\mathrm{S}^k(X) \times \mathrm{S}^n(Y))\end{array} \\
\Big\downarrow{\scriptstyle \mathrm{S}^k(X)\times\sigma^{m,n}} & & \Big\downarrow{\scriptstyle \rho^{k,m}\times\rho^{k,n}} \\
\mathrm{S}^k(X) \times \mathrm{S}^{m+n}(Y) & & \mathrm{S}^{km}(X \times Y) \times \mathrm{S}^{kn}(X \times Y) \\
\Big\downarrow{\scriptstyle \rho^{k,m+n}} & & \Big\downarrow{\scriptstyle \sigma^{km,kn}} \\
\mathrm{S}^{k(m+n)}(X \times Y) & = & \mathrm{S}^{k(m+n)}(X \times Y)
\end{array}$$

- *Functoriality of multiplication:*

$$\begin{array}{ccc}
\mathrm{S}^k(\mathrm{S}^l(X)) \times \mathrm{S}^m(\mathrm{S}^n(Y)) & \xrightarrow{\tau^{k,l}\times\tau^{m,n}} & \mathrm{S}^{kl}(X) \times \mathrm{S}^{mn}(Y) \\
\Big\downarrow{\scriptstyle \rho^{k,m}} & & \Big\downarrow \\
\mathrm{S}^{km}(\mathrm{S}^l(X) \times \mathrm{S}^n(Y)) & & {\scriptstyle \rho^{kl,mn}} \\
\Big\downarrow{\scriptstyle \mathrm{S}^{km}(\rho^{l,n})} & & \Big\downarrow \\
\mathrm{S}^{km}(\mathrm{S}^{ln}(X \times Y)) & \xrightarrow{\tau^{km,ln}} & \mathrm{S}^{klmn}(X \times Y).
\end{array}$$

*Proof.* This can be checked in a tedious but straightforward way similar to the proof of Lemma 3.3.8: all properties correspond via the inclusions of Convention 1.1.2 by the definition of $\tau^{m,n}$, $\sigma^{m,n}$ and $\rho^{m,n}$ to commuting diagrams of fibre products. □

**Remark 3.3.10.** Let $\mathcal{C}$ be a category that admits finite coproducts as well as finite products. Furthermore assume that it admits quotients by actions of finite groups and which satisfy $(X \times Y)/G \cong X \times (Y/G)$ for all $X, Y \in \mathcal{C}$ and finite groups $G$ acting on $Y$.

Then the proofs of this section work mutatis mutandis for the objects $\mathrm{S}^n(X) := X^n/\,\mathrm{S}_n$. We recover our more specific setting by letting $\mathcal{C}$ be the category of separated algebraic spaces flat over a base $S$, but have to be careful: only the symmetric products, not quotients in general, preserve the flatness that guarantees the isomorphism $(X \times Y)/G \cong X \times (Y/G)$ of Lemma 3.2.3. This, however, was sufficient for our arguments. Alternatively, one could by Lemmas 3.2.3 and 3.2.4 use the category of separated schemes which are flat and AF over a separated scheme $S$.

## 3.4 Multivalued morphisms

This section is based on Appendix 1 of [Ayo14b], but see Section 10 of [Ryd08b] for a closely related yet different approach.

The following formalizes the idea of a map $X \multimap Y$ that associates to every $x \in X$ a finite set of elements of $Y$:

**Definition 3.4.1.** Let $S$ be scheme and let $X$, $Y$ be schemes over $S$. An *effective multivalued morphism* $X \multimap Y$ over $S$ is a morphism

$$X \to \coprod_{d=0}^{\infty} \mathrm{S}^d(Y|S)$$

of schemes over $S$. We call it *homogeneous of degree $d$* if it is simply a morphism $X \to \mathrm{S}^d(Y|S)$.

We denote the set of all effective multivalued morphisms $X \multimap Y$ by $\mathrm{Multi}_S^{\mathrm{eff}}(X, Y)$.

**Remark 3.4.2.** Note that we assumed the basic objects to be schemes, yet the resulting symmetric products $\mathrm{S}^d(Y|S)$ might only exist as algebraic spaces. One could without any other change define $\mathrm{Multi}_S^{\mathrm{eff}}$ for algebraic spaces which are flat and separated over an algebraic space $S$. We will further explore such generalities in Section 3.9.

**Definition 3.4.3** (Addition of multivalued morphisms)**.** Let $S$ be a scheme and let $X$, $Y$ be schemes which are flat over $S$. Let $\alpha, \beta \colon X \multimap Y$ be effective multivalued morphisms over $S$ given by morphisms $f \colon X \to \coprod_{m=0}^{\infty} \mathrm{S}^m(Y|S)$ and $g \colon X \to \coprod_{n=0}^{\infty} \mathrm{S}^n(Y|S)$.

Using the morphisms $\sigma^{m,n}$ of Proposition 3.3.4, we define their *sum* $\alpha + \beta \colon X \multimap Y$ as the composition

$$\begin{array}{ccc}
X & \xrightarrow{(f,g)_S} & \coprod_{m=0}^{\infty} \mathrm{S}^m(Y|S) \times_S \coprod_{n=0}^{\infty} \mathrm{S}^n(Y|S) \\
\vdots & & \| \\
\downarrow & & \\
\coprod_{k=0}^{\infty} \mathrm{S}^k(Y|S) & \xleftarrow{\coprod_{m,n=0}^{\infty} \sigma^{m,n}} & \coprod_{m,n=0}^{\infty} \mathrm{S}^m(Y|S) \times_S \mathrm{S}^n(Y|S)
\end{array}$$

**Lemma 3.4.4.** *Let $S$ be a scheme and let $X$, $Y$ be schemes which are flat over $S$.*

*Then addition of multivalued morphisms turns the set of effective multivalued morphisms $X \multimap Y$ into an abelian monoid.*

*Proof.* A neutral element is given by the structure morphism $X \to S \cong \mathrm{S}^0(Y|S)$. Commutativity and associativity follow immediately from the aptly named compatibilities of Proposition 3.3.9.

This is also found in [Ayo14b], Appendix A. □

**Proposition 3.4.5.** *Let $S$ be a scheme. Let $I$ and $J$ be arbitrary sets. Also let schemes $X_i$, $i \in I$, and $Y_j$, $j \in J$, which are flat over $S$ be given.*

*Then there is a natural additive decomposition*

$$\mathrm{Multi}_S^{\mathrm{eff}}\left(\coprod_{i\in I} X_i, \coprod_{j\in J} Y_j\right) \cong \prod_{i\in I}\bigoplus_{j\in J} \mathrm{Multi}_S^{\mathrm{eff}}(X_i, Y_j).$$

*Proof.* This is part of Appendix A of [Ayo14b], but we give a quick reminder of their short proofs:

From the universal property of disjoint unions, i.e. coproducts, we trivially have an additive bijection

$$\mathrm{Multi}_S^{\mathrm{eff}}\left(\coprod_{i\in I} X_i, \coprod_{j\in J} Y_j\right) \cong \prod_{i\in I} \mathrm{Multi}_S^{\mathrm{eff}}\left(X_i, \coprod_{j\in J} Y_j\right).$$

We are thus left to show that additively

$$\mathrm{Multi}_S^{\mathrm{eff}}\left(X_i, \coprod_{j\in J} Y_j\right) = \bigoplus_{j\in J} \mathrm{Multi}_S^{\mathrm{eff}}(X_i, Y_j).$$

By a simple argument, decomposing and recomposing the $X_i$ into connected components, we may assume that the $X_i$ are connected.

We have from Proposition 3.3.5 a natural isomorphism

$$\mathrm{S}^d\left(\coprod_{j\in J} Y_j|S\right) \cong \coprod_{\substack{(d_j)\in(\mathbb{N}_0)^J \\ \sum_{j\in J} d_j = d}} \prod_{\substack{j\in J \\ d_j\neq 0}} \mathrm{S}^{d_j}(Y_j|S).$$

Hence a multivalued morphism $X_i \multimap \coprod_{j\in J} Y_j$ is the same as a morphism

$$X_i \to \coprod_{\substack{(d_j)\in(\mathbb{N}_0)^J \\ \sum_{j\in J} d_j < \infty}} \prod_{\substack{j\in J \\ d_j\neq 0}} \mathrm{S}^{d_j}(Y_j|S).$$

But we assumed $X_i$ to be connected, hence it gets mapped to a single

$$\prod_{\substack{j\in J\\ d_j\neq 0}} \mathrm{S}^{d_j}(Y_j|S).$$

The required bijection then becomes the universal property of fibre products over $S$. The additivity is immediate from the definition of the monoid structure. $\square$

**Remark 3.4.6.** Lemma 3.4.4 states that every effective multivalued morphism $\alpha : X \multimap Y$, where $X$ is connected, has a unique decomposition as a finite sum of homogeneous multivalued morphisms $\alpha_j : X \multimap Y_j$. Here $Y_j$ runs through the different connected components of $Y$.

On arbitrary $X$ this then holds individually for each connected component.

**Definition 3.4.7** (Composition of multivalued morphisms)**.** Let $S$ be a scheme. Let $\alpha\colon X \multimap Y$ and $\beta\colon Y \multimap Z$ be effective multivalued morphisms between flat schemes over $S$.

If $\alpha$ and $\beta$ are homogeneous of respective degrees $m$ and $n$, given by morphisms $f\colon X \to \mathrm{S}^m(Y|S)$ and $g\colon Y \to \mathrm{S}^n(Z|S)$, then we use the morphisms $\tau^{m,n}$ of Proposition 3.3.3 to define their *composition* $\beta\circ\alpha$ as the morphism

$$\tau^{m,n}\circ \mathrm{S}^m(g|S)\circ f\colon X \to \mathrm{S}^m(Y|S) \to \mathrm{S}^m(\mathrm{S}^n(Z|S)|S) \to \mathrm{S}^{mn}(Z|S).$$

It is in particular homogeneous of degree $mn$.

For general $\alpha$ and $\beta$ we extend this additively by Remark 3.4.6.

**Remark 3.4.8.** A direct definition of composition, not using Remark 3.4.6, can be given: Proposition 3.3.5 and the compatibilities compose to a chain

$$\begin{array}{ccccc}
X & \xrightarrow{\quad f \quad} & \coprod_{m=0}^{\infty} \mathrm{S}^m(Y|S) & \xrightarrow{\coprod \mathrm{S}^\bullet(g|S)} & \coprod_{m=0}^{\infty} \mathrm{S}^m\left(\coprod_{n=0}^{\infty} \mathrm{S}^n(Z|S)|S\right) \\
\vdots & & & & \downarrow \cong \\
\coprod_{k=0}^{\infty} \mathrm{S}^k(Z|S) & \xleftarrow{\ \sigma\ } & \coprod_{(d_i)} \prod_{\substack{i\in\mathbb{N}_0\\ d_i\neq 0}} \mathrm{S}^{i\cdot d_i}(Z|S) & \xleftarrow{\ \tau\ } & \coprod_{(d_i)} \prod_{\substack{i\in\mathbb{N}_0\\ d_i\neq 0}} \mathrm{S}^{d_i}\left(\mathrm{S}^i(Z|S)|S\right)
\end{array}$$

where we used $(d_i)$ to denote all sequences $d_0, d_1, d_2, \ldots \in \mathbb{N}_0$ with finite sum.

**Definition 3.4.9** (Tensor product of multivalued morphisms)**.** Let $S$ be a scheme. Let $\alpha_i\colon X_i \multimap Y_i$, $i\in\{1,2\}$, be effective multivalued morphisms between flat schemes over $S$, given by morphisms $f_i\colon X_i \to \coprod_{n_i=0}^{\infty} \mathrm{S}^{n_i}(Y_i|S)$.

Using the morphisms $\rho^{m,n}$ of Proposition 3.3.7, we define their *tensor product* $\alpha_1 \otimes \alpha_2 \colon X_1 \times_S X_2 \multimap Y_1 \times_S Y_2$ as the composition

$$\begin{array}{ccc}
X_1 \times_S X_2 & \xrightarrow{f_1 \times_S f_2} & \coprod_{n_1=0}^{\infty} \mathrm{S}^{n_2}(Y_1|S) \times_S \coprod_{n_2=0}^{\infty} \mathrm{S}^{n_2}(Y_2|S) \\
\Big\downarrow & & \Big\| \\
\coprod_{k=0}^{\infty} \mathrm{S}^k(Y_1 \times_S Y_2|S) & \xleftarrow{\coprod_{n_1,n_2=0}^{\infty} \rho^{n_1,n_2}} & \coprod_{n_1,n_2=0}^{\infty} \mathrm{S}^{n_1}(Y_1|S) \times_S \mathrm{S}^{n_2}(Y_2|S).
\end{array}$$

**Definition 3.4.10.** Let $S$ be a scheme. We define the category $\mathrm{Multi}_S^{\mathrm{eff}}$ of *effective multivalued morphisms* over $S$ as follows:

- its objects are the separated schemes that are flat over $S$,
- its morphisms are the multivalued morphisms $X \multimap Y$ of Definition 3.4.1,
- composition is given by Definition 3.4.7.

This indeed gives a category, and even more:

**Theorem 3.4.11.** *Let $S$ be a scheme.*

*Then composition (Definition 3.4.7), addition (Definition 3.4.3) and tensor product (Definition 3.4.9) turn the category $\mathrm{Multi}_S^{\mathrm{eff}}$ of effective multivalued morphisms over $S$ into a symmetric monoidal category enriched over abelian monoids.*

**Remark 3.4.12.** Note that all three structures, composition, addition and tensor product, preserve homogeneity. Conversely, Remark 3.4.6 shows that it is sufficient to define them in this case. Recalling the definitions, they amount to simple applications of $\sigma^{m,n}$, $\tau^{m,n}$ and $\rho^{m,n}$.

*Proof of Theorem 3.4.11.* By Remark 3.4.12 it suffices to check everything for homogeneous multivalued morphisms. There, the theorem follows readily from the compatibilities of Proposition 3.3.9. $\square$

All but the tensor structure can also already be found in [Ayo14b] as Proposition A.3 and Proposition A.5.

**Remark 3.4.13.** Similar to Remark 4.5.4 we have a natural embedding $\mathrm{Sch}_S \to \mathrm{Multi}_S^{\mathrm{eff}}$:

On objects it is just the identity, and a morphism $f \colon X \to Y$ can via $\mathrm{S}^1(Y|S) \cong (Y|S)^1 \cong Y$ be interpreted as a homogeneous multivalued morphism of degree 1. Under the identification $\mathrm{S}^1(Y|S) \cong Y$ one immediately finds that $\tau^{1,1}$ and $\rho^{1,1}$ are the respective identities. Hence this is indeed a functor that is furthermore compatible with the product structures. It is also clearly an embedding.

Due to the abelian monoid structure the following makes sense:

**Definition 3.4.14.** Let $S$ be a scheme and let $\Lambda$ be a ring.

The category $\mathrm{Multi}_{S,\Lambda} = \mathrm{Multi}(S, \Lambda)$ of $\Lambda$-*multivalued morphisms* has the same objects as $\mathrm{Multi}_S^{\mathrm{eff}}$ and the $\Lambda$-linear extensions

$$\mathrm{Multi}_S^{\mathrm{eff}}(X, Y) \otimes_{\mathbb{N}} \Lambda$$

as morphisms.

If $\Lambda = \mathbb{Z}$ we will call it the category of *multivalued morphisms* and denote it by $\mathrm{Multi}_S$.

We hence get from Theorem 3.4.11, generalizing and extending Appendix A of [Ayo14b]:

**Proposition 3.4.15.** *Let $S$ be a scheme and let $\Lambda$ be a ring.*

*Then $\mathrm{Multi}_{S,\Lambda}$ is a $\Lambda$-linear symmetric tensor category.*

We also can easily extend Proposition 3.4.5 to:

**Proposition 3.4.16.** *Let $S$ be a scheme. Let $I$ and $J$ be arbitrary sets. Also let schemes $X_i$, $i \in I$, and $Y_j$, $j \in J$, which are flat over $S$ be given.*

*Then there is a natural additive decomposition*

$$\mathrm{Multi}_{S,\Lambda}\left(\coprod_{i\in I} X_i, \coprod_{j\in J} Y_j\right) \cong \prod_{i\in I}\bigoplus_{j\in J} \mathrm{Multi}_{S,\Lambda}(X_i, Y_j).$$

**Remark 3.4.17.** Also compare this result to Proposition 1.6.4.

## 3.5 Symmetrization

Assume that $f\colon Y \to X$ is a map of sets where each $x \in X$ has the same number $n$ of preimages. Then this defines a canonical map $\mathrm{sym}(f)\colon X \to Y^n/\mathrm{S}_n$, sending $x \in X$ to the unordered $n$-tuple $\mathrm{sym}(f)(x)$ consisting of its preimages under $f$.

The scheme-theoretic analogue of such a map $f$ would be a finite flat morphism of constant degree $n$. We now show how the above idea works in this setting and generalize it to non-flat morphisms. This is possible, at least over normal schemes, due to generic flatness. It is also helpful that we have, unlike the set-theoretic setting, an intrinsic notion of multiplicity.

**Remark 3.5.1.** The reader not interested in the higher generalities offered by algebraic spaces may always assume to work in the category $\mathrm{Sch}_S^{\mathrm{AF}}$ of schemes AF over a fixed base scheme $S$. By Lemma 3.1.12 (a) every morphism in this category is separated and hence AF by Lemma 3.1.12 (d). By Lemma 3.2.4 it is closed under fibre products and symmetric products.

We opted for higher generality and hence sometimes get algebraic spaces instead of schemes when looking at symmetric products. One should be wary that the theory of cycles in the sense of Section 1.4 has only been defined over noetherian schemes, which forces us to work within these constraints. But this is never a relevant problem.

The following construction is taken from Section 6 of [SV96]. The ring-theoretic version is much older and can for example be found in [Rob63] and [Fer98].

**Lemma 3.5.2.** *Let $X$ be a normal and connected scheme. Let $\pi\colon \Gamma \to X$ be a finite and surjective morphism of degree $n$ between schemes.*

*Then there exists a splitting* $\mathrm{sym}(\pi)\colon X \to \mathrm{S}^n(\Gamma|X)$ *of the morphism* $\mathrm{S}^n(\pi)\colon \mathrm{S}^n(\Gamma|X) \to \mathrm{S}^n(X|X) \cong X$.

*Proof.* As every finite morphism is AF, the symmetric product $\mathrm{S}^n(\Gamma|X)$ exists as a scheme by Lemma 3.2.4. Let $\xi$ be the generic point of the irreducible scheme $X$ and let $K(X) = \mathcal{O}_{X,\xi}$ be the corresponding function field. Then the fibre $\pi_\xi\colon \Gamma_\xi \to \xi$ is by assumption a morphism of degree $n$. By Definition 2.3.4 this induces a morphism $\vartheta_\xi\colon\ \mathrm{S}_n(\mathcal{O}_{\Gamma_\xi}|K(X)) \to K(X)$.

Now let $U \hookrightarrow X$ be an affine open subset of $X$ and set $W = \pi^{-1}(U)$. Then $W$ is affine because $\pi$ is finite. We consider the diagram

$$\begin{array}{ccccc}
(\mathcal{O}_W|\mathcal{O}_U)^{\otimes n} & \hookleftarrow & \mathrm{S}_n(\mathcal{O}_W|\mathcal{O}_U) & \overset{\vartheta_U}{\dashrightarrow} & \mathcal{O}_U \\
\downarrow & & \downarrow & & \downarrow \\
(\mathcal{O}_{\Gamma_\xi}|K(X))^{\otimes n} & \hookleftarrow & \mathrm{S}_n(\mathcal{O}_{\Gamma_\xi}|K(X)) & \xrightarrow{\vartheta} & K(X),
\end{array}$$

where we want to construct the dotted arrow. The left side commutes by definition. Now $W \to U$ is finite, hence $\mathcal{O}_W$ is integral over $\mathcal{O}_U$ and thus the same is true for $(\mathcal{O}_W|\mathcal{O}_U)^{\otimes n}$ and its subring $\mathrm{S}_n(\mathcal{O}_W|\mathcal{O}_U)$. Therefore the image of $\mathrm{S}_n(\mathcal{O}_W|\mathcal{O}_U)$ under $\vartheta$ is integral over $\mathcal{O}_U$ and contained in $K(X)$, hence it is contained in $\mathcal{O}_U$ because $X$ is normal. This gives the dashed arrow, amounting to a morphism

$$\theta_U\colon U \to \mathrm{S}^n(W|U) \hookrightarrow \mathrm{S}^n(\Gamma|X)$$

by Proposition 3.1.4 (e).

If $V \subseteq U$ is open and affine, then $\theta_V = (\theta_U)_{|V}$ by construction. Hence the $\theta_U$ glue to a morphism $\mathrm{sym}(\pi)\colon X \to \mathrm{S}^n(\Gamma|X)$. $\square$

**Definition 3.5.3** (Symmetrizations)**.** Let $X$ be a noetherian, connected and normal scheme.

Let $\pi\colon \Gamma \to X$ be a finite and surjective morphism of degree $n$. Then the construction of Lemma 3.5.2 gives us the *symmetrization*

$$\mathrm{sym}(\pi) = \mathrm{sym}(\Gamma|X)\colon X \to \mathrm{S}^n(\Gamma|X).$$

Let $Y \to X$ be a morphism of finite type. Let $\alpha \in c^{\mathrm{nai}}(Y|X)$ be a basic cycle, i.e. by Definition 1.4.3 the generic point of a closed integral subscheme $i\colon \operatorname{supp}(\alpha) \hookrightarrow Y$ finite and surjective over $X$. It induces a morphism

$$\operatorname{sym}(\alpha) := \mathrm{S}^n(i|X) \circ \operatorname{sym}(\pi)\colon X \to \mathrm{S}^n(\operatorname{supp}(\alpha)|X) \to \mathrm{S}^n(Y|X)$$

which we again call *symmetrization.*

Let now $\alpha = \sum_{i=1}^r \alpha_i \in c^{\mathrm{eff},\mathrm{nai}}(Y|X)$ be an effective naive cycle of degree $n = \sum_{i=1}^r n_i$. Here, the $\alpha_i$ are basic cycles of degree $n_i \in \mathbb{N}$. We also assume $Y \to X$ to be flat. Then, using the addition morphism $\sigma^{(n_i)} = \sigma^{n_1,\dots,n_r}$ of Proposition 3.3.4, we define the *symmetrization* $\operatorname{sym}(\alpha)\colon X \to \mathrm{S}^n(Y|X)$ as the composition

$$X \xrightarrow{\prod_{i=1}^r \operatorname{sym}(\alpha_i)} \prod_{i=1}^r \mathrm{S}^{n_i}(Y|X) \xrightarrow{\sigma^{(n_i)}} \mathrm{S}^n(Y|X).$$

**Definition 3.5.4** (Homogeneous symmetrization)**.** Let $Y \to X$ be a flat morphism of finite type between noetherian schemes and assume that $X$ is normal. Let $\alpha \in c^{\mathrm{eff},\mathrm{nai}}(Y|X)$ be of constant degree $n$. Let $X = \coprod_{i=1}^r X_i$ be the decomposition into connected components and let $\alpha_i = (X_i \hookrightarrow X)^{\circledast}\alpha$ be the restrictions to the connected components.

Then we extend Definition 3.5.3 by setting

$$\operatorname{sym}(\alpha) = \coprod_{i=1}^r \operatorname{sym}(\alpha_i)\colon X = \coprod_{i=1}^r X_i \to \mathrm{S}^n(Y|X).$$

**Remark 3.5.5.** Definition 3.5.3 should be understood as an additive map

$$\operatorname{sym}\colon c^{\mathrm{eff},\mathrm{nai}}(Y|X) \to \mathrm{Multi}_X(X,Y) = \mathrm{Sch}_X(X, \coprod_{d=0}^{\infty} \mathrm{S}^d(Y|X)),$$

the additivity being immediate from the definition.

The compatibility of addition with composition and the tensor structure on both sides allows to check some properties on a set of additive generators. We will often use this to reduce from arbitrary effective cycles to basic ones. This, however, requires that our property of interest is defined on all naive cycles $c^{\mathrm{eff},\mathrm{nai}}(Y|X)$, not just the subset of genuine relative cycles $c^{\mathrm{eff}}(Y|X)$.

**Remark 3.5.6.** Let $Y \to X$ be a flat morphism of finite type between noetherian schemes and let $\alpha \in c^{\mathrm{nai},\mathrm{eff}}(Y|X)$ be of degree $n$. Assume that $X$ is normal and connected with generic point $g\colon \xi \to X$.

Recall the naive pullback of Definition 1.4.11. From the construction in Lemma 3.5.2 and Definition 3.5.3 we get the equality

$$\operatorname{sym}(\alpha) \circ g = \mathrm{S}^n(g \times_X Y|g) \circ \operatorname{sym}(g^*\alpha)$$

of morphisms $\xi \to \mathrm{S}^n(Y|X)$. This is by Proposition 1.5.4 a special case of Proposition 3.6.6 below. Note also that all our morphisms are assumed to be separated and that $X$ is clearly reduced, thus by Proposition 9.7 ii) of [GW10] a morphism $X \to Z$ is made unique by fixing a morphism $\xi \to Z$.

Hence when showing properties of sym(−), we may by its construction via the generic point $\xi \to X$ often assume that the base $X$ is the spectrum of a field.

The two main advantages of this reduction are:

- $X$ and hence any $\Gamma$ finite over it become affine, allowing us to use the results of Chapter 2,

- everything becomes flat over $X$, allowing the usage of Proposition 3.3.3 and Proposition 3.3.4.

Both aspects will be exploited several times. A prime example of this technique will be the proof of Lemma 3.6.1.

## 3.6 Structural compatibilities

We show the compatibility of symmetrization with the fundamental operations on cycles, which are: decomposition, pushforward, pullback, the correspondence map and the exterior product. The most difficult one is the pullback, whose definition we have to unravel. As soon as it is out of the way one could, as an alternative to the approach presented in this section, proceed similar to Section 7 of [Ryd08b]. This requires the formalism of op. cit., which we will only meet in Section 3.9. Given this formalism, the results of this section can also be extracted from [Ryd08d], of which we only learned after significant parts of this section had been completed. It is noteworthy that the approach in op. cit. is opposite to ours, mapping multivalued morphisms to cycles, the latter not necessarily being relative ones.

**Lemma 3.6.1.** *Let*

$$\begin{array}{ccc} X & \longleftarrow & \Gamma \\ \big\downarrow f & & \big\downarrow g \\ X' & \longleftarrow & \Gamma' \end{array}$$

*be a cartesian square of noetherian schemes. Assume that $X$ and $X'$ are connected and normal and that the lower horizontal morphism is flat, surjective and finite of degree $n$.*

*Then the upper horizontal morphism is also flat, surjective and finite of*

*the same degree $n$ and the resulting diagram*

$$\begin{array}{ccc} X & \xrightarrow{\operatorname{sym}(\Gamma|X)} & \mathrm{S}^n(\Gamma|X) \\ \downarrow{\scriptstyle f} & & \downarrow{\scriptstyle \mathrm{S}^n(g|f)} \\ X' & \xrightarrow{\operatorname{sym}(\Gamma'|X')} & \mathrm{S}^n(\Gamma'|X') \end{array}$$

*commutes.*

*Proof.* Flatness, surjectivity and finiteness of the upper morphism are satisfied as it is a pullback of such a morphism.

By Remark 3.5.6 we may assume that $X = \{\xi\}$ consists only of its generic point and thus is the spectrum of a field $\lambda$. Similarly we then may replace $X'$ with the spectrum of the local ring $A = \mathcal{O}_{X',f(\xi)}$. As the horizontal morphisms are finite we get that $\Gamma$ and $\Gamma'$ are affine. This reduces the lemma to the following:

Let $A$ be a local integral domain with residue field $\kappa$ and field of quotients $K$, let $B$ be a finite, flat and therefore free $A$-algebra of rank $n$ and let $\lambda$ be a field extension of $\kappa$. Then we have to show that the horizontal morphisms in the diagram

$$\begin{array}{ccc} \lambda & \hookrightarrow & B \otimes_A \lambda \\ \uparrow & & \uparrow \\ A & \hookrightarrow & B \\ \downarrow & & \downarrow \\ K & \hookrightarrow & B \otimes_A K \end{array}$$

are free of rank $n$ and induce a commutative diagram

$$\begin{array}{ccc} \lambda & \xleftarrow{\vartheta_{B\otimes_A \lambda|\lambda}} & \mathrm{S}_n(B \otimes_A \lambda|\lambda) \\ \uparrow & & \uparrow \\ A & \xleftarrow{\vartheta_{B|A}} & \mathrm{S}_n(B|A) \\ \downarrow & & \downarrow \\ K & \xleftarrow{\vartheta_{B\otimes_A A|K}} & \mathrm{S}_n(B \otimes_A K|K). \end{array}$$

This is just a twofold application of Lemma 2.3.6. □

**Lemma 3.6.2.** *Let*

$$\Theta \xrightarrow{\rho} \Gamma \xrightarrow{\pi} X$$

*be finite surjective morphisms between noetherian schemes. Assume that $X$ and $\Gamma$ are both normal and connected. Let $m$ and $n$ be the degrees of $\pi$ and $\rho$, respectively. Lastly assume that $\pi$ and $\pi \circ \rho$ are flat.*

*Then $\pi \circ \rho$ is finite and surjective of degree $mn$. Combined with the morphism $\tau^{m,n}$ of Proposition 3.3.3 we have:*

*(a)* $\mathrm{sym}(\pi \circ \rho) = \tau^{m,n} \circ \mathrm{S}^m(\mathrm{sym}(\rho)|X) \circ \mathrm{sym}(\pi)$, *i.e. the diagram*

$$\begin{array}{ccc} X & \xrightarrow{\mathrm{sym}(\pi)} & \mathrm{S}^m(\Gamma|X) \\ {\scriptstyle \mathrm{sym}(\pi\circ\rho)}\downarrow & & \downarrow{\scriptstyle \mathrm{S}^m(\mathrm{sym}(\rho)|X)} \\ \mathrm{S}^{mn}(\Theta|X) & \xleftarrow{\tau^{m,n}} & \mathrm{S}^m(\mathrm{S}^n(\Theta|\Gamma)|X) \end{array}$$

*commutes.*

*(b)* $\sigma^{\overbrace{m,m,\ldots,m}^{n \text{ times}}} \circ (\mathrm{sym}(\pi)|X)^n = \mathrm{S}^{mn}(\rho|X) \circ \mathrm{sym}(\pi \circ \rho)$, *i.e. the diagram*

$$\begin{array}{ccc} X & \xrightarrow{(\mathrm{sym}(\pi)|X)^n} & (\mathrm{S}^m(\Gamma|X)|X)^n \\ {\scriptstyle \mathrm{sym}(\pi\circ\rho)}\downarrow & & \downarrow{\scriptstyle \sigma^{\overbrace{m,m,\ldots,m}^{n \text{ times}}}} \\ \mathrm{S}^{mn}(\Theta|X) & \xrightarrow{\mathrm{S}^{mn}(\rho|X)} & \mathrm{S}^{mn}(\Gamma|X) \end{array}$$

*commutes.*

*Proof.* The finiteness and surjectivity of $\pi \circ \rho$ is trivial. It suffices by Remark 3.5.6 to check the equalities over the generic point of $X$, i.e. we may assume that $X = \mathrm{Spec}(k)$ is the spectrum of a field and that all three schemes are therefore affine and flat over $X$.

(a) By definition of $\mathrm{sym}(-)$ it suffices to show that the maps

$$\vartheta_{\mathcal{O}_\Theta|k} = \vartheta_{\mathcal{O}_\Gamma|k} \circ \mathrm{S}_m(\vartheta_{\mathcal{O}_\Theta|\mathcal{O}_\Gamma}|k) \circ \tau_{m,n}$$

from $\mathrm{S}_{mn}(\Theta|k)$ to $k$ agree. But this is Proposition 2.4.4.

(b) Note that this can restated as the equality

$$n \cdot \mathrm{sym}(\pi) = \rho \circ \mathrm{sym}(\pi \circ \rho)$$

of multivalued morphisms over $X$. We use the shorthand

$$\mathrm{sym}'(\rho) := \mathrm{S}^n(\Theta|\pi) \circ \mathrm{sym}(\rho) \colon \Gamma \to \mathrm{S}^n(\Gamma|X)$$

and observe that it defines a multivalued endomorphism of $\Gamma$ over $X$. Furthermore, part (a) states that

$$\mathrm{sym}(\pi \circ \rho) = \mathrm{sym}'(\rho) \circ \mathrm{sym}(\pi)$$

as multivalued morphisms.

If $X = \Gamma$ and $\pi = \mathrm{id}_X$ then the result is trivially true. This translates into $n \cdot \mathrm{id}_\Gamma = \rho \circ \mathrm{sym}'(\rho)$ as multivalued morphisms over $X$. Hence we conclude from the linearity of composition of multivalued morphisms that indeed

$$n \cdot \mathrm{sym}(\pi) = \rho \circ \mathrm{sym}'(\rho) \circ \mathrm{sym}(\pi) = \rho \circ \mathrm{sym}(\pi \circ \rho).$$

□

**Corollary 3.6.3.** *Let*

$$\Theta \xrightarrow{\ \rho\ } \Gamma \xrightarrow{\ \pi\ } X$$

*be morphisms of finite type between noetherian schemes. Assume that $X$ and $\Gamma$ are both normal and connected. Also assume that $\pi$ is finite and surjective of degree $m$. Lastly assume that $\pi$ and $\pi \circ \rho$ are flat.*

*Let $\beta \in c^{\mathrm{nai,eff}}(\Theta|\Gamma)$ be of degree $n$ and let $\pi_\# \beta \in c^{\mathrm{nai,eff}}(\Theta|X)$ be the same naive cycle interpreted as one over $X$ as in Definition 1.4.19.*

*Then the diagram*

$$\begin{array}{ccc} X & \xrightarrow{\ \mathrm{sym}(\pi)\ } & \mathrm{S}^m(\Gamma|X) \\ \Big\downarrow{\scriptstyle \mathrm{sym}(\pi_\#\beta)} & & \Big\downarrow{\scriptstyle \mathrm{S}^m(\mathrm{sym}(\beta))} \\ \mathrm{S}^{mn}(\Theta|X) & \xleftarrow{\ \tau^{m,n}\ } & \mathrm{S}^m(\mathrm{S}^n(\Theta|\Gamma)|X) \end{array}$$

*commutes.*

*Proof.* Due to Remark 3.5.5 we may assume that $\beta$ is basic and by Remark 3.5.6 we may assume that $X$ is the spectrum of a field. Hence the tower

$$\mathrm{supp}(\beta) \longrightarrow \Gamma \xrightarrow{\ \pi\ } X$$

satisfies the conditions of Lemma 3.6.2, which immediately implies the result. □

**Proposition 3.6.4** (Symmetrization commutes with decomposition into cycles)**.** *Let $Y \to X$ be a flat morphism of finite type between noetherian schemes where $X$ is connected and normal. Let $\Gamma$ be a closed subscheme of $Y$ such that the resulting morphism $\pi\colon \Gamma \to X$ is finite and surjective of degree $n$. Also recall the decomposition $\mathrm{cycl}_{Y|X}(\Gamma) \in c^{\mathrm{nai,eff}}(Y|X)$ of $\Gamma$ into an effective naive cycle from Definition 1.4.7.*

*Then*

$$\mathrm{sym}(\pi) = \mathrm{sym}(\mathrm{cycl}_{Y|X}(\Gamma))$$

*as morphisms $X \to \mathrm{S}^n(Y|X)$.*

*Proof.* The statement holds by definition if $\Gamma$ is integral. We show the general case in three steps:

(a) By Remark 3.5.6 it suffices to check this property over the generic point of $X$, i.e. we may assume that $X = \mathrm{Spec}(k)$ is the spectrum of a field $k$ and thus that $\Gamma$ is affine. As $\Gamma$ is now flat over $X$, we may by the naturality of $\sigma^-$ (cf. Proposition 3.3.4) in Definition 3.5.3 assume that $Y = \Gamma$.

(b) If $\Gamma$ is irreducible, we let $d$ be the degree of the $k$-algebra $\mathcal{O}_\Gamma$. Now $\mathrm{cycl}_{Y|X}(\Gamma) = m \cdot \Gamma^{\mathrm{red}}$ by definition, where $m$ is the length of $\mathcal{O}_\Gamma$. Thus the theorem follows from Proposition 2.4.3.

(c) If $\Gamma$ is arbitrary, we write $\Gamma = \bigcup_{i=1}^r \Gamma_i$ as union of its irreducible components with individual degrees $n_i$. Note that $\Gamma$ is finite over a field, hence this is actually the decomposition into its points. Then

$$\mathrm{cycl}_{Y|X}(\Gamma) = \sum_{i=1}^{r} \mathrm{cycl}_{Y|X}(\Gamma_i),$$

so by the previous case we are left to show that

$$\mathrm{sym}\left(\mathrm{cycl}_{Y|X}(\Gamma)\right) = \sigma^{(n_1,n_2,\dots,n_r)} \circ \prod_{i=1}^{r} \mathrm{sym}\left(\mathrm{cycl}_{Y|X}(\Gamma_i)\right).$$

Switching to global sections, this amounts to

$$\vartheta_{\mathcal{O}_\Gamma|k} = \bigotimes_{i=1}^{r} \vartheta_{\mathcal{O}_{\Gamma_i}|k} \circ \sigma_{(n_1,n_2,\dots,n_r)},$$

which follows from Lemma 2.4.1 by induction on $r$. □

**Proposition 3.6.5** (Symmetrization commutes with pushforwards)**.** *Let*

$$Z \xrightarrow{\ g\ } Y \xrightarrow{\ f\ } X$$

*be morphisms of finite type between noetherian schemes, where $X$ is normal. Assume that $f$ and $f \circ g$ are flat. Let $\alpha \in c^{\mathrm{nai,eff}}(Z|X)$ be an effective naive cycle of constant degree $n$.*

*Then $g_*\alpha$ has the same constant degree $n$ and we get an equality*

$$\mathrm{sym}(g_*\alpha) = \mathrm{S}^n(g|X) \circ \mathrm{sym}(\alpha)$$

*of morphisms $X \to \mathrm{S}^n(Y|X)$.*

*Proof.* The equality $\deg(g_*\alpha) = \deg(\alpha)$, including the degree being constant, is Lemma 1.8.5. A pullback to a connected component reduces us to the case where $X$ is connected.

Both sides of the equation are by definition linear in $\alpha$ as described in Remark 3.5.5. Hence we may assume that $\alpha$ is basic.

By Lemma 1.4 of [MVW06] we have an equality $g(\operatorname{supp}(\alpha)) = \operatorname{supp}(g_*\alpha)$ of sets. Equipping the supports with the reduced induced structure we let $m$ be the degree of the induced finite and surjective morphism $\operatorname{supp}(\alpha) \to \operatorname{supp}(g_*\alpha)$. Hence by definition of the pushforward we have $g_*\alpha = m \cdot g(\alpha)$.

We let $d$ be the degree of $\operatorname{supp}(g_*\alpha)$ over $X$, thus we have $n = dm$. Then the commutative diagram

$$\begin{array}{ccccc}
 & & \overset{g_{|\operatorname{supp}(\alpha)}}{\curvearrowright} & & \\
\operatorname{supp}(\alpha) & \xrightarrow[(f\circ g)_{|\operatorname{supp}(\alpha)}]{} & X & \xleftarrow[f_{|\operatorname{supp}(g_*\alpha)}]{} & \operatorname{supp}(g_*\alpha) \\
\downarrow & & & & \downarrow \\
Z & & \xrightarrow[g]{} & & Y
\end{array}$$

induces a diagram

$$\begin{array}{ccccc}
 & & \overset{\mathrm{S}^n(g_{|\operatorname{supp}(\alpha)})}{\curvearrowright} & & \\
\mathrm{S}^n(\operatorname{supp}(\alpha)) & \xleftarrow{\operatorname{sym}((f\circ g)_{|\operatorname{supp}(\alpha)})} & X & \xrightarrow{\operatorname{sym}(d\cdot f_{|\operatorname{supp}(g_*\alpha)})} & \mathrm{S}^n(\operatorname{supp}(g_*\alpha)) \\
\downarrow & \swarrow{\scriptstyle \operatorname{sym}(\alpha)} & & \searrow{\scriptstyle \operatorname{sym}(g_*(\alpha))} & \downarrow \\
\mathrm{S}^n(Z) & & \xrightarrow[\mathrm{S}^n(g)]{} & & \mathrm{S}^n(Y).
\end{array}$$

We used the suggestive notation $\operatorname{sym}(d \cdot f_{|\operatorname{supp}(g_*\alpha)})$ to denote $\sigma^{\overbrace{1,1,\dots,1}^{d \text{ times}}} \circ (\operatorname{sym}(f_{|\operatorname{supp}(g_*\alpha)})|X)^d$.

The upper area commutes by Lemma 3.6.2, part b), and the exterior commutes by functoriality of $\mathrm{S}^n$. The left and right triangles are simply the definition of $\operatorname{sym}(-)$. Thus the bottom triangle commutes, proving the proposition. $\square$

**Proposition 3.6.6** (Symmetrization commutes with pullbacks)**.** *Let*

$$\begin{array}{ccc} X' & \longleftarrow & Y' \\ \downarrow{\scriptstyle f} & & \downarrow{\scriptstyle f'} \\ X & \longleftarrow & Y \end{array}$$

*be a cartesian square of noetherian schemes. Assume that the horizontal morphisms are flat and of finite type. Also assume that $X$ and $X'$ are normal. Let $\alpha \in c^{\mathrm{eff}}(Y|X)$ be effective of constant degree $n$.*

*Then $f^{\circledast}\alpha \in c^{\mathrm{eff}}(Y'|X')$ has the same constant degree $n$ and we have*

$$\mathrm{S}^n(f'|f) \circ \mathrm{sym}(f^{\circledast}\alpha) = \mathrm{sym}(\alpha) \circ f,$$

*i.e. the diagram*

$$\begin{array}{ccc} X' & \xrightarrow{\mathrm{sym}(f^{\circledast}\alpha)} & \mathrm{S}^n(Y'|X') \\ \downarrow{\scriptstyle f} & & \downarrow{\scriptstyle \mathrm{S}^n(f'|f)} \\ X & \xrightarrow{\mathrm{sym}(\alpha)} & \mathrm{S}^n(Y|X) \end{array}$$

*commutes.*

*Proof.* The statement about the degree is Lemma 1.8.7. We may also assume that $X$ and $X'$ are connected.

We write $\alpha = \sum_{i=1}^{r} \alpha_i$ as a formal sum of points. The points give us integral subschemes $\gamma_i\colon \mathrm{supp}(\alpha_i) \hookrightarrow Y$ which are finite and surjective over $X$.

We assume for a moment that all $\mathrm{supp}(\alpha_i)$ are flat over $X$. Then by Proposition 1.5.6 each $\alpha_i$ is a relative cycle over $X$ and its pullback $f^{\circledast}\alpha_i$ is given by $\mathrm{cycl}_{X'}(\mathrm{supp}(\alpha_i) \times_X X')$. By linearity (cf. Remark 3.5.5) it suffices to consider each of the individual $\alpha_i$.

Applying Proposition 3.6.4 we have to show that

$$\mathrm{S}^{\deg(\alpha_i)}(f'|p) \circ \mathrm{sym}(\gamma_i \times_Y Y') = \mathrm{sym}(\alpha) \circ f,$$

which is Lemma 3.6.1. This finishes the flat case.

For the general case we closely follow the Definition 1.4.17 of $f^{\circledast}\alpha$:

By Remark 3.5.6 and linearity (again cf. Remark 3.5.5) we may, due to the third step of the construction of $f^{\circledast}\alpha$, replace $X'$ by its generic point $\xi'$. Combined with Lemma 3.6.1 we may, according to the second step in the definition of $f^{\circledast}\alpha$, even assume that $\xi'$ is a point of $X$ via the morphism $f$.

Proposition 1.4.14 gives us a finite field extension $\lambda\colon u := \mathrm{Spec}(L) \to \xi'$ of the residue field $\kappa(\xi')$ such that the resulting morphism $u \to \xi' \to X$ has a good factorization

$$p \circ i\colon \mathrm{Spec}(L) \to \widetilde{X} \to X$$

with respect to $\alpha$ (cf. Definition 1.4.13). This includes that the individual supports $\mathrm{supp}(\widetilde{X} \times_X \alpha_i)$ are flat over $\widetilde{X}$.

We may also assume that $\widetilde{X}$ is normal by replacing it with its normalization.

By Definition 1.4.15 of relative cycles we find a naive cycle

$$\alpha_{|\xi'} \in c^{\mathrm{eff,nai}}(Y \times_X \xi' | \xi')$$

such that

$$i^* p^* \alpha = \lambda^* \alpha_{|\xi'}.$$

Then the pullback $p^{\circledast}\alpha$ was defined to be this cycle.

We look at the following diagram:

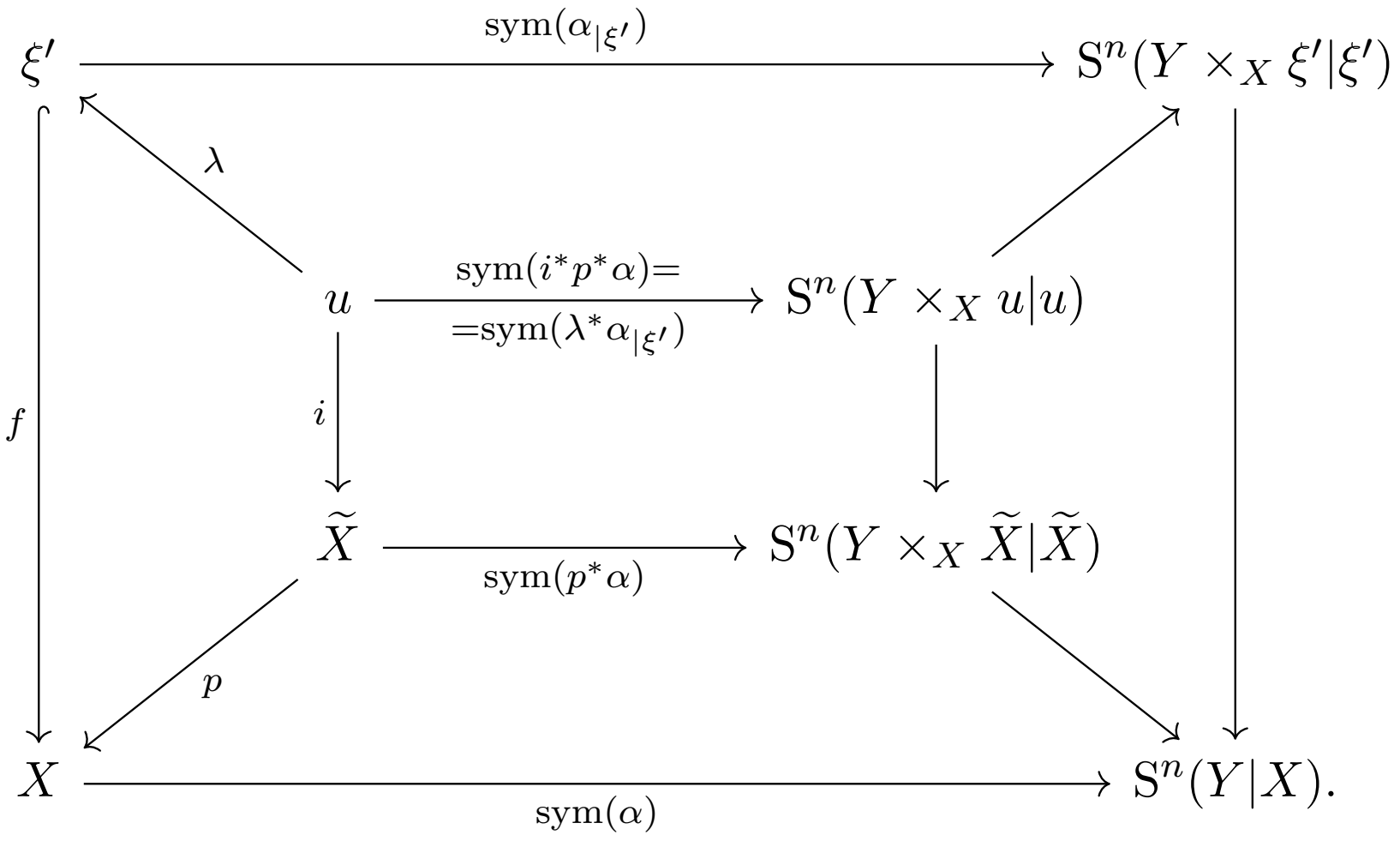

The inner square and the top square commute by the previously considered flat case. The left square commutes by construction and the right one by functoriality of $\mathrm{S}^n(-|-)$.

The bottom square commutes on an open dense subscheme where the abstract blow-up $p\colon \widetilde{X} \to X$ is an isomorphism. Hence by $\widetilde{X}$ being reduced and $\mathrm{S}^n(Y|X)$ being separated over $X$ it commutes at all of $\widetilde{X}$ by Proposition 9.7 ii) of [GW10].

We note that $u \to \xi'$ is an epimorphism in the category of schemes, e.g. by being surjective and $\xi'$ being reduced. We then conclude by a simple diagram chase that the outside of the diagram commutes, finishing the proof. □

For an alternative proof see Proposition 10.8 of [Ryd08d].

**Corollary 3.6.7.** *Under the assumptions of Proposition 3.6.6 we get a*

*commutative triangle*

$$\begin{array}{ccc} X' & \xrightarrow{\operatorname{sym}(p^{\circledast}\alpha)} & \mathrm{S}^n(Y'|X') \\ & {\scriptstyle X'\times_X \operatorname{sym}(\alpha)} \searrow \qquad \swarrow {\scriptstyle (\mathrm{S}^n(p'|p),p)_X} & \\ & \mathrm{S}^n(Y|X) \times_X X' & \end{array}$$

**Proposition 3.6.8** (Symmetrization commutes with Cor)**.** *Let*

$$Z \xrightarrow{\ g\ } Y \xrightarrow{\ f\ } X$$

*be flat morphisms of finite type between noetherian schemes, where $X$ and $Y$ are normal. Let $\alpha \in c^{\mathrm{eff}}(Y|X)$ and $\beta \in c^{\mathrm{eff}}(Y|X)$ be effective relative cycles of constant degrees $m$ and $n$, respectively.*

*Then $\operatorname{Cor}(\beta, \alpha) \in c^{\mathrm{eff}}(Z|X)$ has constant degree $mn$ and we have the equality*

$$\operatorname{sym}(\operatorname{Cor}(\beta, \alpha)) = \tau^{m,n} \circ \mathrm{S}^m(\operatorname{sym}(\beta)) \circ \operatorname{sym}(\alpha),$$

*i.e. the diagram*

$$\begin{array}{ccc} X & \xrightarrow{\operatorname{sym}(\alpha)} & \mathrm{S}^m(Y) \\ {\scriptstyle \operatorname{sym}(\operatorname{Cor}(\beta,\alpha))}\downarrow & & \downarrow{\scriptstyle \mathrm{S}^m(\operatorname{sym}(\beta))} \\ \mathrm{S}^{mn}(Z) & \xleftarrow[\tau^{m,n}]{} & \mathrm{S}^m(\mathrm{S}^n(Z)) \end{array}$$

*commutes.*

*Proof.* By Lemma 1.8.9 we indeed have constant degree $\deg(\operatorname{Cor}(\beta, \alpha)) = \deg(\beta)\deg(\alpha) = mn$. We may assume that $X$ is connected. Even more, by Remark 3.5.6 and Proposition 1.5.9 we may pullback to the generic point of $X$, i.e. assume that $X$ is the spectrum of a field. In particular, every naive cycle is already a relative cycle.

Both sides are linear in the sense of Remark 3.5.5, hence we may assume that $\alpha$ is basic. Then $\Gamma := \operatorname{supp}(\alpha)$ is an integral scheme and finite over $X$, hence itself the spectrum of a field. We can replace $Y$ by the connected component which contains $\Gamma$. By definition we have

$$\operatorname{Cor}(\beta, \alpha) = (i \times_Y Z)_*(f \circ i)_\# i^{\circledast}\beta,$$

where $i \colon \Gamma \hookrightarrow Y$ is the inclusion and

$$(f \circ i)_\# \colon c^{\mathrm{nai,eff}}(\operatorname{supp}(\alpha) \times_Y Z|X) \to c^{\mathrm{nai,eff}}(\operatorname{supp}(\alpha)|X)$$

is the interpretation of cycles over $\Gamma$ as cycles over $X$ (cf. Definition 1.4.19).

Then $\Gamma$ is an integral scheme that is finite over $X$, hence is itself the spectrum of a field. In particular we may assume that $\Gamma$ is normal and connected. We look at the following diagram:

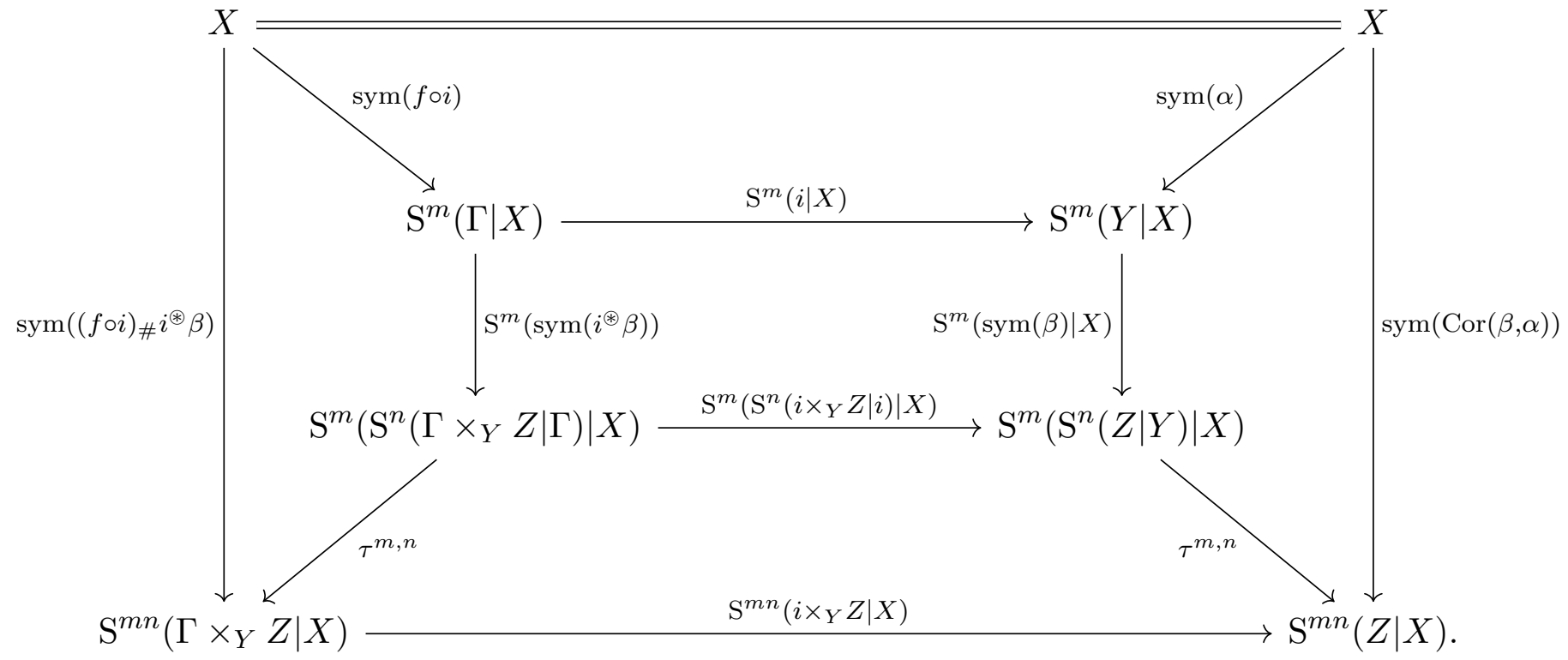

The proposition is equivalent to the commutativity of the right square. Hence by an easy diagram chase it suffices to show the commutativity of the other squares:

- top square: Definition 3.5.3 of $\mathrm{sym}(\alpha)$,
- bottom square: naturality of $\tau^{m,n}$, i.e. Proposition 3.3.3,
- outer square: compatibility of $\mathrm{sym}(-)$ with the pushforward $(i \times_Y Z)_*$, i.e. Proposition 3.6.5,
- inner square: apply $\mathrm{S}^m(-|X)$ to the compatibility of $\mathrm{sym}(-)$ with the pullback $i^\circledast$, i.e. Proposition 3.6.6,
- left square: this is Corollary 3.6.3. □

Also compare this to Proposition 7.8 of [Ryd08b].

**Proposition 3.6.9** (Symmetrization is compatible with the exterior product)**.** *Let $S$ be a noetherian scheme and let $Z_1, Z_2$ be flat schemes of finite type over $S$. Let $\alpha_1 \in c^{\mathrm{eff}}(Z_1|S)$ and $\alpha_2 \in c^{\mathrm{eff}}(Z_2|S)$ be effective relative cycles of constant degrees $n_1$ and $n_2$, respectively.*

*Then $\alpha_1 \otimes_S \alpha_2 \in c^{\mathrm{eff}}(Z_1 \times_S Z_2|S)$ has constant degree $n_1 n_2$ and there is an equality*

$$\mathrm{sym}(\alpha_1 \otimes \alpha_2) = \rho^{n_1,n_2} \circ (\mathrm{sym}(\alpha_1) \times_S \mathrm{sym}(\alpha_2))$$

*of morphisms $S \to \mathrm{S}^{n_1 n_2}(Z_1 Z_2)$.*

*Proof.* The statement about the degree is Lemma 1.8.10.

We let $z_1 \colon Z_1 \to S$ be the structure morphism and set $a = \operatorname{sym}(\alpha_2)$. We get the diagram

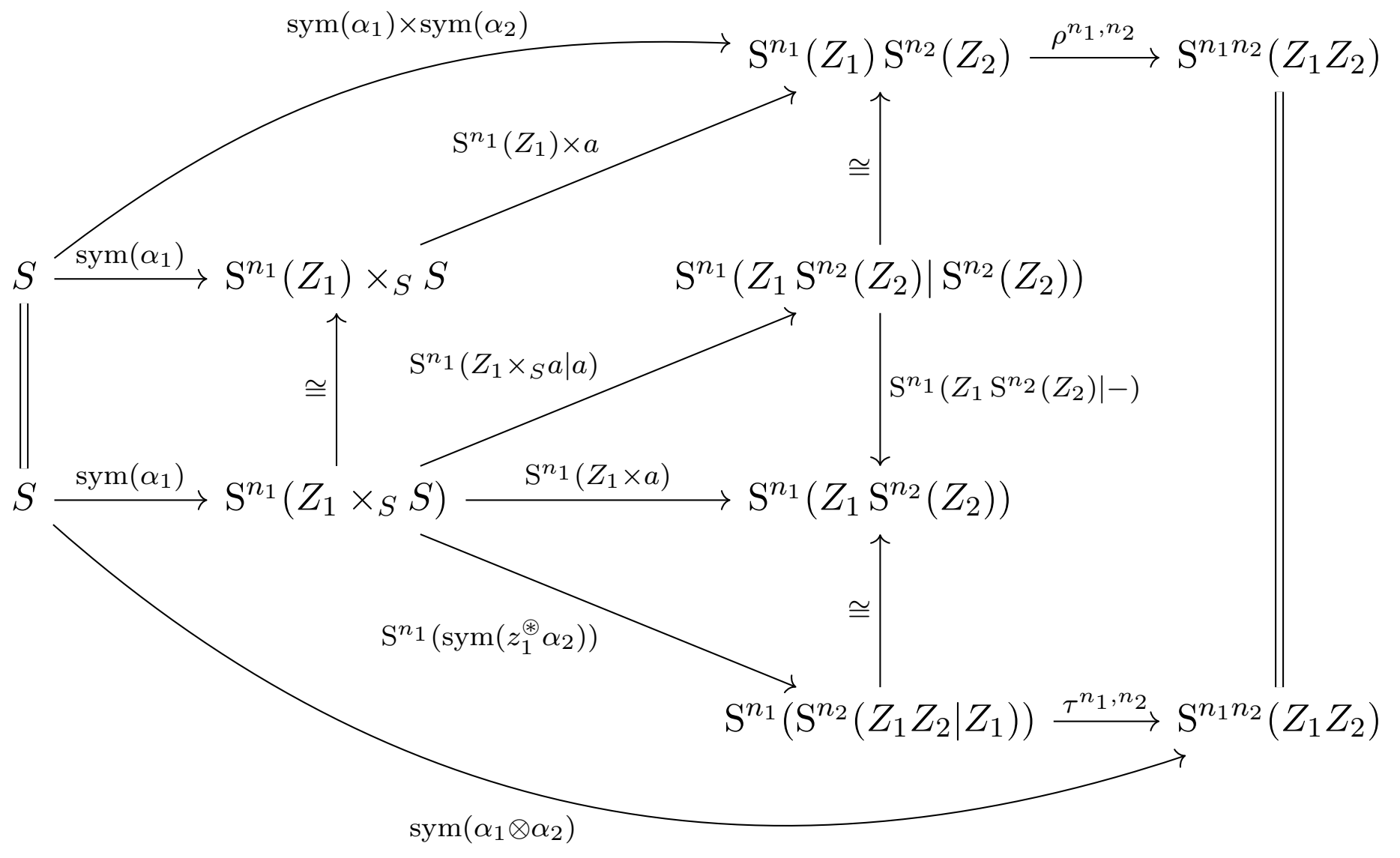

where we suppressed the base $S$ whenever possible and used our abbreviations for products. The isomorphisms are that of Lemma 3.3.2, where the required flatnesses are either given or follows from Lemma 3.2.3. We have to show that the outside commutes, i.e. after a simple diagram chase want to show that all the areas on the inside commute. The reasons are as follows:

- leftmost square: trivial,
- top left area: trivial,
- bottom area: recall that

  $$\alpha_1 \otimes \alpha_2 = \operatorname{Cor}(z_1^{\circledast}\alpha_2, \alpha_1)$$

  by Definition 1.7.1, hence this is just an instance of the compatibility of $\operatorname{sym}(-)$ and $\operatorname{Cor}(-,-)$, i.e. Proposition 3.6.8,
- inner parallelogram: naturality of the isomorphism of Lemma 3.3.2,
- rectangle on the right: this is the alternative description of $\rho^{n_1,n_2}$ of Lemma 3.3.8,
- lower inner triangle: functoriality of $\mathrm{S}^{n_1}$ applied to the triangle of Corollary 3.6.7,
- upper inner triangle: functoriality of $\mathrm{S}^{n_1}(-|-)$. □

## 3.7 Translation into multivalued morphisms

**Definition 3.7.1** (Symmetrization of a finite correspondence)**.** Let $S$ be a noetherian scheme. Let $X$ and $Y$ be schemes of finite type over $S$. Assume that $Y$ is flat over $S$ and that $X$ is normal. Also let $\alpha \in c^{\mathrm{nai,eff}}(X \times_S Y|X)$.

We define its *symmetrization* $\mathrm{sym}_S(\alpha)\colon X \to \mathrm{S}^n(Y|S)$ as the composition

$$X \xrightarrow{\mathrm{sym}(\alpha)} \mathrm{S}^n(X \times_S Y|X) \xrightarrow{\mathrm{S}^n(\mathrm{pr}_Y|X\to S)} \mathrm{S}^n(Y|S).$$

**Proposition 3.7.2** (Symmetrization is compatible with composition)**.** *Let $S$ be a noetherian scheme. Let $X$, $Y$ and $Z$ be schemes over $S$. Assume furthermore that $X$ and $Y$ are normal and that $Y$ and $Z$ are flat and of finite type over $S$. Let $\alpha\colon X \rightsquigarrow Y$ and $\beta\colon Y \rightsquigarrow Z$ be effective finite correspondences of constant degrees $m$ and $n$, respectively.*

*Then the composition $\beta\circ\alpha$ is effective of constant degree $mn$ and the symmetrization morphisms of Definition 3.7.1 are compatible with composition, i.e. the diagram*

$$\begin{array}{ccc}
X & \xrightarrow{\mathrm{sym}_S(\alpha)} & \mathrm{S}^m(Y|S) \\
\downarrow{\scriptstyle \mathrm{sym}_S(\beta\circ\alpha)} & & \downarrow{\scriptstyle \mathrm{S}^m(\mathrm{sym}_S(\beta)|S)} \\
\mathrm{S}^{mn}(Z|S) & \xleftarrow{\tau^{m,n}} & \mathrm{S}^m(\mathrm{S}^n(Z|S)|S)
\end{array}$$

*commutes.*

*Proof.* The effectivity of $\beta \circ \alpha$ is clear from Remark 1.4.25. The statement about the degree is Lemma 1.8.13.

Decomposing $X$ into connected, hence by normality irreducible, components and Remark 3.5.6 allow us to assume that $X$ is the spectrum of a field. We can then by linearity (cf. Remark 3.5.5) assume that $\alpha \in c^{\mathrm{nai,eff}}(X \times_S Y|X) = c^{\mathrm{eff}}(X \times_S Y|X)$ is basic, the equality by Proposition 1.4.18. Then $\Gamma := \mathrm{supp}(\alpha)$ is also the spectrum of a field. Let $i\colon \Gamma \hookrightarrow Y$ be the inclusion.

Combining Definition 1.4.22 with Definition 1.6.5 we see that

$$\beta \circ \alpha = (i \times_S Z)_* \operatorname{Cor}(i^{\circledast}\beta, \alpha).$$

We look at the diagram

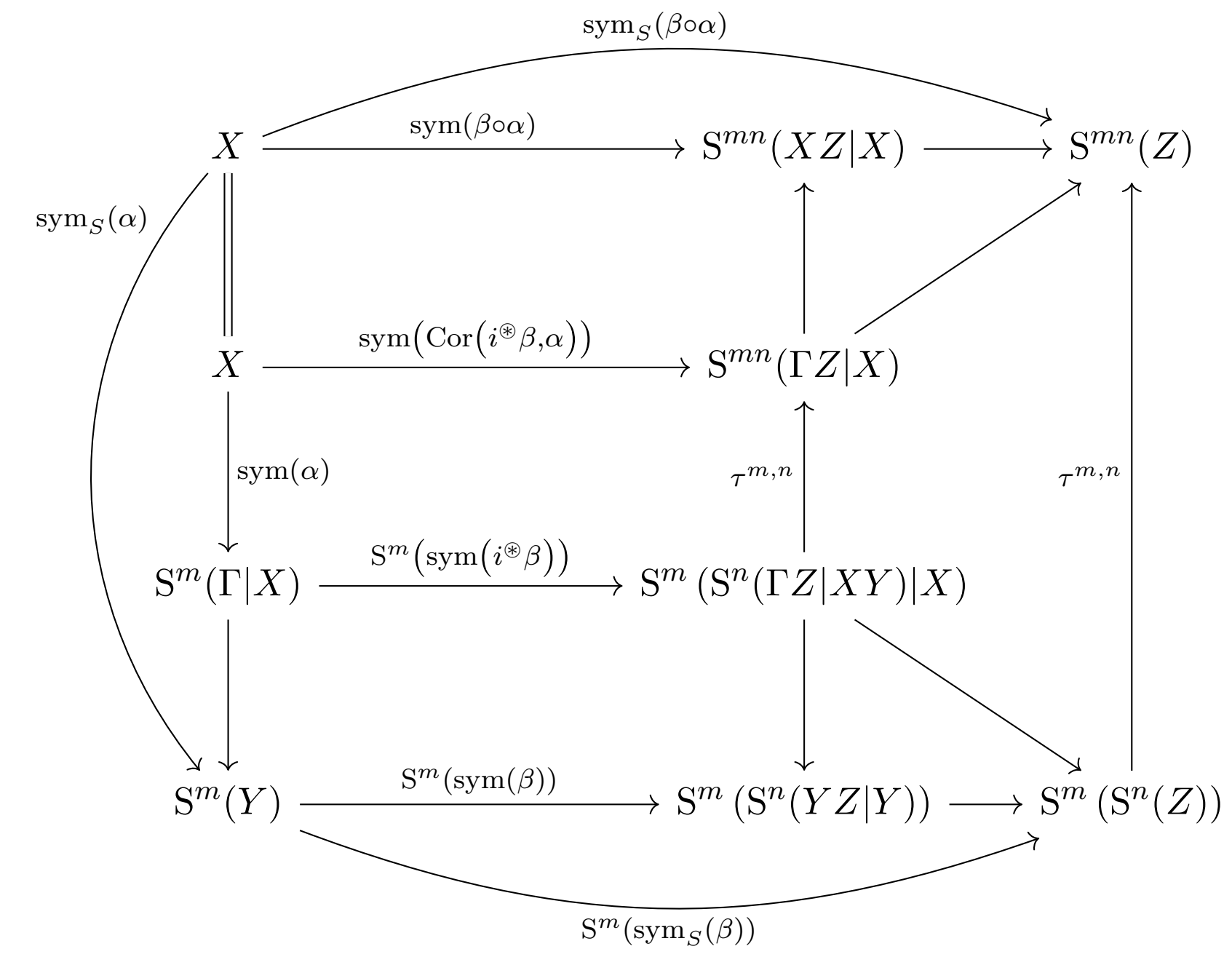

where we suppressed the base $S$, used our usual abbreviations for products and refrained from labelling the morphisms that are only the functoriality of symmetric products. We have to show that the boundary commutes, hence by a diagram chase want to argue that all the areas on the inside commute:

- areas inside the three arcs: these are the Definition 3.7.1 of $\mathrm{sym}_S(-)$,
- two triangles: functoriality of $\mathrm{S}^m(-)$ and $\mathrm{S}^n(-)$,
- trapezoid: naturality of $\tau^{m,n}$ (Proposition 3.3.3),
- top square: this is the compatibility of $\mathrm{sym}(-)$ with pushforward (Proposition 3.6.5),
- middle square: compatibility of $\mathrm{sym}(-)$ with $\mathrm{Cor}(-,-)$ (Proposition 3.6.8),
- lower square: this is the functoriality of $\mathrm{S}^m(-)$ applied to the compatibility of $\mathrm{sym}(-)$ with pullbacks (Proposition 3.6.6). □

**Remark 3.7.3.** For smooth varieties over a field $k$ of characteristic 0 there exists a short argument involving generic étaleness that can be found in [Ayo14b], Appendix A.

Proposition 3.7.2 has been stated over arbitrary fields in [BV08]. We found their proof hard to follow, partially because they make the incorrect assertion that the algebra $\mathrm{S}_n(B|A)$ is generated by the pure tensor powers $b^{\otimes n}$, where $b$ runs over the elements of $B$. Indeed, counterexamples already exist when $B|A$ is a finite extension of perfect fields of positive characteristic, see for example Remark 2.1.20 and recall that $\rho_3(b) = b^{\otimes 3}$.

**Theorem 3.7.4** (Symmetrization is compatible with tensor products)**.** *Let $S$ be a noetherian normal scheme and let $X_1, X_2, Y_1, Y_2$ be schemes which are smooth and of finite type over $S$. Let*

$$\alpha_1\colon X_1 \rightsquigarrow Y_1$$
$$\alpha_2\colon X_2 \rightsquigarrow Y_2$$

*be effective finite correspondences over $S$ of constant degrees $n_1$ and $n_2$, respectively.*

*Then the exterior product*

$$\alpha_1 \otimes \alpha_2\colon X_1 \times_S X_2 \rightsquigarrow Y_1 \times_S Y_2$$

*is effective of constant degree $n_1 n_2$. Furthermore, the symmetrization of finite correspondences as in Definition 3.5.3 is multiplicative, i.e. the diagram*

$$\begin{array}{ccc} X_1 \times_S X_2 & \xrightarrow{\mathrm{sym}_S(\alpha_1)\times \mathrm{sym}_S(\alpha_2)} & \mathrm{S}^{n_1}(Y_1|S) \times_S \mathrm{S}^{n}(Y_2|S) \\ & {\scriptstyle \mathrm{sym}_S(\alpha_1\otimes\alpha_2)} \searrow \qquad \swarrow {\scriptstyle \rho^{n_1,n_2}} & \\ & \mathrm{S}^{n_1 n_2}(Y_1 \times_S Y_2|S) & \end{array}$$

*commutes.*

*Proof.* The effectivity of $\alpha_1 \otimes \alpha_2$ is clear from the definition of the exterior product. It has constant degree $n_1 n_2$ by Lemma 1.8.15.

From now on we suppress the base $S$ whenever possible. Corollary 3.6.7 gives us two commutative triangles

$$\begin{array}{ccc} X_1X_2 & \xrightarrow{\mathrm{sym}\left((\mathrm{pr}^{X_1X_2}_{X_1})^{\circledast}\alpha_1\right)} & \mathrm{S}^n(X_1Y_1X_2|X_1X_2) \\ & {\scriptstyle \mathrm{sym}(\alpha_1)\times X_2} \searrow \qquad \swarrow {\scriptstyle \cong} & \\ & \mathrm{S}^n(X_1Y_1|X_1) \times X_2 & \end{array}$$

$$\begin{array}{ccc} X_1X_2 & \xrightarrow{\mathrm{sym}\left((\mathrm{pr}^{X_1X_2}_{X_2})^{\circledast}\alpha_2\right)} & \mathrm{S}^n(X_1X_2Y_2|X_1X_2) \\ & {\scriptstyle X_1\times\mathrm{sym}(\alpha_2)} \searrow \qquad \swarrow {\scriptstyle \cong} & \\ & X_1 \times \mathrm{S}^n(X_2Y_2|X_2) & \end{array}$$

where the isomorphisms are those of Lemma 3.3.2. Taking their tensor product over $X_1 \times X_2$ gives us a single triangle

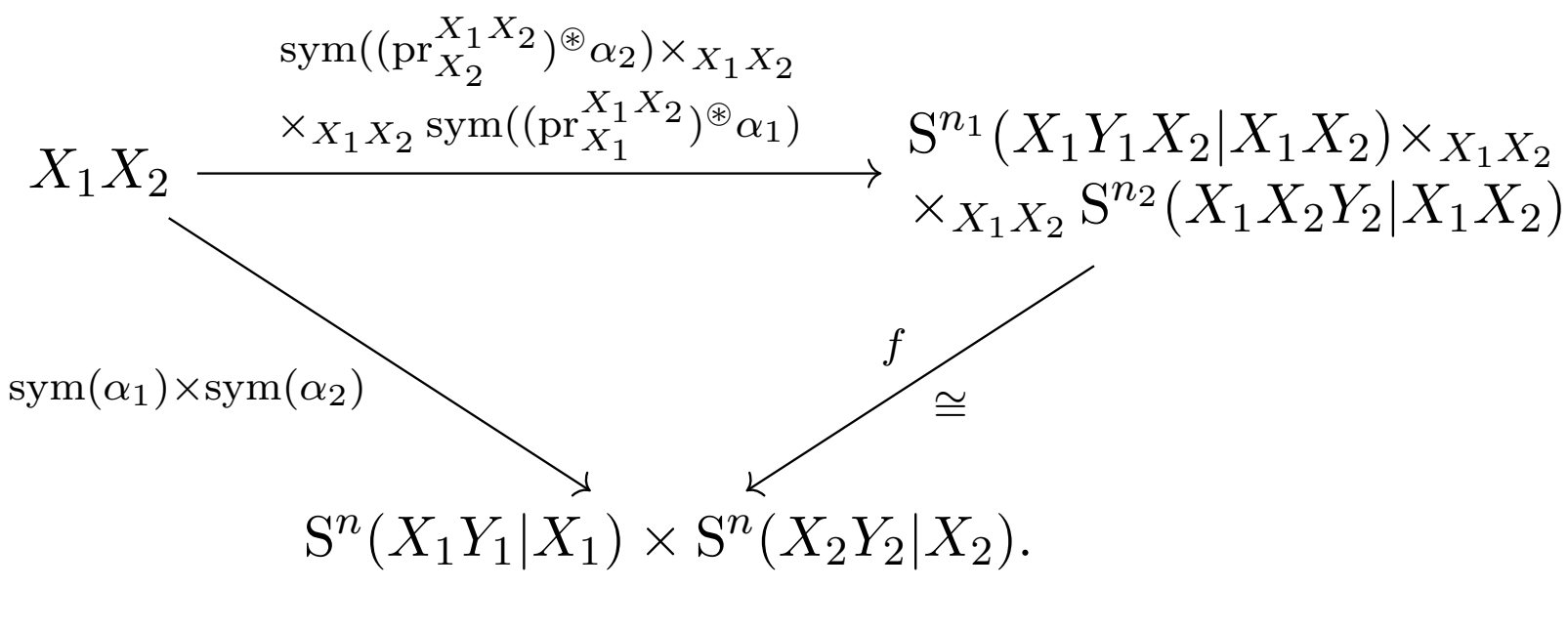

We extend it to a diagram

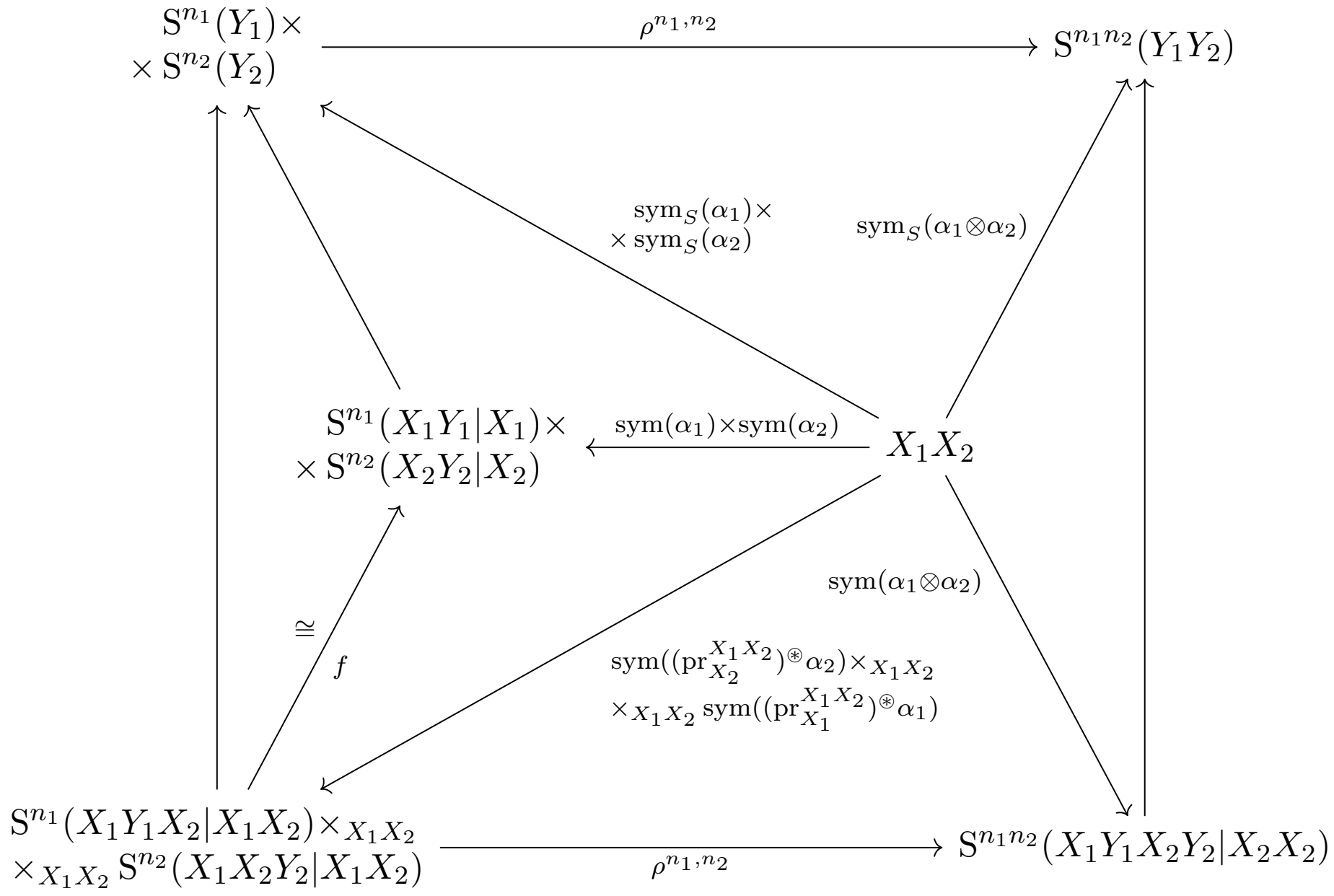

where we refrained from labelling the morphisms that are just the functoriality of symmetric products. We have to show that the top triangle commutes, i.e. want to check that all other areas commutes:

- outside: functoriality of $\rho^{n_1,n_2}$ (Proposition 3.3.7),
- leftmost triangle: functoriality of $\mathrm{S}^{n_1}(-)$, $\mathrm{S}^{n_2}(-)$ and fibre products,
- rightmost triangle: Definition 3.7.1 of $\mathrm{sym}_S(-)$,
- upper inner triangle: Definition 3.7.1 of $\mathrm{sym}_S(-)$,
- lower inner triangle: by the construction above,

- bottom triangle: by the Definition 1.7.4 of the tensor product of finite correspondences and the compatibility of sym(−) with the exterior product (Proposition 3.6.9). □

We are now able to prove the main theorem of this chapter, showing that the interpretation of effective relative cycles as multivalued maps is functorial.

**Theorem 3.7.5.** *Let $S$ be a noetherian normal scheme.*

*Then symmetrization induces an additive, strict monoidal, graded and faithful functor*

$$\mathrm{sym}^{\mathrm{eff}} = \mathrm{sym}_{\mathbb{N}} \colon \mathrm{SmCor}_S^{\mathrm{eff}} \to \mathrm{Multi}_S^{\mathrm{eff}}$$

*between symmetric monoidal $\mathbb{N}_0$-graded categories enriched over abelian monoids. Additionally, this functor is compatible with the embeddings of* $\mathrm{Sm}_S$ *(cf. Remark 4.5.4 and Remark 3.4.13) into both sides.*

*If $S$ is purely of characteristic* 0*, then this functor is full.*

*Proof.* On objects we let $\mathrm{sym}^{\mathrm{eff}}$ to be the identity. Let $X, Y \in \mathrm{Sm}_S$ and $\alpha \in \mathrm{SmCor}_S^{\mathrm{eff}}(X, Y)$. If $X$ is connected, then we get a morphism $\mathrm{sym}_S(\alpha)\colon X \to \mathrm{S}^{\deg(\alpha)}(Y|S)$. Thus if in general $X = \coprod_{i=1}^r X_i$ is its decomposition into connected components and $\alpha_i$ is the pullback of $\alpha$ to $X_i$, we get a multivalued morphism

$$\coprod_{i=1}^r \mathrm{sym}_S(\alpha_i)\colon \coprod_{i=1}^r X \to \coprod_{i=1}^r \mathrm{S}^{\deg(\alpha_i)}(Y|S)$$

from $X$ to $Y$ as desired.

It is by construction enough to check the remaining properties for a connected domain $X$. Functoriality is immediate from Proposition 3.7.2. The preservation of degree and additivity are clear from the construction. The compatibility with the tensor structures is Theorem 3.7.4. The compatibility with the embeddings from $\mathrm{Sm}_S$ is immediate from the definitions.

Next we show the injectivity on hom-sets. Let $\xi = \mathrm{Spec}(k)$ be the generic point of $X$ and let $\zeta = \mathrm{Spec}(L)$ be the spectrum of an algebraically closed field $L|k$. Then we have by Proposition 3.6.6 a commutative square

$$\begin{array}{ccc}
\mathrm{SmCor}_S^{\mathrm{eff}}(X,Y) & \xrightarrow{\mathrm{sym}^{\mathrm{eff}}} & \mathrm{Multi}_S(X,Y) \\
\big\downarrow f^{\circledast} & & \big\downarrow -\circ f \\
c^{\mathrm{eff}}(\zeta \times_S Y|\zeta) & \xrightarrow{a} & \mathrm{Hom}_S(\zeta, \coprod_{d=0}^{\infty} \mathrm{S}^d(Y|S)),
\end{array}$$

where the lower morphism $a$ is induced by symmetrization. Note that $f^{\circledast}$ is injective by Proposition 1.5.5. It hence suffices to show that $a$ is injective.

We note that every basic cycle $\alpha$ over $\zeta$ is the spectrum of a finite field extension of $L$, thus $\alpha = \mathrm{supp}(\alpha) \cong \zeta$. This identifies $c^{\mathrm{eff}}(\zeta \times_S Y|\zeta)$ with the additive monoid $\mathbb{N}_0[Y(L)]$ of formal sums $\sum_{i=1}^{r} g_i$ of not necessarily distinct $L$-valued points $g_i\colon \zeta \to Y$ over $S$.

On the other hand, we note that a morphism $\zeta \to \mathrm{S}^d(Y|S)$ over $S$ is the same as a morphism

$$\zeta \to \zeta \times_S \mathrm{S}^d(Y|S) \cong \mathrm{S}^d(\zeta \times_S Y|\zeta)$$

over $\zeta$, where the isomorphism is that of Proposition 3.1.9. Such a morphism, i.e. a closed point, is by Proposition 3.1.4 (a) nothing else than an unordered $d$-tuple of morphisms $\zeta \to Y$ over $S$, i.e. a formal sum of $d$ of these points. Therefore $\mathrm{Hom}_S(\zeta, \coprod_{d=0}^{\infty} \mathrm{S}^d(Y|S))$ can also be identified with $\mathbb{N}_0[Y(L)]$.

Hence we have a chain

$$c^{\mathrm{eff}}\left(\zeta \times_S Y|\zeta\right) \cong \mathbb{N}_0[Y(L)] \cong \mathrm{Hom}_S\left(\zeta, \coprod_{d=0}^{\infty} \mathrm{S}^d(Y|S)\right)$$

of additive bijections. On the $L$-points themselves, i.e. for actual morphisms, we already know this identification to be $\mathrm{sym}^{\mathrm{eff}}$, and hence by additivity it is $\mathrm{sym}^{\mathrm{eff}}$ everywhere. This shows the required injectivity.

If $S$ is purely of characteristic 0, then the fullness of $\mathrm{sym}^{\mathrm{eff}}$ follows from the arguments in Appendix A of [Ayo14b]. □

**Remark 3.7.6.** Following Definitions 8.4 and 8.5 of [Ryd08d] we have the abelian submonoid

$$Q_S^{\mathrm{nai}}(X,Y) \subseteq c^{\mathrm{nai,eff}}(X \times_S Y|X) \otimes_{\mathbb{N}_0} \mathbb{Q}_0^+$$

of *quasi-integral naive cycles*, freely generated by elements $\frac{1}{e}\alpha$. Here $\alpha \in c^{\mathrm{nai,eff}}(X \times_S Y|X)$ is basic and $e$ is its *inseparable discrepancy* (cf. Definition 8.4 of loc. cit.). We let

$$Q(X,Y) := Q_S^{\mathrm{nai}}(X,Y) \cap (c^{\mathrm{eff}}(X \times_S Y|Y) \otimes_{\mathbb{N}_0} \mathbb{Q}_0^+)$$

be the *effective quasi-integral relative cycles.*

Then $\mathrm{sym}_S$ induces by Theorem 10.16 of [Ryd08d] a bijection

$$Q_S(X,Y) \cong \mathrm{Multi}_S^{\mathrm{eff}}(X,Y).$$

This gives us an explicit description of the discrepancy between the hom-sets of both sides in Theorem 3.7.5. In particular, one immediately gets from their definition that the inseparable discrepancies vanish in characteristic 0, giving another proof for the fullness of $\mathrm{sym}_S$ in this case.

**Theorem 3.7.7.** *Let $S$ be a noetherian normal scheme.*

*Then symmetrization induces an additive, strict monoidal, graded and faithful functor*

$$\mathrm{sym}_{\mathbb{Z}}\colon \mathrm{SmCor}_S \to \mathrm{Multi}_S$$

*between preadditive symmetric monoidal $\mathbb{Z}$-graded categories. Additionally, this functor is compatible with the embeddings of* $\mathrm{Sm}_S$ *(cf. Remark 4.5.4 and Remark 3.4.13) into both sides.*

*If $S$ is purely of characteristic* $0$*, then this functor is full.*

*Proof.* This follows from Proposition 1.4.18 (a) by taking the additive group-completions on both sides in Theorem 3.7.5. □

For the reader's interest we also quickly state and prove a version for arbitrary rings of coefficients:

**Theorem 3.7.8.** *Let $S$ be a regular noetherian scheme and let $\Lambda$ be a ring.*

*Then symmetrization induces a $\Lambda$-linear strict monoidal functor*

$$\mathrm{sym}_\Lambda\colon\ \mathrm{SmCor}(S,\Lambda)\to\mathrm{Multi}(S,\Lambda).$$

*Additionally, this functor is compatible with the embeddings of* $\Lambda[\mathrm{Sm}_S]$ *(cf. Remark 4.5.4 and Remark 3.4.13) into both sides.*

*If $\Lambda$ is flat as a $\mathbb{Z}$-module or if $S$ is purely of characteristic* $0$*, then this functor is faithful. In the latter case it is also full.*

*Proof.* Both sides were defined by taking the $\Lambda$-linear extensions of the case with integral coefficients, i.e. by applying $-\otimes_\mathbb{Z}\Lambda$ to all hom-sets. Hence the existence of the functor $\mathrm{sym}_\Lambda$ and its compatibility with the embeddings from $\Lambda[\mathrm{Sm}_S]$ follow from Theorem 3.7.7.

If $\Lambda$ is a flat $\mathbb{Z}$-module, then $-\otimes_\mathbb{Z}\Lambda$ preserves injectivity, proving faithfulness. If $S$ is purely of characteristic 0, then by Theorem 3.7.7 we find $\mathrm{sym}_\mathbb{Z}$ to be an isomorphism between the integral hom-sets, which is preserved by the base change. Hence the functor is fully faithful in this case. □

**Remark 3.7.9.** The regularity condition in Theorem 3.7.8 is necessary for the reasons explained in Remark 1.6.2.

**Definition 3.7.10.** For a scheme $S$ we let

$$\mathrm{char}(S)=\{\mathrm{char}(\kappa(s))\mid s\in S\}\backslash\{0\}$$

be the set of non-zero characteristics occurring at residue fields.

**Remark 3.7.11.** Using Remark 3.7.6 and unravelling the definitions mentioned there, one can show that the functors $\mathrm{sym}_\mathbb{N}$, $\mathrm{sym}_\mathbb{Z}$ and $\mathrm{sym}_\Lambda$ of Theorems 3.7.5, 3.7.7 and 3.7.8 are full whenever $\mathrm{char}(S)\subseteq\Lambda^\times$.

## 3.8 Commutative group schemes with transfers

**Theorem 3.8.1.** *Let $S$ be a noetherian normal scheme. Let $\mathrm{Sm}_S$ be the category of schemes that are smooth, separated and of finite type over $S$ and let $\mathrm{AlgSp}^\flat_S$ be the category of algebraic spaces that are flat and separated over $S$. Let $\mathcal{G}$ be an abelian group object in $\mathrm{AlgSp}^\flat_S$.*

*Then the corresponding presheaf*

$$\mathcal{G}(-) = \mathrm{AlgSp}^\flat_S(-, \mathcal{G})\colon\ \mathrm{Sm}_S \hookrightarrow \mathrm{AlgSp}^\flat_S \to \mathrm{Ab}$$

*admits transfers, i.e. it extends to a presheaf*

$$\widetilde{\mathcal{G}}\colon\ \mathrm{SmCor}(S, \mathbb{Z}) \to \mathrm{Ab}\,.$$

**Remark 3.8.2.** Theorem 3.8.1 generalizes a result by Spieß and Szamuely (cf. Lemma 3.2 of [SS03]), who consider the special case where the base $S$ is the spectrum of a field and where the group object $\mathcal{G}$ is represented by a scheme. Note that the functoriality of the extension was not checked in [SS03].

Also compare Theorem 3.8.1 to the results of Ancona, Huber and Pepin Lehalleur (cf. Theorem 2.8 of [AHP15]), who give a similar result over an excellent scheme, but require rational coefficients.

*Proof of Theorem 3.8.1.* Due to Theorem 3.7.7 it suffices to extend $\mathcal{G}(-)$ to a presheaf

$$\widetilde{\mathcal{G}}\colon\ \mathrm{Multi}(S, \mathbb{Z}) \to \mathrm{Ab}\,.$$

On objects we get nothing new, hence we only have to functorially define $\widetilde{\mathcal{G}}(\alpha)$ for all multivalued morphisms $\alpha\colon X \multimap Y$ between schemes $X, Y \in \mathrm{Sm}_S$. Assume for now that $\alpha$ is effective and that $X$ is connected. Hence $\alpha$ is homogeneous and we denote its degree by $n$. We also suppress writing the base $S$.

Now consider the morphism $\Delta_n\colon \mathcal{G}^n \to \mathcal{G}$ induced by the group structure of $\mathcal{G}$. Because $\mathcal{G}$ is commutative, $\Delta_n$ is unaffected by permuting the factors of $\mathcal{G}^n$. It thus descends, due to the universal property of quotients by a group action, to a morphism $\delta_n\colon\ \mathrm{S}^n(\mathcal{G}) \to \mathcal{G}$.

Let $f\colon Y \to \mathcal{G}$ be an element of $\mathcal{G}(Y)$. The functoriality of $\mathrm{S}^n(-)$ then gives us a composite morphism

$$X \xrightarrow{\ \alpha\ } \mathrm{S}^n(Y) \xrightarrow{\mathrm{S}^n(f)} \mathrm{S}^n(\mathcal{G}) \xrightarrow{\ \delta_n\ } \mathcal{G}$$

which defines us an element $\widetilde{\mathcal{G}}(\alpha)(f) \in \mathcal{G}(X)$. This is additive in $\alpha$, hence extends uniquely to the non-effective multivalued morphisms and similarly to non-connected $X$.

We still have to show functoriality. It suffices to check it for effective homogeneous multivalued morphisms $\alpha\colon X \multimap Y$ and $\beta\colon Y \multimap Z$ of degrees

$m$ and $n$, respectively. For this we pick any morphism $f\colon Z \to \mathcal{G}$ and consider the diagram

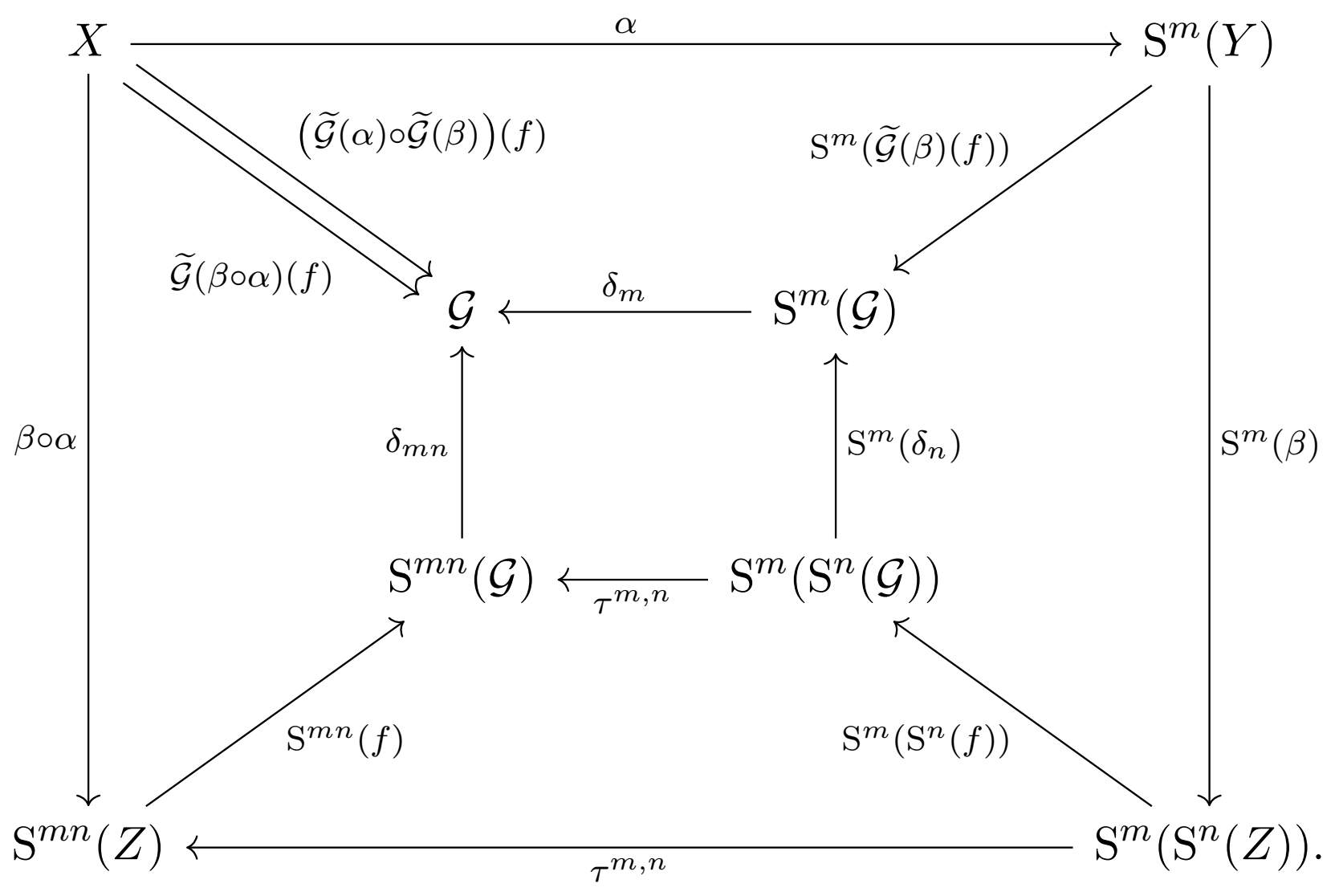

We observe that

- the inner square commutes because both paths are quotients of the multiplication morphism $\mathcal{G}^{mn} \to \mathcal{G}$,
- the outside is the Definition 3.4.7 of composition of multivalued morphisms,
- the bottom area is the functoriality of $\tau^{m,n}$ (Proposition 3.3.3),
- the remaining three areas are the construction of $\widetilde{\mathcal{G}}(-)$, where we also applied $\mathrm{S}^m(-)$ to the one on the right.

Hence a simple diagram chase shows that indeed

$$\widetilde{\mathcal{G}}(\beta \circ \alpha)(f) = \Big(\widetilde{\mathcal{G}}(\alpha) \circ \widetilde{\mathcal{G}}(\beta)\Big)(f),$$

i.e. that $\widetilde{\mathcal{G}}$ is a contravariant functor $\mathrm{Multi}_S \to \mathrm{Ab}$. $\square$

**Remark 3.8.3.** The assumption that $\mathcal{G}$ is flat over $S$ was only needed for the functorial morphism $\tau^{m,n}$ to exist. This assumption can be removed by using divided powers as defined in the next Section 3.9, Theorem 3.9.10 in particular. Note that symmetric products and divided powers agree on flat objects by Proposition 3.9.12, hence this change does not influence any other part of the proof.

## 3.9 Divided powers

Let us explain the more general theory of divided powers as developed by Rydh in [Ryd08a], [Ryd08b], [Ryd08c] and [Ryd08d]. Just like Section 2.6, this section will not be used anywhere else in this thesis.

We have already met the affine versions in Section 2.6. The general definition is as follows:

**Definition 3.9.1.** Let $f\colon X \to S$ be an affine morphism of schemes. Then, as described in more details in Section 1.4 of [Ryd08a], the scheme $\Gamma^d(X|S)$ of *divided powers* is defined by glueing $\mathrm{Spec}(\Gamma_d(\mathcal{O}_{X\times_S U}|\mathcal{O}_U))$ where $U$ runs over all affine open subschemes $U \hookrightarrow S$. In other words, one can describe it as $\Gamma^d(X|S) = \mathrm{Spec}(\Gamma_d(f_*\mathcal{O}_X|\mathcal{O}_S))$.

By loc. cit. this extends via étale descent to an algebraic space $\Gamma^d(X|S)$ whenever $X \to S$ is an affine morphism of algebraic spaces.

Definitions 3.1.1 and 3.1.3 of [Ryd08a] extend this to non-affine morphisms as follows:

**Definition 3.9.2.** Let $X \to S$ be a separated morphism of algebraic spaces.

A *family of zero-cycles of degree* $d$ on $X|S$ consists of:

- a closed subscheme $Z \hookrightarrow X$ such that the resulting morphism $Z \to S$ is integral, in particular affine,
- a morphism $\alpha\colon S \to \Gamma^d(Z|S)$.

Two families $(Z_1, \alpha_1)$ and $(Z_2, \alpha_2)$ of zero-cycles are called *equivalent* if:

- they have the same degree $d$,
- there exists a family of zero-cycles $(Z, \alpha)$ of degree $d$, where $Z \hookrightarrow Z_1 \cap Z_2$ is a common closed subscheme,
- we have factorizations $\alpha_k = \Gamma^d(i_k|S)\circ\alpha$, where $k \in \{1, 2\}$ and $i_k\colon Z \hookrightarrow Z_k$ is the inclusion.

Let $T$ be another scheme over $S$. We define the functor

$$\underline{\Gamma}^d_{X|S}\colon \mathrm{Sch}_S \to \mathrm{Set}$$

by sending $T$ to the set of equivalence classes of families of zero cycles of degree $d$ on $X \times_S T|T$.

Lastly, by [Ryd08a], Theorem 3.4.1, this functor is represented by an algebraic space $\Gamma^d(X|S)$ of *divided powers*, which agrees with the algebraic space of Definition 3.9.1 if $X \to S$ is affine.

**Remark 3.9.3.** In the case $d = 1$ we always have a natural isomorphism $\Gamma^1(X|S) \cong X$ by Remark 3.1.5 of [Ryd08a].

**Remark 3.9.4.** The existence of the algebraic space $\Gamma^d(X|S)$ is far from trivial and is one of the main results of [Ryd08a]. It becomes significantly easier in the following case:

**Proposition 3.9.5.** *Let the scheme $X$ be AF over a separated scheme $S$. Then $\Gamma^d(X|S)$ is a scheme and AF over $S$.*

*Proof.* This is Theorem 3.1.11 of [Ryd08a]. □

**Remark 3.9.6.** Divided powers behaves very well with respect to our functorialities: all the results of Sections 3.2 and 3.3 remain true when replacing $\mathrm{S}^\bullet(-|-)$ by $\Gamma^\bullet(-|-)$, even if we remove all assumptions of flatness. Let us now elaborate on these generalities, whose affine counterparts we saw in Section 2.6.

**Lemma 3.9.7.** *Let $X$, $Y$ be an algebraic spaces separated over an algebraic space $S$ and let $d$ be a non-negative integer.*

*Then there is a natural isomorphism*

$$\Gamma^d(X \times_S Y|Y) \cong \Gamma^d(X|S) \times_S Y.$$

*Proof.* One can directly check that the functor $\underline{\Gamma}^d_{X\times_S Y|Y}$ is represented by $\Gamma^d(X|S) \times_S Y$. □

**Lemma 3.9.8.** *Let $d$ be a non-negative integer and let $\mathrm{AlgSp}^{\to,\mathrm{sep}}$ be the category of separated morphisms between algebraic spaces.*

*Then $\Gamma^d(-|-)$ constitutes a base-preserving endofunctor*

$$\mathrm{AlgSp}^{\to,\mathrm{sep}} \to \mathrm{AlgSp}^{\to,\mathrm{sep}}.$$

*Proof.* For every morphism $f\colon X \to Y$ of algebraic spaces separated over an algebraic space $S$, a pushforward $f_*\colon \Gamma^d(X|S) \to \Gamma^d(Y|S)$ is constructed in Definition 3.3.1 of [Ryd08a] by taking scheme-theoretic images.

Let now $(Y|T) \to (X|S)$ be any morphism in $\mathrm{AlgSp}^{\to,\mathrm{sep}}$. This induces a morphism $a\colon Y \to X \times_S T$. We get the desired natural morphism $\Gamma^d(Y|T) \to \Gamma^d(X|S)$ as the composition

$$\Gamma^d(Y|T) \xrightarrow{a_*} \Gamma^d(X \times_S T|T) \cong \Gamma^d(X|S) \times_S T \xrightarrow{\mathrm{pr}_1} \Gamma^d(X|S),$$

the isomorphism being that of Lemma 3.9.7. □

**Proposition 3.9.9.** *Let $X$, $Y$ be algebraic spaces separated over an algebraic space $S$ and let $k$ be a non-negative integer.*

*There exists a natural decomposition*

$$\Gamma^k(X \sqcup Y|S) \cong \coprod_{m'+n'=k} \Gamma^{m'}(X|S) \times_S \Gamma^{n'}(Y|S).$$

*Proof.* On the represented functors this is Proposition 3.1.8 of [Ryd08a], whence it also holds on the representing algebraic spaces. □

**Theorem 3.9.10.** *Let $S$ be an algebraic space.*

1. *For all separated morphisms $X \to T \to S$ of algebraic spaces and non-negative integers $m, n$ there exists a natural morphism*
$$\tau^{m,n} \colon \Gamma^m(\Gamma^n(X|T)|S) \to \Gamma^{mn}(X|S).$$
2. *For all separated morphism $X \to S$ of algebraic spaces and non-negative integers $n_1, \dots, n_r$ there exists a natural morphism*
$$\sigma^{(n_i)} \colon \prod_{i=1}^{r} \Gamma^{n_i}(X|S) \to \Gamma^{\sum_{i=1}^{r} n_i}(X|S).$$
3. *For all algebraic spaces $X_1, \dots, X_r$ separated over $S$ and non-negative integers $n_1, \dots, n_r$ there exists a natural morphism*
$$\rho^{(n_i)} \colon \prod_{i=1}^{r} \Gamma^{n_i}(X_i|S) \to \Gamma^{\prod_{i=1}^{r} n_i}\left(\prod_{i=1}^{r} X_i|S\right).$$

*Proof.*

1. The morphism $\tau^{m,n}$ can be found as Definition 7.2 of [Ryd08b], where it follows from the affine variant, which again is induced by the homogeneous multiplicative polynomial law
$$\begin{aligned} A &\mapsto \Gamma^m(\Gamma^n(B|A)|A) \\ x &\mapsto \gamma_m(\gamma_n(x)) \end{aligned}$$
of degree $mn$ (cf. Proposition 2.6.9).
2. This is Definition-Proposition 4.1.1 of [Ryd08a]: The morphism $\sigma^{m,n}$ can be constructed from Proposition 3.9.9 by means of Remark 3.3.6, i.e. by setting $X = Y$ and using the obvious projection $X \sqcup X \to X$. The variant with three or more integers $n_i$ follows by iteration.
3. Similarly it suffices to construct the morphism $\rho^{m,n}$. It is defined in Section 8 of [Ryd08b], but we give a different yet related construction: the morphism can be obtained by mimicking Lemma 3.3.8, i.e. using Lemma 3.9.7 to get a diagram
$$\begin{array}{ccc} \Gamma^m\big(X \times_S \Gamma^n(Y|S)|S\big) & \xleftarrow{\cong} & \Gamma^m\big(\Gamma^n(X \times_S Y|X)|S\big) \\ \uparrow & & \\ \Gamma^m\big(X \times_S \Gamma^n(Y|S)|\Gamma^n(Y|S)\big) & & \downarrow \tau^{m,n} \\ \downarrow \cong & & \\ \Gamma^m(X|S) \times_S \Gamma^n(Y|S) & \xrightarrow{\rho^{m,n}} & \Gamma^{mn}(X \times_S Y|S). \end{array}$$

defining the desired morphism. □

We state without an explicit proof:

**Proposition 3.9.11.** *Let $S$ be an algebraic space, let $X, Y, Z$ be algebraic spaces separated over $S$ and let $k, l, m, n$ be positive integers.*

*Then the three natural morphisms of Theorem 3.9.10 satisfy the commutative diagrams of Proposition 3.3.9, but with all occurrences of $\mathrm{S}^\bullet$ replaced by $\Gamma^\bullet$.*

We can also compare divided powers to symmetric products, which turns the results of Sections 3.2 and 3.3 into special cases of the above:

**Proposition 3.9.12.** *Let $X \to S$ be a separated morphism of algebraic spaces and let $n$ be a non-negative integer.*

*Then there exists a natural morphism*

$$\mathrm{SG}^n_{X|S} \colon \mathrm{S}^n(X|S) \to \Gamma^n(X|S)$$

*which is a universal homeomorphism with trivial residue field extensions.*

*It is an isomorphism if $X \to S$ is flat or $S$ is purely of characteristic $0$. This isomorphism is compatible with the functorialities $\sigma^{m,n}$, $\tau^{m,n}$ and $\rho^{m,n}$ in the obvious sense.*

*Proof.* If $n = 1$ then both sides are naturally isomorphic to $X$, hence this is trivially true. The general existence of the morphism is found as Proposition 4.1.5 i) and Remark 4.2.2 of [Ryd08a]: there is a canonical $S_n$-invariant morphism

$$\Psi_X = \sigma^{\overbrace{1,1,\ldots,1}^{n \text{ times}}} \circ (\mathrm{SG}^1_{X|S}|S)^n \colon (X|S)^n \cong (\Gamma^1(X|S)|S)^n \to \Gamma^n(X|S),$$

which descends to the desired morphism due to the universal property of quotients by group actions.

Corollary 4.2.5 of op. cit. shows that $\mathrm{SG}^n_{X|S}$ is an isomorphism when the additional assumptions are satisfied. The compatibility with $\sigma^{m,n}$ is almost by construction, and that with $\tau^{m,n}$ and $\rho^{m,n}$ is a lengthy but straightforward check. We omit the details. □

**Remark 3.9.13.** Apart from results regarding relative cycles, we only used the aforementioned results in sections 3.4 to 3.7. Therefore one can now mutatis mutandis transfer everything into the formalism of divided powers, while omitting assumptions of flatness where not needed, i.e. when they were only required to assure the functorial properties of symmetric products. Many of these results already appear in [Ryd08d].

We finish this section by giving some important examples, but often omit the proof.

**Definition 3.9.14.** Let $S$ be an algebraic space and let $X$, $Y$ be algebraic spaces separated over $S$. An *effective multivalued morphism* $X \multimap Y$ over $S$ is a morphism

$$X \to \coprod_{d=0}^{\infty} \Gamma^d(Y|S)$$

of algebraic spaces over $S$. We call it *homogeneous of degree* $d$ if it is simply a morphism $X \to \Gamma^d(Y|S)$.

We denote the set of all multivalued morphisms $X \multimap Y$ by $\widetilde{\mathrm{Multi}}_S^{\mathrm{eff}}(X, Y)$.

*Composition*, *sum* and *tensor* product of multivalued morphisms are then defined completely analogous to Definitions 3.4.3, 3.4.7 and 3.4.9.

Note that this agrees with our previous definition whenever it applies.

**Proposition 3.9.15.** *Let $S$ be an algebraic space. Let $I$ and $J$ be arbitrary sets. Also let algebraic spaces $X_i$, $i \in I$, and $Y_j$, $j \in J$, separated over $S$ be given.*

*Then there is a natural additive decomposition*

$$\widetilde{\mathrm{Multi}}_S^{\mathrm{eff}}\left(\coprod_{i\in I} X_i, \coprod_{j\in J} Y_j\right) \cong \prod_{i\in I}\bigoplus_{j\in J}\widetilde{\mathrm{Multi}}_S^{\mathrm{eff}}(X_i, Y_j).$$

**Definition 3.9.16.** Let $S$ be an algebraic space. We define the category $\widetilde{\mathrm{Multi}}_S^{\mathrm{eff}}$ of *effective extended multivalued morphisms* over $S$ as follows:

- its objects are the algebraic spaces that are separated over $S$,
- its morphisms are the multivalued morphisms $X \multimap Y$ of Definition 3.9.14,
- composition is as described in Definition 3.9.14.

Almost identical to Theorem 3.4.11 we get, this time from Proposition 3.9.11:

**Theorem 3.9.17.** *Let $S$ be an algebraic space.*

*Then composition, addition and tensor product turn the category $\widetilde{\mathrm{Multi}}_S^{\mathrm{eff}}$ of effective extended multivalued morphisms over $S$ into a symmetric monoidal category enriched over abelian monoids.*

**Remark 3.9.18.** The natural isomorphism $\Gamma^1(Y|S) \cong Y$ induces, analogous to Remark 3.4.13, a natural embedding of the category of algebraic spaces separated over $S$ into $\mathrm{Multi}_S^{\mathrm{eff}}$.

**Definition 3.9.19.** Let $S$ be an algebraic space and let $\Lambda$ be a ring.

The category $\widetilde{\mathrm{Multi}}_{S,\Lambda}$ of *extended $\Lambda$-multivalued morphisms* has the same objects as $\widetilde{\mathrm{Multi}}_S^{\mathrm{eff}}$ and the $\Lambda$-linear extensions

$$\widetilde{\mathrm{Multi}}_S^{\mathrm{eff}}(X,Y) \otimes_{\mathbb{N}} \Lambda$$

as morphisms.

If $\Lambda = \mathbb{Z}$ we will call it the category of *extended multivalued morphisms* and denote it by $\widetilde{\mathrm{Multi}}_S$.

**Proposition 3.9.20.** *Let $S$ be an algebraic space and let $\Lambda$ be a ring.*

*Then $\widetilde{\mathrm{Multi}}_{S,\Lambda}$ is a $\Lambda$-linear symmetric tensor category.*

The results of Sections 3.5 to 3.7 now translate, which in particular gives us:

**Theorem 3.9.21.** *Let $S$ be a noetherian scheme. Let $\mathrm{Nm}_S$ be the category of normal schemes which are separated and of finite type over $S$. Let $\mathrm{NmCor}_S^{\mathrm{eff}}$ be the full subcategory of $\mathrm{SchCor}_S^{\mathrm{eff}}$ consisting of normal schemes which are separated and of finite type over $S$.*

*Then symmetrization induces an additive, graded and faithful functor*

$$\widetilde{\mathrm{sym}}^{\mathrm{eff}} = \mathrm{sym}_{\mathbb{N}} \colon \mathrm{NmCor}_S^{\mathrm{eff}} \to \widetilde{\mathrm{Multi}}_S^{\mathrm{eff}}$$

*between preadditive $\mathbb{N}_0$-graded categories. This functor is compatible with the embeddings of $\mathrm{Nm}_S$ into both sides.*

*Sketch of proof.* This follows from the results of [Ryd08b] and [Ryd08d]. Let us give a quick description how one argues based on our results:

For any finite surjective $\pi \colon \Gamma \to X$ of degree $n$ into a normal and connected scheme $X$, a morphism $\widetilde{\mathrm{sym}}(\pi) \colon X \to \Gamma^n(\Gamma|X)$ is given by $\widetilde{\mathrm{sym}}(\pi) = \mathrm{SG}^d_{X|S} \circ \mathrm{sym}(\pi)$ (cf. Corollary 2.6.14 for the relation to the universal property in the affine case). It extends as in Definitions 3.5.3 and 3.5.4.

The compatibilities analogous to Section 3.6 follow by similar proofs. Often they can, using Remark 3.5.6, be reduced to spectra of fields, where everything becomes flat. Hence in those cases the required properties follow from Proposition 3.9.12 and the already established variants with symmetric products. More directly, they can be extracted from Section 7 of [Ryd08b] and Section 10 of [Ryd08d].

After this, the theorem follows from the same arguments as found in the proofs of Proposition 3.7.2 and Theorem 3.7.5. □

**Remark 3.9.22.** The other functors of Section 3.7 have similar analogues.

**Remark 3.9.23.** Note that we cannot simply transfer the tensor structure: a fibre product of normal schemes is not always normal. Restricting to smooth schemes over a normal base will by Proposition 3.9.12 give us nothing new beyond Theorem 3.7.5.

# Chapter 4

# Yoga of Good Filtrations

Working with Nori motives often involves the so-called 'yoga of good pairs'. Its main idea is to find certain algebraic filtrations $\mathcal{F}$ on varieties that behave like skeleta of CW-structures on cohomology (cf. Definition 4.1.1) to get (very) good pairs, hence the name. A lot of the required groundwork was already laid by Nori, Huber, Müller-Stach and others.

We first repeat the basic definitions and properties as found in [HM16], in particular Nori's cellular complexes $\mathrm{C}^\bullet_{\mathcal{F}}(X) \in C^b(\mathcal{MM}^{\mathrm{eff}}_{\mathrm{Nori}})$ (cf. Definition 4.2.1) that will not depend on $\mathcal{F}$ as an object of $D^b(\mathcal{MM}^{\mathrm{eff}}_{\mathrm{Nori}})$ (cf. Lemma 4.2.4). Their existence is a consequence of Nori's Basic Lemma 4.3.1.

Following Nori's original construction, we associate to every finite correspondence $\alpha\colon X \to Y$ between smooth affine varieties a morphism

$$\mathrm{C}^\bullet_{\mathcal{F}}(\alpha)\colon\ \mathrm{C}^\bullet_{\mathcal{F}}(Y) \to \mathrm{C}^\bullet_{\mathcal{F}}(X)$$

(cf. Definition 4.6.2). As a connected smooth variety is automatically irreducible and normal, we will be able to apply the results of Chapter 3, which allows us to replace finite correspondences by multivalued morphisms.

Thereafter, the primary goal of this section will be to prove Theorem 4.8.1. It roughly states that assigning to a smooth affine variety $X$ the pro-object

$$\varprojlim_{\mathcal{F}} \mathrm{C}^\bullet_{\mathcal{F}}(X) \in \text{pro-}C^b(\mathcal{MM}^{\mathrm{eff}}_{\mathrm{Nori}})$$

is functorial with respect to arbitrary finite correspondences. By the nature of this functor and the $\mathrm{C}^\bullet_{\mathcal{F}}(X)$ it will then follow quite formally that it extends to a functor

$$C^b(\mathrm{SmCor}^{\mathrm{aff}}) \to \text{pro-}C^b(\mathcal{MM}^{\mathrm{eff}}_{\mathrm{Nori}})$$

such that the individual objects of the projective system become naturally isomorphic in $D^b(\mathcal{MM}^{\mathrm{eff}}_{\mathrm{Nori}})$ (cf. Theorem 7.2.4). This completes Step 3. of our overarching Strategy 0.3.1.

Similarities with the arguments in Section 9.2 of [HM16] are not coincidental: on morphisms, they are simply identical as explained in Remark 4.5.4. But to get functoriality for finite correspondences, not just morphisms, one needs to add several additional arguments. Most prominently this involves Nori's equivariant version (cf. Theorem 4.4.3) of his Basic Lemma 4.3.1. The diagrams will, in contrast to the case of morphisms, contain additional objects that serve as stepping stones for constructing morphisms and proving their properties. We thus have to introduce the notion of a filtration on a finite correspondence (cf. Definition 4.5.3) and then that of functorial filtrations (cf. Definition 4.8.2).

**Convention 4.0.1.** For this chapter we assume that $k \subseteq \mathbb{C}$ is a fixed base field. Recall that by a variety we mean a reduced and separated scheme of finite type over $k$. We also fix a noetherian ring of coefficients $\Lambda$ and will often work in the category $\mathrm{SmCor}(k, \Lambda)$ of finite correspondences with coefficients in $\Lambda$ between smooth varieties over $k$ as given in Definition 1.9.1. We will also use the category $\mathrm{SmCor}^{\mathrm{eff}}(k) = \mathrm{SmCor}(k, \mathbb{N})$ of effective finite correspondences defined there.

An important idea when working with Nori motives is to work in-between the algebraic and the analytic topology. All topological notions are, unless stated otherwise, within the context of varieties, i.e. in the Zariski topology. In contrast we will consider sheaves to be on the analytic side:

**Convention 4.0.2.** Throughout this chapter, a *sheaf* on a variety $X \in \mathrm{Var}_k$ is by definition a sheaf on the analytification $X^{\mathrm{an}} = X(\mathbb{C})$ with the analytic topology. Consequently, if $f\colon X \to Y$ is a morphism of varieties over $k$, then the pushforward $f_*$, the proper pushforward $f_!$ and the pullback $f^*$ denote their analytic counterparts. Therefore $\Lambda_X$ will be the constant sheaf $\Lambda_{X^{\mathrm{an}}}$ on $X^{\mathrm{an}}$. Constructible sheaves (cf. Definition 4.3.5) follow suit, but we require all stratifications to come from ones on varieties.

If now $F$ is a sheaf on $X$, i.e. by our convention a sheaf on $X^{\mathrm{an}}$, then we use the shortcut $H^n(X, F)$ instead of $H^n(X^{\mathrm{an}}, F)$ for sheaf cohomology. If $Y$ is a closed subset of $X$, we denote by $H^n_{\mathrm{sing}}(X, Y, \Lambda)$ the singular cohomology $H^n_{\mathrm{sing}}(X^{\mathrm{an}}, Y^{\mathrm{an}}, \Lambda)$ with coefficients in $\Lambda$. In particular we set

$$H^n_{\mathrm{sing}}(X, \Lambda) = H^n_{\mathrm{sing}}(X^{\mathrm{an}}, \Lambda).$$

## 4.1 Cellular filtrations

To formalize our results and techniques we will need some notation, most of which are standard. Also recall Definitions 1.2.2, 1.2.5 and 1.2.9.

Based on Nori's work and [HM16] we define:

**Definition 4.1.1** (Cellular filtrations)**.**

- A *standard filtration* $\mathcal{F}$ *on a scheme* $X$ is an increasing filtration $X_\bullet = \mathcal{F}_\bullet X$ of $X$ by closed subsets such that $\dim X_i \leq i$. Here we set $\dim \emptyset = -\infty$ and thus require $X_i = \emptyset$ for all $i < 0$. Unless stated otherwise we understand the $X_i$ as closed subschemes of $X$ with the reduced induced structure.

- If $\mathcal{F}$ and $\mathcal{G}$ are two standard filtrations on $X$, we say that $\mathcal{F}$ is *finer* than $\mathcal{G}$, or that $\mathcal{G}$ is *coarser* than $\mathcal{F}$, if $\mathcal{F}_i X \subseteq \mathcal{G}_i X$ for all $i$.

- We call standard filtrations $\mathcal{F}_\bullet X$ and $\mathcal{F}_\bullet Y$ on schemes $X$ and $Y$ *compatible with a morphism* $f \colon X \to Y$ if $f(\mathcal{F}_i X) \subseteq \mathcal{F}_i Y$ for all $i$.

- A *standard filtration on a diagram of schemes* consists of standard filtrations on each object. We do not require them to be compatible. We call a standard filtration *finer* (respectively *coarser*) than another one if it is finer (respectively coarser) on each object.

- Let $D$ be a diagram of schemes and let $\mathcal{F}$ be a standard filtration on $D$. Then its *restriction to a subdiagram* $D' \subseteq D$ is the filtration $\mathcal{F}_{|D'}$ on $D'$ obtained by forgetting the objects outside of $D'$. We say that $\mathcal{F}$ *extends* a standard filtration $\mathcal{G}$ on $D'$ if $\mathcal{F}_{|D'} = \mathcal{G}$.

- A *(cohomological)* $\Lambda$*-cellular filtration on a variety* $X$ over $k$ is a standard filtration $X_\bullet = \mathcal{F}_\bullet X$ on $X$ such that:

  - $X_i \backslash X_{i-1}$ is smooth for all $i \in \mathbb{Z}$,
  - $X_i = X$ for $i \geq \dim(X)$,
  - $H^j_{\mathrm{sing}}(X_i, X_{i-1}, \Lambda) = 0$ for all integers $i \neq j$,
  - $H^i_{\mathrm{sing}}(X_i, X_{i-1}, \Lambda)$ is a finitely generated projective $\Lambda$-module for all integers $i$.

- A $\Lambda$*-cellular filtration on a diagram of varieties* is a standard filtration on the diagram which is a $\Lambda$-cellular filtration on each object and compatible with all morphisms.

**Remark 4.1.2.** Note that we do not require standard filtrations on diagrams to be compatible with the morphisms, in contrast to the case of $\Lambda$-cellular filtrations.

**Remark 4.1.3.** The dual notion of a *homological* $\Lambda$*-cellular filtration* exists, defined by using singular homology instead of cohomology. We have the following interesting observation:

**Proposition 4.1.4.** *Let* $X$ *be a variety over a field* $k \subseteq \mathbb{C}$*. The following are equivalent for a standard filtration* $\mathcal{F}_\bullet$ *on* $X$*:*

(a) $\mathcal{F}_\bullet$ *is cohomologically* $\mathbb{Z}$*-cellular,*

(b) $\mathcal{F}_\bullet$ *is cohomologically* $\Lambda$*-cellular for all rings* $\Lambda$*,*

(c) $\mathcal{F}_\bullet$ *is homologically* $\mathbb{Z}$*-cellular,*

(d) $\mathcal{F}_\bullet$ *is homologically* $\Lambda$*-cellular for all rings* $\Lambda$*.*

*Proof.* Trivially (b) implies (a) and (d) implies (c). The universal coefficient theorems for homology and cohomology show that property (c) implies the others. It hence suffices to show that (a) implies (c).

Assume that $\mathcal{F}_\bullet$ is cohomologically $\mathbb{Z}$-cellular and let $i, j \in \mathbb{Z}$. As singular homology of varieties, all the abelian groups $H_j^{\mathrm{sing}}(X_i^{\mathrm{an}}, X_{i-1}^{\mathrm{an}})$ are finitely generated. Thus there exist non-canonical decompositions

$$H_j^{\mathrm{sing}}(X_i^{\mathrm{an}}, X_{i-1}^{\mathrm{an}}) \cong F_{i,j} \oplus T_{i,j}$$

into a free and a torsion part. We furthermore have non-canonical isomorphisms

$$\mathrm{Hom}(H_j^{\mathrm{sing}}(X_i^{\mathrm{an}}, X_{i-1}^{\mathrm{an}}), \mathbb{Z}) \cong F_{i,j}$$

and

$$\mathrm{Ext}^1(H_j^{\mathrm{sing}}(X_i^{\mathrm{an}}, X_{i-1}^{\mathrm{an}}), \mathbb{Z}) \cong T_{i,j}.$$

Therefore the universal coefficient theorem (cf. [Hat02], Theorem 3.2) for singular cohomology turns into the exactness of

$$0 \longrightarrow T_{i,j-1} \longrightarrow H^j_{\mathrm{sing}}(X_i, X_{i-1}) \longrightarrow F_{i,j} \longrightarrow 0.$$

The middle term is by assumption free for all $i, j$, thus all $T_{i,j-1}$ are trivial. Hence, again non-canonically, $H^j_{\mathrm{sing}}(X_i^{\mathrm{an}}, X_{i-1}^{\mathrm{an}}) \cong F_{i,j} \cong H_j^{\mathrm{sing}}(X_i, X_{i-1})$, proving that $\mathcal{F}_\bullet$ is homologically $\mathbb{Z}$-cellular. □

## 4.2 Cellular complexes

Recall Definition 1.3.5.

**Definition 4.2.1** (Cellular complex)**.** Let $X \in \mathrm{Var}_k$ be an affine variety with a $\Lambda$-cellular filtration $X_\bullet = \mathcal{F}_\bullet X$. Then the $(X_i, X_{i-1}, i)$ are very good pairs (cf. Definition 1.3.7) and we define *Nori's cellular complex* $\mathrm{C}^\bullet_{\mathcal{F}}(X) = \mathrm{C}^\bullet_{\mathcal{F}}(X, \Lambda)$ as the complex

$$\{\ldots \to H^i_{\mathrm{Nori}}(X_i, X_{i-1}, \Lambda) \to H^{i+1}_{\mathrm{Nori}}(X_{i+1}, X_i, \Lambda) \to \ldots\} \in \mathcal{MM}^{\mathrm{eff}}_{\mathrm{Nori}}(k, \Lambda).$$

Here we put $H^i_{\mathrm{Nori}}(X_i, X_{i-1}, \Lambda)$ in degree $i$ and let

$$H^i_{\mathrm{Nori}}(X_i, X_{i-1}, \Lambda) \to H^{i+1}_{\mathrm{Nori}}(X_{i+1}, X_i, \Lambda)$$

be the boundary morphisms of $\mathcal{MM}^{\mathrm{eff}}_{\mathrm{Nori}}(k, \Lambda)$ corresponding to the triple $(X_{i+1}, X_i, X_{i-1})$.

**Definition 4.2.2.** If $f\colon X \to Y$ is a morphism of affine varieties over $k$ and $\mathcal{F}$ is a $\Lambda$-cellular filtration on each of $X$ and $Y$ compatible with $f$, then the pullback edges of Definition 1.3.4 induce morphisms

$$\mathrm{C}^i_{\mathcal{F}}(f)\colon H^i_{\mathrm{Nori}}(\mathcal{F}_iY, \mathcal{F}_{i-1}Y, \Lambda) \to H^i_{\mathrm{Nori}}(\mathcal{F}_iX, \mathcal{F}_{i-1}X, \Lambda)$$

in $\mathcal{MM}^{\mathrm{eff}}_{\mathrm{Nori}}(k, \Lambda)$.

They assemble to a morphism

$$\mathrm{C}^\bullet_{\mathcal{F}}(f) = \mathrm{C}^\bullet_{\mathcal{F}}(f, \Lambda)\colon\ \mathrm{C}^\bullet_{\mathcal{F}}(Y, \Lambda) \to \mathrm{C}^\bullet_{\mathcal{F}}(X, \Lambda).$$

Note that $\mathrm{C}^\bullet_{\mathcal{F}}(X)$ is a complex and that $\mathrm{C}^\bullet_{\mathcal{F}}(f)$ is indeed a morphism of complexes. This can, for example, be seen by applying the faithful forgetful functor $\omega_{\mathrm{sing}}\colon \mathcal{MM}^{\mathrm{eff}}_{\mathrm{Nori}}(k, \Lambda) \to \Lambda\text{-}\mathrm{Mod}$.

**Proposition 4.2.3.** *Let $X \in \mathrm{Var}^{\mathrm{aff}}_k$ be an affine variety and let $\mathcal{F}$ be a $\Lambda$-cellular filtration on $X$.*

*Then $\mathrm{C}^\bullet_{\mathcal{F}}(X, \Lambda)$ naturally calculates, after application of the forgetful functor $\omega_{\mathrm{sing}}\colon \mathcal{MM}^{\mathrm{eff}}_{\mathrm{Nori}}(k, \Lambda) \to \Lambda\text{-}\mathrm{Mod}$, the singular cohomology of $X$: there is a natural isomorphism*

$$H^n_{\mathrm{sing}}(X, \Lambda) \cong \omega_{\mathrm{sing}}(H^n(\mathrm{C}^\bullet_{\mathcal{F}}(X, \Lambda))).$$

*Proof.* As $\omega_{\mathrm{sing}}$ is exact and faithful it is sufficient to show that the complex

$$\ldots \to H^i_{\mathrm{sing}}(\mathcal{F}_iX, \mathcal{F}_{i-1}X, \Lambda) \to H^{i+1}_{\mathrm{sing}}(\mathcal{F}_{i+1}X, \mathcal{F}_iX, \Lambda) \to \ldots$$

of boundary morphisms computes singular cohomology. This is shown as in the standard proofs that cellular cohomology computes singular cohomology of CW-complexes. If we use the shortcut $X_i := \mathcal{F}_i(X)$ this is for example the dual of [Hat02], Theorem 2.35. □

All maps in the construction of $\mathrm{C}^\bullet_{\mathcal{F}}(X)$ are the natural ones on singular cohomology. Thus the following is immediate from Proposition 4.2.3 and the faithfulness of $\omega_{\mathrm{sing}}\colon \mathcal{MM}^{\mathrm{eff}}_{\mathrm{Nori}}(k, \Lambda) \to \Lambda\text{-}\mathrm{Mod}$.

**Lemma 4.2.4.** *Let $k \subseteq \mathbb{C}$ be a field and let $\Lambda$ be a ring. Then:*

*(a) Let $X \in \mathrm{Var}_k$ be an affine variety. Let $\mathcal{F}$ and $\mathcal{G}$ be two $\Lambda$-cellular filtrations on $X$ with $\mathcal{F}$ finer than $\mathcal{G}$. Then the inclusion induces a natural quasi-isomorphism $\mathrm{C}^\bullet_{\mathcal{G}}(X) \to \mathrm{C}^\bullet_{\mathcal{F}}(X)$.*

*(b) Let $f\colon X \to Y$ be a morphism of affine varieties over $k$. Let $\mathcal{F}$ and $\mathcal{G}$ be two $\Lambda$-cellular filtrations on the corresponding diagram $X \to Y$, i.e. on both objects and compatible with $f$. Assume that $\mathcal{F}$ is finer than $\mathcal{G}$.*

*Then the diagram*

$$\begin{array}{ccc} \mathrm{C}^\bullet_{\mathcal{G}}(Y) & \xrightarrow{\mathrm{C}^\bullet_{\mathcal{G}}(f)} & \mathrm{C}^\bullet_{\mathcal{G}}(X) \\ \downarrow{\scriptstyle q.\text{-is.}} & & \downarrow{\scriptstyle q.\text{-is.}} \\ \mathrm{C}^\bullet_{\mathcal{F}}(Y) & \xrightarrow{\mathrm{C}^\bullet_{\mathcal{F}}(f)} & \mathrm{C}^\bullet_{\mathcal{F}}(X) \end{array}$$

*commutes.*

## 4.3 Nori's Basic Lemma

The existence of $\Lambda$-cellular filtrations on an affine variety is the content of Nori's Basic Lemma which appeared in [Nor02]:

**Theorem 4.3.1** (Nori's Basic Lemma, classical version)**.** *Let $X \in \mathrm{Var}_k^{\mathrm{aff}}$ be an affine variety of dimension $n$ over $k \subseteq \mathbb{C}$ and let $\Lambda$ be a ring. Let $Y \subset X$ be a closed subset of dimension at most $n-1$.*

*Then there exists a closed subset $Z \subset X$ containing $Y$ such that:*

- $H^i_{\mathrm{sing}}(X, Z, \Lambda) = 0$ *for all integers $i \neq n$,*
- $H^n_{\mathrm{sing}}(X, Z, \Lambda)$ *is a finitely generated free $\Lambda$-module.*

**Remark 4.3.2.** For $\Lambda = \mathbb{Z}$ this is Theorem 2.5.2 of [HM16]. The general version follows easily from Nori's more general sheaf-theoretic version of Theorem 4.3.8 below. We point to [HM16], after the statement of Theorem 2.5.7, for an explicit reduction.

**Remark 4.3.3.** We cannot expect Nori's Basic Lemma to hold for non-affine varieties. Indeed, it is wrong for any projective variety of positive dimension.

**Remark 4.3.4.** By enlarging $Y$ to contain the singular locus of $X$ we can also achieve that $X \backslash Z$ is smooth.

Recall that by Convention 4.0.2 all sheaves are in the analytic topology.

**Definition 4.3.5** (Stratification, constructible)**.** Let $\Lambda$ be a ring and let $X$ be a scheme.

A *stratification* of $X$ is a disjoint decomposition $\mathcal{S} = \{S_i | i \in I\}$ of the Zariski-topological space $X$ into finitely many locally closed subsets such that the *frontier condition*

$$\overline{S_i} = \bigcup_{S_j \in \overline{S_i}} S_j$$

holds.

A sheaf $F$ of $\Lambda$-modules on $X$ is called *weakly constructible* if there is a stratification $\mathcal{S} = \{S_i | i \in I\}$ such that the restrictions $F_{|S_i}$ are locally constant. It is then called *constructible* if all stalks of $F$ are finitely generated $\Lambda$-modules.

We list their most important properties in the form of two lemmas:

**Lemma 4.3.6.** *Constructibility and weak constructibility are stable under finite limits and finite colimits of sheaves. Explicitly this means for a variety $X$ over a field $k \subseteq \mathbb{C}$ and a finite diagram $D\colon \mathcal{C} \to \mathrm{Shv}(X)$ of (weakly) constructible sheaves on $X$:*

*(a)* $\lim_{c\in\mathcal{C}} D(c)$ *is (weakly) constructible,*

*(b)* $\operatorname{colim}_{c\in\mathcal{C}} D(c)$ *is (weakly) constructible.*

*Proof.* By taking intersections of the individual stratifications we may assume that $\mathcal{S} = \{S_i | i \in I\}$ is a stratification of $X$ such that each of the $D(c)$, $c \in \mathcal{C}$, is locally constant on each $S_i$. As pullbacks, restrictions in particular, commute with finite (co)limits, we can restrict to a single $S_i$. Thus we need to show that locally constant sheaves are stable under finite (co)limits. This property is clear for constant sheaves, to which it can be reduced by taking sufficiently fine open covers.

If additionally the stalks of all the $D(c)$, $c \in \mathcal{C}$, are finitely generated, then the same is clearly true for the (co)limit over the finite diagram $D$. □

**Lemma 4.3.7.** *Constructibility and weak constructibility are stable under arbitrary pullbacks, pushforwards along finite morphisms and the extension-by-0 functor. Explicitly this means for a morphism $f\colon X \to Y$ of varieties over a field $k \subseteq \mathbb{C}$:*

*(a) If $\mathcal{G}$ is a (weakly) constructible sheaf on $Y$, then $f^*(\mathcal{G})$ is a (weakly) constructible sheaf on $X$.*

*(b) If $f$ is finite and $\mathcal{F}$ is a (weakly) constructible sheaf on $X$, then $f_*(\mathcal{F})$ is a (weakly) constructible sheaf on $Y$.*

*(c) If $f$ is an open immersion and $\mathcal{F}$ is a (weakly) constructible sheaf on $X$, then $f_!(\mathcal{F})$ is a (weakly) constructible sheaf on $Y$.*

*Proof.* Part (c) is easy and together with part (b) given as Lemma 2.5.10 of [HM16]. Part (a) is almost trivial:

Let $\mathcal{S} = \{S_i | i \in I\}$ be a stratification of $Y$ such that the restrictions $\mathcal{G}_{|S_i}$ are locally constant. Then the restriction of $f^*(\mathcal{G})$ to the locally closed $f^{-1}(S_i)$ is the pullback $f^*(\mathcal{G}_{|S_i})$. Using an open cover of $S_i$ on which $S_i$ is constant then shows that $f^*(\mathcal{G}_{|S_i})$ is locally constant.

As taking stalks commutes with pullbacks, we conclude that constructibility is preserved as well. □

This allows us to formulate the following more precise version of Theorem 4.3.1, that was independently shown in [Nor02] and, with a slight twist, in [Beĭ87] as part of Proof 3.3.1. Both arguments are elaborated in Section 2.5 of [HM16].

**Theorem 4.3.8** (Nori's Basic Lemma, sheaf-theoretic version)**.** *Let $k \subseteq \mathbb{C}$ be a field and let $\Lambda$ be a ring. Let $X$ be an affine variety of dimension $n$ over $k$ and let $Y \subset X$ be closed subset of dimension at most $n-1$. Assume that $\mathcal{F}$ is a weakly constructible sheaf of $\Lambda$-modules on $X$.*

*Then there exists a closed subset $Z \subset X$ of dimension $n-1$, containing $Y$ and with open complement $j\colon U \hookrightarrow X$ such that:*

*(a) $H^m(X, j_! j^* \mathcal{F}) = 0$ for $m \neq n$,*

*(b) $H^n(X, j_! j^* \mathcal{F})$ is a finite sum of stalks $\mathcal{F}_x$.*

*Proof.* This is a variant of Theorem 2.5.7 of [HM16] mentioned in Subsection 2.5.2 of op. cit. at the beginning of Nori's proof of their Theorem 2.5.7. They work over $\Lambda = \mathbb{Z}$, but the proof as given works for an arbitrary ring because all morphisms used (i.e. pushforwards and pullbacks) restrict to $\Lambda$-modules. □

Two easy inductions (cf. Theorem 4.4.5, Lemma 4.5.5 and Theorem 4.8.4 for generalizations and explicit proofs) using Nori's Basic Lemma 4.3.1 show:

**Lemma 4.3.9.** *Let $k \subseteq \mathbb{C}$ be a field and let $\Lambda$ be a ring.*

*(a) Every affine variety $X \in \mathrm{Var}_k^{\mathrm{aff}}$ admits a $\Lambda$-cellular filtration coarser than a given standard filtration on $X$.*

*(b) Let $D\colon \mathcal{C} \to \mathrm{Var}_k^{\mathrm{aff}}$ be a finite acyclic diagram of affine varieties and let $\mathcal{F}$ be a standard filtration on $D$. Then $D$ admits a $\Lambda$-cellular filtration coarser than $\mathcal{F}$.*

## 4.4 The equivariant version of Nori's Basic Lemma

Recall that our sheaves were defined on the analytifications. The following lemma is surely well-known, but we have not been able to find a reference:

**Lemma 4.4.1.** *Let $\pi\colon X \to Y$ be a finite morphism of varieties over $k \subseteq \mathbb{C}$.*

*Then the pushforward $\pi_*\colon \mathrm{Shv}(X) \to \mathrm{Shv}(Y)$ is exact. Thus for any sheaf $\mathcal{F}$ on $X$ we have a canonical isomorphism $H^n(X, \mathcal{F}) \cong H^n(Y, \pi_* \mathcal{F})$.*

*Proof.* By the Leray spectral sequence it suffices to show that the right derived sheaves $(R^n \pi_*)\mathcal{F}$ vanish for $n > 0$. Looking at the stalk $X_y^{\mathrm{an}}$ over $y \in Y^{\mathrm{an}}$ the stalk-wise version of proper base change states that

$$((R^n \pi_*)\mathcal{F})_y \cong H^n(X_y^{\mathrm{an}}, \mathcal{F}_{|X_y^{\mathrm{an}}}).$$

But $\pi$ is quasi-finite, thus $X_y^{\mathrm{an}}$ is a finite union of points. Therefore all cohomology in degree $n > 0$ vanishes as required. □

**Lemma 4.4.2.** *Let $A$ be a finite set with a transitive action by a finite group $G$. Let $\Lambda$ be a ring.*

*Consider the morphism*

$$q\colon \Lambda^A \to \Lambda^{A\times G},$$
$$(\lambda_a)_{a\in A} \mapsto (\lambda_a - \lambda_{ga})_{a\in A,\ g\in G}$$

*and the $G$-action on $\Lambda^A$ induced by that on $A$.*

*Then the kernel* $\ker(q)$ *is the $\Lambda$-submodule of $G$-invariant elements of $\Lambda^A$ and is isomorphic to $\Lambda$ via the diagonal embedding $\Delta\colon \Lambda \to \Lambda^A$. The image* $\mathrm{im}(q)$ *and cokernel* $\mathrm{coker}(q)$ *are free $\Lambda$-modules.*

*Proof.* An element $(\lambda_a)_{a\in A} \in \Lambda^A$ is annihilated by $q$ if and only if $\lambda_a = \lambda_{ga}$ for all $a \in A$ and $g \in G$. In other words, if it is $G$-invariant. By transitivity this means that all $\lambda_a$ are equal.

Let us now deal with the image and the cokernel. We understand $\Lambda^A$ as the free $\Lambda$-module with basis $A$ and $\Lambda^{A\times G}$ as the free $\Lambda$-module with basis $A \times G$.

Fix an element $b \in A$ and for each $a \in A$ a $g_a \in G$ such that $g_a b = a$, possible by transitivity. We claim that a $\Lambda$-basis of $\mathrm{coker}(q)$ is represented by $B := (A \times G)\backslash\{(b, g_a) \mid a \in A,\ a \neq b\}$. We furthermore claim that a $\Lambda$-basis $C$ of $\mathrm{im}(q)$ is given by the $q(a)$, $a \in A\backslash\{b\}$.

The latter is easily verified: $\Lambda^A$ is freely generated by $\sum_{a\in A} a \in \ker(q) \cong \Lambda$ and the $a \in A\backslash\{b\}$, hence the image of $q$ is freely generated by $C$. It now suffices to check that $B \sqcup C$ forms a $\Lambda$-basis of $\Lambda^{A\times G}$.

Let $M$ be the $\Lambda$-submodule of $\Lambda^{A\times G}$ generated by $B$. Note that

$$q(a) = \sum_{g\in G}(a, g) - \sum_{g\in G}(g^{-1}a, g).$$

For $a \in A$, $a \neq b$, we find that $q(a) + (b, g_a) \in M$. Indeed:

The first sum for $q(a)$ contains only pairs $(a, g) \in B$. Let now $c \in A$ be such that the summand $(b, g_c)$ appears in the second sum. This means that $g_c^{-1}a = b$, i.e. $a = g_c b = c$. Hence $(b, g_a)$ is the only summand not in $B$.

This shows that $B \cup C$ is a generating system of $\Lambda^{A\times G}$. Observing $\#(B \cup C) \leq \#B + \#C = \#A \cdot \#G$ we conclude that $B \sqcup C = B \cup C$ is a basis. □

The next Theorems 4.4.3 and 4.4.5 mimic Nori's ideas on how to work with quotients by finite groups. As before we give a cohomological instead of his original homological version.

**Theorem 4.4.3** (Nori's Basic Lemma, equivariant version)**.** *Let $k \subseteq \mathbb{C}$ be a field and let $\Lambda$ be a ring. Let $X$ be an affine variety of dimension $n$ over $k$ and let $G$ be a finite group acting on $X$ via morphisms over $k$.*

*Then there exists a $G$-invariant closed subset $D \subset X$ of dimension $n-1$ such that:*

- $H^n_{\text{sing}}(X, D, \Lambda)$ *and* $H^n_{\text{sing}}(X/G, D/G)$ *are finitely generated projective* $\Lambda$*-modules,*
- $H^i_{\text{sing}}(X/G, D/G, \Lambda) = 0$ *for* $i \neq n$,
- $H^i_{\text{sing}}(X, D, \Lambda) = 0$ *for* $i \neq n$,
- *The projection* $\pi\colon X \to X/G$ *induces an isomorphism*

$$H^n_{\text{sing}}(X/G, D/G, \Lambda) \cong H^n_{\text{sing}}(X, D, \Lambda)^G.$$

*Additionally, $D$ can be chosen to contain an arbitrary closed subset of $X$ of dimension at most $n-1$.*

*Proof.* Let $K$ be an arbitrary closed subscheme of $X$ of dimension at most $n-1$. Then $\pi(K)$ has dimension at most $n-1$, hence there is a closed subset $D'$ of $X/G$ of dimension $n-1$ containing it. Let $U = (X/G)\backslash D'$ be its open complement and let $i\colon U \hookrightarrow X/G$ and $\widetilde{i}\colon \pi^{-1}(U) \to X$ be the inclusions.

By proper base change we have an isomorphism $\pi^* i_! \Lambda_U \cong \widetilde{i}_! \pi^* \Lambda_U \cong \widetilde{i}_! \Lambda_{\pi^{-1}(U)}$. Its adjoint gives the middle morphism of a sequence

$$0 \longrightarrow i_! \Lambda_U \longrightarrow \pi_* \widetilde{i}_! \Lambda_{\pi^{-1}(U)} \xrightarrow{\;q\;} \prod_{g \in G} \pi_* \widetilde{i}_! \Lambda_{\pi^{-1}(U)} \tag{4.1}$$

of sheaves on $X/G$. The last morphism $q$ is the product over the morphisms

$$\mathrm{id} - \pi_* \widetilde{i}_! g_* \colon \pi_* \widetilde{i}_! \Lambda_{\pi^{-1}(U)} \to \pi_* \widetilde{i}_! \Lambda_{\pi^{-1}(U)}$$

induced by the adjoint action $g_* \colon \Lambda_{\pi^{-1}(U)} \to \Lambda_{\pi^{-1}(U)}$ of $g \in G$ on the constant sheaf $\Lambda_{\pi^{-1}(U)} \cong g^* \Lambda_{\pi^{-1}(U)}$. This sequence is exact:

If $x \in D'(\mathbb{C})$, then all stalks vanish. If $x \in U(\mathbb{C})$, pick any preimage $y \in X(\mathbb{C})$. Then $\pi^{-1}(x)$ is the orbit $Gy$ (cf. Proposition 3.1.4 (a)). Thus, by the stalk-wise version of proper base change, the sequence at the stalk at $x$ reads as

$$0 \longrightarrow \Lambda \longrightarrow \Lambda^{Gy} \longrightarrow \Lambda^{Gy \times G},$$

where the middle morphism is the diagonal one and the last morphism sends $(\lambda_z)_{z \in Gy}$ to $(\lambda_z - \lambda_{gz})_{z \in Gy,\, g \in G}$. As the action of $G$ on the orbit $Gy$ is transitive we conclude the exactness from Lemma 4.4.2. It also tells us that the stalks of $\mathrm{im}(q)$ and $\mathrm{coker}(q)$ are finitely generated free $\Lambda$-modules.

We set

$$\mathcal{F} := i_!\Lambda_U \oplus \pi_*\widetilde{i}_!\Lambda_{\pi^{-1}(U)} \oplus \mathrm{im}(q) \oplus \mathrm{coker}(q).$$

Each of the summands is constructible by Lemmas 4.3.6 and 4.3.7, so by the sheaf-theoretic version of Nori's Basic Lemma 4.3.8 there is an open dense subset $j\colon V \hookrightarrow X/G$ contained in $U$ and such that $H^m(X/G, j_!j^*\mathcal{F})$ is a free $\Lambda$-module of finite rank for $m = n$ and 0 otherwise. We claim that the obviously $G$-invariant complement $D$ of $\pi^{-1}(V)$ satisfies the conditions of the lemma.

Let $\widetilde{j}\colon \pi^{-1}(V) \to X$ and $u\colon V \hookrightarrow U$ be the inclusions. We note that

$$j^*i_!\Lambda_U \cong u^*i^*i_!\Lambda_U \cong u^*\Lambda_U \cong \Lambda_V \tag{4.2}$$

and analogously $\widetilde{j}^*\widetilde{i}_!\Lambda_{\pi^{-1}(U)} \cong \Lambda_{\pi^{-1}(V)}$. The second isomorphism, $\pi_* = \pi_!$ and proper base change give us

$$j_!j^*\pi_*\widetilde{i}_!\Lambda_{\pi^{-1}(U)} \cong j_!\pi_*\widetilde{j}^*\widetilde{i}_!\Lambda_{\pi^{-1}(U)} \cong \pi_*\widetilde{j}_!\Lambda_{\pi^{-1}(V)}. \tag{4.3}$$

On cohomology we therefore get

$$H^m(X/G, j_!j^*i_!\Lambda_U) \overset{(4.2)}{\cong} H^m(X/G, j_!\Lambda_V) \cong H^m_{\mathrm{sing}}(X/G, D/G, \Lambda)$$

as well as

$$\begin{aligned} H^m(X/G, j_!j^*\pi_*\widetilde{i}_!\Lambda_{\pi^{-1}(U)}) &\overset{(4.3)}{\cong} H^m(X/G, \pi_*\widetilde{j}_!\Lambda_{\pi^{-1}(V)}) \cong \\ &\cong H^m(X, \widetilde{j}_!\Lambda_{\pi^{-1}(V)}) \cong H^m_{\mathrm{sing}}(X, D, \Lambda), \end{aligned}$$

the second isomorphism by Lemma 4.4.1.

Both cohomologies are summands of $H^m(X/G, j_!j^*\mathcal{F})$, therefore trivial if $m \neq n$ and projective if $m = n$. As summands, in particular quotients, of a finitely generated $\Lambda$-module we find them to be finitely generated as well. Thus we have shown everything except the last point.

By definition we can split (4.1) into two short exact sequences

$$0 \longrightarrow i_!\Lambda_U \longrightarrow \pi_*\widetilde{i}_!\Lambda_{\pi^{-1}(U)} \longrightarrow \mathrm{im}(q) \longrightarrow 0,$$

$$0 \longrightarrow \mathrm{im}(q) \longrightarrow \prod_{g\in G} \pi_*\widetilde{i}_!\Lambda_{\pi^{-1}(U)} \longrightarrow \mathrm{coker}(q) \longrightarrow 0.$$

Now we apply $H^n(X/G, j_!j^*-)$ to them. The long exact sequence of sheaf cohomology together with the vanishing of $H^m(X/G, j_!j^*-)$ for $m \neq n$ and the above isomorphisms turn them into short exact sequences

$$0 \to H^m_{\mathrm{sing}}(X/G, D/G, \Lambda) \to H^m_{\mathrm{sing}}(X, D, \Lambda) \to H^n(X/G, j_!j^*\,\mathrm{im}(q)) \to 0,$$

$$0 \to H^m(X/G, j_! j^* \operatorname{im}(q)) \to \prod_{g \in G} H^m_{\text{sing}}(X, D, \Lambda) \to H^n(X/G, j_! j^* \operatorname{coker}(q)) \to 0.$$

Recombining them shows that $H^m_{\text{sing}}(X/G, D/G, \Lambda)$ is the kernel of

$$H^m_{\text{sing}}(X, D, \Lambda) \longrightarrow \prod_{g \in G} H^m_{\text{sing}}(X, D, \Lambda)$$

induced by the action of $G$ on $X$, i.e. that indeed

$$H^n_{\text{sing}}(X/G, D/G, \Lambda) \cong H^n_{\text{sing}}(X, D)^G.$$

□

**Definition 4.4.4.** Let $X$ be an affine variety with an action of a finite group $G$. A standard filtration $\mathcal{F}$ on $X$ is called *$G$-invariant* if each $\mathcal{F}_i X$ is. If this is satisfied we get a *quotient filtration* $\mathcal{F}/G$ on $X/G$ defined by $(\mathcal{F}/G)_i(X/G) = (\mathcal{F}_i X)/G$.

Note that $X$ being affine assures the existence of the quotients as varieties (cf. Remark 3.1.5).

**Theorem 4.4.5.** *Let $G$ be a finite group acting on an affine variety $X$ over $k \subseteq \mathbb{C}$ via morphisms over $k$. Let $\pi\colon X \to X/G$ be the corresponding quotient morphism and let $\mathcal{G}$ be a standard filtration on $X$ (cf. Definition 4.1.1). Let $\Lambda$ be a noetherian ring.*

*Then there exists a standard filtration $\mathcal{F}$ on $X$ coarser than $\mathcal{G}$ such that:*

- *$\mathcal{F}$ is $G$-invariant,*
- *$\mathcal{F}$ is $\Lambda$-cellular,*
- *the quotient filtration $\mathcal{F}/G$ on $X/G$ is $\Lambda$-cellular,*
- *the quotient morphism $\pi\colon X \twoheadrightarrow X/G$ induces an isomorphism*

$$\mathrm{C}^\bullet_{\mathcal{F}/G}(X/G) \cong \mathrm{C}^\bullet_{\mathcal{F}}(X)^G$$

*of complexes over $\mathcal{MM}^{\text{eff}}_{\text{Nori}}(k, \Lambda)$.*

*Proof.* Let $X_n = X$. By using Theorem 4.4.3 we can inductively choose closed subsets $X_{n-1}, X_{n-2}, \ldots, X_0, X_{-1} = \emptyset$ of $X$ such that:

- $X_j \supseteq \mathcal{G}_j(X)$ is invariant under $G$ and has dimension $j$
- $\mathcal{F}_i(X) := X_i$ and $(\mathcal{F}/G)_i(X/G) = X_i/G$ are $\Lambda$-cellular filtrations,
- $H^i_{\text{sing}}(X_i/G, X_{i-1}/G) \cong H^i_{\text{sing}}(X_i, X_{i-1})^G$.

The filtrations on $X$ and $X/G$ are in particular compatible with $\pi$. The natural map $\mathrm{C}^i_{\mathcal{F}/G}(X/G) \to \mathrm{C}^i_{\mathcal{F}}(X)$ becomes by construction an isomorphism onto $\mathrm{C}^i_{\mathcal{F}}(X)^G$ when considered in singular cohomology, i.e. after applying the faithful exact forgetful functor $\omega_{\mathrm{sing}} \colon \mathcal{MM}^{\mathrm{eff}}_{\mathrm{Nori}}(k,\Lambda) \to \Lambda\text{-}\mathrm{Mod}$. Thus $\pi$ induces an isomorphism

$$\mathrm{C}^\bullet_{\mathcal{F}/G}(X/G) \to \mathrm{C}^\bullet_{\mathcal{F}}(X)^G.$$

□

**Definition 4.4.6.** Let $D$ be a diagram of affine varieties over $k \subseteq \mathbb{C}$ and let $\pi \colon X \to X/G$ be a morphism in $D$ that is a quotient by the action of a finite group. We call a $\Lambda$-cellular filtration $\mathcal{F}_\bullet$ on $D$ *nice at* $\pi$ if $\mathcal{F}_\bullet(X)$ is $G$-invariant and $\pi$ induces an isomorphism $\mathrm{C}^\bullet_{\mathcal{F}}(X/G) \cong \mathrm{C}^\bullet_{\mathcal{F}}(X)^G$ of complexes over $\mathcal{MM}^{\mathrm{eff}}_{\mathrm{Nori}}(k,\Lambda)$.

It is easy to find $\Lambda$-cellular filtrations on finite acyclic diagrams of affine schemes. But if we want them to be nice at group quotients $X \to X/G$ we require a more technical obstruction than acyclicity. The reason is that the filtrations of $X$ and $X/G$ need to be chosen at the same time.

**Lemma 4.4.7.** *Let $D \colon \mathcal{C} \to \mathrm{Var}^{\mathrm{aff}}_k$ be a finite acyclic diagram of affine varieties over $k \subseteq \mathbb{C}$ and let $\Lambda$ be a noetherian ring. Let $S$ be a set of arrows in $\mathcal{C}$ that map to finite group quotients $X \to X/G$. We assume for all arrows $s \colon x \to y$ in $S$:*

- *no other arrow in $S$ has $x$ or $y$ as one if its endpoints,*
- *there is no other arrow $x \to y$ in $\mathcal{C}$,*
- *there is no third object $z \neq x, y$ in $\mathcal{C}$ such that there are arrows $x \to z$ and $z \to y$ in $\mathcal{C}$.*

*Then there exists a $\Lambda$-cellular filtration $\mathcal{F}$ on $D$ coarser than a given standard filtration on $D$ and nice at every arrow of $S$. If furthermore we are already given such a $\Lambda$-cellular filtration on an initial segment $I$ of $D$ that for each $s \in S$ contains either none or both of its endpoints, then $\mathcal{F}$ can be chosen to extend it.*

*Proof.* By the finiteness and acyclicity of $D$ there is a terminal object $y \in \mathcal{C}$, i.e. one that has no arrow into any other object of $\mathcal{C}$. There is then at most one object $x \in \mathcal{C}$ mapping to $y$ via a morphism $\pi \in S$; we set $x = y$ and $\pi = \mathrm{id}_y$ if no such object and morphism exist. Then by assumption the subdiagram on $\mathcal{C}\backslash\{x, y\}$ obtained by removing $x$, $y$ and all involved arrows is an initial segment of $\mathcal{C}$.

Let $\mathcal{G}$ be a standard filtration on $D$. By induction we may assume that we have found a $\Lambda$-cellular filtration $\mathcal{F}$ coarser than $\mathcal{G}$ on the initial segment $I = \mathcal{C}\backslash\{x,y\} \to \mathrm{Var}_k^{\mathrm{aff}}$. Let $f$ be an arrow from any $z \in \mathcal{C}\backslash\{x,y\}$ to $x$. As our varieties are of finite type over a field we get

$$\dim\left(\overline{f(\mathcal{F}_i D(z))}\right) \leq \dim\left(\mathcal{F}_i D(z)\right) \leq i$$

and thus may assume that $\overline{f(\mathcal{F}_i D(z))} \subseteq \mathcal{G}_i D(x)$. We do this for all the finitely many such $f$ and analogously for all arrows with target $y$.

Because $D(\pi)$ is a finite morphism, we have

$$\dim\left(D(\pi)^{-1}(\mathcal{G}_i D(y))\right) \leq \dim\left(\mathcal{G}_i D(y)\right) \leq i$$

and can therefore also assume that $D(\pi)^{-1}(\mathcal{G}_i D(y)) \subseteq \mathcal{G}_i D(x)$. We apply Theorem 4.4.5 to the finite group quotient $D(\pi)\colon D(x) \to D(y)$ and the given standard filtration on $\mathcal{G}_\bullet D(x)$, thereby extending $\mathcal{F}$ to all of $D$. By construction it is coarser than $\mathcal{G}$, extends the given filtration on $I$, is a $\Lambda$-cellular filtration, compatible with all morphisms and nice at all elements of $S$. □

## 4.5 Cellular filtrations on finite correspondences

In order to work with finite correspondences it is not enough to use $\Lambda$-cellular filtrations on the objects. One has to use filtrations on the morphisms as well, a notion we now want to make precise:

**Definition 4.5.1** (Edges). Let $\alpha\colon X \rightsquigarrow Y$ be an effective finite correspondence of constant degree $n$ (cf. Definition 1.8.11) between smooth affine varieties over a field $k \subseteq \mathbb{C}$. We define the *(effective) edge* $\widetilde{E}(\alpha)$ as the diagram

$$X \xrightarrow{\mathrm{sym}_{\mathbb{N}}(\alpha)} \mathrm{S}^n(Y) \xleftarrow{\quad\pi\quad} Y^n \overset{\mathrm{pr}_i}{\rightrightarrows} Y.$$

Here, the left arrow is the symmetrization $\mathrm{sym}_{\mathbb{N}}(\alpha) = \mathrm{sym}_{\mathrm{Spec}(k)}(\alpha)$ (see Definition 3.7.1 and Theorem 3.7.5), the middle arrow is the quotient by the symmetric group $S_n$ acting via permutation of the factors and the rightmost morphisms are the different projections onto the factors.

Assume now that $\alpha = \sum_{i=1}^r a_i\alpha_i\colon X \rightsquigarrow Y$ is a finite correspondence in $\mathrm{SmCor}^{\mathrm{aff}}(k,\Lambda)$ with $a_i \in \Lambda$ and basic finite correspondences $\alpha_i\colon X_i \rightsquigarrow Y_i$, where $X_i$ and $Y_i$ are the connected components of $X$ and $Y$, respectively, over which $\alpha_i$ is defined (cf. Proposition 1.6.4). We define the *edge* $E(\alpha)$ as the disjoint union $\coprod_{i=1}^r \widetilde{E}(\alpha_i)$ of the edges of the individual $\alpha_i$.

By abuse of notation we will call a $\Lambda$-cellular filtration on an edge $\widetilde{E}(\alpha)$ or $E(\alpha)$ *nice* if it is nice at the finite group quotients $Y^n \to \mathrm{S}^n(Y)$.

**Remark 4.5.2.** There is some slight overlap between the two definitions of edges: every effective finite correspondence $\alpha\colon X \rightsquigarrow Y$ can be understood as a finite $\Lambda$-correspondence $\alpha \otimes \Lambda$, amounting to a functor

$$- \otimes \Lambda\colon\ \mathrm{SmCor}^{\mathrm{eff}}(k) \to \mathrm{SmCor}(k, \Lambda).$$

But the edges will not be the same, as the effective version $\widetilde{E}$ creates a single connected diagram while the general version $E$ creates one for each basic summand. In theory, this is never a problem if one does not suppress the functor $- \otimes \Lambda$. We prefer to be cautious and choose to have different symbols $E$ and $\widetilde{E}$.

**Definition 4.5.3** (Filtrations on finite correspondences)**.** Let $k \subseteq \mathbb{C}$ be a field and let $\alpha = \sum_{i=1}^{r} a_i \alpha_i \colon X \rightsquigarrow Y$ be a finite correspondence in $\mathrm{SmCor}^{\mathrm{aff}}(k, \Lambda)$ as in Definition 4.5.1.

- A *standard filtration on* $\alpha$ is a standard filtration $\mathcal{F}$ on the edge $E(\alpha)$. We call it a $\Lambda$-*cellular filtration on* $\alpha$ if it is a nice $\Lambda$-cellular filtration. The notions of being *finer* or *coarser* extend in the obvious way.
- Let $\mathcal{F}_\bullet X$, $\mathcal{G}_\bullet Y$ and $\mathcal{H}_\bullet \alpha$ be standard filtrations on the two objects and the finite correspondence. We call them *compatible* if
  - the filtration $\mathcal{H}_\bullet \alpha$ on $E(\alpha)$ is compatible at each morphism in $E(\alpha)$,
  - the filtrations $\mathcal{F}_\bullet$ and $\mathcal{H}_\bullet$ agree on the connected components of $X$ appearing in the individual edges $\widetilde{E}(\alpha_i)$,
  - the filtrations $\mathcal{G}_\bullet$ and $\mathcal{H}_\bullet$ agree on the connected components of $Y$ appearing in the individual edges $\widetilde{E}(\alpha_i)$.
- A *standard filtration on a diagram in* $\mathrm{SmCor}(k, \Lambda)$ consists of a standard filtration on each object and each finite correspondence. The notions of being *finer* or *coarser* transfer accordingly.
- A $\Lambda$-*cellular filtration on a diagram* $D$ *in* $\mathrm{SmCor}(k, \Lambda)$ is a standard filtration on $D$ such that the standard filtrations on the objects and the finite correspondences are $\Lambda$-cellular as well as compatible with each other.

**Remark 4.5.4.** There is the embedding $\mathrm{Sm}_k \hookrightarrow \mathrm{SmCor}(k, \Lambda)$ of Definition 1.6.9, mapping a morphism $f\colon X \to Y$ to its graph $\Gamma_f$. In this sense every standard (respectively $\Lambda$-cellular) filtration on a diagram of smooth affine varieties induces a standard (respectively $\Lambda$-cellular) filtration on the induced diagram in $\mathrm{SmCor}^{\mathrm{aff}}(k, \Lambda)$ via the identification $Y \cong Y^1 \cong \mathrm{S}^1(Y)$. This obviously preserves coarseness, thus our constructions are fully compatible with those obtained in [HM16] for morphisms of smooth affine varieties.

Once more we find enough $\Lambda$-cellular filtrations:

**Lemma 4.5.5.** *Let $k \subseteq \mathbb{C}$ be a field and let $\Lambda$ be a noetherian ring. Let $D\colon \mathcal{C} \to \mathrm{SmCor}^{\mathrm{aff}}(k, \Lambda)$ be a finite acyclic diagram.*

*Then there exists a $\Lambda$-cellular filtration $\mathcal{F}$ on $D$ coarser than a given standard filtration on $D$. In particular they form a directed system. If we are already given such a $\Lambda$-cellular filtration on an initial segment $I$ of $D$, then $\mathcal{F}$ can be chosen to extend it.*

*Proof.* As in the proof of Lemma 4.4.7 it suffices to consider the case where $I$ contains all but one terminal object $D(z)$ of $D$. For each arrow $\alpha\colon y \to z$ in $\mathcal{C}$ we find, using Lemma 4.4.7, a $\Lambda$-cellular filtration on $D(\alpha)$ which extends the given one on $D(y)$. In the end we also pick a $\Lambda$-cellular filtration on $D(z)$ coarser than those on the individual $D(\alpha)$ mapping to $D(z)$. ☐

## 4.6 Mapping finite correspondences to Nori motives

Let us state the following trivial observation:

**Lemma 4.6.1.** *Let $\mathcal{A}$ be an abelian category and let $G$ be a finite group acting on $B \in \mathcal{A}$. Let $f\colon A \to B$ be a morphism in $\mathcal{A}$.*

*Then the orbit sum*

$$\sum_{f' \in G \circ f} f'\colon A \to B$$

*factors over the $G$-invariants*

$$B^G := \ker\left(\prod_{g\in G} (\mathrm{id}_B - g)\colon B \to \prod_{g \in G} B\right).$$

**Definition 4.6.2.** Let $\alpha\colon X \rightsquigarrow Y$ be an effective finite correspondence of constant degree $n$ between smooth affine varieties over a field $k \subseteq \mathbb{C}$ and let $\mathcal{F}$ be a nice $\Lambda$-cellular filtration on the edge $\widetilde{E}(\alpha)$. Let

$$\Sigma_Y^n := \sum_{i=1}^n \mathrm{C}^\bullet_{\mathcal{F}}(\mathrm{pr}_i)\colon\ \mathrm{C}^\bullet_{\mathcal{F}}(Y) \stackrel{\Sigma}{\to} \mathrm{C}^\bullet_{\mathcal{F}}(Y^n)$$

be the sum over the $\mathrm{C}^\bullet_{\mathcal{F}}(-)$ (cf. Definition 4.2.1) of the projections $\mathrm{pr}_i$ onto the $i$-th factor.

By Lemma 4.6.1 the image $\mathrm{im}(\Sigma_Y^n)$ is contained in $\mathrm{C}^\bullet_{\mathcal{F}}(Y^n)^{S_n}$. Thus we can define a morphism

$$\widetilde{\mathrm{C}}^\bullet_{\mathcal{F}}(\alpha)\colon\ \mathrm{C}^\bullet_{\mathcal{F}}(Y) \to \mathrm{C}^\bullet_{\mathcal{F}}(X)$$

between the respective cellular complexes as the composition

$$\mathrm{C}^\bullet_{\mathcal{F}}(Y) \xrightarrow{\Sigma^n_Y} \mathrm{C}^\bullet_{\mathcal{F}}(Y^n)^{S_n} \cong \mathrm{C}^\bullet_{\mathcal{F}}(\mathrm{S}^n(Y)) \xrightarrow{\mathrm{C}^\bullet_{\mathcal{F}}(\mathrm{sym}_{\mathbb{N}}(\alpha))} \mathrm{C}^\bullet_{\mathcal{F}}(X)$$

along the edge $\widetilde{E}(\alpha)$, the isomorphism due to the Definition 4.4.6 of niceness.

Assume now that $\alpha = \sum_{i=1}^r a_i\alpha_i \colon X \rightsquigarrow Y$ is a finite correspondence in $\mathrm{SmCor}^{\mathrm{aff}}(k,\Lambda)$ with basic finite correspondences $\alpha_i$ and coefficients $a_i \in \Lambda$. Let $\mathcal{F}$ be a $\Lambda$-cellular filtration on $E(\alpha)$, thus satisfying the above assumption on the individual effective edges $\widetilde{E}(\alpha_i)$. We define the *transfer map*

$$\mathrm{C}^\bullet_{\mathcal{F}}(\alpha) := \sum_{i=1}^r a_i \widetilde{\mathrm{C}}^\bullet_{\mathcal{F}}(\alpha_i) \colon \mathrm{C}^\bullet_{\mathcal{F}}(Y) \to \mathrm{C}^\bullet_{\mathcal{F}}(X)$$

by $\Lambda$-linear continuation from the basic finite correspondences.

By Definition 4.5.3 a refinement is one on each object in all relevant edges. Thus Lemma 4.2.4 implies the following generalization of itself:

**Lemma 4.6.3.** *Let $k \subseteq \mathbb{C}$ be a field and let $\Lambda$ be a noetherian ring. Let $\alpha \colon X \rightsquigarrow Y$ be a finite correspondence in $\mathrm{SmCor}^{\mathrm{aff}}(k,\Lambda)$. Let $\mathcal{F}$ and $\mathcal{G}$ be two $\Lambda$-cellular filtrations on this diagram such that $\mathcal{F}$ is finer than $\mathcal{G}$.*

*Then the diagram*

$$\begin{array}{ccc}
\mathrm{C}^\bullet_{\mathcal{G}}(Y) & \xrightarrow{\mathrm{C}^\bullet_{\mathcal{G}}(\alpha)} & \mathrm{C}^\bullet_{\mathcal{G}}(X) \\
\downarrow{\scriptstyle \text{q.}-\text{is.}} & & \downarrow{\scriptstyle \text{q.}-\text{is.}} \\
\mathrm{C}^\bullet_{\mathcal{F}}(Y) & \xrightarrow{\mathrm{C}^\bullet_{\mathcal{F}}(\alpha)} & \mathrm{C}^\bullet_{\mathcal{F}}(X)
\end{array}$$

*commutes.*

## 4.7 Additivity of $\mathrm{C}^\bullet_{\mathcal{F}}$

Note that $\widetilde{\mathrm{C}}^\bullet_{\mathcal{F}}(\alpha)$ is defined for all effective finite correspondences, not just basic ones. But this gives nothing new:

**Proposition 4.7.1.** *Let $\Lambda$ be a noetherian ring. Let $\alpha_1, \alpha_2 \colon X \rightsquigarrow Y$ be effective finite correspondences of constant degree between smooth affine varieties over a field $k \subseteq \mathbb{C}$.*

*Then for every sufficiently coarse $\Lambda$-cellular filtration $\mathcal{F}$ we have*

$$\widetilde{\mathrm{C}}^\bullet_{\mathcal{F}}(\alpha_1 + \alpha_2) = \widetilde{\mathrm{C}}^\bullet_{\mathcal{F}}(\alpha_1) + \widetilde{\mathrm{C}}^\bullet_{\mathcal{F}}(\alpha_2) \tag{4.4}$$

*as morphisms $\mathrm{C}^\bullet_{\mathcal{F}}(Y) \to \mathrm{C}^\bullet_{\mathcal{F}}(X)$ between Nori's cellular complexes.*

*Precisely this means: for every given $\Lambda$-cellular filtration $\mathcal{H}$ on $X$ there exists a nice $\Lambda$-cellular filtration $\mathcal{G}$ on the effective edges $\widetilde{E}(\alpha_1)$, $\widetilde{E}(\alpha_2)$ and $\widetilde{E}(\alpha_1+\alpha_2)$ such that*

- *the filtrations on the edges agree with $\mathcal{H}$ on $X$,*
- *the filtrations on the edges agree on $Y$,*
- *if $\mathcal{F}$ is a nice $\Lambda$-cellular filtration on the edges, coarser than $\mathcal{G}$ and satisfying the previous two properties, then equation* (4.4) *holds.*

*Proof.* Let $n_1$ and $n_2$ be the degrees of $\alpha_1$ and $\alpha_2$, respectively. The morphism

$$\mathrm{sym}_{\mathbb{N}}(\alpha_1+\alpha_2)\colon X \to \mathrm{S}^{n_1+n_2}(Y)$$

factors by Definitions 3.5.3 and 3.5.4 through

$$\alpha = (\mathrm{sym}_{\mathbb{N}}(\alpha_1), \mathrm{sym}_{\mathbb{N}}(\alpha_2))\colon X \to \mathrm{S}^{n_1}(Y) \times \mathrm{S}^{n_2}(Y)$$

via the morphism

$$\sigma_{n_1,n_2}\colon \mathrm{S}^{n_1}(Y) \times \mathrm{S}^{n_2}(Y) \to \mathrm{S}^{n_1+n_2}(Y)$$

of Proposition 3.3.4. The other two symmetrization morphisms $\mathrm{sym}_{\mathbb{N}}(\alpha_1)$ and $\mathrm{sym}_{\mathbb{N}}(\alpha_2)$ also factor through $\alpha$, this time via the projections $\mathrm{p}_1$ and $\mathrm{p}_2$ onto the factors. Thus we get a commutative diagram

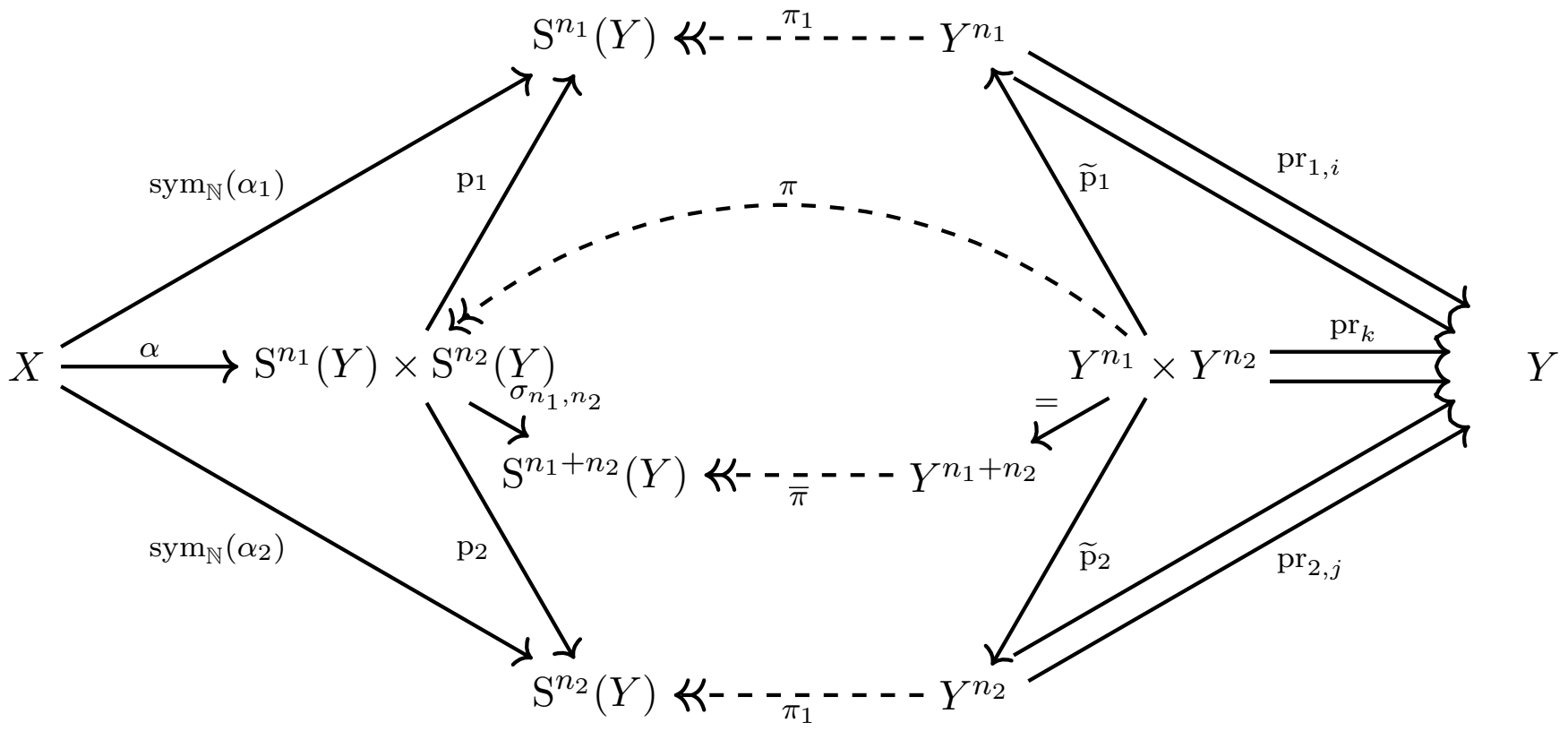

where we omitted to draw $\mathrm{sym}_{\mathbb{N}}(\alpha_1+\alpha_2)$ and the projections $\overline{\mathrm{pr}}_k\colon Y^{n_1+n_2} \to Y$ for clarity. Invoking Lemma 4.4.7 gives us a $\Lambda$-cellular filtration $\mathcal{G}$ on this diagram which is nice at the four dotted arrows and is the given one on $X$.

For a morphism $f$ we now use the shorthand $f^* := \mathrm{C}^\bullet_{\mathcal{G}}(f)$.

We note that $\pi_1^*$, $\pi_2^*$, $\pi^*$ and $\widetilde{\pi}^*$ are, by niceness of $\mathcal{G}$ at the dotted arrows, isomorphisms onto their images, which are certain invariant complexes. We will freely use their inverses $(\pi_1^*)^{-1}$, $(\pi_2^*)^{-1}$, $(\pi^*)^{-1}$ and $(\widetilde{\pi}^*)^{-1}$ as morphisms starting at those images. In all cases where they appear it is easily verified, using Lemma 4.6.1, that we are inside the respective images.

Then the following calculation shows that the lemma is true for $\mathcal{G}$:

$$\begin{aligned}
&\widetilde{\mathrm{C}}^{\bullet}_{\mathcal{G}}(\alpha_1) + \widetilde{\mathrm{C}}^{\bullet}_{\mathcal{G}}(\alpha_2) = \\
&= \mathrm{sym}_{\mathbb{N}}(\alpha_1)^*(\pi_1^*)^{-1} \sum_{i=1}^{n_1} \mathrm{pr}_{1,i}^* + \mathrm{sym}_{\mathbb{N}}(\alpha_2)^*(\pi_2^*)^{-1} \sum_{j=1}^{n_2} \mathrm{pr}_{2,j}^* = \\
&= \alpha^* \mathrm{p}_1^*(\pi_1^*)^{-1} \sum_{i=1}^{n_1} \mathrm{pr}_{1,i}^* + \alpha^* \mathrm{p}_2^*(\pi_2^*)^{-1} \sum_{j=1}^{n_2} \mathrm{pr}_{2,j}^* = \\
&= \alpha^*(\pi^*)^{-1} \widetilde{\mathrm{p}}_1^* \sum_{i=1}^{n_1} \mathrm{pr}_{1,i}^* + \alpha^*(\pi^*)^{-1} \widetilde{\mathrm{p}}_2^* \sum_{j=1}^{n_2} \mathrm{pr}_{2,j}^* = \\
&= \alpha^*(\pi^*)^{-1} \left( \sum_{i=1}^{n_1} \widetilde{\mathrm{p}}_1^* \, \mathrm{pr}_{1,i}^* + \sum_{j=1}^{n_2} \widetilde{\mathrm{p}}_2^* \, \mathrm{pr}_{2,j}^* \right) = \\
&= \alpha^*(\pi^*)^{-1} \sum_{k=1}^{n_1+n_2} \mathrm{pr}_k^* = \alpha^* \sigma_{n_1,n_2}^* (\overline{\pi})^{-1} \sum_{k=1}^{n_1+n_2} \overline{\mathrm{pr}}_k^* = \\
&= \mathrm{sym}_{\mathbb{N}}(\alpha_1 + \alpha_2)^* (\overline{\pi})^{-1} \sum_{k=1}^{n_1+n_2} \overline{\mathrm{pr}}_k^* = \widetilde{\mathrm{C}}^{\bullet}_{\mathcal{G}}(\alpha_1 + \alpha_2)
\end{aligned}$$

Now consider a general $\mathcal{F}$ as in the precise formulation. Lemma 4.6.3 gives us a factorization $\mathrm{C}^{\bullet}_{\mathcal{F}}(Y) \to \mathrm{C}^{\bullet}_{\mathcal{G}}(Y) \to \mathrm{C}^{\bullet}_{\mathcal{G}}(X) = \mathrm{C}^{\bullet}_{\mathcal{F}}(X)$ for each of the three maps $\widetilde{\mathrm{C}}^{\bullet}_{\mathcal{F}}(\alpha_1)$, $\widetilde{\mathrm{C}}^{\bullet}_{\mathcal{F}}(\alpha_2)$ and $\widetilde{\mathrm{C}}^{\bullet}_{\mathcal{F}}(\alpha_1 + \alpha_2)$, proving the lemma. □

**Proposition 4.7.2.** *Let $\alpha\colon X \rightsquigarrow Y$ be an effective finite correspondence of constant degree between smooth affine varieties over a field $k \subseteq \mathbb{C}$. Let $\Lambda$ be a noetherian ring.*

*Then for every sufficiently coarse $\Lambda$-cellular filtration $\mathcal{F}$ we have the equality*

$$\mathrm{C}^{\bullet}_{\mathcal{F}}(\alpha \otimes \Lambda) = \widetilde{\mathrm{C}}^{\bullet}_{\mathcal{F}}(\alpha) \tag{4.5}$$

*of morphisms $\mathrm{C}^{\bullet}_{\mathcal{F}}(Y) \to \mathrm{C}^{\bullet}_{\mathcal{F}}(X)$ between Nori's cellular complexes.*

*The precise meaning of this is as follows:*

*Given a $\Lambda$-cellular filtration $\mathcal{H}$ on $X$ there exists a nice $\Lambda$-cellular filtration $\mathcal{G}$ on the edges $\widetilde{E}(\alpha)$ and $E(\alpha \otimes \Lambda)$ such that:*

- *the filtrations on the edges agree with $\mathcal{H}$ on $X$,*
- *the filtrations on the edges agree on $Y$,*
- *if $\mathcal{F}$ is a nice $\Lambda$-cellular filtration on the two edges, coarser than $\mathcal{G}$ and satisfying the previous two properties, then equation (4.5) holds.*

*Proof.* Let $X = \coprod_{i\in I} X_i$ and $Y = \coprod_{j\in J} Y_j$ be the decompositions into connected components and let $\alpha_{i,j}\colon X_i \rightsquigarrow Y_j$ be the resulting effective finite correspondences as in Proposition 1.6.4.

Note a finite correspondence over the connected, thus by smoothness irreducible, $X_i$ is automatically of constant degree. Therefore a simple induction on the number of basic summands based on Proposition 4.7.1 gives us a filtration $\mathcal{G}_{i,j}$ on the respective edges, agreeing with $\mathcal{H}$ on $X_i$, such that

$$\mathrm{C}^{\bullet}_{\mathcal{G}_{i,j}}(\alpha_{i,j}\otimes\Lambda) = \widetilde{\mathrm{C}}^{\bullet}_{\mathcal{G}_{i,j}}(\alpha_{i,j}).$$

Let us write $\alpha'_{i,j}$ for $\alpha_{i,j}$ interpreted as a finite correspondence $X_i \rightsquigarrow Y$ via postcomposition with the embedding $Y_j \hookrightarrow Y$. This does clearly not change the degree. Note that the objects $Y_j$, $Y_j^{\deg(\alpha_{i,j})}$ and $\mathrm{S}^{\deg(\alpha_{i,j})}(Y_j)$ occurring in the effective edge $\widetilde{E}(\alpha_{i,j})$ consist of unions of connected components of the objects $Y$, $Y^{\deg(\alpha_{i,j})}$ and $\mathrm{S}^{\deg(\alpha_{i,j})}(Y)$ occurring in the effective edge $\widetilde{E}(\alpha'_{i,j})$, the third one by Proposition 3.3.5. We can thus extend $\mathcal{G}_{i,j}$ to $\widetilde{E}(\alpha'_{i,j})$ and get

$$\mathrm{C}^{\bullet}_{\mathcal{G}_{i,j}}(\alpha_{i,j}\otimes\Lambda) = \widetilde{\mathrm{C}}^{\bullet}_{\mathcal{G}_{i,j}}(\alpha'_{i,j}).$$

Another application of Proposition 4.7.1 shows that

$$\widetilde{\mathrm{C}}^{\bullet}_{\mathcal{G}_i}\left(\sum_{j\in J}\alpha'_{i,j}\right) = \sum_{j\in J}\widetilde{\mathrm{C}}^{\bullet}_{\mathcal{G}_i}(\alpha'_{i,j})$$

for a sufficiently coarse $\mathcal{G}_i$, still agreeing with $\mathcal{H}$ on $X_i$ and coarser than the individual $\mathcal{G}_{i,j}$.

Now we note that the degree $n$ of $\sum_{j\in J}\alpha_{i,j}$ does not depend on $i$ by assumption. We choose a common coarsening $\mathcal{G}$ of the $\mathcal{G}_i$ on $Y$, $Y^n$ and $\mathrm{S}^n(Y)$ with $\mathcal{G} = \mathcal{H}$ on $X$. This is by definition a nice $\Lambda$-cellular filtration on the edge $\widetilde{E}(\alpha)$ and satisfies

$$\widetilde{\mathrm{C}}^{\bullet}_{\mathcal{G}}\left(\sum_{\substack{i\in I\\ j\in J}}\alpha'_{i,j}\right) = \sum_{i\in I}\widetilde{\mathrm{C}}^{\bullet}_{\mathcal{G}}\left(\sum_{j\in J}\alpha'_{i,j}\right).$$

Combining the equations we find that indeed

$$\mathrm{C}^{\bullet}_{\mathcal{F}}(\alpha\otimes\Lambda) = \mathrm{C}^{\bullet}_{\mathcal{F}}\left(\sum_{\substack{i\in I\\ j\in J}}\alpha'_{i,j}\otimes\Lambda\right) = \sum_{\substack{i\in I\\ j\in J}}\mathrm{C}^{\bullet}_{\mathcal{F}}(\alpha'_{i,j}\otimes\Lambda) =$$

$$= \sum_{\substack{i\in I\\ j\in J}}\widetilde{\mathrm{C}}^{\bullet}_{\mathcal{F}}(\alpha'_{i,j}) = \sum_{i\in I}\widetilde{\mathrm{C}}^{\bullet}_{\mathcal{F}}\left(\sum_{j\in J}\alpha'_{i,j}\right) = \widetilde{\mathrm{C}}^{\bullet}_{\mathcal{F}}\left(\sum_{\substack{i\in I\\ j\in J}}\alpha'_{i,j}\right) = \widetilde{\mathrm{C}}^{\bullet}_{\mathcal{F}}(\alpha)$$

for every $\Lambda$-cellular filtration $\mathcal{F}$ coarser than the ones constructed above, yet agreeing with $\mathcal{H}$ on $X$. $\square$

**Remark 4.7.3.** Proposition 4.7.1 and Proposition 4.7.2 demonstrate the behaviour of $\Lambda$-cellular filtrations on finite correspondences: it is no longer enough to have sufficiently coarse filtrations on the objects, but we get constraints from the morphisms as well. This will become even more apparent in Theorems 4.8.1 and 4.9.6. From the viewpoint of ind-systems it means that certain properties are only eventually true.

**Remark 4.7.4.** Instead of using finite correspondences it should be possible to directly use multivalued morphisms as defined in Definition 3.4.10. This has the advantage of working over non-smooth varieties. But there is one technical subtlety one has to deal with, especially for more general bases $S$ and non-smooth schemes: there might not be a nicely-behaved set of linearly independent generators of the abelian group $\mathrm{Multi}_S(X, Y)$ of multivalued morphisms, which we have in our case, namely the basic finite correspondences. We will not go deeper into this, but offer a short description how one could resolve this problem:

One has to fix as part of the filtration $\mathcal{F}$ a finite set $S_{\mathcal{F}}$ of effective elements of $\mathrm{Multi}_S(X, Y)$ such that the respective element $\alpha \in \mathrm{Multi}_k(X, Y) \otimes_{\mathbb{N}} \Lambda$ is generated as a $\Lambda$-linear combination from $S_{\mathcal{F}}$. The edge $E_{\mathcal{F}}(\alpha)$ is then the disjoint union of all edges $\widetilde{E}(s)$ with $s \in S_{\mathcal{F}}$. For a filtration $\mathcal{G}$ to be coarser than $\mathcal{F}$ it would then be required that $S_{\mathcal{F}} \subseteq S_{\mathcal{G}}$.

Proposition 4.7.1 then becomes essential to show that $\mathrm{C}^\bullet_{\mathcal{F}}(\alpha)$ is well-defined, as we might have non-trivial relations between the elements of $S_{\mathcal{F}}$.

## 4.8 Functoriality of $\mathrm{C}^\bullet_{\mathcal{F}}$

We are now able to deal with compositions of finite correspondences. As we will see in the proof, it is necessary to further restrict the filtrations.

**Theorem 4.8.1.** *Let $\Lambda$ be a noetherian ring. Let $X$, $Y$ and $Z$ be smooth affine varieties over a field $k \subseteq \mathbb{C}$, and let $\alpha\colon X \rightsquigarrow Y$ and $\beta\colon Y \rightsquigarrow Z$ be finite correspondences in $\mathrm{SmCor}^{\mathrm{aff}}(k, \Lambda)$.*

*Then every sufficiently coarse $\Lambda$-cellular filtration $\mathcal{F}$ on the diagram*

$$X \xrightarrow{\alpha} Y \xrightarrow{\beta} Z$$

*as in Definitions 4.1.1 and 4.5.3 satisfies*

$$\mathrm{C}^\bullet_{\mathcal{F}}(\beta \circ \alpha) = \mathrm{C}^\bullet_{\mathcal{F}}(\alpha) \circ \mathrm{C}^\bullet_{\mathcal{F}}(\beta) \tag{4.6}$$

*as morphisms $\mathrm{C}^\bullet_{\mathcal{F}}(Z) \to \mathrm{C}^\bullet_{\mathcal{F}}(X)$ between Nori's cellular complexes.*

*Rigorously this means:*

*Any given $\Lambda$-cellular filtration $\mathcal{H}$ on the diagram $X \xrightarrow{\alpha} Y$ extends to a $\Lambda$-cellular filtration $\mathcal{G}$ on the diagram $X \xrightarrow{\alpha} Y \xrightarrow{\beta} Z$ with the following property: if $\mathcal{F}$ is an arbitrary $\Lambda$-cellular filtration on the same diagram which extends the restriction of $\mathcal{H}$ to $X$ and which is coarser than $\mathcal{G}$, then equation (4.6) holds.*

*Proof.* By definition and linearity we reduce to the case where $\alpha = \widetilde{\alpha} \otimes \Lambda$ and $\beta = \widetilde{\beta} \otimes \Lambda$ come from basic finite correspondences $\widetilde{\alpha}$ and $\widetilde{\beta}$. In particular, they are each represented by a single irreducible variety, thus their projections to $Y$ each land in a single connected component of $Y$, as also witnessed by Proposition 1.6.4. If these two connected components are different, then $\beta \circ \alpha = 0$ and $\mathrm{C}^\bullet_{\mathcal{F}}(\alpha) \circ \mathrm{C}^\bullet_{\mathcal{F}}(\beta) = 0$ for any $\Lambda$-cellular filtration $\mathcal{F}$ on the diagram, thus the equality is automatically satisfied. Otherwise we may replace $Y$ by the connected component both cycles lie over. We can also replace $X$ and $Z$ by the respective connected components $\widetilde{\alpha}$ and $\widetilde{\beta}$ lie over and observe that by $\Lambda$-linearity $\beta \circ \alpha = (\widetilde{\beta} \circ \widetilde{\alpha}) \otimes \Lambda$.

We have a commutative diagram

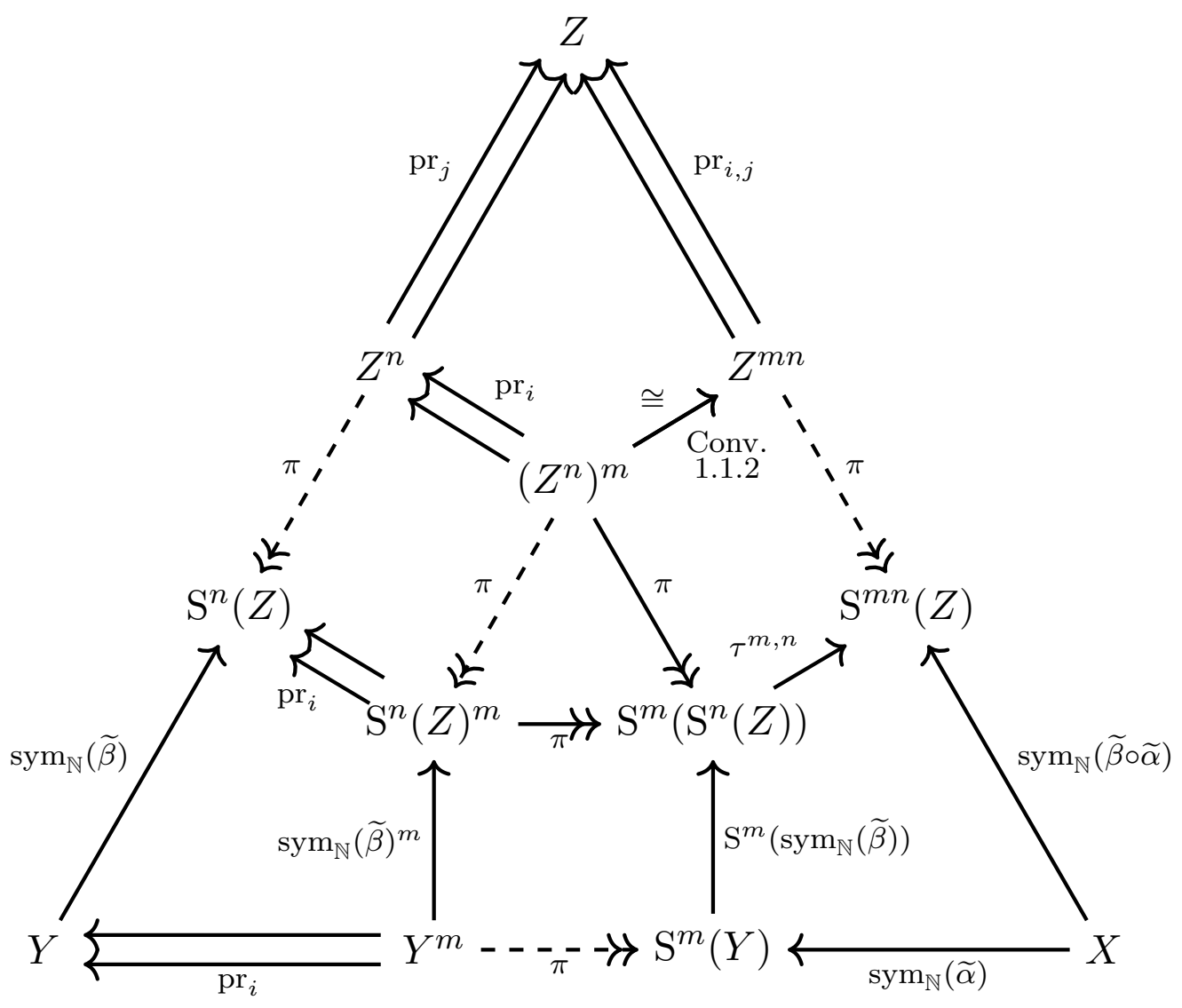

including the edges $E(\alpha) = \widetilde{E}(\widetilde{\alpha})$, $E(\beta) = \widetilde{E}(\widetilde{\beta})$ and $\widetilde{E}(\widetilde{\beta} \circ \widetilde{\alpha})$ as boundary (cf. Remark 4.5.2 for the subtle differences between the edges). Note that the lower right area is the Definition 3.4.7 of composition of multivalued morphisms. Proposition 4.7.2 gives us a $\Lambda$-cellular filtration $\mathcal{G}$ on the edges $E((\widetilde{\beta} \circ \widetilde{\alpha}) \otimes \Lambda) = E(\beta \circ \alpha)$ and $\widetilde{E}(\widetilde{\beta} \circ \widetilde{\alpha})$ such that $\mathrm{C}^\bullet_{\mathcal{G}}(\beta \circ \alpha) = \widetilde{\mathrm{C}}^\bullet_{\mathcal{G}}(\widetilde{\beta} \circ \widetilde{\alpha})$. Hence we are left to show that

$$\widetilde{\mathrm{C}}^\bullet_{\mathcal{G}}(\widetilde{\beta} \circ \widetilde{\alpha}) = \widetilde{\mathrm{C}}^\bullet_{\mathcal{G}}(\widetilde{\alpha}) \circ \widetilde{\mathrm{C}}^\bullet_{\mathcal{G}}(\widetilde{\beta})$$

and are by the precise version of Proposition 4.7.2 allowed to coarsen $\mathcal{G}$ everywhere except on $X$.

By potentially coarsening it on $\widetilde{E}(\widetilde{\beta}\circ\widetilde{\alpha})$, Lemma 4.4.7 allows us to extend $\mathcal{G}$ to a $\Lambda$-cellular filtration on the above diagram that it is nice at the four dotted arrows and extends the given filtration $\mathcal{H}$ on $\widetilde{E}(\alpha)$. In particular, we left the filtration unchanged at $X$ as desired above. We therefore arrive at the diagram which commutes on the three outer areas by the definition of $\widetilde{\mathrm{C}}^\bullet_{\mathcal{G}}(-)$. It commutes on the inner areas because this is already true on the geometric counterparts. Thus a diagram chase using the four isomorphisms shows

$$\widetilde{\mathrm{C}}^\bullet_{\mathcal{G}}(\widetilde{\beta}\circ\widetilde{\alpha}) = \widetilde{\mathrm{C}}^\bullet_{\mathcal{G}}(\widetilde{\alpha})\circ\widetilde{\mathrm{C}}^\bullet_{\mathcal{G}}(\widetilde{\beta}).$$

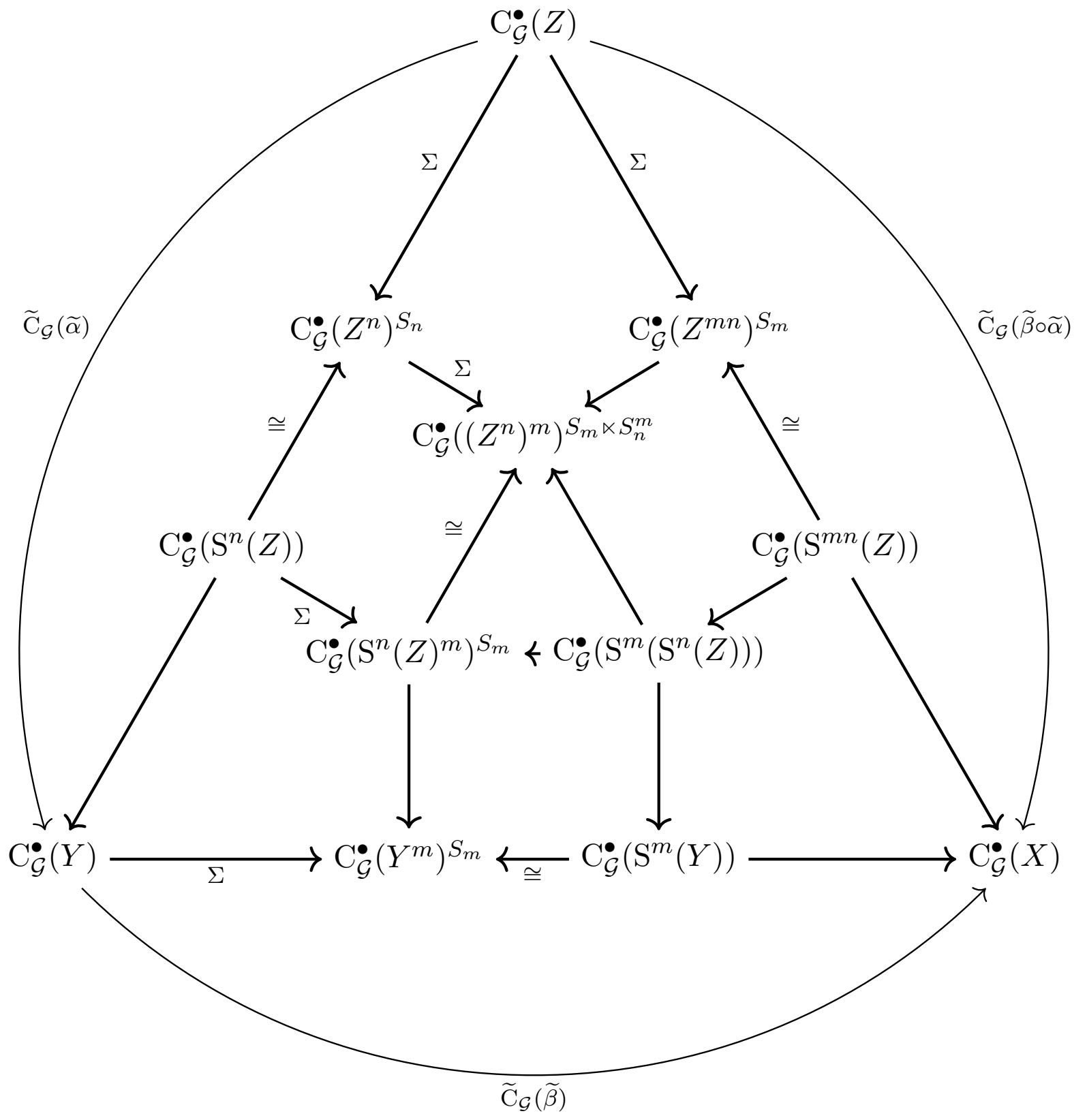

If $\mathcal{F}$ is any $\Lambda$-cellular filtration on $X \xrightarrow{\alpha} Y \xrightarrow{\beta} Z$ coarser than $\mathcal{G}$, but the same on $X$, we argue as in the end of the proof of Proposition 4.7.1: both $\mathrm{C}^\bullet_{\mathcal{F}}(\beta\circ\alpha)$ and $\mathrm{C}^\bullet_{\mathcal{F}}(\alpha)\circ\mathrm{C}^\bullet_{\mathcal{F}}(\beta)$ factor by Lemma 4.6.3 as

$$\mathrm{C}^\bullet_{\mathcal{F}}(Z)\to\mathrm{C}^\bullet_{\mathcal{G}}(Z)\to\mathrm{C}^\bullet_{\mathcal{G}}(X)=\mathrm{C}^\bullet_{\mathcal{F}}(X)$$

over $\mathrm{C}^\bullet_{\mathcal{G}}(\beta\circ\alpha)$ and $\mathrm{C}^\bullet_{\mathcal{G}}(\alpha)\circ\mathrm{C}^\bullet_{\mathcal{F}}(\beta)$ and are thus equal. $\square$

**Definition 4.8.2.** We define a $\Lambda$-cellular filtration $\mathcal{F}$ on a diagram $D$ in $\mathrm{SmCor}^{\mathrm{aff}}(k,\Lambda)$ to be *functorial* if

$$\mathrm{C}^\bullet_{\mathcal{F}}(\beta \circ \alpha) = \mathrm{C}^\bullet_{\mathcal{F}}(\alpha) \circ \mathrm{C}^\bullet_{\mathcal{F}}(\beta)$$

for each pair $(\alpha, \beta)$ of composable finite correspondences in $D$.

**Definition 4.8.3.** An almost acyclic diagram (cf. Definition 1.2.9) $D : \mathcal{C} \to \mathrm{SmCor}(k,\Lambda)$ is called *hom-almost acyclic* if every automorphism $D(a)$, with $a$ an automorphism of $c \in \mathcal{C}$, comes from a scheme-theoretic automorphism of $D(c)$ via the embedding $\mathrm{Sm}_k \hookrightarrow \mathrm{SmCor}(k,\Lambda)$ of Definition 1.6.9.

**Theorem 4.8.4.** *Let $k \subseteq \mathbb{C}$ be a field and let $\Lambda$ be a noetherian ring. Let $D\colon \mathcal{C} \to \mathrm{SmCor}^{\mathrm{aff}}(k,\Lambda)$ be a finite and hom-almost acyclic diagram (cf. Definition 4.8.3).*

*Then there exists a functorial $\Lambda$-cellular filtration $\mathcal{F}$ on $D$ coarser than a given standard filtration. In particular they form a directed system.*

*Proof.* Similar to the proofs of Lemma 4.4.7 and Lemma 4.5.5 this is a straightforward induction using Theorem 4.8.1:

Let $z$ be a terminal object of $\mathcal{C}$, i.e. one that maps to no other object. By induction on the size of $\mathcal{C}$ we may assume to already have a functorial $\Lambda$-cellular filtration on the rest $\mathcal{C}\backslash\{z\}$ of the diagram. Let now $\alpha\colon x \to y$ and $\beta\colon y \to z$ be a pair of composable arrows in $\mathcal{C}$. Then by Theorem 4.8.1 all sufficiently coarse $\Lambda$-cellular filtrations on $D(x) \to D(y) \to D(z)$ extending the given one on $D(x) \to D(y)$ are functorial.

We do this for every such pair $(\alpha, \beta)$. Each time we might have to coarsen the filtration on $D(\beta)$ and thus on $D(z)$, which is clearly allowed by the precise version of Theorem 4.8.1.

Lastly, we consider the finite automorphism group $G := \mathrm{Aut}_D(z) = \mathrm{End}_D(z)$. By Theorem 4.4.5 we can then coarsen the found $\Lambda$-cellular filtration on $D(z)$ to a $G$-invariant one, finishing the proof. $\square$

## 4.9 Tensor structure

If $(X,Y)$ and $(X',Y')$ are good pairs (cf. Definition 1.3.7), then we have the Künneth isomorphism

$$H^{i+j}_{\mathrm{sing}}(X \times X', X \times Y' \cup Y \times X') \cong H^i_{\mathrm{sing}}(X,Y) \otimes_\Lambda H^j_{\mathrm{sing}}(X',Y').$$

Thus expecting

$$H^i_{\mathrm{Nori}}(X,Y) \otimes H^j_{\mathrm{Nori}}(X,Y) \cong H^{i+j}_{\mathrm{Nori}}(X \times X', X \times Y' \cup Y \times X')$$

is natural. Indeed, this motivates the definition of the tensor product on $\mathcal{MM}^{\mathrm{eff}}_{\mathrm{Nori}}(k,\Lambda)$ in [HM16]. Combined with the formal tensor structures on diagrams of op. cit. it turns $\mathcal{MM}^{\mathrm{eff}}_{\mathrm{Nori}}(k,\Lambda)$ into a symmetric monoidal category such that the faithful exact functor $\omega_{\mathrm{sing}}$ to $\Lambda$- Mod preserves the tensor structure.

**Definition 4.9.1.** Let $X$ and $Y$ be affine varieties with standard filtrations $\mathcal{F}$ and $\mathcal{G}$, respectively. The *product filtration* $\mathcal{F}\times\mathcal{G}$ on $X\times Y$ is defined as

$$(\mathcal{F}\times\mathcal{G})_n(X\times Y) := \bigcup_{i+j=n} \mathcal{F}_iX\times\mathcal{G}_jY.$$

**Proposition 4.9.2.** *Let $\Lambda$ be a noetherian ring. Let $X$ and $Y$ be affine varieties over a field $k\subseteq\mathbb{C}$ with $\Lambda$-cellular filtrations $\mathcal{F}$ and $\mathcal{G}$, respectively.*

*Then $\mathcal{F}\times\mathcal{G}$ is a $\Lambda$-cellular filtration on $X\times Y$.*

*Furthermore, the tensor structure of $\mathcal{MM}^{\mathrm{eff}}_{\mathrm{Nori}}(k,\Lambda)$ induces a natural isomorphism*

$$\psi^\bullet_{\mathcal{F},\mathcal{G}}\colon \mathrm{C}^\bullet_{\mathcal{F}\times\mathcal{G}}(X\times Y) \xrightarrow{\cong} \mathrm{Tot}^\bullet(\mathrm{C}^\bullet_{\mathcal{F}}(X)\otimes\mathrm{C}^\bullet_{\mathcal{G}}(Y))$$

*of bounded chain complexes over $\mathcal{MM}^{\mathrm{eff}}_{\mathrm{Nori}}(k,\Lambda)$.*

*Proof.* If $Z$ is any of the objects $X$, $Y$ or $X\times Y$ we use the shorthand $Z_i$ to denote the $i$-th element of the corresponding filtration $\mathcal{F}$, $\mathcal{G}$ or $\mathcal{F}\times\mathcal{G}$, respectively. We also write $H^\bullet_{\mathrm{sing}}(A,B)$ for $H^\bullet_{\mathrm{sing}}(A,B,\Lambda)$.

As each $H^s_{\mathrm{sing}}(X_i,X_{i-1})$ and $H^t_{\mathrm{sing}}(Y_j,Y_{j-1})$ is projective, in particular flat, the Künneth formula for relative cohomology with an arbitrary ring of coefficients states

$$\begin{aligned}&H^n_{\mathrm{sing}}\left(X_i\times Y_j, X_{i-1}\times Y_j\cup X_i\times Y_{j-1}\right)\cong\\ &\cong\bigoplus_{s+t=n} H^s_{\mathrm{sing}}(X_i,X_{i-1})\otimes_\Lambda H^t_{\mathrm{sing}}(Y_j,Y_{j-1}).\end{aligned}$$

The summands at the bottom vanish by assumption unless $s=i$ and $t=j$. Therefore we have an isomorphism

$$H^{i+j}_{\mathrm{sing}}\left(X_i\times Y_j, X_{i-1}\times Y_j\cup X_i\times Y_{j-1}\right)\cong H^i_{\mathrm{sing}}(X_i,X_{i-1})\otimes_\Lambda H^j_{\mathrm{sing}}(Y_j,Y_{j-1})$$

and $H^n_{\mathrm{sing}}\left(X_i\times Y_j, X_{i-1}\times Y_j\cup X_i\times Y_{j-1}\right)=0$ whenever $n\neq i+j$.

Then an iterated application of the Mayer-Vietoris long exact sequence shows

$$\begin{aligned}&H^m_{\mathrm{sing}}\left((X\times Y)_n,(X\times Y)_{n-1}\right)\cong\\ &\cong\bigoplus_{i+j=n} H^m_{\mathrm{sing}}\left(X_i\times Y_j, X_{i-1}\times Y_j\cup X_i\times Y_{j-1}\right).\end{aligned}$$

In particular this is 0 unless $m=n$.

By combining the above isomorphisms we get a natural isomorphism

$$H^n_{\text{sing}}\left((X \times Y)_n, (X \times Y)_{n-1}\right) \cong \bigoplus_{i+j=n} H^i_{\text{sing}}(X_i, X_{i-1}) \otimes_\Lambda H^j_{\text{sing}}(Y_j, Y_{j-1})$$

which we denote as $\varphi^n_{\mathcal{F},\mathcal{G}}$.

All the summands on the left are tensor products of finitely generated projective $\Lambda$-modules and therefore finitely generated and projective. Hence the same is true for the finite sum. Thus $H^m_{\text{sing}}\left((X \times Y)_n, (X \times Y)_{n-1}\right)$ is finitely generated, projective and vanishes if $m \neq n$. Therefore $\mathcal{F} \times \mathcal{G}$ satisfies the defining properties of a $\Lambda$-cellular filtration.

We now use Proposition 9.3.1 of [HM16] to lift this to Nori motives:

The graded multiplicative structure on the quiver $\text{VGood}^{\text{eff}}$ induces a morphism

$$H^n_{\text{Nori}}\left(X_i \times Y_j, X_{i-1} \times Y_j \cup X_i \times Y_{j-1}\right) \to H^i_{\text{Nori}}(X_i, X_{i-1}) \otimes H^j_{\text{Nori}}(Y_j, Y_{j-1})$$

whenever $i + j = n$. The inclusion

$$(X_i \times Y_j, X_{i-1} \times Y_j \cup X_i \times Y_{j-1}) \hookrightarrow ((X \times Y)_n, (X \times Y)_{n-1})$$

of very good pairs gives a morphisms

$$H^n_{\text{Nori}}\left((X \times Y)_n, (X \times Y)_{n-1}\right) \to H^n_{\text{Nori}}\left(X_i \times Y_j, X_{i-1} \times Y_j \cup X_i \times Y_{j-1}\right).$$

Composing them and taking the sum over all $i, j$ with $i + j = n$ gives us a morphism $\psi^n_{\mathcal{F},\mathcal{G}}$:

$$H^n_{\text{Nori}}\left((X \times Y)_n, (X \times Y)_{n-1}\right) \to \bigoplus_{i+j=n} H^i_{\text{Nori}}(X_i, X_{i-1}) \otimes H^j_{\text{Nori}}(Y_j, Y_{j-1})$$

such that, by Proposition 9.3.1 of [HM16] and the construction, we have $\omega_{\text{sing}}(\psi^n_{\mathcal{F},\mathcal{G}}) = \varphi^n_{\mathcal{F},\mathcal{G}}$. Here $\omega_{\text{sing}}$ is the faithful tensor functor $\mathcal{MM}^{\text{eff}}_{\text{Nori}}(k, \Lambda) \to \Lambda\text{-}\mathrm{Mod}$, hence $\psi^n_{\mathcal{F},\mathcal{G}}$ is an isomorphism as well. Due to the graded compatibility of the Künneth formula with boundary morphisms (cf. Proposition 2.4.3 of [HM16]) we thus got the desired isomorphism

$$\psi^\bullet_{\mathcal{F},\mathcal{G}} \colon \mathrm{C}^\bullet_{\mathcal{F}\times\mathcal{G}}(X \times Y) \to \mathrm{Tot}^\bullet(\mathrm{C}^\bullet_{\mathcal{F}}(X) \otimes \mathrm{C}^\bullet_{\mathcal{G}}(Y)).$$

□

**Definition 4.9.3.** Let $X$ and $Y$ be affine varieties over $k \subseteq \mathbb{C}$ with $\Lambda$-cellular filtrations $\mathcal{F}$ and $\mathcal{G}$, respectively. Let $\mathcal{H}$ be a $\Lambda$-cellular filtration on $X \times Y$ coarser than $\mathcal{F} \times \mathcal{G}$. Then we denote the resulting morphism

$$\mathrm{C}^\bullet_{\mathcal{H}}(X \times Y) \to \mathrm{C}^\bullet_{\mathcal{F}\times\mathcal{G}}(X \times Y) \to \mathrm{Tot}^\bullet(\mathrm{C}^\bullet_{\mathcal{F}}(X) \otimes \mathrm{C}^\bullet_{\mathcal{G}}(Y))$$

again by $\psi^\bullet_{\mathcal{F},\mathcal{G}}$.

**Remark 4.9.4.** From Proposition 4.9.2 and Lemma 4.2.4 we see that $\psi^\bullet_{\mathcal{F},\mathcal{G}}$ is a quasi-isomorphism.

**Remark 4.9.5.** The product filtrations $\mathcal{F} \times \mathcal{G}$ on $X \times Y$ do in general not form a final subsystem, neither for standard nor for $\Lambda$-cellular filtrations. For example, the filtration

$$0 \hookrightarrow \Delta \hookrightarrow \mathbb{A}^2_k$$

of $\mathbb{A}^2_k$, where $\Delta$ is the diagonal, is $\Lambda$-cellular, but $\Delta$ is not contained in any $X_0 \times Y_1 \cup X_1 \times Y_0$, even if we allow arbitrary closed subsets $X_i, Y_i \hookrightarrow \mathbb{A}^1$ of dimension at most $i$. For the monoidal structure on motives this means that we cannot expect a strict tensor functor unless we invert the morphisms corresponding to coarsenings of $\Lambda$-cellular filtrations.

**Theorem 4.9.6.** *Let $k \subseteq \mathbb{C}$ be a field and let $\Lambda$ be a noetherian ring. Let*

$$\alpha_1 \colon X_1 \rightsquigarrow Y_1,$$
$$\alpha_2 \colon X_2 \rightsquigarrow Y_1$$

*be finite correspondences in $\mathrm{SmCor}^{\mathrm{aff}}(k, \Lambda)$ and recall the Definition 1.7.4 of their tensor product*

$$\alpha_1 \otimes \alpha_2 \colon X_1 \times X_2 \rightsquigarrow Y_1 \times Y_2.$$

*Then all sufficiently coarse $\Lambda$-cellular filtrations $\mathcal{F}_1$, $\mathcal{F}_2$ and $\mathcal{F}$ on the finite correspondences $\alpha_1$, $\alpha_2$ and $\alpha_1 \otimes \alpha_2$, respectively, induce, together with the quasi-isomorphisms $\psi^\bullet_{\mathcal{F},\mathcal{G}}$ of Definition 4.9.3, a natural commutative square*

$$\begin{array}{ccc}
\mathrm{Tot}^\bullet(\mathrm{C}^\bullet_{\mathcal{F}_1}(X_1) \otimes \mathrm{C}^\bullet_{\mathcal{F}_2}(X_2)) & \xleftarrow{\mathrm{Tot}^\bullet(\mathrm{C}^\bullet_{\mathcal{F}_1}(\alpha_1) \otimes \mathrm{C}^\bullet_{\mathcal{F}_2}(\alpha_2))} & \mathrm{Tot}^\bullet(\mathrm{C}^\bullet_{\mathcal{F}_1}(Y_1) \otimes \mathrm{C}^\bullet_{\mathcal{F}_2}(Y_2)) \\
\uparrow {\scriptstyle \psi^\bullet_{\mathcal{F}_1,\mathcal{F}_2}\ \text{q.-is.}} & & \uparrow {\scriptstyle \psi^\bullet_{\mathcal{F}_1,\mathcal{F}_2}\ \text{q.-is.}} \\
\mathrm{C}^\bullet_{\mathcal{F}}(X_1 \times X_2) & \xleftarrow{\mathrm{C}^\bullet_{\mathcal{F}}(\alpha_1 \otimes \alpha_2)} & \mathrm{C}^\bullet_{\mathcal{F}}(Y_1 \times Y_2).
\end{array} \tag{4.7}$$

*Rigorously this means:*

*Let any $\Lambda$-cellular filtrations $\mathcal{G}_1$ and $\mathcal{G}_2$ on the finite correspondences $\alpha_1$ and $\alpha_2$, respectively, be given. Then there exists a $\Lambda$-cellular filtration $\mathcal{G}$ on the finite correspondence $\alpha_1 \otimes \alpha_2$ which extends $\mathcal{G}_1 \times \mathcal{G}_2$ on $X_1 \times X_2$ and is coarser than $\mathcal{G}_1 \times \mathcal{G}_2$ on $Y_1 \times Y_2$. It is such that, given*

- *a $\Lambda$-cellular filtration $\mathcal{F}_1$ on $\alpha_1$, which is coarser than $\mathcal{G}_1$ and agrees with it on $X_1$,*

- *a $\Lambda$-cellular filtration $\mathcal{F}_2$ on $\alpha_2$, which is coarser than $\mathcal{G}_2$ and agrees with it on $X_2$,*

- *a $\Lambda$-cellular filtration $\mathcal{F}$ on $\alpha_1 \otimes \alpha_2$, which is coarser than $\mathcal{G}$,*

*the diagram 4.7 commutes.*

*Proof.* By definition and linearity we may assume that $\alpha_1$ and $\alpha_2$ are basic of degrees $n_1$ and $n_2$, respectively. We then replace $X_1$, $Y_1$, $X_2$ and $Y_2$ by their respective connected components over which the correspondences are defined. By Remark 1.4.25 and Lemma 1.8.15 the tensor product $\alpha_1 \otimes \alpha_2$ is represented by an effective finite correspondence $\widetilde{\alpha}\colon X_1 \times X_2 \rightsquigarrow Y_1 \times Y_2$ of constant degree $n_1 n_2$.

We consider the diagram

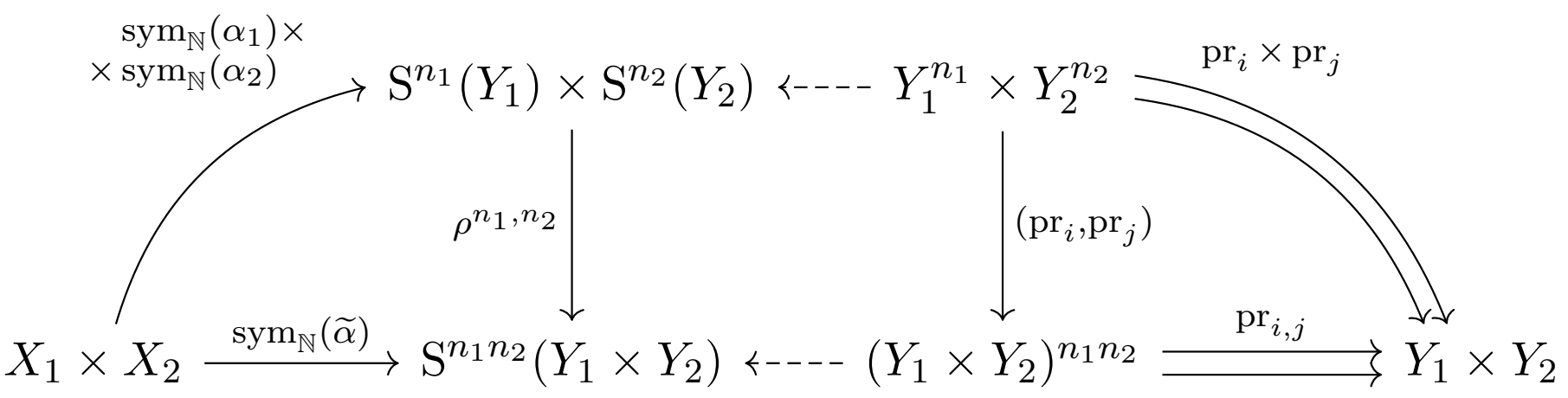

where the dashed arrows are the respective quotients by the group actions. The middle square is the defining one of Proposition 3.3.7 and therefore commutes. The area on the right commutes trivially. Finally, the area on the left commutes by Theorem 3.7.4. We also observe that the bottom row is the effective edge $\widetilde{E}(\widetilde{\alpha})$.

By Lemma 4.4.7 and Proposition 4.7.2 we find a $\Lambda$-cellular filtration $\mathcal{G}$ on the entire diagram and the edge $E(\widetilde{\alpha} \otimes \Lambda)$ which

- agrees with $\mathcal{G}_1 \times \mathcal{G}_2$ on $X_1 \times X_2$,

- is coarser than the product filtrations $\mathcal{G}_1 \times \mathcal{G}_2$ on $\mathrm{S}^{n_1}(Y_1) \times \mathrm{S}^{n_2}(Y_2)$, $Y_1^{n_1} \times Y_2^{n_2}$, and $Y_1 \times Y_2$,

- is nice at the dotted arrows,

- satisfies
$$\mathrm{C}^{\bullet}_{\mathcal{G}}(\widetilde{\alpha} \otimes \Lambda) = \widetilde{\mathrm{C}}^{\bullet}_{\mathcal{G}}(\widetilde{\alpha}).$$

Hence we produced a commutative diagram

$$
\begin{array}{ccccccc}
\multicolumn{7}{c}{\xleftarrow{\quad \mathrm{Tot}^\bullet(\mathrm{C}^\bullet_{\mathcal{G}_1}(\alpha_1)\otimes \mathrm{C}^\bullet_{\mathcal{G}_2}(\alpha_2)) \quad}} \\
\mathrm{Tot}^\bullet\begin{pmatrix} \mathrm{C}^\bullet_{\mathcal{G}_1}(X_1)\otimes \\ \otimes\,\mathrm{C}^\bullet_{\mathcal{G}_2}(X_2)\end{pmatrix} & \longleftarrow & \mathrm{Tot}^\bullet\begin{pmatrix} \mathrm{C}^\bullet_{\mathcal{G}_1}(\mathrm{S}^{n_1}(Y_1))\otimes \\ \otimes\,\mathrm{C}^\bullet_{\mathcal{G}_2}(\mathrm{S}^{n_2}(Y_2))\end{pmatrix} & \xrightarrow{\cong} & \mathrm{Tot}^\bullet\begin{pmatrix} \mathrm{C}^\bullet_{\mathcal{G}_1}(Y_1^{n_1})^{S_{n_1}}\otimes \\ \otimes\,\mathrm{C}^\bullet_{\mathcal{G}_2}(Y_2^{n_2})^{S_{n_2}}\end{pmatrix} & \longleftarrow & \mathrm{Tot}^\bullet\begin{pmatrix} \mathrm{C}^\bullet_{\mathcal{G}_1}(Y_1)\otimes \\ \otimes\,\mathrm{C}^\bullet_{\mathcal{G}_2}(Y_2)\end{pmatrix} \\
\uparrow & & \uparrow{\scriptstyle \psi^\bullet_{\mathcal{G}_1,\mathcal{G}_2}} & & \uparrow{\scriptstyle \psi^\bullet_{\mathcal{G}_1,\mathcal{G}_2}} & & \uparrow \\
{\scriptstyle \psi^\bullet_{\mathcal{G}_1,\mathcal{G}_2}} & & \mathrm{C}^\bullet_{\mathcal{G}}(\mathrm{S}^{n_1}(Y_1)\times \mathrm{S}^{n_2}(Y_2)) & \xrightarrow{\cong} & \mathrm{C}^\bullet_{\mathcal{G}}(Y_1^{n_1}\times Y_2^{n_2})^{S_{n_1}\times S_{n_2}} & & {\scriptstyle \psi^\bullet_{\mathcal{G}_1,\mathcal{G}_2}} \\
\uparrow & \swarrow & \uparrow & & \uparrow & \nwarrow & \uparrow \\
\mathrm{C}^\bullet_{\mathcal{G}_1\times\mathcal{G}_2}(X_1\times X_2) & \longleftarrow & \mathrm{C}^\bullet_{\mathcal{G}}(\mathrm{S}^{n_1n_2}(Y_1\times Y_2)) & \xrightarrow{\cong} & \mathrm{C}^\bullet_{\mathcal{G}}(Y_1^{n_1n_2}\times Y_2^{n_1n_2})^{S_{n_1n_2}} & \longleftarrow & \mathrm{C}^\bullet_{\mathcal{G}}(Y_1\times Y_2). \\
\multicolumn{7}{c}{\xleftarrow[\quad \mathrm{C}^\bullet_{\mathcal{G}}(\widetilde{\alpha}\otimes\Lambda)=\widetilde{\mathrm{C}}^\bullet_{\mathcal{G}}(\widetilde{\alpha}) \quad]{}}
\end{array}
$$

This shows the result for the constructed filtrations. If $\mathcal{F}_1$, $\mathcal{F}_2$ and $\mathcal{F}$ are as in the precise version, then the diagram (4.7) commutes by Lemma 4.6.3 due to a similar argument as seen at the end of the proof of Theorem 4.8.1. $\square$

There is again a version on diagrams:

**Theorem 4.9.7.** *Let $k \subseteq \mathbb{C}$ be a field and let $\Lambda$ be a noetherian ring. Let*

$$
\begin{aligned}
D_1 &\colon \mathcal{C}_1 \to \mathrm{SmCor}^{\mathrm{aff}}(k,\Lambda),\\
D_2 &\colon \mathcal{C}_2 \to \mathrm{SmCor}^{\mathrm{aff}}(k,\Lambda)
\end{aligned}
$$

*be two finite hom-almost acyclic diagrams (cf. Definition 4.8.3).*

*Then the product diagram (cf. Definition 1.2.11)*

$$D_1 \otimes D_2 \colon \mathcal{C}_1 \times \mathcal{C}_2 \to \mathrm{SmCor}^{\mathrm{aff}}(k,\Lambda)$$

*induced by the tensor structure of $\mathrm{SmCor}^{\mathrm{aff}}(k,\Lambda)$ is finite and hom-almost acyclic. Furthermore, there exist functorial $\Lambda$-cellular filtrations $\mathcal{F}_1$, $\mathcal{F}_2$ and $\mathcal{F}$ on $D_1$, $D_2$ and $D_1 \otimes D_2$, respectively, which are finer than given standard filtrations on them and compatible with tensor products in the following sense:*

*For any arrows $a_i\colon x_i \to y_i$ in $\mathcal{C}_i$, where $i \in \{1,2\}$, we set $X_i = D_i(x_i)$, $Y_i = D_i(y_i)$ and $\alpha_i = D_i(a_i)$. Then $\mathcal{F}$ is coarser than $\mathcal{F}_1 \times \mathcal{F}_2$ on each $X_1 \times X_2$ and the diagrams*

$$
\begin{array}{ccc}
\mathrm{Tot}^\bullet(\mathrm{C}^\bullet_{\mathcal{F}_1}(X_1)\otimes \mathrm{C}^\bullet_{\mathcal{F}_2}(X_2)) & \xleftarrow{\mathrm{Tot}^\bullet(\mathrm{C}^\bullet_{\mathcal{F}_1}(\alpha_1)\otimes \mathrm{C}^\bullet_{\mathcal{F}_2}(\alpha_2))} & \mathrm{Tot}^\bullet(\mathrm{C}^\bullet_{\mathcal{F}_1}(Y_1)\otimes \mathrm{C}^\bullet_{\mathcal{F}_2}(Y_2)) \\
{\scriptstyle \psi^\bullet_{\mathcal{F}_1,\mathcal{F}_2}}\uparrow{\scriptstyle \text{q.-is.}} & & {\scriptstyle \psi^\bullet_{\mathcal{F}_1,\mathcal{F}_2}}\uparrow{\scriptstyle \text{q.-is.}} \\
\mathrm{C}^\bullet_{\mathcal{F}}(X_1\times X_2) & \xleftarrow{\quad \mathrm{C}^\bullet_{\mathcal{F}}(\alpha_1\otimes\alpha_2) \quad} & \mathrm{C}^\bullet_{\mathcal{F}}(Y_1\times Y_2).
\end{array}
$$

*involving the morphisms $\psi^\bullet_{\mathcal{F}_1,\mathcal{F}_2}$ of Definition 4.9.3 commute.*

*Proof.* The finiteness and hom-almost acyclicity of $D_1 \otimes D_2$ are trivial. We can choose functorial $\Lambda$-cellular filtrations on $D_1$ and $D_2$ by Theorem 4.8.4. Afterwards one uses Theorem 4.9.6 to do a straightforward induction on $D_1 \otimes D_2$. □

# Chapter 5

# Covers, Bridges and Čech Complexes

As seen in Chapter 4, $\Lambda$-cellular filtrations functorially exist on affine varieties and by extension complexes of them. But, as explained in Remark 4.3.3, we cannot expect to extend this to non-affine varieties.

One might attempt to resolve this problem by choosing an open affine cover $\mathcal{U} = \{U_i\} \twoheadrightarrow X$ and then to work with its (reduced) Čech complex (cf. Definitions 5.4.1 and 5.4.2). We will see in Theorem 7.2.4 that this indeed gives us an adequate object of $D^b(\mathcal{MM}^{\mathrm{eff}}_{\mathrm{Nori}})$ and Section 9.2 of [HM16] explains how to deal with morphisms of varieties. This does, however, not allow us to extend the functor to arbitrary finite correspondences $X \rightsquigarrow Y$.

It is central to this problem that Zariski covers do not induce a pretopology on $\mathrm{SmCor}(S)$. The solution is to instead use the finer Nisnevich covers. Consequently, we have to deal with several technicalities not encountered in [HM16]. The main ideas are, given a regular scheme $S$ of finite dimension:

1. Every object of $\mathrm{Sm}_S$ becomes in $\mathrm{DM}^{\mathrm{eff}}(S, \Lambda)$ isomorphic to the (reduced) Čech complex of appropriate Nisnevich covers as witnessed by Theorem 10.3.3 of [CD12] and by Theorem A.0.10. Hence it is natural to work with them.

2. Finite correspondences lift via so-called *bridges* (cf. Definition 5.2.1) to Nisnevich or even more general covers and their Čech complexes (Definition 6.4.1 and Proposition 6.4.6).

3. Morphisms, and consequently finite correspondences that come from bridges, can be made unique by choosing *rigidifications* (Definition 6.1.1, Lemma 6.2.6, and Theorem 6.4.9).

4. Every finite acyclic diagram in $\mathrm{SmCor}(S)$ admits enough Nisnevich covers satisfying the previous properties (Theorem 6.5.3).

5. The diagrams resulting from the infinite Čech complexes can be reinterpreted to finite diagrams by identifying schemes isomorphic over the respective bases (Proposition 5.5.7). This finiteness then extends to the finite correspondences between them (Proposition 6.6.3).

This chapter will focus on the notions of covers, bridges and Čech complexes. The following Chapter 6 then uses those techniques to define and use rigidifications. Altogether, this will, as desired, ultimately allow us to apply the results of Chapter 4 to finite acyclic diagrams in $\mathrm{SmCor}(k, \Lambda)$.

We will also briefly talk about the new notion of *unifibrant* Nisnevich covers. They are interesting on their own because they behave similar to the pretopology induced by Nisnevich cd-squares and because they satisfy a strong finiteness result on which we elaborate in Appendix A.

Lastly, we note that most of our results, with the exceptions of those involving the finitistic covers of Definition A.0.1, have analogues for étale covers, which we will neither need nor prove.

## 5.1 Generalities on covers

Fix a category $\mathcal{C}$. The following definitions are standard:

**Definition 5.1.1.** Let $X \in \mathcal{C}$. A *pseudocover* of $X$ is a set of morphisms $f_i \colon \mathcal{U}_i \to X$, $i \in I$, in $\mathcal{C}$. If the category $\mathcal{C}$ permits coproducts over $I$ we will often understand this as a single morphism $\mathcal{U} = \coprod_{i \in I} \mathcal{U}_i \to X$ together with the fixed decomposition $\mathcal{U} = \coprod_{i \in I} \mathcal{U}_i$. We call $X$ its *base* and $I$ its *indexing set*. We write

$$f \colon (\mathcal{U}|I) \twoheadrightarrow X$$

to denote a pseudocover, or simply $f \colon \mathcal{U} \twoheadrightarrow X$ if the indexing set or the decomposition do not matter.

We call a pseudocover *finitely indexed* if its indexing set is a finite set. We will respectively use the words *piecewise* and *jointly* to differentiate between the properties of the individual morphism $\mathcal{U}_i \to X$ and the total morphism $\mathcal{U} \to X$, assuming that the latter exists. We call the pseudocover *ordered* if its indexing set is equipped with a total order.

**Remark 5.1.2.** We will only work in categories where the relevant coproducts exist.

**Convention 5.1.3.** We assume all pseudocovers of schemes to be jointly of finite type unless explicitly stated otherwise.

**Definition 5.1.4.** Let $f \colon (\mathcal{U}|I) \twoheadrightarrow X$ be a pseudocover. If we are additionally given a collection $g = \{g_i \colon (\mathcal{V}_i|J_i) \twoheadrightarrow U_i \mid i \in I\}$ of pseudocovers of the individual $U_i$, then the *composite pseudocover* is defined as follows:

Its indexing set is

$$J_I := \{(i,j) \mid i \in I, j \in J_i\},$$

which is nothing else than the formal construction of the disjoint union $\bigsqcup_{i \in I} J_i$ of sets. The individual morphism for $(i,j) \in J_I$ is

$$(f \circ g)_{(i,j)} := f_i \circ g_j \colon \mathcal{V}_j \to X.$$

If $I$ and the $J_i$ are ordered, then so is $J_I$ by the *lexicographic order*: $(i,j) \le (i',j')$ if either $i = i'$ and $j \le j'$ or $i < i'$. Thus the composite of ordered pseudocovers is an ordered pseudocover and it is easily checked that this operation is associative.

**Definition 5.1.5.** A *morphism of pseudocovers*

$$f \colon (\mathcal{U}|I) \twoheadrightarrow X,$$
$$g \colon (\mathcal{V}|J) \twoheadrightarrow Y$$

consists of a map $\delta \colon I \to J$ and morphisms $s_i \colon \mathcal{U}_i \to \mathcal{V}_{\delta(i)}$ for all $i \in I$. If the coproducts exist we interpret this as a single morphism $s \colon \mathcal{U} \to \mathcal{V}$ preserving the decompositions.

This morphism is called *compatible* with a given morphism $t \colon X \to Y$ if the obvious squares commute: $t \circ f_i = g_{\delta(i)} \circ s_i$ for all $i \in I$.

A *morphism of pseudocovers over* $X$, or a *refinement of pseudocovers*, is a morphism of pseudocovers of $X$ that is a morphism over $X$, i.e. compatible with the identity on $X$.

**Definition 5.1.6.** Assume that $\mathcal{C}$ admits fibre products.

The *pullback* of a pseudocover $f \colon (\mathcal{U}|I) \twoheadrightarrow X$ along a morphism $s \colon Y \to X$ is the pseudocover

$$f \times_X Y \colon (\mathcal{U} \times_X Y|I) \twoheadrightarrow Y$$

with the same indexing set $I$ and the individual pullbacks

$$(f \times_X Y)_i = f_i \times_X Y \colon \mathcal{V}_i \times_X Y \to Y.$$

If $f_1 \colon (\mathcal{U}_1|I_1) \twoheadrightarrow X$ and $f_2 \colon (\mathcal{U}_2|I_2) \twoheadrightarrow X$ are two pseudocovers, then their *fibre product* is the pseudocover

$$f_1 \times_X f_2 \colon (\mathcal{U}_1 \times_X \mathcal{U}_2|I_1 \times I_2) \twoheadrightarrow X$$

with indexing set $I_1 \times I_2$ and individual morphisms

$$(f_1 \times_X f_2)_{(i_1,i_2)} = f_{1,i_1} \times_X f_{2,i_2} \colon (\mathcal{U}_1 \times_X \mathcal{U}_2)_{i_1,i_2} = \mathcal{U}_{1,i_1} \times_X \mathcal{U}_{2,i_2} \to X.$$

If both indexing sets are ordered we equip $I_1 \times I_2$ with the lexicographic order.

## 5.2 Bridges

The notion of 'Bridges' introduced in this section offers a method to encode finite correspondences via mere morphisms, but in a quite different manner than the multivalued morphisms of Chapter 3. They are well-behaved under composition, refinement and fibre products. We will see in Section 6.4 that *rigidifications* can be added to them to make them unique, at least if we work with Nisnevich covers.

While we aim to use their properties only for Nisnevich covers of smooth varieties, we nonetheless state more general versions for pseudocovers as they come with no additional complexity.

**Definition 5.2.1.** Let $\alpha\colon X \rightsquigarrow Y$ be a finite correspondence over a noetherian scheme $S$.

- A *pylon* over $\alpha$ is a pair $(\gamma, \widetilde{\alpha})$ consisting of:
  - a morphism $\gamma\colon \Gamma \to X \times_S Y$ of schemes over $S$ such that the composition $\mathrm{pr}_X \circ \gamma\colon \Gamma \to X$ is finite,
  - a relative cycle $\widetilde{\alpha} \in c(\Gamma|X)$ with pushforward $\gamma_*(\widetilde{\alpha}) = \alpha$.
- Let $f\colon \mathcal{U} \twoheadrightarrow X$ and $g\colon \mathcal{V} \twoheadrightarrow Y$ be pseudocovers. A *bridging* over a pylon $(\gamma\colon \Gamma \to X \times_S Y, \widetilde{\alpha})$ over $\alpha$ between these pseudocovers is a morphism

  $$b\colon \mathcal{U} \times_X \Gamma \to \Gamma \times_Y V$$

  of schemes over $\Gamma$.
- A *bridge* over $\alpha\colon X \rightsquigarrow Y$ is a tuple

  $$(f, g, \gamma, \widetilde{\alpha}, b)$$

  consisting of:
  - two pseudocovers $f$ and $g$ of $X$ and $Y$, respectively,
  - a pylon $(\gamma, \widetilde{\alpha})$ over $\alpha$,
  - a bridging $b$ over this pylon between the two pseudocovers.

**Remark 5.2.2.** The relation between the given data is best understood in the following commutative diagram:

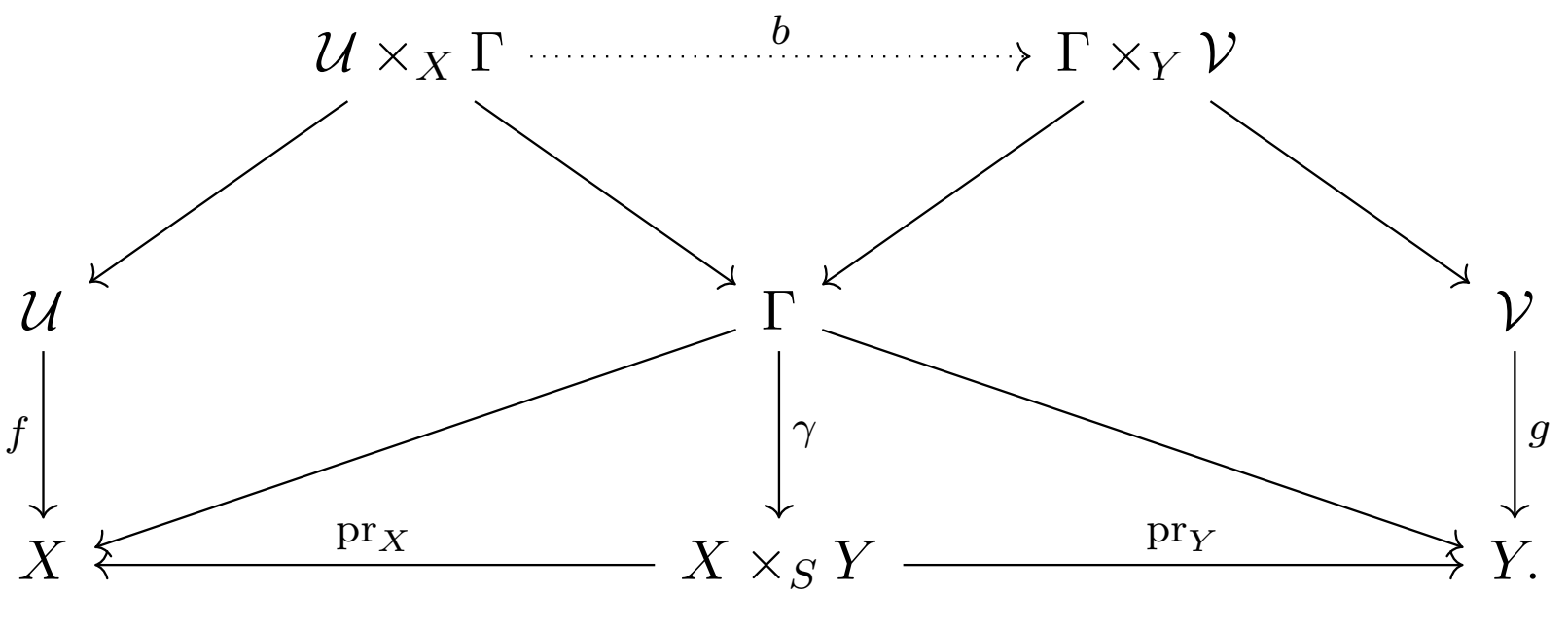

One should understand a pylon as a sort of support that is not required to be inside $X \times_S Y$. This allows more degrees of freedom, which will allow us to define composition directly at the level of bridges.

**Remark 5.2.3.** We do not require that the bridging is compatible with the given decompositions on the pullbacks of the pseudocovers to $\Gamma$. This is not without reason, as it is often impossible to achieve such a strong condition when constructing a bridge, for example in Proposition 6.4.6.

We now explain how bridges relate to finite correspondences.

**Definition 5.2.4.** Let $\alpha\colon X \rightsquigarrow Y$ be a finite correspondence over a noetherian scheme $S$ and let

$$\begin{aligned} f &\colon \mathcal{U} \twoheadrightarrow X, \\ g &\colon \mathcal{V} \twoheadrightarrow Y \end{aligned}$$

be pseudocovers. Also let

$$\mathcal{B} = (f, g, \gamma\colon \Gamma \to X \times_S Y, \widetilde{\alpha}, b)$$

be a bridge over $\alpha$. The projections $\mathrm{pr}_{\mathcal{U}}\colon \mathcal{U} \times_X \Gamma \to \mathcal{U}$ and $\mathrm{pr}_{\mathcal{V}}\colon \Gamma \times_Y \mathcal{V} \to \mathcal{V}$ give us a morphism

$$(\mathrm{pr}_{\mathcal{U}}, \mathrm{pr}_{\mathcal{V}} \circ b) = (\mathrm{pr}_{\mathcal{U}}, \mathrm{pr}_{\mathcal{V}} \circ b)_S : \mathcal{U} \times_X \Gamma \to \mathcal{U} \times_S \mathcal{V}.$$

Using this morphism one defines the *associated finite correspondence*

$$\alpha_{\mathcal{B}} := (\mathrm{pr}_{\mathcal{U}}, \mathrm{pr}_{\mathcal{V}} \circ b)_* f^{\circledast}(\widetilde{\alpha}) \in c(\mathcal{U} \times_S \mathcal{V} | \mathcal{U})$$

of $\mathcal{B}$ over $\alpha$.

We call a finite correspondence $\mathcal{U} \rightsquigarrow \mathcal{V}$ *bridgeable* if it is the associated finite correspondence of a bridge $\mathcal{B}$ over $\alpha$.

**Remark 5.2.5.** The morphism $(\mathrm{pr}_{\mathcal{U}}, \mathrm{pr}_{\mathcal{V}} \circ b) : \mathcal{U} \times_X \Gamma \to \mathcal{U} \times_S \mathcal{V}$ of Definition 5.2.4 has a different yet longer description: it can be written as the composition

$$\begin{aligned} \mathcal{U} \times_X \Gamma \hookrightarrow (\mathcal{U} \times_X \Gamma) \times_\Gamma (\Gamma \times_Y \mathcal{V}) \cong \mathcal{U} \times_X \Gamma \times_Y \mathcal{V} \to \\ \to \mathcal{U} \times_X (X \times_S Y) \times_Y \mathcal{V} \cong \mathcal{U} \times_S \mathcal{V}, \end{aligned}$$

the first morphism being the graph of $b$ over $\Gamma$ and the second one being the obvious base change of $\gamma\colon \Gamma \to X \times_S Y$. In particular it is finite, and if $\gamma$ is a closed immersion, then so is $(\mathrm{pr}_{\mathcal{U}}, \mathrm{pr}_{\mathcal{V}} \circ b)$.

**Proposition 5.2.6.** *Let $\alpha\colon X \rightsquigarrow Y$ be a finite correspondence over a noetherian scheme $S$. Let*

$$\mathcal{B} = (f\colon \mathcal{U} \twoheadrightarrow X, g\colon \mathcal{V} \twoheadrightarrow Y, \gamma, \widetilde{\alpha}, b)$$

*be a bridge over $\alpha$.*

*Then the associated finite correspondence $\alpha_{\mathcal{B}}$ completes the commutative square*

$$\begin{array}{ccc} \mathcal{U} & \overset{\alpha_B}{\dashrightarrow} & \mathcal{V} \\ {\scriptstyle f}\downarrow & & \downarrow{\scriptstyle g} \\ X & \overset{\alpha}{\rightsquigarrow} & Y, \end{array}$$

*i.e. we have*

$$\alpha \circ f = g \circ \alpha_{\mathcal{B}}.$$

*Proof.* Let $(\gamma : \Gamma \to X \times_S Y, \widetilde{\alpha})$ be the pylon of the bridge. The diagram

$$\begin{array}{ccc} \mathcal{U} \times_X \Gamma & \xrightarrow{(\mathrm{pr}_U, \mathrm{pr}_V \circ b)} & \mathcal{U} \times_S \mathcal{V} \\ \downarrow{\scriptstyle \mathcal{U} \times_X \gamma} & & \downarrow{\scriptstyle \mathcal{U} \times_S g} \\ \mathcal{U} \times_X (X \times_S Y) & = & \mathcal{U} \times_S Y \end{array}$$

commutes, thus we get from Proposition 1.6.10 and Proposition 1.5.3 that

$$\begin{aligned} g \circ \alpha_B &= (\mathcal{U} \times_S g)_* \, (\mathrm{pr}_U, \mathrm{pr}_V \circ b)_* \, f^{\circledast}(\widetilde{\alpha}) = \\ &= (\mathcal{U} \times_X \gamma)_* \, f^{\circledast}(\widetilde{\alpha}) = \\ &= f^{\circledast}\gamma_*(\widetilde{\alpha}) = f^{\circledast}\alpha = \alpha \circ f \end{aligned}$$

as claimed. □

**Definition 5.2.7** (Refinement of bridges)**.** Let $\alpha\colon X \rightsquigarrow Y$ be a finite correspondence and let

$$\begin{aligned} \mathcal{B}_1 &= (f_1\colon \mathcal{U}_1 \twoheadrightarrow X, g_1\colon \mathcal{V}_1 \twoheadrightarrow Y, \gamma_1\colon \mathcal{U}_1 \times_X \Gamma_1 \to X \times_S Y, \widetilde{\alpha}_1, b_1) \\ \mathcal{B}_2 &= (f_2\colon \mathcal{U}_2 \twoheadrightarrow X, g_2\colon \mathcal{V}_2 \twoheadrightarrow Y, \gamma_2\colon \mathcal{U}_2 \times_X \Gamma_2 \to X \times_S Y, \widetilde{\alpha}_2, b_2) \end{aligned}$$

be bridges over $\alpha$. A *morphism $\mathcal{B}_2 \to \mathcal{B}_1$ of bridges* over $\alpha$, or a *refinement* of $\mathcal{B}_1$ into $\mathcal{B}_2$, is a triple $(u, v, \delta)$ of morphisms

$$u\colon \mathcal{U}_2 \to \mathcal{U}_1, \qquad v\colon \mathcal{V}_2 \to \mathcal{V}_1, \qquad \delta\colon \Gamma_2 \to \Gamma_1$$

such that:

- $u$ and $v$ are morphisms of pseudocovers over $X$ and $Y$, respectively,

- $\gamma_2 = \gamma_1 \circ \delta$,
- $\delta_* (\widetilde{\alpha}_2) = \widetilde{\alpha}_1$,
- the diagram

$$\begin{array}{ccc} \mathcal{U}_2 \times_X \Gamma_2 & \xrightarrow{b_2} & \Gamma_2 \times_Y \mathcal{V}_2 \\ \big\downarrow{\scriptstyle u \times_X \delta} & & \big\downarrow{\scriptstyle \delta \times_Y v} \\ \mathcal{U}_1 \times_X \Gamma_1 & \xrightarrow{b_1} & \Gamma_1 \times_Y \mathcal{V}_1 \end{array}$$

  commutes.

Composition of refinements of bridges is defined in the obvious way.

**Remark 5.2.8.** We illustrate the setting with the following commutative diagram:

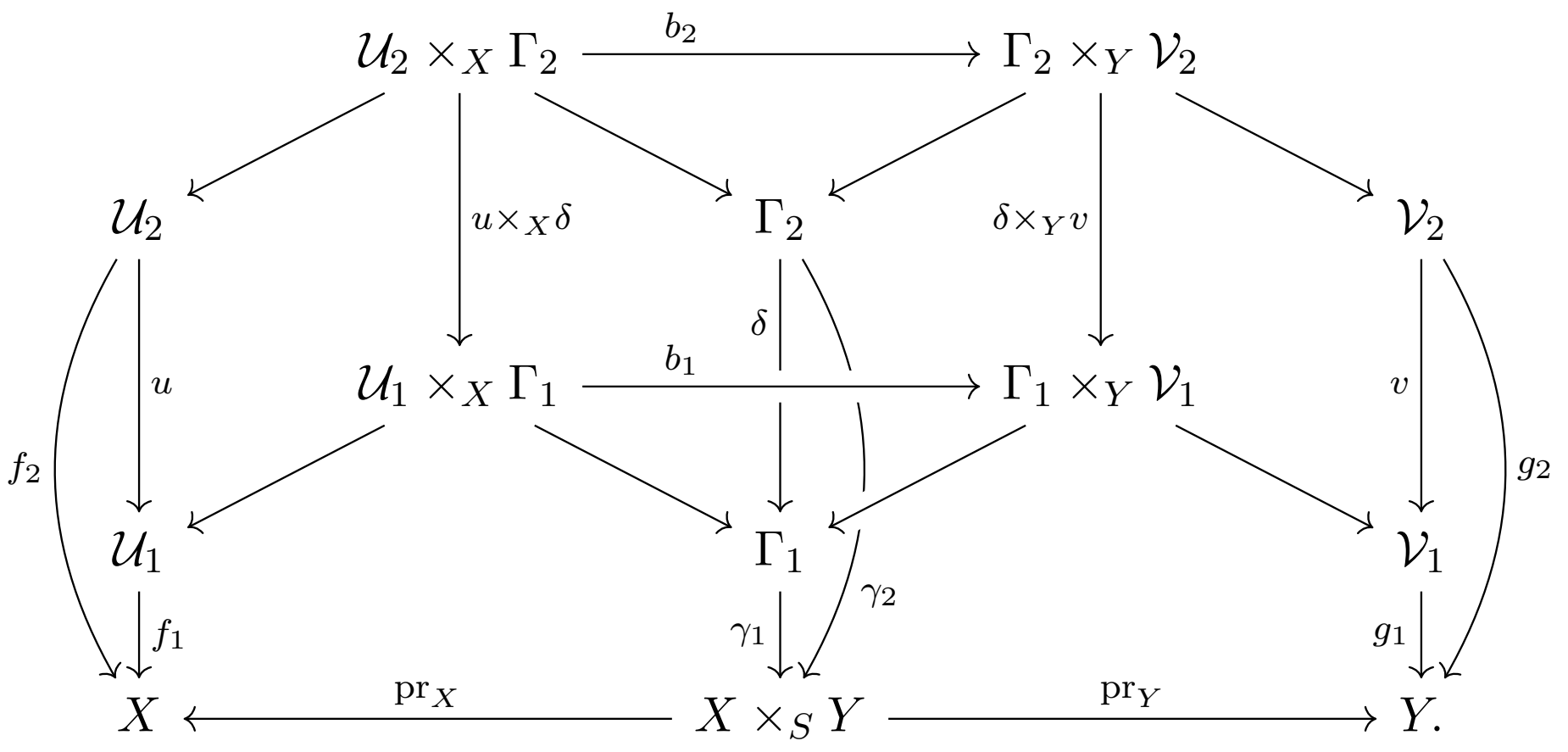

**Remark 5.2.9.** The definition of a refinement almost splits into two parts, a refinement of pseudocovers and a refinement of pylons, if defined accordingly. The only thing connecting them is the commutative diagram at the end of Definition 5.2.7. We will take this decomposition into two parts to good use in the proof of Proposition 5.2.10.

**Proposition 5.2.10.** *Let everything be as in Definition 5.2.7. Then*

$$\alpha_{\mathcal{B}_1} \circ u = v \circ \alpha_{\mathcal{B}_2},$$

*i.e the commutative diagram of Proposition 5.2.6 extends to a commutative diagram*

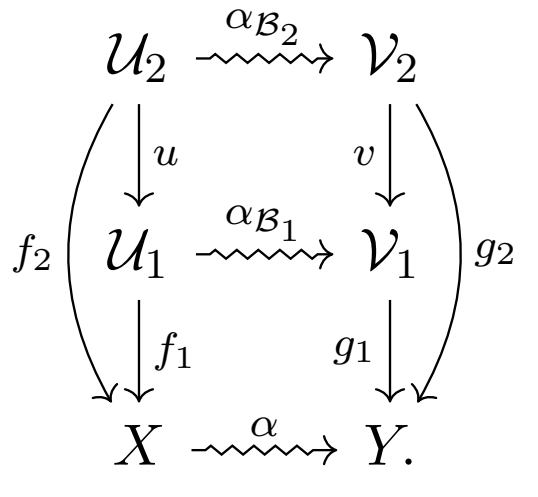

*Proof.* The diagram

$$\begin{array}{ccc}
\mathcal{U}_2 \times_X \Gamma_2 & \xrightarrow{b_2} & \Gamma_2 \times_Y \mathcal{V}_2 \\
{\scriptstyle u\times_X\Gamma_2}\downarrow & & \downarrow{\scriptstyle \Gamma_2\times_Y v} \\
\mathcal{U}_1 \times_X \Gamma_2 & \xrightarrow{b_1\times_{\Gamma_1}\Gamma_2} & \Gamma_2 \times_Y \mathcal{V}_1 \\
{\scriptstyle \mathcal{U}_1\times_X\delta}\downarrow & & \downarrow{\scriptstyle \delta\times_Y\mathcal{V}_1} \\
\mathcal{U}_1 \times_X \Gamma_1 & \xrightarrow{b_1} & \Gamma_1 \times_Y \mathcal{V}_1
\end{array}$$

commutes, which follows easily from the definition of a morphism of bridges. Thus there is an intermediate bridge

$$\mathcal{B}_{1,2} := (f_1, g_1, \gamma_2, \widetilde{\alpha}_2, b_1 \times_{\Gamma_1} \delta)$$

which corresponds to the factorization $(u, v, \delta) = (\mathrm{id}_{\mathcal{U}_1}, \mathrm{id}_{\mathcal{V}_1}, \delta) \circ (u, v, \mathrm{id}_{\Gamma_2})$ of refinements of bridges. Hence it suffices to consider the following two special cases:

(a) The covers are identical, i.e.

$$\begin{aligned}
\mathcal{U} &:= \mathcal{U}_1 = \mathcal{U}_2, & \mathcal{V} &:= \mathcal{V}_1 = \mathcal{V}_2, \\
u &= \mathrm{id}_{\mathcal{U}}, & v &= \mathrm{id}_{\mathcal{V}}, \\
f &:= f_1 = f_2, & g &:= g_1 = g_2.
\end{aligned}$$

We have $b_1 \circ (\mathcal{U} \times_X \delta) = (\delta \times_Y \mathcal{V}) \circ b_2$. Using this and Proposition 1.5.3 we get

$$\begin{aligned}
\alpha_{B_1} \circ u = \alpha_{B_1} &\overset{\text{Def.}}{=} (\mathrm{pr}_{\mathcal{U}}, \mathrm{pr}_{\mathcal{V}} \circ b_1)_* f^{\circledast} \delta_* \widetilde{\alpha}_2 = \\
&= (\mathrm{pr}_{\mathcal{U}}, \mathrm{pr}_{\mathcal{V}} \circ b_1)_* (\mathcal{U} \times_X \delta)_* f^{\circledast} \widetilde{\alpha}_2 = \\
&= (\mathrm{pr}_{\mathcal{U}}, \mathrm{pr}_{\mathcal{V}} \circ (\delta \times_Y \mathcal{V}) \circ b_2)_* f^{\circledast} \widetilde{\alpha}_2 = \\
&= (\mathrm{pr}_{\mathcal{U}}, \mathrm{pr}_{\mathcal{V}} \circ b_2)_* f^{\circledast} \widetilde{\alpha}_2 \overset{\text{Def.}}{=} \alpha_{B_2} = v \circ \alpha_{B_2}.
\end{aligned}$$

(b) The pylons are identical, i.e.

$$\begin{aligned}
\Gamma &:= \Gamma_1 = \Gamma_2, \\
\gamma &:= \gamma_1 = \gamma_2, \\
\delta &= \mathrm{id}_{\Gamma}, \\
\widetilde{\alpha} &:= \widetilde{\alpha}_1 = \widetilde{\alpha}_2.
\end{aligned}$$

The conditions readily imply $\mathrm{pr}_{\mathcal{V}_1} \circ b_1 \circ (u \times_X \Gamma) = v \circ \mathrm{pr}_{\mathcal{V}_2} \circ b_2$. Using this, Proposition 1.6.10 and Proposition 1.5.3 we get

$$\begin{aligned}\alpha_{B_1} \circ u = u^{\circledast}\alpha_{B_1} &\overset{\text{Def.}}{=} u^{\circledast}(\mathrm{pr}_{\mathcal{U}_1}, \mathrm{pr}_{\mathcal{V}_1} \circ b_1)_* f_1^{\circledast}\widetilde{\alpha} = \\ &= (\mathcal{U}_2 \times_{\mathcal{U}_1} (\mathrm{pr}_{\mathcal{U}_1}, \mathrm{pr}_{\mathcal{V}_1} \circ b_1))_* u^{\circledast} f_1^{\circledast}\widetilde{\alpha} = \\ &= (\mathrm{pr}_{\mathcal{U}_2}, v \circ \mathrm{pr}_{\mathcal{V}_2} \circ b_2)_*(f_1 \circ u)^{\circledast}\widetilde{\alpha} = \\ &= (\mathcal{U}_2 \times_S v)_*(\mathrm{pr}_{\mathcal{U}_2}, \mathrm{pr}_{\mathcal{V}_2} \circ b_2)_* f_2^{\circledast}\widetilde{\alpha} = v \circ \alpha_{B_2}.\end{aligned}$$

□

**Definition 5.2.11** (Fibre product of bridges)**.** Let $\alpha\colon X \rightsquigarrow Y$ be a finite correspondence over a noetherian scheme $S$. Furthermore let

$$\begin{aligned}\mathcal{B}_1 &= (f_1\colon \mathcal{U}_1 \to X, g_1\colon \mathcal{V}_1 \to Y, \gamma, \widetilde{\alpha}, b_1)\,, \\ \mathcal{B}_2 &= (f_2\colon \mathcal{U}_2 \to X, g_2\colon \mathcal{V}_2 \to Y, \gamma, \widetilde{\alpha}, b_2)\end{aligned}$$

be bridges over $\alpha$ with the same pylon $(\gamma\colon \Gamma \to X \times_S Y, \widetilde{\alpha})$.

The *fibre product* $\mathcal{B}_1 \times \mathcal{B}_2$ of the bridges is defined by taking fibre products over the respective bases:

$\mathcal{B}_1 \times \mathcal{B}_2$ is the bridge over the same pylon $(\gamma, \widetilde{\alpha})$ between the fibre products

$$\begin{aligned}f_1 \times_X f_2 &\colon \mathcal{U}_1 \times_X \mathcal{U}_2 \twoheadrightarrow X, \\ g_1 \times_Y g_2 &\colon \mathcal{V}_1 \times_Y \mathcal{V}_2 \twoheadrightarrow Y\end{aligned}$$

of pseudocovers, the bridging being

$$\begin{aligned}b_1 \times_\Gamma b_2\colon\ \mathcal{U}_1 \times_X \mathcal{U}_2 \times_X \Gamma \cong (\mathcal{U}_1 \times_X \Gamma) \times_\Gamma (\mathcal{U}_2 \times_X \Gamma) \longrightarrow \\ \longrightarrow (\Gamma \times_Y \mathcal{V}_1) \times_\Gamma (\Gamma \times_Y \mathcal{V}_2) \cong \Gamma \times_Y \mathcal{V}_1 \times_Y \mathcal{V}_2.\end{aligned}$$

The fibre product of bridges is clearly associative and commutative in the obvious sense. It is also immediate that the fibre product $\mathcal{B}_1 \times \mathcal{B}_2$ refines $\mathcal{B}_1$ and $\mathcal{B}_2$, thus we immediately get from Proposition 5.2.10:

**Corollary 5.2.12.** *Let everything be as in Definition 5.2.11.*

*Then the diagram*

$$\begin{array}{ccc} \mathcal{U}_1 & \overset{\alpha_{\mathcal{B}_1}}{\rightsquigarrow} & \mathcal{V}_1 \\ {\scriptstyle \mathrm{pr}_{\mathcal{U}_1}}\uparrow & & \uparrow{\scriptstyle \mathrm{pr}_{\mathcal{V}_1}} \\ \mathcal{U}_1 \times_X \mathcal{U}_2 & \overset{\alpha_{\mathcal{B}_1 \times_\Gamma \mathcal{B}_2}}{\rightsquigarrow} & \mathcal{V}_1 \times_Y \mathcal{V}_2 \\ {\scriptstyle \mathrm{pr}_{\mathcal{U}_2}}\downarrow & & \downarrow{\scriptstyle \mathrm{pr}_{\mathcal{V}_2}} \\ \mathcal{U}_2 & \overset{\alpha_{\mathcal{B}_2}}{\rightsquigarrow} & \mathcal{V}_2. \end{array}$$

*commutes.*

**Remark 5.2.13.** A very general notion of pullbacks of relative cycles is introduced in chapter 8 of [CD12]. Our definition of the fibre product of bridges should be understood as a different approach to this. Note that we require a good notion of rigidifications later on, for which bridges will prove suitable as well.

**Definition 5.2.14.** Let $\alpha\colon X \rightsquigarrow Y$ and $\beta\colon Y \rightsquigarrow Z$ be finite correspondences over $S$ and let

$$\begin{aligned} f\colon\ & \mathcal{U} \twoheadrightarrow X, \\ g\colon\ & \mathcal{V} \twoheadrightarrow Y, \\ h\colon\ & \mathcal{W} \twoheadrightarrow Z \end{aligned}$$

be pseudocovers. Furthermore let

$$\begin{aligned} \mathcal{A} &= (f, g, \gamma\colon \Gamma \to X \times_S Y, \widetilde{\alpha}, a), \\ \mathcal{B} &= (g, h, \delta\colon \Delta \to Y \times_S Z, \widetilde{\beta}, b) \end{aligned}$$

be bridges between those pseudocovers over $\alpha$ and $\beta$, respectively.

We define their *composition* $\mathcal{B} \circ \mathcal{A}$ as the tuple

$$(f,\ h,\ \epsilon\colon \Gamma \times_Y \Delta \to X \times_S Z,\ \widetilde{\beta} \odot \widetilde{\alpha},\ b \boxtimes a),$$

where:

- $\epsilon = \mathrm{pr}^{XYZ}_{XZ} \circ (\gamma \times_Y \delta)\colon \Gamma \times_Y \Delta \to X \times_S Z$,
- $b \boxtimes a$ is the composition

$$\mathcal{U} \times_X \Gamma \times_Y \Delta \xrightarrow{a \times_Y \Delta} \Gamma \times_Y \mathcal{V} \times_Y \Delta \xrightarrow{\Gamma \times_Y b} \Gamma \times_Y \Delta \times_Z \mathcal{W},$$

- $\widetilde{\beta} \odot \widetilde{\alpha} = \mathrm{Cor}((\mathrm{pr}^{XY}_{Y} \circ \gamma)^{\circledast} \widetilde{\beta}, \widetilde{\alpha}) \in c(\Gamma \times_Y \Delta | X)$ (cf. Definition 1.4.22).

We now check that this is actually a bridge between $\mathcal{U}$ and $\mathcal{W}$ over $\beta \circ \alpha$.

**Proposition 5.2.15.** *Let everything be as in Definition 5.2.14. Then:*

*(a) The tuple $\mathcal{B} \circ \mathcal{A} = (f, h, \epsilon, \widetilde{\beta} \odot \widetilde{\alpha}, b \boxtimes a)$ as defined in Definition 5.2.14 is a bridge over $\beta \circ \alpha$ between $\mathcal{U}$ and $\mathcal{W}$,*

*(b) We have the equality $(\beta \circ \alpha)_{\mathcal{B} \circ \mathcal{A}} = \beta_{\mathcal{B}} \circ \alpha_{\mathcal{A}}$ of associated finite correspondences $\mathcal{U} \rightsquigarrow \mathcal{W}$.*

*Proof.*

(a) As all morphisms are between the correct schemes, we only have to check that $(\epsilon, \widetilde{\alpha} \odot \widetilde{\beta})$ is indeed a pylon over $\beta \circ \alpha$. By construction the diagram

$$\begin{array}{ccccccc}
\Gamma \times_Y \Delta & & & \xrightarrow{\qquad\qquad\qquad\epsilon\qquad\qquad\qquad} & & & X \times_S Z \\
\downarrow{\scriptstyle \Gamma \times_Y \delta} & & & & & & \\
\Gamma \times_Y Y \times_S Z & = & \Gamma \times_S Z & \xrightarrow{\gamma \times_S Z} & X \times_S Y \times_S Z & \xrightarrow{\mathrm{pr}_{XZ}} & X \times_S Z \\
\downarrow{\scriptstyle \Gamma \times_Y \mathrm{pr}_Y} & & \downarrow{\scriptstyle \mathrm{pr}_\Gamma} & & \downarrow{\scriptstyle \mathrm{pr}_{XY}} & & \downarrow{\scriptstyle \mathrm{pr}_X} \\
\Gamma \times_Y Y & = & \Gamma & \xrightarrow{\gamma} & X \times_S Y & \xrightarrow{\mathrm{pr}_X} & X
\end{array}$$

commutes. Thus $\mathrm{pr}_X^{XZ} \circ \epsilon = (\mathrm{pr}_X^{XY} \circ \gamma) \circ (\Gamma \times_Y (\mathrm{pr}_Y^{YZ} \circ \delta))$ is finite as composition of such morphisms.

We are left with the calculation

$$\begin{aligned}
\epsilon_*(\widetilde{\beta} \odot \widetilde{\alpha}) &\overset{a)}{=} \epsilon_* \operatorname{Cor}\Big((\mathrm{pr}_Y^{XY} \circ \gamma)^{\circledast} \widetilde{\beta},\ \widetilde{\alpha}\Big) = \\
&\overset{b)}{=} (\mathrm{pr}_{XZ}^{XYZ})_* (X \times_S \delta)_* (\gamma \times_Y \Delta)_* \operatorname{Cor}\Big(\gamma^{\circledast} (\mathrm{pr}_Y^{XY})^{\circledast} \widetilde{\beta},\ \widetilde{\alpha}\Big) = \\
&\overset{c)}{=} (\mathrm{pr}_{XZ}^{XYZ})_* (X \times_S \delta)_* \operatorname{Cor}\Big((\mathrm{pr}_Y^{XY})^{\circledast} \widetilde{\beta},\ \gamma_* \widetilde{\alpha}\Big) = \\
&\overset{d)}{=} (\mathrm{pr}_{XZ}^{XYZ})_* \operatorname{Cor}\Big((X \times_S \delta)_* (\mathrm{pr}_Y^{XY})^{\circledast} \widetilde{\beta},\ \gamma_* \widetilde{\alpha}\Big) = \\
&\overset{e)}{=} (\mathrm{pr}_{XZ}^{XYZ})_* \operatorname{Cor}\Big((\mathrm{pr}_Y^{XY})^{\circledast} \delta_* \widetilde{\beta},\ \gamma_* \widetilde{\alpha}\Big) = \\
&\overset{f)}{=} (\mathrm{pr}_{XZ}^{XYZ})_* \operatorname{Cor}\big((\mathrm{pr}_Y^{XY})^{\circledast} \beta,\ \alpha\big) \overset{g)}{=} \beta \circ \alpha,
\end{aligned}$$

where we used the following:

a) Definition 5.2.14 of $\widetilde{\beta} \odot \widetilde{\alpha}$,

b) functoriality of pushforward and pullback (Proposition 1.5.1 and Proposition 1.5.2),

c) compatibility of $\operatorname{Cor}(-,-)$ with pushforward on the right (Proposition 1.5.7),

d) compatibility of $\operatorname{Cor}(-,-)$ with pushforward on the left (Proposition 1.5.8),

e) compatibility of pushforward with pullback (Proposition 1.5.3),

f) Definition 5.2.1 of $\widetilde{\alpha}$ and $\widetilde{\beta}$,

g) Definition 1.6.5 of composition of finite correspondences.

(b) We let

$$\begin{aligned}
\mathrm{pr}_{\mathcal{U}}\colon\ & \mathcal{U}\times_X\Gamma\to\mathcal{U},\\
\mathrm{pr}_{\mathcal{V}}\colon\ & \Gamma\times_Y\mathcal{V}\to\mathcal{V},\\
\mathrm{pr}'_{\mathcal{V}}\colon\ & \mathcal{V}\times_Y\Delta\to\mathcal{V},\\
\mathrm{pr}'_{\mathcal{W}}\colon\ & \Delta\times_Z\mathcal{W}\to\mathcal{W}
\end{aligned}$$

be the projections. For simplicity and readability and use the shorthand

$$\begin{aligned}
\hat{a} &:= (\mathrm{pr}_{\mathcal{U}}, \mathrm{pr}_{\mathcal{V}}\circ a)\colon \mathcal{U}\times_X\Gamma\to\mathcal{U}\times_S\mathcal{V},\\
\hat{b} &:= (\mathrm{pr}'_{\mathcal{V}}, \mathrm{pr}'_{\mathcal{W}}\circ b)\colon \mathcal{V}\times_Y\Delta\to\mathcal{V}\times_S\mathcal{W}
\end{aligned}$$

and

$$\hat{c} := \left(\mathrm{pr}^{\mathcal{U}\times_X\Gamma\times_Y\Delta}_{\mathcal{U}\times_X\Gamma}, \mathrm{pr}'_{\mathcal{W}}\circ(b\boxtimes a)\right)_S : \mathcal{U}\times_X\Gamma\times_Y\Delta\to\mathcal{U}\times_X\Gamma\times_S\mathcal{W}.$$

Note that latter is the pullback of $\hat{b}$ along $\mathrm{pr}_{\mathcal{V}}\circ a\colon \mathcal{U}\times_X\Gamma\to\mathcal{V}$. By examining the individual coordinates it is also clear that

$$(\mathrm{pr}_{\mathcal{U}}\times_S\mathcal{W})\circ\hat{c} = \left(\mathrm{pr}_{\mathcal{U}}, \mathrm{pr}'_{\mathcal{W}}\circ(b\boxtimes a)\right).$$

Using these we again do a lengthy calculation:

$$\begin{aligned}
&\beta_B\circ\alpha_A =\\
&\overset{\text{a)}}{=} (\mathrm{pr}^{\mathcal{UVW}}_{\mathcal{UW}})_*\operatorname{Cor}\left((\mathrm{pr}^{\mathcal{UV}}_{\mathcal{V}})^{\circledast}\beta_B,\ \alpha_A\right) =\\
&\overset{\text{b)}}{=} (\mathrm{pr}^{\mathcal{UVW}}_{\mathcal{UW}})_*\operatorname{Cor}\left((\mathrm{pr}^{\mathcal{UV}}_{\mathcal{V}})^{\circledast}\hat{b}_*g^{\circledast}\widetilde{\beta},\ \hat{a}_*f^{\circledast}\widetilde{\alpha}\right) =\\
&\overset{\text{c)}}{=} (\mathrm{pr}^{\mathcal{UVW}}_{\mathcal{UW}})_*(\hat{a}\times_S\mathcal{W})_*\operatorname{Cor}\left(\hat{a}^{\circledast}(\mathrm{pr}^{\mathcal{UV}}_{\mathcal{V}})^{\circledast}\hat{b}_*g^{\circledast}\widetilde{\beta},\ f^{\circledast}\widetilde{\alpha}\right) =\\
&\overset{\text{d)}}{=} (\mathrm{pr}_{\mathcal{U}}\times_S\mathcal{W})_*\operatorname{Cor}\left((\mathrm{pr}_{\mathcal{V}}\circ a)^{\circledast}\hat{b}_*g^{\circledast}\widetilde{\beta},\ f^{\circledast}\widetilde{\alpha}\right) =\\
&\overset{\text{e)}}{=} (\mathrm{pr}_{\mathcal{U}}\times_S\mathcal{W})_*\operatorname{Cor}\left(\hat{c}_*(\mathrm{pr}_{\mathcal{V}}\circ a)^{\circledast}g^{\circledast}\widetilde{\beta},\ f^{\circledast}\widetilde{\alpha}\right) =\\
&\overset{\text{f)}}{=} (\mathrm{pr}_{\mathcal{U}}\times_S\mathcal{W})_*\hat{c}_*\operatorname{Cor}\left((\mathrm{pr}_{\mathcal{V}}\circ a)^{\circledast}g^{\circledast}\widetilde{\beta},\ f^{\circledast}\widetilde{\alpha}\right) =\\
&\overset{\text{g)}}{=} (\mathrm{pr}_{\mathcal{U}}, \mathrm{pr}'_{\mathcal{W}}\circ(b\boxtimes a))_*\operatorname{Cor}\left((f\times_X\Gamma)^{\circledast}(\mathrm{pr}^{XY}_Y\circ\gamma)^{\circledast}\widetilde{\beta},\ f^{\circledast}\widetilde{\alpha}\right) =\\
&\overset{\text{h)}}{=} (\mathrm{pr}_{\mathcal{U}}, \mathrm{pr}'_{\mathcal{W}}\circ(b\boxtimes a))_*f^{\circledast}\operatorname{Cor}\left((\mathrm{pr}^{XY}_Y\circ\gamma)^{\circledast}\widetilde{\beta},\ \widetilde{\alpha}\right) =\\
&\overset{\text{i)}}{=} (\beta\circ\alpha)_{B\circ A}.
\end{aligned}$$

Here we used:

a) Definition 1.6.5 of composition of finite correspondences,

b) Definition 5.2.1 of $\widetilde{\alpha}$ and $\widetilde{\beta}$,

c) compatibility of $\mathrm{Cor}(-,-)$ with pushforward on the right (Proposition 1.5.7),

d) functoriality of pushforward and pullback (Proposition 1.5.1 and Proposition 1.5.2),

e) compatibility of pushforward with pullback (Proposition 1.5.3),

f) compatibility of $\mathrm{Cor}(-,-)$ with pushforward on the left (Proposition 1.5.8),

g) functoriality of pushforward (Proposition 1.5.1) and pullback (Proposition 1.5.2),

h) compatibility of $\mathrm{Cor}(-,-)$ with pullback (Proposition 1.5.9),

i) Definition 5.2.14 of $B \circ A$. □

**Proposition 5.2.16.** *Let $\alpha\colon X \rightsquigarrow Y$ and $\beta\colon Y \rightsquigarrow Z$ be finite correspondences over a noetherian scheme $S$ and let*

$$\begin{aligned} f_1\colon\ \mathcal{U}_1 &\twoheadrightarrow X, \\ f_2\colon\ \mathcal{U}_2 &\twoheadrightarrow X, \\ g_1\colon\ \mathcal{V}_1 &\twoheadrightarrow Y, \\ g_2\colon\ \mathcal{V}_2 &\twoheadrightarrow Y, \\ h_1\colon \mathcal{W}_1 &\twoheadrightarrow Z, \\ h_2\colon \mathcal{W}_2 &\twoheadrightarrow Z \end{aligned}$$

*be pseudocovers. Let $(\gamma\colon \Gamma \to X \times_S Y, \widetilde{\alpha})$ be a pylon over $\alpha$ and let $(\delta\colon \Delta \to Y \times_S Z, \widetilde{\beta})$ by a pylon over $\beta$. Furthermore let*

$$\begin{aligned} \mathcal{A}_1 &= (f_1, g_1, \gamma, \widetilde{\alpha}, a_1), \\ \mathcal{A}_2 &= (f_2, g_2, \gamma, \widetilde{\alpha}, a_2), \\ \mathcal{B}_1 &= (g_1, h_1, \delta, \widetilde{\beta}, b_1), \\ \mathcal{B}_2 &= (g_2, h_2, \delta, \widetilde{\beta}, b_2) \end{aligned}$$

*be two sets of bridges between those pseudocovers over the two chosen pylons, respectively.*

*Then we have an equality*

$$(\mathcal{B}_1 \circ \mathcal{A}_1) \times (\mathcal{B}_2 \circ \mathcal{A}_2) = (\mathcal{B}_1 \times \mathcal{B}_2) \circ (\mathcal{A}_1 \times \mathcal{A}_2)$$

*of bridges between $\mathcal{U}$ and $\mathcal{W}$ over $\beta \circ \alpha$.*

*Proof.* Trivially, both bridges in question involve the same pseudocovers. They have the same pylons because fibre products do not change them. Hence the equality of the bridges reduces to that of the bridgings: we want

$$\begin{aligned}((a_1 \times_Y \Delta) \circ (\Gamma \times_Y b_1)) \times_{\Gamma \times_Y \Delta} ((a_2 \times_Y \Delta) \circ (\Gamma \times_Y b_2)) = \\ = (a_1 \times_\Gamma a_2 \times_Y \Delta) \circ (\Gamma \times_Y b_1 \times_\Delta b_2)\end{aligned}$$

as morphisms

$$\mathcal{U}_1 \times_X \mathcal{U}_2 \times_X \Gamma \times_Y \Delta \to \Gamma \times_Y \Delta \times_Z \mathcal{W}_1 \times_Z \mathcal{W}_2.$$

This is a lengthy but straightforward equality of fibre products of morphisms. □

**Definition 5.2.17** (Tensor product of bridges)**.** Let $\alpha_1 \colon X_1 \rightsquigarrow Y_1$ and $\alpha_2 \colon X_2 \rightsquigarrow Y_2$ be finite correspondences over a noetherian scheme $S$. Furthermore let

$$\begin{aligned}\mathcal{B}_1 &= (f_1 \colon \mathcal{U}_1 \to X_1, g_1 \colon \mathcal{V}_1 \to Y_1, \gamma_1 \colon \Gamma_1 \to X_1 \times_S Y_1, \widetilde{\alpha}_1, b_1)\,, \\ \mathcal{B}_2 &= (f_2 \colon \mathcal{U}_2 \to X_2, g_2 \colon \mathcal{V}_2 \to Y_2, \gamma_2 \colon \Gamma_2 \to X_2 \times_S Y_2, \widetilde{\alpha}_2, b_2)\end{aligned}$$

be bridges over them.

Their *tensor product* $\mathcal{B}_1 \otimes \mathcal{B}_2$ is defined as follows:

We have a pylon

$$\gamma_1 \times_S \gamma_2 \colon \Gamma_1 \times_S \Gamma_2 \to X_1 \times_S Y_1 \times_S X_2 \times_S Y_2 \cong (X_1 \times_S X_2) \times_S (Y_1 \times_S Y_2),$$

the cycle being

$$\widetilde{\alpha}_1 \otimes \widetilde{\alpha}_2 := (\mathrm{pr}^{X_1X_2}_{X_1})^\circledast \widetilde{\alpha}_1 \otimes_{X_1 \times_S X_2} (\mathrm{pr}^{X_1X_2}_{X_2})^\circledast \widetilde{\alpha}_2.$$

The bridging is given by $b_1 \times_S b_2$.

This again behaves as expected:

**Proposition 5.2.18.** *Let everything be as in Definition 5.2.17. Then:*

*(a) The tuple*

$$\mathcal{B}_1 \otimes \mathcal{B}_2 = (f_1 \times_S f_2, g_1 \times_S g_2, \gamma_1 \times_S \gamma_2, \widetilde{\alpha}_1 \otimes \widetilde{\alpha}_2, b_1 \times_S b_2)$$

*as defined in Definition 5.2.17 is a bridge over $\alpha_1 \otimes \alpha_2$ between the pseudocovers*

$$\begin{aligned}f_1 \times_S f_2 \colon \mathcal{U}_1 \times_S \mathcal{U}_1 \twoheadrightarrow X_1 \times_S X_2, \\ g_1 \times_S g_2 \colon \mathcal{V}_1 \times_S \mathcal{V}_1 \twoheadrightarrow Y_1 \times_S Y_2.\end{aligned}$$

(b) *We have the equality* $(\alpha_1)_{\mathcal{B}_1} \otimes (\alpha_2)_{\mathcal{B}_2} = (\alpha_1 \otimes \alpha_2)_{\mathcal{B}_1 \otimes \mathcal{B}_2}$ *of associated finite correspondences* $\mathcal{U}_1 \times_S \mathcal{U}_2 \rightsquigarrow \mathcal{V}_1 \times_S \mathcal{V}_2$.

*Proof.* By Remark 1.6.3 we identify $\widetilde{\alpha}_1$, $\widetilde{\alpha}_2$ and $\widetilde{\alpha}_1 \otimes \widetilde{\alpha}_2$ with finite correspondences $X_1 \rightsquigarrow \Gamma_1$, $X_2 \rightsquigarrow \Gamma_2$ and $X_1 \times_S X_2 \rightsquigarrow \Gamma_1 \times_S \Gamma_2$, respectively. In this sense, the notation $\widetilde{\alpha}_1 \otimes \widetilde{\alpha}_2$ is suggestive, as it becomes identified with this tensor product taken in $\mathrm{SchCor}_S$. Indeed, this follows from Remark 1.7.5 and the results of Section 1.5 by a straightforward calculation.

(a) By our prior identifications and Proposition 1.6.10 we have

$$\begin{aligned} &(\gamma_1 \times_S \gamma_2)_*(\widetilde{\alpha}_1 \otimes \widetilde{\alpha}_2) = \\ =&[(\mathrm{pr}_{Y_1} \circ \gamma_1) \times_S (\mathrm{pr}_{Y_2} \circ \gamma_2)] \circ (\widetilde{\alpha}_1 \otimes \widetilde{\alpha}_2) = \\ =&\big([\mathrm{pr}_{Y_1} \circ \gamma_1] \otimes [\mathrm{pr}_{Y_2} \circ \gamma_2]\big) \circ (\widetilde{\alpha}_1 \otimes \widetilde{\alpha}_2) = \\ =&(\mathrm{pr}_{Y_1} \circ \gamma_1 \circ \widetilde{\alpha}_1) \otimes (\mathrm{pr}_{Y_2} \circ \gamma_2 \circ \widetilde{\alpha}_2) = \\ =&(\gamma_1)_*\widetilde{\alpha}_1 \otimes (\gamma_2)_*\widetilde{\alpha}_2 = \alpha_1 \otimes \alpha_2. \end{aligned}$$

(b) We interpret $f_1^{\circledast}\widetilde{\alpha}_1$ as a finite correspondence $\mathcal{U}_1 \rightsquigarrow \mathcal{U}_1 \times_{X_1} \Gamma_1$ over $S$. Then it is easily seen that

$$(\alpha_1)_{\mathcal{B}_1} = [\mathrm{pr}_{\mathcal{V}_1}] \circ [b_1] \circ f_1^{\circledast}\widetilde{\alpha}_1.$$

Analogous statements hold for $\widetilde{\alpha}_2$ and $\widetilde{\alpha}_1 \otimes \widetilde{\alpha}_2$. From Remark 1.7.5 and Proposition 1.5.9 one concludes that

$$f_1^{\circledast}\widetilde{\alpha}_1 \otimes f_2^{\circledast}\widetilde{\alpha}_2 = (f_1 \times_S f_2)^{\circledast}(\widetilde{\alpha}_1 \otimes \widetilde{\alpha}_2).$$

Now the claimed equality follows by another calculation involving the functoriality of the tensor product:

$$\begin{aligned} &(\alpha_1)_{\mathcal{B}_1} \otimes (\alpha_2)_{\mathcal{B}_2} = \\ =&\big([\mathrm{pr}_{\mathcal{V}_1}] \circ [b_1] \circ f_1^{\circledast}\widetilde{\alpha}_1\big) \otimes \big([\mathrm{pr}_{\mathcal{V}_2}] \circ [b_2] \circ f_2^{\circledast}\widetilde{\alpha}_2\big) = \\ =&\big([\mathrm{pr}_{\mathcal{V}_1}] \otimes [\mathrm{pr}_{\mathcal{V}_2}]\big) \circ ([b_1] \otimes [b_2]) \circ \big((f_1^{\circledast}\widetilde{\alpha}_1) \otimes (f_2^{\circledast}\widetilde{\alpha}_2)\big) = \\ =&[\mathrm{pr}_{\mathcal{V}_1 \times_S \mathcal{V}_2}] \circ [b_1 \times_S b_2] \circ (f_1 \times_S f_2)^{\circledast}(\widetilde{\alpha}_1 \otimes \widetilde{\alpha}_2) = \\ =&(\alpha_1 \otimes \alpha_2)_{\mathcal{B}_1 \otimes \mathcal{B}_2}. \end{aligned}$$

□

**Remark 5.2.19.** The three structures on bridges, composition, fibre product and tensor product, also satisfy the other functorialities and compatibilities one would expected from them. For example, all three are associative, and the two products are symmetric. The proofs of these statements are all of a similar nature than those found in this section. We, however, do not need them in what is to come and therefore omit them.

## 5.3 Standard covers

We now give a brief reminder on some standard covers and describe a new one. We will only work with noetherian schemes, hence the reader can verify that the implicit assumptions of finiteness we make are of no concern.

**Definition 5.3.1** (Covers)**.** Let $f\colon (\mathcal{U}|I) \twoheadrightarrow X$ be a finitely indexed and jointly surjective pseudocover (cf. Definition 5.1.1) in the category of separated noetherian schemes:

- We call it a *Zariski cover* if each $f_i\colon \mathcal{U}_i \to X$ is an open immersion.

- We call it an *étale cover* if $f\colon \mathcal{U} \to X$ is étale and of finite type.

- We call it a *Nisnevich cover* if it is an étale cover and for every, possibly non-closed, $x \in X$ there is a $u \in \mathcal{U}$ above $x$ where $f\colon \mathcal{U} \to X$ induces an isomorphism $\kappa(u) \cong \kappa(x)$ of residue fields.

- A *unifibrant* Nisnevich cover is an étale cover such that for every, possibly non-closed, $x \in X$ there is an $i \in I$ such that $f_i^{-1}(x)$ consists of a single point $u \in \mathcal{U}_i$ where $f_i$ induces an isomorphism $\kappa(x) \cong \kappa(u)$ of residue fields.

- A *Nisnevich cd-square* is a pullback square

$$\begin{array}{ccc} U \times_X V & \hookrightarrow & V \\ \downarrow & & \downarrow g \\ U & \overset{f}{\hookrightarrow} & X \end{array}$$

  such that $f$ is an open immersion, $g$ is étale and $g$ induces an isomorphism $g^{-1}(X\backslash f(U)) \cong X\backslash f(U)$ between the reduced induced structures on those closed subsets.

  It is immediate that $\{U, V\} \twoheadrightarrow X$ is a unifibrant Nisnevich cover. We call such covers *basic cd-Nisnevich covers.* A *cd-Nisnevich cover* is any pseudocover obtained as a composition of basic cd-Nisnevich covers and isomorphisms.

**Remark 5.3.2.** All the notions of Definition 5.3.1 are standard, except that of unifibrant Nisnevich covers. Conceptually, one could understand them as a slight generalization of cd-Nisnevich covers offering a non-recursive definition.

Theorem A.0.10 shows that unifibrancy has an interesting second description in terms of the reduced Čech complex of Definition 5.4.2. It is this property which will make it very useful in Chapter 7 and which is what originally drew our attention to them.

**Remark 5.3.3.** It is well known that a composition of Zariski, Nisnevich or étale covers is again of the respective type. By definition this is also true for the Nisnevich cd-covers and it is easily verified that unifibrant Nisnevich covers also satisfy this property. It is also easily checked that a pullback of any of these five types of covers results again in such a cover.

Thus they form pretopologies on suitable categories, for example the noetherian separated schemes.

**Definition 5.3.4.** We use Zar, cd-Nis, uni-Nis, Nis and Ét to denote the respective pretopologies of Zariski, cd-Nisnevich, unifibrant Nisnevich, Nisnevich and étale covers.

**Remark 5.3.5.** Every unifibrant Nisnevich cover is clearly a Nisnevich one and every Nisnevich cover is by definition étale. Furthermore, by simple inductions, every Zariski cover is a cd-Nisnevich cover and every cd-Nisnevich cover is a unifibrant Nisnevich cover. Thus we have shown the following chain of inclusions, i.e. refinements, of pretopologies:

$$\text{Zar} \subseteq \text{cd-Nis} \subseteq \text{uni-Nis} \subseteq \text{Nis} \subseteq \text{Ét} \tag{5.1}$$

**Example 5.3.6.** Let us demonstrate that the above inclusions are in general strict. We give three tautological examples followed by four more complicated ones.

Let $k$ be a field:

- Nis $\subsetneq$ Ét: take any non-trivial finite separable extensions $L|k$ to get an étale cover $\mathrm{Spec}(L) \to \mathrm{Spec}(k)$ that is not a Nisnevich cover.

- uni-Nis $\subsetneq$ Nis: cover $\mathrm{Spec}(k)$ by the disjoint union $\mathcal{U} = \mathcal{U}_1 := \mathrm{Spec}(k) \sqcup \mathrm{Spec}(k)$ understood as a single scheme.

- cd-Nis $\subsetneq$ uni-Nis: let $X = a_1 \sqcup a_2$, where the $a_i$ are copies of $\mathrm{Spec}(k)$. Let $\mathcal{U}_1 = a_1 \sqcup a_1 \sqcup a_2$ and $\mathcal{U}_2 = a_1 \sqcup a_2 \sqcup a_2$ and map every copy of $a_i$ to itself. Then this clearly constitutes a unifibrant Nisnevich cover, but it cannot be a cd-Nisnevich cover as neither $\mathcal{U}_1$ nor $\mathcal{U}_2$ is openly immersed into $X$.

Let us now give more serious examples over $\mathbb{A}^1_k = \mathrm{Spec}(k[x])$, where we assume that the characteristic of $k$ is different from 2. Choose a non-zero $a \in k$:

Let $f\colon U = \mathbb{A}^1_k \backslash \{a\} \hookrightarrow \mathbb{A}^1_k$ be the open immersion as a subset. Let also $g\colon V = \mathbb{A}^1_k \backslash \{0\} \to \mathbb{A}^1_k$ be the morphism induced by $x \mapsto x^2$. We consider the pseudocover $f \sqcup g\colon U \sqcup V \to \mathbb{A}^1_k$ consisting of two components:

It is clear that $f$ and $g$ are jointly surjective. The morphism $g$ is not an open immersion, hence $U \sqcup V$ is not a Zariski cover. But $g$ is étale, thus $U \sqcup V$ forms an étale cover.

The fibre $g^{-1}(a)$ at $a$ is isomorphic to $\operatorname{Spec}(k[x]/(x^2-a))$, therefore $U \sqcup V$ is a Nisnevich cover if and only if $a$ is a square in $k$, showing the strictness of the rightmost inclusion of (equation (5.1)). If now $a$ is indeed a square in $k$, then $U \sqcup V$ is still not a unifibrant Nisnevich cover as the fibre has two points $\widetilde{a}_1, \widetilde{a}_2$, proving once more that the third inclusion of (equation (5.1)) is a proper one.

Putting $V' = V \backslash \{\widetilde{a}_1\}$ we get from the induced morphism $g' \colon V' \hookrightarrow V \to \mathbb{A}^1_k$ a basic cd-Nisnevich cover $\{U, V'\}$. As $g'$ is still not an open immersion, we see that the leftmost inclusion is strict as well.

Lastly, we pick two non-zero $a_1, a_2 \in k$. We consider the cover given by the open immersion $f \colon U = \mathbb{A}^1_k \backslash \{a_1^2, a_2^2\}$ and the two morphisms $g_i \colon V_i := \mathbb{A}^1_k \backslash \{0, a_i\} \to \mathbb{A}^1_k$, $i \in \{1, 2\}$, induced by $x \mapsto x^2$. By the same arguments as before we see that they together form a Nisnevich cover. Even more, $g_i^{-1}(a_i^2)$ consists of a single point at which $g_i$ induces an isomorphism, hence they form a unifibrant Nisnevich cover. If this cover were a cd-Nisnevich cover, then the same would be true for its pullback along the closed immersion $\{a_1, a_2\} \hookrightarrow \mathbb{A}^1_k$, but this results in the tautological counterexample we gave earlier.

Despite Example 5.3.6 there exists an almost-converse for the middle inclusions of (equation (5.1)):

**Lemma 5.3.7.** *Every Nisnevich cover has a refinement into a cd-Nisnevich cover.*

*Proof.* This is Proposition 2.16 of [Voe10b]. □

Note that we will mimic the proof of loc. cit. to show a slightly stronger version stated below as Lemma 6.2.5.

## 5.4 Čech complexes

For each pseudocover in a category admitting fibre products one has the usual Čech nerve, but there is also a variant without repetition. If the pseudocover is ordered, a second distinction between increasing and arbitrary tuples emerges. We only need the largest and smallest versions, but the interested reader should find it easy to adapt the definitions. As we do not need the Čech nerves as simplicial objects, but only the associated complexes, we directly skip to the latter.

We fix a category $\mathcal{C}$ that admits fibre products.

**Definition 5.4.1** (Full Čech complex)**.** The *(full) Čech complex* of a pseudocover $f \colon (\mathcal{U}|I) \twoheadrightarrow X$ in $\mathcal{C}$ is the complex

$$\check{C}^\bullet(\mathcal{U} \twoheadrightarrow X) := \left\{ \ldots \longrightarrow \coprod_{i,j,k} \mathcal{U}_i \times_X \mathcal{U}_j \times_X \mathcal{U}_k \longrightarrow \coprod_{i,j} \mathcal{U}_i \times_X \mathcal{U}_j \longrightarrow \coprod_i \mathcal{U}_i \right\}$$

in $C^-(\mathbb{Z}[\mathcal{C}])$. Here

$$\check{\mathrm{C}}^{-n}(\mathcal{U}) = (\mathcal{U}|X)^{n+1} = \coprod_{(i_0,\ldots,i_n)\in I^n} \mathcal{U}_{i_0} \times_X \ldots \times_X \mathcal{U}_{i_n}$$

for $n \in \mathbb{N}_0$ and the boundary maps are the signed sums

$$d\colon\ \check{\mathrm{C}}^{-n}(\mathcal{U}) \to \check{\mathrm{C}}^{-n+1}(\mathcal{U}) = \sum_{m=0}^{n} (-1)^m \operatorname{pr}_m$$

of the projections forgetting the $m$-th factor.

We often simply write $\check{\mathrm{C}}^\bullet(\mathcal{U} \twoheadrightarrow X)$, or even just $\check{\mathrm{C}}^\bullet(\mathcal{U})$ if $X$ is clear from the context.

The defining morphism $f\colon\ \check{\mathrm{C}}^0(\mathcal{U} \twoheadrightarrow X) = \mathcal{U} \to X$ induces a morphism $f[0]\colon\ \check{\mathrm{C}}^\bullet(\mathcal{U}) \to X[0]$, i.e. to $X$ placed in degree 0. Its cone is the *augmented Čech complex* $\underline{\check{\mathrm{C}}}^\bullet(\mathcal{U})$.

**Definition 5.4.2** (Reduced Čech complex)**.** The *reduced Čech complex* of a pseudocover $f\colon (\mathcal{U}|I) \twoheadrightarrow X$ ordered by $\leq$ on $I$ is the subcomplex of $\check{\mathrm{C}}^\bullet(\mathcal{U})$ corresponding to strictly increasing tuples, i.e. given by

$$\check{\mathrm{c}}^\bullet_\leq(\mathcal{U} \twoheadrightarrow X) := \left\{ \ldots \to \coprod_{i<j<k} \mathcal{U}_i \times_X \mathcal{U}_j \times_X \mathcal{U}_k \longrightarrow \coprod_{i<j} \mathcal{U}_i \times_X \mathcal{U}_j \longrightarrow \coprod_i \mathcal{U}_i \right\}.$$

We often simply write $\check{\mathrm{c}}^\bullet(\mathcal{U} \twoheadrightarrow X)$, or even just $\check{\mathrm{c}}^\bullet(\mathcal{U})$, if the ordering or $X$ are clear from the context.

The defining morphism $f\colon\ \check{\mathrm{c}}^0(\mathcal{U} \twoheadrightarrow X) = \mathcal{U} \to X$ induces a morphism $f[0]\colon\ \check{\mathrm{c}}^\bullet(\mathcal{U}) \to X[0]$. Its cone is the *augmented reduced Čech complex* $\underline{\check{\mathrm{c}}}^\bullet(\mathcal{U})$.

**Lemma 5.4.3.** *Let $X \in \mathcal{C}$ and let $f\colon (\mathcal{U}|I) \twoheadrightarrow X$ be a pseudocover. Let $\leq$ and $\preceq$ be two orders on $I$.*

*Then there is a natural isomorphism of complexes*

$$\check{\mathrm{c}}^\bullet_\leq(\mathcal{U} \twoheadrightarrow X) \cong \check{\mathrm{c}}^\bullet_\preceq(\mathcal{U} \twoheadrightarrow X)$$

*in $C^-(\mathbb{Z}[\mathcal{C}])$.*

*Proof.* Let $J \subseteq I$ be a finite set with elements $j_0 < j_1 < \ldots < j_s$. Let $\sigma\colon J \to J$ be the permutation reordering $J$ such that $\sigma(j_m) \prec \sigma(j_n)$ for all $m < n$. We consider the morphism

$$\mathcal{U}_{j_0} \times_X \ldots \times_X \mathcal{U}_{j_s} \to \mathcal{U}_{\sigma(j_0)} \times_X \ldots \times_X \mathcal{U}_{\sigma(j_s)}$$

that is $\operatorname{sgn}(\sigma)$ times the obvious re-ordering of the factors. Varying $J$ now defines the desired isomorphism, as is easily checked. □

By Proposition 1.6.10 we can understand the Čech complexes of pseudocovers in $\mathrm{Sch}_S$ as elements of $C^-(\mathrm{SchCor}(S))$. One of the central goals of bridges was to enable us the following construction:

**Definition 5.4.4** (Finite correspondences between Čech complexes)**.** Let $\alpha\colon X \rightsquigarrow Y$ be a finite correspondence over a noetherian scheme $S$. Let

$$\begin{aligned} f\colon (\mathcal{U}|I) &\twoheadrightarrow X, \\ g\colon (\mathcal{V}|J) &\twoheadrightarrow Y \end{aligned}$$

be pseudocovers and let $\mathcal{B}$ be a bridge over $\alpha$ between those pseudocovers. By Definition 5.2.11 we get bridges

$$\mathcal{B}^n := \underbrace{\mathcal{B} \times \ldots \times \mathcal{B}}_{n \text{ times}}$$

over $\alpha$ between $(\mathcal{U}|X)^n$ and $(\mathcal{V}|Y)^n$.

Hence we get from Definition 5.2.4 an associated finite correspondences

$$\check{\mathrm{C}}^{-n}(\mathcal{B}) := \alpha_{\mathcal{B}^{n+1}}\colon \check{\mathrm{C}}^{-n}(\mathcal{U} \twoheadrightarrow X) \rightsquigarrow \check{\mathrm{C}}^{-n}(\mathcal{V} \twoheadrightarrow Y)$$

which by Corollary 5.2.12 combine to a finite correspondence

$$\check{\mathrm{C}}^{\bullet}(\mathcal{B})\colon \check{\mathrm{C}}^{\bullet}(\mathcal{U} \twoheadrightarrow X) \rightsquigarrow \check{\mathrm{C}}^{\bullet}(\mathcal{V} \twoheadrightarrow Y)$$

of Čech complexes. By Proposition 5.2.10 it extends to a finite correspondence

$$\underline{\check{\mathrm{C}}}^{\bullet}(\mathcal{B})\colon \underline{\check{\mathrm{C}}}^{\bullet}(\mathcal{U} \twoheadrightarrow X) \rightsquigarrow \underline{\check{\mathrm{C}}}^{\bullet}(\mathcal{V} \twoheadrightarrow Y)$$

between the augmented Čech complexes.

**Proposition 5.4.5** (Čech complexes are functorial)**.** *Let $\alpha\colon X \rightsquigarrow Y$ and $\beta\colon Y \rightsquigarrow Z$ be a finite correspondence over a noetherian scheme $S$. Let*

$$\begin{aligned} f\colon &\quad (\mathcal{U}|I) \twoheadrightarrow X, \\ g\colon &\quad (\mathcal{V}|J) \twoheadrightarrow Y, \\ h\colon &(\mathcal{W}|K) \twoheadrightarrow Z \end{aligned}$$

*be pseudocovers. Also let $\mathcal{A}$ and $\mathcal{B}$ be bridges over $\alpha$ and $\beta$, respectively, between the appropriate pseudocovers.*

*Then*

$$\underline{\check{\mathrm{C}}}^{\bullet}(\mathcal{B}) \circ \underline{\check{\mathrm{C}}}^{\bullet}(\mathcal{A}) = \underline{\check{\mathrm{C}}}^{\bullet}(\mathcal{B} \circ \mathcal{A})$$

*as finite correspondences $\underline{\check{\mathrm{C}}}^{\bullet}(\mathcal{U} \twoheadrightarrow X) \rightsquigarrow \underline{\check{\mathrm{C}}}^{\bullet}(\mathcal{W} \twoheadrightarrow Z)$.*

*Proof.* This is immediate from Proposition 5.2.16. □

**Remark 5.4.6.** The functoriality of Čech complexes is clearly defined in the general setting of a category $\mathcal{C}$ which admits fibre products. In fact, this is even easier as it amounts to simple facts on simplicial objects. We used the somewhat convoluted ad hoc definitions and proofs above because we actually do not have such a structure in our case of interest: $\mathrm{SmCor}(S, \Lambda)$ does in general not admit fibre products, hence we need to restrict to genuine morphisms of schemes when defining the Čech complexes.

**Remark 5.4.7.** The reduced Čech complex is functorial as well, assuming that we are given a compatible morphism between the pseudocovers. This condition, while irrelevant for full Čech complexes, is crucial. As explained in Remark 5.2.3, we do not assume it for bridges, hence we cannot properly talk about their functoriality regarding reduced Čech complexes.

## 5.5 Čech diagrams of étale morphisms

**Lemma 5.5.1.** *Let $f : \mathcal{U} \to \mathcal{V}$ be an étale and separated morphism of irreducible schemes. Assume that the induced morphism between their generic points is an isomorphism.*

*Then $f$ is an open immersion.*

*Proof.* This is a special case of [Stacks, Tag 09NQ]. □

**Corollary 5.5.2.** *Let $\mathcal{U} \to X$ be an étale and separated morphism of finite type between schemes and assume that $\mathcal{U}$ is irreducible. Let $f : \mathcal{U} \to \mathcal{U}$ be a morphism over $X$.*

*Then $f$ is an automorphism of $\mathcal{U}$.*

*Proof.* We find that $f$ is étale. By Lemma 5.5.1 it is thus an open immersion. Note that $\mathcal{U} \to X$ is by assumption quasi-finite, therefore a simple dimension comparison shows that $f$ is fibre-wise surjective, hence an isomorphism. □

**Lemma 5.5.3.** *Let $f : \mathcal{U} \to X$ be an étale and separated morphism of noetherian schemes. Let $\widetilde{X} \to X$ be the normalization of $X$.*

*If the pullback $f \times_X \widetilde{X} : \mathcal{U} \times_X \widetilde{X} \to \widetilde{X}$ is an isomorphism, then $f$ is an isomorphism as well.*

*Proof.* Looking at generic points, Lemma 5.5.1 shows that $f$ is an open immersion. If $x \in X$ is any point, then there is a $y \in \widetilde{X}$ above it, and hence by assumption a $v \in \mathcal{U} \times_X \widetilde{X}$ above $y$. Its image in $\mathcal{U}$ hence lies over $x$, thus $f$ is surjective and therefore an isomorphism. □

**Lemma 5.5.4.** *Let $\mathcal{U}$ and $\mathcal{V}$ be schemes which are étale and of finite type over a noetherian scheme $X$.*

*Then there are only finitely many morphism $\mathcal{U} \to \mathcal{V}$ of schemes over $X$.*

*Proof.* It is sufficient to check this for each of the finitely many connected components of $\mathcal{U}$, i.e. we may assume that $\mathcal{U}$ is connected. We observe that $\mathcal{V}$ is quasi-finite over $X$. Fixing a point $u \in \mathcal{U}$, there are thus only finitely many morphisms $u \to \mathcal{V}$ over $X$. Each of them comes from at most one morphism $\mathcal{U} \to \mathcal{V}$ over $X$ by Corollary 3.13 in Chapter I of [Mil80], proving the result. □

We need the following variant of König's Lemma:

**Lemma 5.5.5.** *Let $G$ be a directed graph with infinitely many vertices. Assume that every vertex has finite indegree, i.e. only finitely many edges end in it. Also assume that there is a vertex $v_0$ which can be reached from infinitely many vertices of $G$ via the directed edges.*

*Then there exists an infinite path $\ldots \to v_2 \to v_1 \to v_0$ in $G$ consisting of pairwise distinct vertices.*

*Proof.* We construct this chain inductively:

Assume that we have already found a chain $v_n \to \ldots \to v_1 \to v_0$ such that infinitely many vertices of $G$ can reach $v_n$. For $n = 0$ this is clearly satisfied. Consider the edges $v \to v_n$. By assumption, at least one of these finitely many $v$ is reached from infinitely many vertices of $G$, hence we can pick $v_{n+1} = v$. □

**Definition 5.5.6.** Let $f\colon \mathcal{U} \to X$ be an étale morphism of schemes. Its *Čech diagram* $\check{\mathrm{D}}(f)$ is defined as follows:

- The objects are the different isomorphism classes of schemes over $X$, formally given by fixing a representative, that occur as connected components of a $(\mathcal{U}|X)^n$, $n \in \mathbb{N}_0$.
- The morphisms between those objects are all morphisms of schemes over $X$ between them.

Hence, informally, the Čech diagram is obtained by identifying all duplicates among the connected components occurring in the Čech complex.

**Proposition 5.5.7.** *Let $f\colon \mathcal{U} \to X$ be an étale and separated morphism of finite type between noetherian schemes. Assume that $X$ is normal and of finite dimension.*

*Then the Čech diagram $\check{\mathrm{D}}(f)$ is finite and almost acyclic (cf. Definition 1.2.9). In particular, there are only finitely many isomorphism classes of schemes over $X$ that occur as connected components in the Čech complex $\check{\mathrm{C}}(f)$.*

*Proof.* By Lemma 5.5.4 there are only finitely many morphisms between any two objects of $\check{\mathrm{D}}(f)$. Combined with Corollary 5.5.2 this shows the almost acyclicity. Due to the aforementioned finiteness of all hom-sets it is now enough to show that $\check{\mathrm{D}}(f)$ has only finitely many objects. We may assume that $X$ is connected, hence irreducible by normality.

If $f\colon \mathcal{U} \to X$ is a finite étale Galois covering with automorphism group $G$, then $\mathcal{U} \times_X \mathcal{U} \cong \mathcal{U} \times G$, hence the statement holds. If more generally $f$ is finite étale, we choose a Galois cover $\mathcal{V} \to X$ factoring over $f$. Let again $G$ denote its Galois group. Looking at the factorization $\mathcal{V} \times G^{n-1} \cong (\mathcal{V}|X)^n \to (\mathcal{U}|X)^n \to X$ we see that every vertex of $\check{\mathrm{D}}(f)$ is a subcover of $\mathcal{V}$. But the main theorem of Galois covers, e.g found as Theorem 5.3 in Chapter I of [Mil80], states that the set of such isomorphism classes is finite.

In general we proceed by induction on the dimension of $X$:

Instead of $\check{\mathrm{D}}(f)$ we may consider the quiver $Q$ with the same objects as $\check{\mathrm{D}}(f)$, but the arrows being those $W \to V$, $V \in \check{\mathrm{D}}(f)$, that come from restricting the projection $\mathcal{U} \times_X V \to V$ to one of the finitely many connected components $W$ of $\mathcal{U} \times_X V$. Due to it having the same objects as $\check{\mathrm{D}}(f)$, it is enough to show that $Q$ has finitely many vertices.

Assume that $Q$ is infinite and observe that every $U \in \check{\mathrm{D}}(f)$ can reach $X$ via arrows of $Q$. Hence we can apply Lemma 5.5.5 and find an infinite chain

$$\ldots \to U_2 \to U_1 \to X$$

in $Q$. Note that this is a chain of étale morphisms between schemes which are pair-wise non-isomorphic over $X$. They are furthermore connected, hence irreducible by normality. The morphism $f$ is generically finite, say over the non-empty open $W \hookrightarrow X$. By the finite étale case, we find infinitely many members of our chain which are isomorphic over $W$. Deleting all other members and using Lemma 5.5.1 we therefore found an infinite chain

$$\ldots \to V_2 \to V_1 \to V_0$$

of open immersions, none of them isomorphisms.

Let $Y \to X \backslash W$ be the normalization. Then, by induction on the dimension of $X$, we already know that $\check{\mathrm{D}}(Y \times_X f)$ is finite. Hence the pullback of our chain along $Y \to X$ contains duplicates, hence so does the pullback to $X \backslash W$ by Lemma 5.5.3. This means that the corresponding open immersion $\mathcal{V}_j \to \mathcal{V}_i$ is bijective, hence an isomorphism, contradicting our assumption that $\check{\mathrm{D}}(f)$ is infinite. □

# Chapter 6

# Yoga of Nisnevich Covers

We continue the path laid out at the beginning of the last Chapter 5: we wish to rigidify Nisnevich covers. The general idea is to enforce an additional structure on the Nisnevich covers that makes morphisms, bridges and finite correspondences between them unique. We then use rigidifications and our results to properly deal with Nisnevich covers on diagrams. Indeed, one has many choices when lifting a morphism to given covers, while one surely wants to get the same morphisms when lifting along two compositions $a \circ b = c \circ d$.

We point to [HM16], Section 9.2, for the case of Zariski covers. An étale variant was described in [Fri82], using geometric points. Both only aimed to rigidify morphisms of schemes, not finite correspondences, hence were able to avoid many of the subtleties and technicalities encountered in this chapter.

## 6.1 Rigidification of Nisnevich covers

**Definition 6.1.1** (Rigidification)**.** Let $f\colon (\mathcal{U}|I) \twoheadrightarrow X$ be a Nisnevich cover.

- A *(Nisnevich) pre-rigidification* of $f$ is a point-wise splitting $r\colon X \to \mathcal{U}$, i.e. for every $x \in X$ a choice of an $r(x) = u \in \mathcal{U}$ mapping to $f(u) = x$ where $f$ induces an isomorphism $\kappa(u) \cong \kappa(x)$. We will write $f\colon (\mathcal{U}|I) \twoheadrightarrow_r X$ to add the chosen rigidification into our notation for pseudocovers.

- A pre-rigidification $r$ of $f$ induces a set-theoretic map $I_r\colon X \to I$ sending $x$ to the $i \in I$ with $r(x) \in \mathcal{U}_i$, which we call its *indexing*.

- We call a pre-rigidification *constructible* if for every $x \in X$ there exists a locally closed subset $E \subseteq \mathcal{U}$ containing $r(x)$ such that $r(f(e)) = e$ for all $e \in E$.

- A *rigidification* of $f$ is a pre-rigidification $r$ meeting every connected component of $\mathcal{U}$. This can be stronger than meeting all the given components $\mathcal{U}_i$. We call it *strong* if it even meets every irreducible component of $\mathcal{U}$.

- A *(pre)-rigidified* Nisnevich cover is a Nisnevich cover with a chosen (pre)-rigidification, respectively. If it is furthermore constructible, then we call it *constructably (pre-)rigidified* or *c-(pre-)rigidified* for short. *Strongly (c-)rigidified* Nisnevich covers are defined accordingly.

- A *unifibrantly* pre-rigidified Nisnevich cover is a unifibrant Nisnevich cover with a pre-rigidification such that for each $x \in X$ the corresponding $r(x)$ is the single point in $\mathcal{U}_{I_r(x)}$ lying over $x$.

- To every pre-rigidified Nisnevich cover $f\colon (\mathcal{U}|I) \twoheadrightarrow_r X$ we associate a rigidification $f^{\mathrm{rig}}$ by removing those connected components from $\mathcal{U}$ which are not met by $r$.

**Remark 6.1.2.** Note that if $X$ is normal, then every rigidification of a Nisnevich cover is already strong. Indeed, a scheme étale over a normal one is itself normal, hence its irreducible components are the connected components.

An alternative approach to c-pre-rigidified cover comes from [Voe10b], Definition 2.14:

**Definition 6.1.3.** Let $f\colon (\mathcal{U}|I) \twoheadrightarrow X$ be a Nisnevich cover.

A *rigid splitting sequence* for $f$ consists of a sequence

$$\emptyset = Z_{m+1} \hookrightarrow Z_m \hookrightarrow \ldots \hookrightarrow Z_1 \hookrightarrow Z_0 = X$$

of closed embeddings together with morphisms $s_k\colon Z_k \backslash Z_{k+1} \to \mathcal{U}$ which are sections of $f$, i.e. such that $f \circ s_k = \mathrm{id}_{Z_k \backslash Z_{i+k}}$.

The *associated pre-rigidification* of this splitting sequence is the pre-rigidification of $f$ defined by $r(x) = s_i(x)$ where $x \in Z_i \backslash Z_{i+1}$.

**Remark 6.1.4.** As elaborated in [Ivo07], every scheme $X$ has a *henselization*

$$X^{\mathfrak{h}} := \coprod_{x \in X} X^{\mathfrak{h}}_x$$

defined as the disjoint union of the point-wise henselizations. It comes with a canonical morphism $\mathfrak{i}^{\mathfrak{h}}_X\colon X^{\mathfrak{h}} \to X$. Then pre-rigidifications of a Nisnevich cover $f\colon (\mathcal{U}|I) \twoheadrightarrow X$ can also be defined as factorizations of $\mathfrak{i}^{\mathfrak{h}}_X$ over $f$. Standard results for henselizations show that both definitions are equivalent.

In this light, several of our results should be understood as finite versions of the Nisnevich localization techniques of [Ivo07].

We will use the following well-known result (cf. [EGA4-1], Théorème 1.8.4):

**Lemma 6.1.5** (Chevalley's Theorem)**.** *Let $f\colon X \to Y$ be a morphism of finite presentation and let $\eta \in X$.*

*Then the image $f(X)$ contains a locally closed subset of $Y$ containing $f(\eta)$.*

**Lemma 6.1.6.** *Let $f\colon (\mathcal{U}|I) \twoheadrightarrow X$ be a Nisnevich cover of a noetherian scheme $X$.*

*Then the pre-rigidification associated to a rigid splitting sequence for $f$ is constructible. Conversely, every constructible pre-rigidification of $f$ is the pre-rigidification associated to a rigid splitting sequence.*

*Proof.* The first part is simply Lemma 6.1.5. Now let $f\colon (\mathcal{U}|I) \twoheadrightarrow_r X$ be a c-pre-rigidified Nisnevich cover. We construct a sequence $\emptyset = A_0 \subseteq A_1 \subseteq A_2 \subseteq \dots$ of open subsets as follows:

Assume $A_i$ is already constructed. Unless we already have $A_i = X$ we pick an arbitrary generic point $\eta \in X \backslash A_i$ of maximal dimension. Giving them the induced reduced structures, the isomorphism $\kappa(r(\eta)) \cong \kappa(\eta)$ amounts to a birational morphism $\overline{r(\eta)} \to \overline{\eta}$. Thus there is a non-empty open $\widetilde{E}_i \subseteq \overline{r(\eta)}$ where $f$ induces an isomorphism onto an open subscheme of $\overline{\eta}$.

The constructibility of $r$ allows us to shrink the locally closed $\widetilde{E}_i \subseteq \mathcal{U}$ to get $r(f(u)) = u$ for all $u \in \widetilde{E}_i$. By Lemma 6.1.5 the image $f(\widetilde{E}_i)$ in $\overline{\eta}$ contains a locally closed subset containing $\eta$, hence a non-empty open subset $E_i$ of $\overline{\eta}$. By construction we have an isomorphism $s_i\colon E_i \to \widetilde{E}_i \hookrightarrow \mathcal{U}$ such that $r(x) = s_i(x)$ for all $x \in E_i$.

We let $A_{i+1} = A_i \cup E_i$. Note that it is open, hence its complement is closed. Therefore this procedure terminates because $X$ is noetherian, say $A_m = X$. We got a decomposition $X = \bigcup_{i=0}^m E_i$ into pairwise disjoint subsets and set $Z_k := \bigcup_{i=k}^m E_i$, which is closed in $X$. Giving $Z_i$ the induced reduced structure gives us by consruction a section $s_i\colon E_i = Z_i \backslash Z_{i+1} \to \mathcal{U}$ as required. □

**Corollary 6.1.7.** *Every Nisnevich cover $f\colon (\mathcal{U}|I) \twoheadrightarrow X$ of a noetherian scheme has a constructible pre-rigidification.*

*Proof.* It is enough to find a rigid splitting sequence for $f$. This follows by the same construction as in the proof of Lemma 6.1.6 by ignoring the condition that $r(f(u)) = u$ for all $u \in \widetilde{E}_i$. □

## 6.2 Rigid morphisms

**Definition 6.2.1** (Rigid morphisms)**.** Let

$$f\colon (\mathcal{U}|I) \twoheadrightarrow_r X,$$
$$g\colon (\mathcal{V}|J) \twoheadrightarrow_s Y$$

be pre-rigidified Nisnevich covers. Let $a\colon X \to Y$ be a morphism of schemes.

A morphism $b\colon \mathcal{U} \to \mathcal{V}$ is *rigid* if it is compatible with the rigidifications, i.e. if $s \circ a = b \circ r$ as maps $X \to \mathcal{V}$. If $Y = X$ and $a = \mathrm{id}_X$, then we call it a *rigid refinement*.

**Remark 6.2.2.** Note that the rigid morphism $b$ restricts to a morphism between the rigidifications $f^{\mathrm{rig}}$ and $g^{\mathrm{rig}}$. Indeed, every connected component $U$ of $\mathcal{U}$ is mapped to a single connected component $V$ of $\mathcal{V}$, and if $r$ meets $U$, then by definition $s$ meets $V$.

**Definition 6.2.3** (Rigid pullbacks)**.** Let $f\colon (\mathcal{U}|I) \twoheadrightarrow_r X$ be a pre-rigidified Nisnevich cover.

- If $Z \to X$ is a morphism of schemes, then the *pullback* $r \times_X Z$ of $r$ to the Nisnevich cover $(\mathcal{U} \times_X Z|I) \twoheadrightarrow Z$ is the map sending a point $z \in Z$ over $x \in X$ to $r(x) \times_x z \in \mathcal{U} \times_X Z$. This clearly preserves constructibility. Removing all connected components not met by $r \times_X Z$ gives us the *rigidified pullback* $f \times_X^{\mathrm{rig}} Z$ which is a rigidified Nisnevich cover of $Z$.

- Let $f'\colon (\mathcal{U}'|I') \twoheadrightarrow_{r'} X$ be a second pre-rigidified Nisnevich cover of $X$. Then *fibre product* of the pre-rigidifications $r$ and $r'$ is the pre-rigidification $r \times_X r'$ on the fibre product $f \times_X f'\colon (\mathcal{U} \times_X \mathcal{U}'|I \times I') \twoheadrightarrow X$ defined by
$$x \mapsto r(x) \times_x r'(x) = r(x) \times_X r'(x).$$
We define their *rigidified fibre product* $r \times_X^{\mathrm{rig}} r'$ as the rigidification $(r \times_X r')^{\mathrm{rig}}$, i.e. by removing all unmet connected components.

Fibre products preserve the respective structures, thus we have shown:

**Lemma 6.2.4.** *Let $X$ be a noetherian scheme.*

*Then any syntactically correct combination of affine, rigidified, c-, pre-, cd-, unifibrant(ly), Zariski and Nisnevich causes its corresponding set of such covers to be directed.*

**Lemma 6.2.5.** *Let $X$ be a noetherian scheme.*

*Every c-pre-rigidified Nisnevich cover of $X$ has a refinement into a unifibrantly c-pre-rigidified cd-Nisnevich cover.*

*Proof.* We repeat the argument of [Voe10b], Proposition 2.16:

Let $f\colon (\mathcal{U}|I) \twoheadrightarrow_r X$ be a c-rigidified Nisnevich cover. By Lemma 6.1.6 its rigidification comes from a rigid splitting sequence
$$\emptyset = Z_{m+1} \hookrightarrow Z_m \hookrightarrow Z_{m-1} \hookrightarrow \ldots \hookrightarrow Z_1 \hookrightarrow Z_0 = X$$
with sections $s_i\colon Z_i \backslash Z_{i+1} \to \mathcal{U}$. We may assume that its length $m$ is minimal and that the statement is already proven for all c-pre-rigidified Nisnevich covers that admit a rigid splitting sequence of smaller length.

Because the pullback $Z_m \times_X \mathcal{U} \to Z_m$ is étale, the image of the section $s_m \colon Z_m \to \mathcal{U}$ is open in $Z_m \times_X \mathcal{U}$ and thus its complement $Y$ is closed. Now $Z_m \times_X \mathcal{U} \hookrightarrow \mathcal{U}$ is closed, thus $Y$ is also closed in $\mathcal{U}$.

Set $A = X \backslash Z_m$. Then the pair $\{A, \mathcal{U}\backslash Y\}$ forms a unifibrantly c-pre-rigidified cd-Nisnevich cover of $X$. Furthermore, the pullback $\mathcal{U} \times_X A \twoheadrightarrow_{r \times_X A} A$ is constructible and induced by the splitting sequence

$$\emptyset \hookrightarrow Z_{m-1} \times_X A \hookrightarrow \ldots \hookrightarrow Z_1 \times_X A \hookrightarrow Z_0 \times_X A = A$$

with sections $s_i \times_X A$. As this sequence has length smaller than $m$ we can refine $\mathcal{U} \times_X A \twoheadrightarrow A$ into a unifibrantly c-pre-rigidified Nisnevich cover. Combining it with $\mathcal{U}\backslash Y$ shows the result. □

**Lemma 6.2.6.** *Let $f \colon (\mathcal{U}|I) \twoheadrightarrow_r X$ be a rigidified Nisnevich cover and let $g \colon (\mathcal{V}|J) \twoheadrightarrow_s Y$ be a pre-rigidified Nisnevich cover. Also let a morphism $a \colon X \to Y$ be given.*

*Then there exists at most one rigid morphism $b \colon \mathcal{U} \to \mathcal{V}$ between the Nisnevich covers.*

*Proof.* We can assume that at least one such $b$ exists. Replacing $g \colon (\mathcal{V}|J) \twoheadrightarrow_{r'} Y$ by its rigid pullback along $a$ we may assume that $X = Y$ and $a = \mathrm{id}_X$. Let $U$ be any connected component of $\mathcal{U}$ and choose an $x \in X$ with $r(x) \in U$. Because $g$ is étale, we find that $b$ is the unique morphism $U \to \mathcal{V}$ on the connected $U$ inducing the isomorphism $\kappa(r(x)) \cong \kappa(s(x))$ (cf. [Mil80], chapter I, Corollary 3.13). Varying $U$ shows that $b$ is fixed on all of $\mathcal{U}$. □

We will later need the following to properly deal with bridges:

**Lemma 6.2.7.** *Let $\pi \colon \Gamma \to X$ be a universally closed and pseudo-dominant (cf. Definition 1.4.1) morphism of noetherian schemes. Let $f \colon (\mathcal{U}|I) \twoheadrightarrow_r X$ be a strongly rigidified Nisnevich cover.*

*Then the pullback $r \times_X \pi$ of $r$ along $\pi$ meets every irreducible component of $\mathcal{U} \times_X \Gamma$. In other words: the pullback of a strongly rigidified Nisnevich cover along a universally closed pseudo-dominant morphism is again strongly rigidified.*

*Proof.* Let $V$ be any irreducible component of $\mathcal{U} \times_X \Gamma$. The projection $\mathrm{pr}_{\mathcal{U}} \colon \mathcal{U} \times_X \Gamma \to \mathcal{U}$ is closed and by Lemma 1.4.2 also pseudo-dominant. Hence the image $\mathrm{pr}_{\mathcal{U}}(V)$ is closed and contains a generic point of $\mathcal{U}$. Thus $\mathrm{pr}_{\mathcal{U}}(V)$ contains an irreducible component $W$ of $\mathcal{U}$.

Then by assumption there is an $x \in X$ with $r(x) \in W \subseteq \mathrm{pr}_{\mathcal{U}}(V)$. Choose a point $v \in V$ mapping to $r(x)$ and let $\gamma$ be its image in $\Gamma$. Then $(r \times_X \Gamma)(\gamma) = r(x) \times_x \gamma = v \in V$ as required. □

## 6.3 Nisnevich-locality

Recall the following:

**Definition 6.3.1** (Pointed morphism)**.** Let $X$, $Y$ be schemes and fix $n$ points $x_1, \dots, x_n \in X$ as well as $n$ points $y_1, \dots, y_n \in Y$. An *$n$-pointed morphism* from $X$ to $Y$ with respect to the chosen points is a morphism $f\colon X \to Y$ of schemes sending $x_i$ to $y_i$. We denote this simply by $f\colon (X, (x_i)) \to (Y, (y_i))$.

A *pointed morphism* is a 1-pointed morphism.

**Definition 6.3.2** (Nisnevich neighbourhood)**.** Let $X$ be a scheme. A *Nisnevich neighbourhood* of $x \in X$ is a pointed étale morphism $\varphi\colon (U, u) \to (X, x)$ such that $\varphi$ induces an isomorphism $\kappa(u) \cong \kappa(x)$ of residue fields.

More generally, a *joint Nisnevich neighbourhood* of $n$ points $x_1, \dots, x_n \in X$ is an $n$-pointed étale morphism $\varphi\colon (U, (u_i)) \to (X, (x_i))$ inducing isomorphisms of residue fields at the $u_i$.

The following is folklore, but we give a proof for the reader's convenience:

**Lemma 6.3.3.** *Let $\pi\colon \Gamma \to X$ be a finite morphism of schemes, let $x \in X$ be any point and let $\gamma_1, \dots, \gamma_n \in \Gamma$ be the distinct points of $\Gamma$ above it.*

*Then this induces a natural isomorphism*

$$X_x^{\mathfrak{h}} \times_X \Gamma \cong \coprod_{i=1}^{n} \Gamma_{\gamma_i}^{\mathfrak{h}}$$

*between henselizations of local schemes.*

*Proof.* Let $y$ be the closed point of $X_x^{\mathfrak{h}}$, which is also the unique point over $x$. There we have an isomorphism $\kappa(y) \cong \kappa(x)$ of residue fields, hence the closed points of $X_x^{\mathfrak{h}} \times_X \Gamma$ are in 1-to-1-correspondence with the $\gamma_i$. Let $\gamma_i' \in X_x^{\mathfrak{h}} \times_X \Gamma$ be the unique point above $\gamma_i$.

From Proposition 18.5.10 of [EGA4-4] we get that $X_x^{\mathfrak{h}} \times_X \Gamma$ is a finite disjoint union of henselian local schemes. Furthermore the $\gamma_i'$ constitute exactly the distinct closed points of $X_x^{\mathfrak{h}} \times_X \Gamma$ by the above, hence we have an isomorphism

$$X_x^{\mathfrak{h}} \times_X \Gamma \cong \coprod_i (X_x^{\mathfrak{h}} \times_X \Gamma)_{\gamma_i'}.$$

By [Stacks, Tag 05WP] we also have an isomorphism

$$(X_x^{\mathfrak{h}} \times_X \Gamma_{\gamma_i})_{\gamma_i'} \cong \Gamma_{\gamma_i}^{\mathfrak{h}}.$$

Putting things together we got a chain

$$X_x^{\mathfrak{h}} \times_X \Gamma \cong \coprod_{i=1}^{n} (X_x^{\mathfrak{h}} \times_X \Gamma)_{\gamma_i'} \cong \coprod_{i=1}^{n} (X_x^{\mathfrak{h}} \times_X \Gamma_{\gamma_i})_{\gamma_i'} \cong \coprod_{i=1}^{n} \Gamma_{\gamma_i}^{\mathfrak{h}}$$

of isomorphisms. □

We now prove a rigid version of *Nisnevich-locality* for finite morphisms.

**Proposition 6.3.4.** *Let $\pi\colon \Gamma \to X$ be a finite morphism of noetherian schemes. Let $x \in X$ be an arbitrary point and let $\gamma_1, \dots, \gamma_n \in \Gamma$ be the distinct points of $\Gamma$ above it. Let $g\colon (\mathcal{V}, (v_i)) \to (\Gamma, (\gamma_i))$ be a joint Nisnevich neighbourhood.*

*Then there exists a pointed Nisnevich neighbourhood $f\colon (\mathcal{U}, u) \to (X, x)$ and a morphism $b\colon U \times_X \Gamma \to \mathcal{V}$ of schemes over $\Gamma$ sending the point $u \times_x \gamma_i$ to $v_i$ for all $i \in \{1, 2, \dots, n\}$.*

*Proof.* For every $i \in \{1, 2, \dots, n\}$ we let $W_i = \{(\mathcal{W}_i, w_i) \to (\Gamma, \gamma_i)\}$ be the directed system of pointed Nisnevich neighbourhoods of $\gamma_i \in \Gamma$. We also let $W = \{(\mathcal{W}, (w_i)) \to (\Gamma, (\gamma_i))\}$ be the directed system of joint Nisnevich neighbourhoods of $\gamma_1, \dots, \gamma_n \in \Gamma$.

Then we have maps

$$\begin{aligned}
F\colon \qquad \prod_{i=1}^{n} W_i &\to W, \\
((\mathcal{W}_i, w_i)) &\mapsto (\coprod_{i=1}^{n} \mathcal{W}_i, (w_i)), \\
G\colon \qquad W &\to \prod_{i=1}^{n} W_i, \\
(\mathcal{W}, (w_i)) &\mapsto ((\mathcal{W}^i, w_i))
\end{aligned}$$

of directed systems, where we set $\mathcal{W}^i$ to be the connected component of $\mathcal{W}$ containing $w_i$. It is immediately checked that both $G \circ F$ and $F \circ G$ map a given Nisnevich neighbourhood to a finite union of some of its closed components, possibly containing duplicates. Thus both compositions are refinements, which gives us an isomorphism

$$\coprod_{i=1}^{n} \varprojlim_{(\mathcal{W}_i, w_i) \to (\Gamma, \gamma_i)} (\mathcal{W}_i, w_i) \cong \varprojlim_{(\mathcal{W}, (w_i)) \to (\Gamma, (\gamma_i))} (\mathcal{W}_i, (w_i))$$

of projective limits of ($n$-)pointed morphisms.

With this, by Lemma 6.3.3 and from [Stacks, Tag 04HY] we get $n$-pointed

isomorphisms

$$\begin{aligned}
&\varprojlim_{(\mathcal{U},u)\to(X,x)} (\mathcal{U}\times_X \Gamma, (u\times_x \gamma_i)) \cong \\
\cong &\varprojlim_{(\mathcal{U},u)\to(X,x)} ((\mathcal{U},u)\times_{(X,x)} (\Gamma,(\gamma_i)) \cong \\
\cong &(\varprojlim_{(\mathcal{U},u)\to(X,x)} (\mathcal{U},u))\times_{(X,x)} (\Gamma,(\gamma_i)) \cong \\
\cong &X_x^{\mathfrak{h}} \times_{(X,x)} (\Gamma,(\gamma_i)) \cong \coprod_{i=1}^{n} \Gamma_{\gamma_i}^{\mathfrak{h}} \cong \\
\cong &\coprod_{i=1}^{n} \varprojlim_{(\mathcal{W}_i,w_i)\to(\Gamma,\gamma_i)} (\mathcal{W}_i, w_i) \cong \\
\cong &\varprojlim_{(\mathcal{W},(w_i))\to(\Gamma,(\gamma_i))} (\mathcal{W},(w_i))
\end{aligned}$$

of projective limits running over the respective directed systems of pointed Nisnevich neighbourhoods.

Now $(\mathcal{V},(v_i))$ is one of the Nisnevich neighbourhoods of $(\Gamma,(\gamma_i))$ in the last limit. The scheme $\mathcal{V}$ is finitely presented over the noetherian $X$, thus the isomorphism of projective limits gives us an $n$-pointed morphism $\mathcal{U}\times_X \Gamma \to \mathcal{V}$ for some sufficiently fine Nisnevich neighbourhood $(\mathcal{U},u)\to(X,x)$ as desired. $\square$

## 6.4 Rigid bridges

**Definition 6.4.1** (Rigid bridges)**.** Let $\alpha\colon X \rightsquigarrow Y$ be a finite correspondence over a noetherian scheme $S$. Let $f\colon (\mathcal{U}|I) \twoheadrightarrow_r X$ and $g\colon (\mathcal{V}|J) \twoheadrightarrow_s Y$ be pre-rigidified Nisnevich covers.

A *rigid bridge* over $\alpha$ between the given covers is a bridge

$$\mathcal{B} = (f, g, \gamma\colon \Gamma \to X\times_S Y, \widetilde{\alpha}, b\colon \mathcal{U}\times_X \Gamma \to \Gamma\times_Y \mathcal{V})$$

where $b\colon \mathcal{U}\times_X \Gamma \to \Gamma\times_Y \mathcal{V}$ is a rigid refinement between the pre-rigidified Nisnevich covers on $\Gamma$ induced by pullback.

We call a finite correspondence $\widetilde{\alpha}\colon \mathcal{U}\to\mathcal{V}$ between the covers *rigidly bridgeable over* $\alpha$ if it is the associated finite correspondence (cf. Definition 5.2.4) of a rigid bridge.

Hence this means that $b(r(x)\times_X \gamma) = \gamma\times_Y s(y)$ for every point $\gamma\in\Gamma$ over points $x\in X$ and $y\in Y$. This notion clearly generalizes that of rigid morphisms in the light of Definition 1.6.9 and Proposition 1.6.10.

**Remark 6.4.2.** Similar to Remark 6.2.2, we note that the rigid bridge $b$ restricts to one between the rigidifications $f^{\mathrm{rig}}$ and $g^{\mathrm{rig}}$.

**Proposition 6.4.3.** *Being a rigid bridge is preserved by composition, fibre product and tensor product of bridges whenever they are defined.*

*Proof.* Let us exemplary check this for composition, the other two being equally simple:

Rigidity is clearly preserved by composition of morphisms and by pullbacks. Hence the claimed property follows immediately from the Definition 5.2.14 of composition of bridges. □

**Theorem 6.4.4.** *Let $S$ be a noetherian scheme. Let $\alpha\colon X \rightsquigarrow Y$ and $\beta\colon Y \rightsquigarrow Z$ be finite correspondences over $S$ and let*

$$\mathcal{U} \twoheadrightarrow_r X, \quad \mathcal{V} \twoheadrightarrow_s Y, \quad \mathcal{W} \twoheadrightarrow_t Z$$

*be pre-rigidified Nisnevich covers. Furthermore let $\widetilde{\alpha}\colon \mathcal{U} \rightsquigarrow \mathcal{V}$ and $\widetilde{\beta}\colon \mathcal{V} \rightsquigarrow \mathcal{W}$ be rigidly bridgeable finite correspondences between those Nisnevich covers over $\alpha$ and $\beta$, respectively.*

*Then $\widetilde{\beta} \circ \widetilde{\alpha}\colon \mathcal{U} \rightsquigarrow \mathcal{W}$ is rigidly bridgeable over $\beta \circ \alpha$.*

*Proof.* Choose any bridges $\mathcal{A}$ and $\mathcal{B}$ inducing $\widetilde{\alpha}$ and $\widetilde{\beta}$ over $\alpha$ and $\beta$, respectively. Then the bridge $\mathcal{B} \circ \mathcal{A}$ is rigid by Proposition 6.4.3. Its associated finite correspondence is $\widetilde{\beta} \circ \widetilde{\alpha}$ by Proposition 5.2.15 (b), showing the result. □

**Theorem 6.4.5.** *Let $\alpha_1\colon X_1 \rightsquigarrow Y_1$ and $\alpha_2\colon X_2 \rightsquigarrow Y_2$ be finite correspondences over a noetherian scheme $S$ and let*

$$\begin{aligned} &\mathcal{U}_1 \twoheadrightarrow_{r_1} X_1, && \mathcal{V}_1 \twoheadrightarrow_{s_1} Y_1, \\ &\mathcal{U}_2 \twoheadrightarrow_{r_2} X_2, && \mathcal{V}_2 \twoheadrightarrow_{s_2} Y_2 \end{aligned}$$

*be pre-rigidified Nisnevich covers. Furthermore, let $\widetilde{\alpha}_1\colon \mathcal{U}_1 \rightsquigarrow \mathcal{V}_1$ and $\widetilde{\alpha}_2\colon \mathcal{U}_2 \rightsquigarrow \mathcal{V}_2$ be rigidly bridgeable finite correspondences between those Nisnevich covers over $\alpha_1$ and $\alpha_2$, respectively.*

*Then*

$$\widetilde{\alpha}_1 \otimes \widetilde{\alpha}_2\colon \mathcal{U}_1 \times_S \mathcal{U}_2 \rightsquigarrow \mathcal{V}_1 \times_S \mathcal{V}_2$$

*is rigidly bridgeable over $\alpha_1 \otimes \alpha_2$.*

*Proof.* Choose any bridges $\mathcal{B}_1$ and $\mathcal{B}_2$ inducing $\widetilde{\alpha}_1$ and $\widetilde{\alpha}_2$ over $\alpha_1$ and $\alpha_2$, respectively. Then the bridge $\mathcal{B}_1 \otimes \mathcal{B}_2$ is rigid by Proposition 6.4.3. Its associated finite correspondence is $\widetilde{\alpha}_1 \otimes \widetilde{\alpha}_2$ by Proposition 5.2.18 (b), showing the result. □

The following is a finite and rigid version of the Nisnevich locality found in chapter 2 of [Ivo07]. A special and unpointed case can also be found as Lemma 3.2.2 in [Kel13].

**Proposition 6.4.6** (Building bridges)**.** *Let a finite correspondence $\alpha\colon X \rightsquigarrow Y$ over a noetherian scheme $S$ and a c-pre-rigidified Nisnevich cover*

$$g\colon (\mathcal{V}|J) \twoheadrightarrow_s Y$$

*be given. Also let*

$$(\gamma\colon \Gamma \to X \times_S Y, \widetilde{\alpha})$$

*be a pylon over $\alpha$.*

*Then there exists a c-rigidified Nisnevich cover*

$$f\colon (\mathcal{U}|I) \twoheadrightarrow_r X$$

*and a rigid bridge*

$$\mathcal{B} = (f, g, \gamma, \widetilde{\alpha}, b)$$

*over $\alpha$.*

*Proof.* We proceed by noetherian induction:

Assume we already found an open subset $A \subsetneq X$ over which we have such a bridge. This means we have an étale cover $f_A\colon \mathcal{U}_A \to A$, a constructible splitting $r_A\colon A \to \mathcal{U}_A$ and a morphism $b_A\colon \mathcal{U}_A \times_X \Gamma \to \Gamma \times_Y \mathcal{V}$ of schemes over $\Gamma$ which is rigid, i.e. $b_A(r(a) \times_X w) = w \times_Y s(y)$ for every $w \in \Gamma$ above $a \in A$ and $y \in Y$. Note that the base case $A = \emptyset$ satisfies this tautologically.

Now pick any generic point $x \in X \backslash A$.

Let $s'\colon \Gamma \to \Gamma \times_Y \mathcal{V}$ by the pullback of the pre-rigidification $s\colon Y \to \mathcal{V}$. Let $w_1, \ldots, w_n \in \Gamma$ be the distinct points above $x$ and consider the $n$-pointed Nisnevich cover

$$(\Gamma \times_Y \mathcal{V}, s'(w_i)) \to (\Gamma, (\gamma_i)).$$

Applying Proposition 6.3.4 we find a pointed Nisnevich neighbourhood $f\colon (\mathcal{U}, u) \to (X, x)$ and a morphism

$$b\colon \mathcal{U} \times_X \Gamma \to \Gamma \times_Y \mathcal{V}$$

such that

$$b(u \times_x w_i) = s'(w_i)$$

for all $i \in \{1, 2, \ldots, n\}$.

We interpret $f$ as a birational morphism $\overline{u} \to \overline{x}$, both with the induced reduced structures. On these grounds we pick some non-empty open $E \subseteq \overline{u}$ on which $f$ induces an isomorphism onto its open image $f(E)$. By removing intersections with other irreducible components of $X \backslash A$ we may assume that $f(E)$ is open in $X \backslash A$.

By constructibility of the pullback $s' = \Gamma \times_X s\colon \Gamma \to \Gamma \times_Y \mathcal{V}$ we find for each of the finitely many $w_i \in \Gamma$ over $x$ a locally closed subset $W_i \subseteq \Gamma \times_Y \mathcal{V}$ where $s' \circ (\Gamma \times_Y g)$ is the identity. From Chevalleys's theorem in the version of Lemma 6.1.5 it follows that the images of the $W_i$ in $X$ contain a locally closed subset $Z$ containing $x$. As before we may assume that $Z$ is non-empty and open in $\overline{x}$. We replace $E$ by $E \cap f^{-1}(Z)$ and may thus assume that $f(E) \subseteq Z$.

Therefore we can extend our partial bridge to one over $A \cup f(E)$ by doing the following replacements:

- $\mathcal{U}_A$ by $\mathcal{U}_A \sqcup \mathcal{U}$,
- $b_A$ by $b_A \sqcup b\colon (\mathcal{U}_A \sqcup \mathcal{U}) \times_X \Gamma \to \mathcal{V} \times_Y \Gamma$,
- $r_A$ by $x \mapsto \begin{cases} (f_{|E})^{-1}(x) & \text{if } x \in f(E) \\ r_A(x) & \text{if } x \in A. \end{cases}$

Note that $A \cup f(E)$ is not a disjoint union, hence the slightly more complicated extension of $r_A$. Constructibility of the new $r_A$ is satisfied in $E$ by construction and on $A$ by assumption. All other desired properties of the larger partial bridge follow directly from the construction.

We continue with this larger constructible subset. The construction terminates after finitely many steps because $X$ is noetherian. □

We have just shown existence of rigid bridges and now direct our attention to uniqueness.

**Lemma 6.4.7.** *Let $\alpha\colon X \rightsquigarrow Y$ be a finite correspondence over a noetherian scheme $S$. Let $f\colon (\mathcal{U}|I) \twoheadrightarrow_r X$ be a strongly rigidified Nisnevich cover and let $g\colon (\mathcal{V}|J) \twoheadrightarrow_s Y$ be a pre-rigidified Nisnevich cover. Also let $(\gamma\colon \Gamma \to X \times_S Y, \widetilde{\alpha})$ be a pylon over $\alpha$ and assume that $\mathrm{pr}_X \circ \gamma$ is pseudo-dominant (cf. Definition 1.4.1).*

*Then there is at most one rigid bridge*

$$\mathcal{B} = (f, g, \gamma, \widetilde{\alpha}, b)$$

*over this pylon and between the given covers.*

*Proof.* By Lemma 6.2.7 the pullback of $r$ to $\mathcal{U} \times_X \Gamma$ is strongly rigidified, in particular rigidified. Thus by Lemma 6.2.6 there is at most one rigid bridging $b\colon \mathcal{U} \times_X \Gamma \to \Gamma \times_Y \mathcal{V}$ between the pullbacks. □

**Lemma 6.4.8.** *Let $X \to S$ be a morphism of finite type between noetherian schemes. Let*

$$\begin{aligned} f'\colon & (\mathcal{U}'|J) \twoheadrightarrow_{r'} X, \\ f\colon & (\mathcal{U}|I) \twoheadrightarrow_r X \end{aligned}$$

*be strongly rigidified Nisnevich covers. Let $u\colon \mathcal{U}' \to \mathcal{U}$ be a rigid refinement of Nisnevich covers and let $n \in \mathbb{N}_0$.*

*Then $(u|X)^n\colon (\mathcal{U}'|X)^n \to (\mathcal{U}|X)^n$ is an epimorphism in $\mathrm{SchCor}_S$.*

*Proof.* By Proposition 1.6.11 (b) it suffices to show that $(u|X)^n$ is dominant and pseudo-dominant. It is étale, hence open and in particular pseudo-dominant, leaving us to show that $(u|X)^n$ is dominant.

We start with $n = 1$. Let $V$ be an irreducible component of $\mathcal{U}$. By definition of strong rigidifications there exists an $x \in X$ such that $r(x) \in V$. As $u: \mathcal{U}' \to \mathcal{U}$ is a rigid refinement we have $u(r'(x)) = r(x)$, therefore $V$ is met by $u$. As $u$ is open and because $V$ is irreducible we conclude that the generic point of $V$ is in the image of $u$, which is therefore dominant as $V$ was arbitrary.

In general we conclude by induction: We already know that $u$ and $(u|X)^{n-1}$ are dominant. Hence by [EGA4-2], Proposition 2.3.7, their pullbacks

$$(\mathcal{U}'|X)^n \to (\mathcal{U}'|X)^{n-1} \times_X \mathcal{U},$$
$$(\mathcal{U}'|X)^{n-1} \times_X \mathcal{U} \to (\mathcal{U}|X)^n$$

along the flat morphisms $(\mathcal{U}'|X)^{n-1} \to X$ and $\mathcal{U} \to X$, respectively, are dominant as well. Hence the composition

$$(\mathcal{U}'|X)^n \to (\mathcal{U}'|X)^{n-1} \times_X \mathcal{U} \to (\mathcal{U}|X)^n$$

is dominant as required. □

The main reason for working with rigidified Nisnevich covers is now the following:

**Theorem 6.4.9** (Uniqueness of rigidly bridgeable correspondences)**.** *Let $\alpha: X \rightsquigarrow Y$ be a finite correspondence over a noetherian scheme $S$. Let*

$$f\colon (\mathcal{U}|I) \twoheadrightarrow_r X,$$
$$g\colon (\mathcal{V}|J) \twoheadrightarrow_s Y$$

*be c-rigidified Nisnevich covers. Assume that $X$ is normal.*

*Then there exists at most one rigidly bridgeable finite correspondence $\beta: \mathcal{U} \rightsquigarrow \mathcal{V}$ over $\alpha$. More generally, this is recovered as the degree-0 part of the following:*

*Let $\mathcal{B}_1$ and $\mathcal{B}_2$ be two bridges over $\alpha$ between the same given c-rigidified covers. Then we have the equality*

$$\underline{\check{\mathrm{C}}}^\bullet(\mathcal{B}_1) = \underline{\check{\mathrm{C}}}^\bullet(\mathcal{B}_2)$$

*of finite correspondences between augmented Čech complexes.*

*Proof.* This can, using a limit argument, be extracted from the results of [Ivo07], Section 2.3. We give an independent argument:

Let the bridges be

$$\mathcal{B}_1 = (f, g, \gamma_1 : \Gamma_1 \to X \times_S Y, \widetilde{\alpha}_1, b_1 : \mathcal{U} \times_X \Gamma_1 \to \Gamma_1 \times_Y \mathcal{V}_1),$$
$$\mathcal{B}_2 = (f, g, \gamma_2 : \Gamma_2 \to X \times_S Y, \widetilde{\alpha}_2, b_2 : \mathcal{U} \times_X \Gamma_2 \to \Gamma_2 \times_Y \mathcal{V}_2).$$

We consider the following five pylons:

- for $i \in \{1,2\}$ the given pylons $(\gamma_i \colon \Gamma_i \to X \times_S Y, \widetilde{\alpha}_i)$,

- for $i \in \{1,2\}$ the reduced pylons $(\gamma_i^{\mathrm{red}} \colon \Gamma_i^{\mathrm{red}} \to X \times_S Y, \widetilde{\alpha}_i)$, where we abusively used $\widetilde{\alpha}_i$ to also denote the cycle on the homeomorphic reduced scheme $\Gamma_i^{\mathrm{red}} \hookrightarrow \Gamma_i$ (see Definition 1.4.20),

- the image pylon $(i \colon \Gamma' \hookrightarrow X \times_S Y, \alpha)$, where $\Gamma' := \gamma_1(\Gamma_1) \cup \gamma_2(\Gamma_2)$ with the reduced induced structure, and where we identified $\alpha$ with its counterpart in $\operatorname{supp}(\alpha) \subseteq \Gamma'$ (see again Definition 1.4.20).

Note that the images $\gamma_i(\Gamma_i) \subseteq X \times_S Y$ are indeed closed, e.g. by Lemma 1.4 of [MVW06]. As this does not influence the resulting finite correspondences we may also assume that all five pylons are pseudo-dominant over $X$ by removing superfluous irreducible components.

We then have the following diagram of schemes which are finite and pseudo-dominant over $X$:

$$\begin{array}{ccccccc} & \Gamma_1^{\mathrm{red}} & & & \Gamma_2^{\mathrm{red}} & & \\ \swarrow & & \overset{\gamma_1^{\mathrm{red}}}{\searrow} & & \overset{\gamma_2^{\mathrm{red}}}{\swarrow} & & \searrow \\ \Gamma_1 & & & \Gamma' & & & \Gamma_2. \end{array}$$

We assumed that $X$ is normal, hence every rigidification is strong by Remark 6.1.2. In particular, we see that any two (strongly) c-rigidified covers have a common rigid refinement into a (strongly) c-rigidified Nisnevich cover (cf. Lemma 6.2.4).

Let $u \colon \mathcal{U}' \to \mathcal{U}$ be a rigid refinement of c-rigidified Nisnevich covers. The morphisms $u^n = \check{\mathrm{C}}^{-n+1}(u) \colon \check{\mathrm{C}}^{-n+1}(\mathcal{U}' \twoheadrightarrow X) \to \check{\mathrm{C}}^{-n+1}(\mathcal{U} \twoheadrightarrow X)$ are therefore by Lemma 6.4.8 epimorphisms in $\mathrm{SchCor}_S$. Thus we are allowed to rigidly refine $\mathcal{U}$ to check the equalities.

Thus, using Proposition 6.4.6, we may assume to have a rigid bridging $\widetilde{b} \colon \mathcal{U} \times_X \Gamma' \to \Gamma' \times_Y \mathcal{V}$. It and the two given bridgings $b_1$, $b_2$ now induce the following bridgings on the two reduced pylons for $i \in \{1,2\}$:

- the pullback $b_i^{\mathrm{red}} := b_i \times_{\Gamma_i} \Gamma_i^{\mathrm{red}} \colon \mathcal{U} \times_X \Gamma_i^{\mathrm{red}} \to \Gamma_i^{\mathrm{red}} \times_Y \mathcal{V}$,

- the pullback $\widetilde{b}_i := \widetilde{b} \times_{\Gamma'} \Gamma_i^{\mathrm{red}} \colon \mathcal{U} \times_X \Gamma_i^{\mathrm{red}} \to \Gamma_i^{\mathrm{red}} \times_Y \mathcal{V}$.

As (pre-)rigidified pullbacks of rigid morphisms they are again rigid. Hence Lemma 6.4.7 tells us that $b_i^{\mathrm{red}} = \widetilde{b}_i$. We now have bridges $\widetilde{\mathcal{B}}$, $\mathcal{B}_i$, and $\mathcal{B}_i^{\mathrm{red}} = \widetilde{\mathcal{B}}_i$ resulting from the bridgings $\widetilde{b}$, $b_i$ and $b_i^{\mathrm{red}} = \widetilde{b}_i^{\mathrm{red}}$. We thus have the following diagram of refinements of bridges between the same Nisnevich covers $\mathcal{U} \twoheadrightarrow X$ and $\mathcal{V} \twoheadrightarrow Y$:

$$\begin{array}{ccccccccc} & \mathcal{B}_1^{\mathrm{red}} & = & \widetilde{\mathcal{B}}_1 & & \widetilde{\mathcal{B}}_2 & = & \mathcal{B}_2^{\mathrm{red}} & \\ \swarrow & & & & \searrow \quad \swarrow & & & & \searrow \\ \mathcal{B}_1 & & & & \widetilde{\mathcal{B}} & & & & \mathcal{B}_2. \end{array}$$

Taking fibre products of bridges is directly seen to preserve refinements. Hence Proposition 5.2.10 shows that

$$\alpha_{(\mathcal{B}_1)^n} = \alpha_{(\mathcal{B}_1^{\mathrm{red}})^n} = \alpha_{(\widetilde{\mathcal{B}}_1^{\mathrm{red}})^n} = \alpha_{(\widetilde{\mathcal{B}})^n} = \alpha_{(\widetilde{\mathcal{B}}_2^{\mathrm{red}})^n} = \alpha_{(\mathcal{B}_2^{\mathrm{red}})^n} = \alpha_{(\mathcal{B}_2)^n},$$

where $(-)^n$ was used to denote the $n$-fold fibre product of the respective bridges with themselves. By definition this means that we have an equality $\underline{\check{\mathrm{C}}}^\bullet(\mathcal{B}_1) = \underline{\check{\mathrm{C}}}^\bullet(\mathcal{B}_2)$ of finite correspondences. □

**Remark 6.4.10.** If the base $S$ is regular, then the results of this section generalize to $\mathrm{SmCor}(S, \Lambda)$ by $\Lambda$-linear extension. To be more precise, the existence and uniqueness results are true for the individual basic cycles.

## 6.5 Nisnevich covers on diagrams

**Definition 6.5.1.** Let $S$ be a regular noetherian scheme, let $\Lambda$ be a ring and let $X\colon D \to \mathrm{SmCor}(S)$ be a diagram. A *(rigidified) Nisnevich cover* $f\colon (\mathcal{U}|I) \twoheadrightarrow_r X$ on the diagram $X$ consists of

- a diagram $\mathcal{U}\colon D \to \mathrm{SmCor}(S, \Lambda)$,
- for each $d \in D$ an indexing set $I(d)$ and a decomposition $\mathcal{U}(d) = \coprod_{i \in I(d)} \mathcal{U}_i(d)$,
- a natural transformation $f\colon \mathcal{U} \implies X$ which turns every

  $$f(d)\colon (\mathcal{U}(d)|I(d)) \to X(d),$$

  $d \in D$, into a Nisnevich cover,
- for each $d \in D$ a constructible rigidification $r(d)\colon X(d) \to \mathcal{U}(d)$

such that the finite correspondence $\mathcal{U}(\alpha)$ is rigidly bridgeable over $X(\alpha)$ for every arrow $\alpha$ of $D$. The correspondences are hence unique by Theorem 6.4.9.

We call it a *cd-Nisnevich cover*, a *unifibrant* Nisnevich cover or *unifibrantly rigidified* if the respective properties hold at the individual $d \in D$. We call it *affine* if the individual $\mathcal{U}(d)$, $d \in D$, are affine schemes.

A *Nisnevich pre-cover* of a diagram $X\colon D \to \mathrm{SmCor}(S)$ consists of a choice of constructable pre-rigidified Nisnevich covers $\mathcal{U}(d) \to X(d)$ for all $d \in D$, with no choice of finite correspondences between the $\mathcal{U}(d)$. A *refinement of a Nisnevich pre-cover by a Nisnevich cover* is defined as such a rigid refinement at each object.

**Definition 6.5.2.** Let $S$ be a regular noetherian scheme and let $\Lambda$ be a ring. Let $\mathcal{U}\colon D \to \mathrm{SmCor}(S,\Lambda)$ be a Nisnevich cover of a diagram $X\colon D \to \mathrm{SmCor}(S,\Lambda)$. We define its *(full) Čech complex* as the diagram

$$\check{\mathrm{C}}^\bullet(\mathcal{U} \twoheadrightarrow X)\colon D \to C^-(\mathrm{SmCor}(S,\Lambda))$$

defined as follows:

The objects are the Čech complexes

$$\check{\mathrm{C}}^\bullet(\mathcal{U} \twoheadrightarrow X)(d) = \check{\mathrm{C}}^\bullet(\mathcal{U}(d) \twoheadrightarrow X(d)),$$

$d \in D$. The finite correspondence at an arrow $\alpha\colon c \to d$ of $D$ is

$$\check{\mathrm{C}}^\bullet(\mathcal{B}(\alpha))\colon\ \check{\mathrm{C}}^\bullet(\mathcal{U} \twoheadrightarrow X)(c) \rightsquigarrow \check{\mathrm{C}}^\bullet(\mathcal{U} \twoheadrightarrow X)(d)$$

as in Definition 5.4.4, where $\mathcal{B}(\alpha)$ is any rigid bridge inducing the rigidly bridgeable finite correspondence $\mathcal{U}(\alpha)\colon \mathcal{U}(c) \rightsquigarrow \mathcal{U}(d)$.

Note that the resulting finite correspondences between the Čech complexes do not depend on the choice of the bridge $\mathcal{B}(\alpha)$ by Theorem 6.4.9. Combined with Theorem 6.4.4 and Proposition 5.4.5 this also shows that we indeed got a diagram.

**Theorem 6.5.3.** *Let $S$ be a regular noetherian scheme and let $\Lambda$ be a ring. Let $X\colon D \to \mathrm{SmCor}(S,\Lambda)$ be a finite acyclic diagram of finite correspondences.*

*Then there exists an affine and unifibrantly rigidified cd-Nisnevich cover of $X$ finer than a given Nisnevich pre-cover. In particular such a Nisnevich cover exists and they form a directed system.*

*Proof.* We start with the case where $X \in \mathrm{Sm}_S$ is a single object. If $\mathcal{U} \twoheadrightarrow X$ is an arbitrary c-rigidified Nisnevich (pre-)cover of $X$, then it has by Lemma 6.2.5 a refinement into a unifibrantly c-rigidified cd-Nisnevich cover $f\colon \mathcal{V} \twoheadrightarrow X$. By taking another refinement by an open affine cover we arrive at an affine and unifibrantly c-rigidified cd-Nisnevich cover refining $\mathcal{U}$.

The general case then follows by induction:

By finiteness and acyclicity there is an object $d_{\min} \in D$ into which no other object maps. We may inductively assume that the diagram $X\backslash X(d_{\min})\colon D\backslash d_{\min} \to \mathrm{SmCor}(S,\Lambda)$ has an affine and unifibrantly rigidified cd-Nisnevich cover $\mathcal{V}$ refining a given Nisnevich pre-cover. For each of the finitely many arrows $f_i\colon d_{\min} \to d_i$ from $d_{\min}$ into an object of $D\backslash d_{\min}$ we find by Proposition 6.4.6 a rigidified Nisnevich cover $\mathcal{U}_i \twoheadrightarrow X(d_{\min})$ together with a rigid bridging into the previously constructed Nisnevich cover $\mathcal{V}(d_i) \twoheadrightarrow X(d_i)$.

By Lemma 6.2.4 we can refine the $\mathcal{U}_i$ and the given c-rigidified Nisnevich (pre-)cover on $X(d_{\min})$ into a single c-rigidified Nisnevich cover. From the case of a single object we know that we can refine it further into an affine and unifibrantly c-rigidified cd-Nisnevich cover. The rigid bridgings automatically extend to this refinement. This finishes the proof. □

## 6.6 Čech diagrams of Nisnevich covers on diagrams

**Definition 6.6.1.** Let $f\colon \mathcal{U} \to X$ be a Nisnevich cover of a diagram $X\colon D \to \mathrm{SmCor}(S, \Lambda)$ as in Definition 6.5.1 and recall Definition 6.5.2.

Let $\alpha\colon c \to d$ be an arrow in $D$ and let $n \in \mathbb{N}_0$. Let $U$ and $V$ be connected components of $\check{\mathrm{C}}^{-n}(\mathcal{U}(c))$ and $\check{\mathrm{C}}^{-n}(\mathcal{U}(d))$, respectively. Then the restriction of the finite correspondence

$$\check{\mathrm{C}}^{-n}(\mathcal{U} \twoheadrightarrow X)(\alpha)\colon\ \check{\mathrm{C}}^{-n}(\mathcal{U}(c)) \rightsquigarrow \check{\mathrm{C}}^{-n}(\mathcal{U}(d))$$

to $U$ and $V$ (cf. Proposition 1.6.4) is called a *simple finite correspondence.*

The *Čech diagram* $\check{\mathrm{D}}(\mathcal{U}) = \check{\mathrm{D}}(f)$ of $f$ is the following diagram:

- The set of objects is the disjoint union $\coprod_{d\in D} \check{\mathrm{D}}(f(d))$ of the individual Čech diagrams (cf. Definition 5.5.6).
- The morphisms are all the finite correspondences arising as valid compositions $\alpha_n \circ \ldots \circ \alpha_1$, where each $\alpha_i$ is either a simple finite correspondence as defined above or a morphism inside an individual $\check{\mathrm{D}}(f(d))$, $d \in D$.

The following is immediate from Definitions 4.5.3, 5.5.6 and 6.6.1:

**Theorem 6.6.2.** *Let $S$ be a regular noetherian scheme and let $\Lambda$ be a ring. Furthermore, let $f\colon \mathcal{U} \to X$ be a Nisnevich cover of a diagram $X\colon D \to \mathrm{SmCor}(S, \Lambda)$.*

*Then a $\Lambda$-cellular filtration on the Čech diagram $\check{\mathrm{D}}(f)$ induces a $\Lambda$-cellular filtration on the Čech complex $\check{\mathrm{C}}^{\bullet}(f)$.*

**Proposition 6.6.3.** *Let $S$ be a regular scheme of finite dimension and let $\Lambda$ be a ring. Let $X\colon D \to \mathrm{SmCor}(S, \Lambda)$ be a finite acyclic diagram and let $f\colon \mathcal{U} \to X$ be a Nisnevich cover of $X$.*

*Then the Čech diagram $\check{\mathrm{D}}(f)$ is finite and hom-almost acyclic (cf. Definition 4.8.3).*

*Proof.* The hom-almost acyclicity and the finiteness of the set of objects follows directly from Proposition 5.5.7 and the definitions. Therefore it suffices to prove the finiteness of the hom-sets.

Every arrow in the Čech diagram $\check{\mathrm{D}}(f)$ can be written as a chain

$$U_1 \to V_1 \rightsquigarrow U_2 \to V_2 \rightsquigarrow \ldots \rightsquigarrow U_m \to V_m,$$

where every $U_i \to V_i$ is a morphism inside a $\check{\mathrm{D}}(f(d))$, $d \in D$, and the finite correspondences are simple finite correspondences over a non-identity arrow of $D$. Due to the acyclicity of $D$ the integer $m$ is bounded by the number of objects of $D$. Therefore, and by the finiteness of the $\check{\mathrm{D}}(f(c))$ (cf. Proposition 5.5.7), it suffices to show that there are only finitely many simple finite correspondences $V_i \rightsquigarrow U_{i+1}$. This amounts to the following claim:

Given a bridge

$$\mathcal{B} = (f : \mathcal{U} \twoheadrightarrow X, g : \mathcal{V} \twoheadrightarrow Y, \gamma : \Gamma \to X \times_S Y, \widetilde{\alpha}, b),$$

over a finite correspondence $\alpha : X \rightsquigarrow Y$ there are only finitely many finite correspondences $U \rightsquigarrow V$ that come from restricting an $\alpha_{\mathcal{B}^n} : (\mathcal{U}|X)^n \rightsquigarrow (\mathcal{V}|Y)^n$ to connected components isomorphic to $U$ and $V$ over $X$ and $Y$, respectively. Here we are allowed to vary $n$ freely.

While proving this we may assume that $X$ is connected. We also observe that every such finite correspondence $U \rightsquigarrow (\mathcal{V}|Y)^n$ is bridgeable over the pylon $(\gamma, \widetilde{\alpha})$.

Let $d$ be the degree of the finite morphism $\mathrm{pr}_X \circ \gamma : \Gamma \to X$. Then the number of connected components of $U \times_X \Gamma$ is bounded by $d$. Each of them is mapped by the bridging $b$ to one of the finitely many $\Gamma \times_Y V_j$, where $V \in \check{\mathrm{D}}(f(d))$. Therefore we found that every simple finite correspondence $U \rightsquigarrow V$ comes from restricting a finite correspondence $U \rightsquigarrow \coprod_{i=1}^{d} V_d$ which is bridgeable over the pylon $(\gamma, \widetilde{\alpha})$ to $V = V_1$. Here the $V_j$ are not necessarily distinct elements of $\check{\mathrm{D}}(\mathcal{V} \twoheadrightarrow Y)$.

As $d$ is fixed and $\check{\mathrm{D}}(\mathcal{V} \twoheadrightarrow Y)$ is finite we note that there are only finitely many isomorphism classes of such coproducts $\coprod_{i=1}^{d} V_d$ over $Y$. Furthermore, there are only finitely many bridgings

$$b' : U \times_X \Gamma \to \Gamma \times_Y \coprod_{i=1}^{d} V_d$$

over $\Gamma$ by Lemma 5.5.4. In total this shows the desired finiteness. □

**Remark 6.6.4.** Note that we cannot use the rigidifications for two reasons: we identified many isomorphic schemes and we only have pre-rigidifications on the higher powers $(\mathcal{U}|X)^n$.

# Chapter 7

# Realization between Motives

For the entirety of this chapter we fix a field $k \subseteq \mathbb{C}$ and a noetherian ring $\Lambda$.

We combine the results of Chapters 4 and 6. Let us briefly sketch what is about to happen:

Chapter 6 associates to an object of $\mathrm{SmCor}(k, \Lambda)$ the pro-system of affine (unifibrantly) rigidified Nisnevich covers, with rigidly bridgeable morphisms between them. Taking Čech complexes results in a covariant functor

$$\mathrm{SmCor}(k, \Lambda) \to \text{pro-}C^{-}(\mathrm{SmCor}^{\mathrm{aff}}(k, \Lambda)).$$

On the other hand, the results of Chapter 4 give us a contravariant functor

$$\mathrm{SmCor}^{\mathrm{aff}}(k, \Lambda) \to \text{pro-}C^{b}(\mathcal{MM}^{\mathrm{eff}}_{\mathrm{Nori}}(k, \Lambda)),$$

which is compatible with tensor products.

A careful yoga of pro-ind-systems allows us to combine, refine and extend the above, resulting in a contravariant functor

$$C^{b}(\mathrm{SmCor}^{\mathrm{aff}}(k, \Lambda)) \to \text{pro-ind-}C^{+}(\mathcal{MM}^{\mathrm{eff}}_{\mathrm{Nori}}(k, \Lambda)).$$

Note that the ind-system arises from the contravariance of the second functor.

Afterwards we check that this descends to the desired functor. We show that the constructed functor actually lands in systems of quasi-isomorphisms of bounded complexes, hence it induces a contravariant functor

$$C^{b}(\mathrm{SmCor}^{\mathrm{aff}}(k, \Lambda)) \to D^{b}(\mathcal{MM}^{\mathrm{eff}}_{\mathrm{Nori}}(k, \Lambda)).$$

Then it localizes on the left hand side by a singular cohomology argument to $\mathrm{DM}^{\mathrm{eff}}_{\mathrm{gm}}(k, \Lambda)$.

We then compare the tensor structures. The central observation is that the functor returns complexes of flat objects.

## 7.1 Combining covers and filtrations

The following observation is simple, yet crucial:

**Lemma 7.1.1.** *Let $U \to X$ and $V \to X$ be morphisms of schemes. Assume that $X$ is separated and that $U$ and $V$ are affine.*

*Then $U \times_X V$ is an affine scheme.*

*Proof.* The morphism $X \to \mathrm{Spec}(\mathbb{Z})$ is separated, therefore the diagonal $\Delta\colon X \to X \times_{\mathrm{Spec}(\mathbb{Z})} X$ is a closed immersion. Now $U \times_X V$ is the pullback of the affine scheme $U \times_{\mathrm{Spec}(\mathbb{Z})} V \cong \mathrm{Spec}(\mathcal{O}_U \otimes_{\mathbb{Z}} \mathcal{O}_V)$ along $\Delta$, hence again affine. □

**Definition 7.1.2.** Let $X\colon \mathcal{C} \to \mathrm{SmCor}(k, \Lambda)$ be a diagram. A $\Lambda$-*cellular cover-filtration* on $X$ is a pair $((\mathcal{U}|I) \twoheadrightarrow X, \mathcal{F})$ consisting of:

- an affine Nisnevich cover $(\mathcal{U}|I) \twoheadrightarrow X$ (cf. Definition 6.5.1),
- a functorial $\Lambda$-cellular filtration (cf. Definition 4.8.2) $\mathcal{F}$ on the Čech complex $\check{\mathrm{C}}^\bullet(\mathcal{U} \twoheadrightarrow X)$ (cf. Definition 6.5.2).

**Theorem 7.1.3.** *Let $k \subseteq \mathbb{C}$ be a field and let $\Lambda$ be a noetherian ring.*

*Then every finite acyclic diagram $X\colon \mathcal{C} \to \mathrm{SmCor}(k, \Lambda)$ has a $\Lambda$-cellular cover-filtration.*

*Proof.* By Theorem 6.5.3 the diagram $X$ admits an affine Nisnevich cover $\mathcal{U} \twoheadrightarrow X$. The resulting Čech diagram $\check{\mathrm{D}}(\mathcal{U} \twoheadrightarrow X)$ is finite and hom-almost acyclic by Proposition 6.6.3, hence it has a functorial $\Lambda$-cellular filtration due to Theorem 4.8.4. Together with Theorem 6.6.2 this induces the desired $\Lambda$-cellular filtration on the Čech complex $\check{\mathrm{C}}^\bullet(\mathcal{U} \twoheadrightarrow X)$. □

**Remark 7.1.4.** Unlike seen in Theorems 4.8.4 and 6.5.3, $\Lambda$-cellular cover-filtrations do not form a directed system. This is due to a difference in directions: Nisnevich covers are constructed contravariantly, while $\Lambda$-cellular filtrations are constructed covariantly. Hence we cannot simply extend from an initial or terminal segment. This is also the reason why we construct all the Nisnevich covers first, and only then the filtrations.

**Definition 7.1.5.** Let $X^\bullet \in C^b(\mathrm{SmCor}(k, \Lambda))$ and let $\mathfrak{F} = ((\mathcal{U}|I) \twoheadrightarrow X, \mathcal{F})$ be a $\Lambda$-cellular cover-filtration on $X^\bullet$. Then, generalizing the $\mathrm{C}^\bullet_{\mathcal{F}}$ of Definition 4.2.1, its *Nori complex* is defined as the total complex

$$\mathrm{C}^\bullet_{\mathfrak{F}}(X^\bullet) := \mathrm{Tot}^\bullet_{i,-j,-k}(\mathrm{C}^i_{\mathcal{F}}(\check{\mathrm{C}}^j(\mathcal{U}^k \twoheadrightarrow X^k))) \in C^+(\mathcal{MM}^{\mathrm{eff}}_{\mathrm{Nori}}(k, \Lambda)).$$

Note that the minus signs in $\mathrm{Tot}_{i,-j,-k}$ come from the contravariance of $\mathrm{C}^i_{\mathcal{F}}$, which is also the reason why we end up in $C^+(\mathcal{MM}^{\mathrm{eff}}_{\mathrm{Nori}}(k, \Lambda))$.

**Definition 7.1.6.** Let $\alpha^\bullet \colon X^\bullet \rightsquigarrow Y^\bullet$ be a morphism in $\mathrm{SmCor}(k,\Lambda)$ and let $\mathfrak{F}$ be a $\Lambda$-cellular cover-filtration on this diagram. Then the functoriality of the $\Lambda$-cellular filtration induces contravariantly a morphism of the defining triple complexes and hence a morphism

$$\mathrm{C}^\bullet_{\mathfrak{F}}(\alpha)\colon\ \mathrm{C}^\bullet_{\mathfrak{F}}(Y^\bullet) \to \mathrm{C}^\bullet_{\mathfrak{F}}(X^\bullet)$$

in $C^+(\mathcal{MM}^{\mathrm{eff}}_{\mathrm{Nori}}(k,\Lambda))$.

**Remark 7.1.7.** Instead of using the Čech diagram, one could also use a cut-off-construction:

We add an integer $N$ to the data of a $\Lambda$-cellular cover-filtration. Then we only choose $\Lambda$-cellular filtrations on the $\check{\mathrm{C}}^{-n}(\mathcal{U}(d))$ with $n \leq N$ and cut the Nori complex off at degrees above $N$.

The boundedness of singular cohomology then shows that everything we do remains true as long as the degree is sufficiently far away from $N$. One then recovers our results by choosing $N$ sufficiently large. This, however, makes many arguments, the proof of Theorem 7.4.10 in particular, significantly more complicated.

**Lemma 7.1.8.** *The derived forgetful functors*

$$D(\omega_{\mathrm{sing}})\colon D(\mathcal{MM}^{\mathrm{eff}}_{\mathrm{Nori}}(k,\Lambda)) \to D(\Lambda\text{-}\mathrm{Mod})$$
$$D(\omega_{\mathrm{sing}})\colon D(\mathcal{MM}_{\mathrm{Nori}}(k,\Lambda)) \to D(\Lambda\text{-}\mathrm{Mod})$$

*are conservative, i.e. reflect isomorphisms.*

*Proof.* A morphism $f$ in the derived category $D(\mathcal{A})$ of an abelian category $\mathcal{A}$ is an isomorphism if and only if the induced morphisms $H^n(f)$ on cohomology are isomorphisms. The latter property is clearly preserved and reflected by an exact and faithful functor such as $\omega_{\mathrm{sing}}$, implying the result. □

**Proposition 7.1.9.** *Let $X \in \mathrm{Var}_k$ be a variety over a field $k \subseteq \mathbb{C}$ and let $\Lambda$ be a noetherian ring. Let $(\mathcal{U}|I) \twoheadrightarrow X$ be an affine Nisnevich cover of $X$ and let $\mathcal{F}$ be a $\Lambda$-cellular filtration on the Čech complex $\check{\mathrm{C}}^\bullet(\mathcal{U} \twoheadrightarrow X)$.*

*Then there exists a natural isomorphism*

$$\omega_{\mathrm{sing}}(H^n(\mathrm{Tot}^\bullet_{i,-j}(\mathrm{C}^i_{\mathcal{F}}(\check{\mathrm{C}}^j(\mathcal{U} \twoheadrightarrow X)))) \cong H^n_{\mathrm{sing}}(X,\Lambda).$$

*Proof.* The double complex $\mathrm{C}^\bullet_{\mathcal{F}}(\check{\mathrm{C}}^\bullet(\mathcal{U} \twoheadrightarrow X))$ induces a spectral sequence

$$E_2^{p,q} = H^q(H^p(\mathrm{C}^\bullet_{\mathcal{F}}(\check{\mathrm{C}}^\bullet(\mathcal{U} \twoheadrightarrow X)))) \implies H^{p+q}(\mathrm{Tot}^\bullet(\mathrm{C}^\bullet_{\mathcal{F}}(\check{\mathrm{C}}^\bullet(\mathcal{U} \twoheadrightarrow X)))).$$

We apply the exact forgetful functor

$$\omega_{\mathrm{sing}}\colon \mathcal{MM}^{\mathrm{eff}}_{\mathrm{Nori}}(k,\Lambda) \to \Lambda\text{-}\mathrm{Mod},$$

turning the spectral sequence by Proposition 4.2.3 into

$$E_2^{p,q} = H^q(H^p_{\mathrm{sing}}(\check{\mathrm{C}}^\bullet(\mathcal{U} \twoheadrightarrow X))) \implies \omega_{\mathrm{sing}}(H^{p+q}(\mathrm{Tot}^\bullet(\mathrm{C}^i_{\mathcal{F}}(\check{\mathrm{C}}^j(\mathcal{U} \twoheadrightarrow X))))).$$

Using Proposition 4.2.3, we see that the $E_2$-page is that of the Čech spectral sequence, i.e. Theorem B.0.16. A more detailed consideration now shows that we have a morphism between these two spectral sequences, and therefore an isomorphism between the limit terms. □

**Remark 7.1.10.** The same argument can be used to show an analogous result using the repetition-free version $\check{c}^j(\mathcal{U} \twoheadrightarrow X)$ of the Čech complex, assuming that the Nisnevich cover $(\mathcal{U}|I) \twoheadrightarrow X$ is unifibrant.

**Corollary 7.1.11.** *Let $X \in \mathrm{Sm}_k$ be a smooth variety over a field $k \subseteq \mathbb{C}$ and let $\Lambda$ be a noetherian ring. Let $\mathfrak{F}$ be a $\Lambda$-cellular cover-filtration on $X$.*

*Then there exists a natural isomorphism*

$$\omega_{\mathrm{sing}}(H^n(\mathrm{C}^\bullet_{\mathfrak{F}}(X))) \cong H^n_{\mathrm{sing}}(X, \Lambda),$$

*i.e. the Nori complex $\mathrm{C}^\bullet_{\mathfrak{F}}(X)$ calculates singular cohomology after application of the forgetful exact functor $\omega_{\mathrm{sing}} \colon \mathcal{MM}^{\mathrm{eff}}_{\mathrm{Nori}}(k, \Lambda) \to \Lambda\text{-}\mathrm{Mod}$.*

**Corollary 7.1.12.** *Let $X \in \mathrm{Sm}^{\mathrm{aff}}_k$ be a smooth affine variety over a field $k \subseteq \mathbb{C}$ and let $\Lambda$ be a noetherian ring. Let $\mathfrak{F} = (f \colon \mathcal{U} \twoheadrightarrow X, \mathcal{F})$ be a $\Lambda$-cellular cover-filtration on $X$ and let $\mathcal{G}$ be a $\Lambda$-cellular filtration on $X$ such that $f$ is compatible with the given filtrations $\mathcal{F}$ and $\mathcal{G}$ on $\mathcal{U}$ and $X$.*

*Then there exists a natural quasi-isomorphism*

$$\mathrm{C}^\bullet_{\mathcal{G}}(X) \to \mathrm{C}^\bullet_{\mathfrak{F}}(X).$$

*Proof.* Consider the structure morphism $\check{\mathrm{C}}^\bullet(\mathcal{U}) \to X[0]$, which immediately induces the required morphism. It being a quasi-isomorphism can by Lemma 7.1.8 be checked in singular cohomology, where it follows from Proposition 4.2.3 and Corollary 7.1.11. □

## 7.2 Morphisms of cover-filtrations

**Definition 7.2.1.** Let $X \colon \mathcal{C} \to \mathrm{SmCor}(k, \Lambda)$ be a diagram. Let

$$\mathfrak{F} = ((\mathcal{U}|I) \twoheadrightarrow X, \mathcal{F}),$$
$$\mathfrak{G} = ((\mathcal{V}|J) \twoheadrightarrow X, \mathcal{G})$$

be two $\Lambda$-cellular cover-filtrations.

A *functorial morphism* $\mathfrak{G} \to \mathfrak{F}$ between these $\Lambda$-cellular cover-filtrations is a rigid refinement $h \colon \mathcal{V} \to \mathcal{U}$ of Nisnevich covers of $X$ such that the morphisms $\check{\mathrm{C}}^{-n}(\mathcal{V} \twoheadrightarrow X) \to \check{\mathrm{C}}^{-n}(\mathcal{U} \twoheadrightarrow X)$ are compatible with the given $\Lambda$-cellular filtrations and such that the induced diagrams

$$\begin{array}{ccc} \mathrm{C}^\bullet_{\mathcal{F}}(\check{\mathrm{C}}^{-n}(\mathcal{U}(y) \twoheadrightarrow X(y))) & \xrightarrow{\mathrm{C}^\bullet_{\mathcal{F}}(\check{\mathrm{C}}^{-n}(\mathcal{U}\twoheadrightarrow X)(c))} & \mathrm{C}^\bullet_{\mathcal{F}}(\check{\mathrm{C}}^{-n}(\mathcal{U}(x) \twoheadrightarrow X(x))) \\ \downarrow {\scriptstyle \mathrm{C}^\bullet_{\mathcal{F}}(\check{\mathrm{C}}^{-n}(h(y)))} & & \downarrow {\scriptstyle \mathrm{C}^\bullet_{\mathcal{F}}(\check{\mathrm{C}}^{-n}(h(x)))} \\ \mathrm{C}^\bullet_{\mathcal{G}}(\check{\mathrm{C}}^{-n}(\mathcal{V}(y) \twoheadrightarrow X(y))) & \xrightarrow{\mathrm{C}^\bullet_{\mathcal{G}}(\check{\mathrm{C}}^{-n}(\mathcal{U}\twoheadrightarrow X)(c))} & \mathrm{C}^\bullet_{\mathcal{G}}(\check{\mathrm{C}}^{-n}(\mathcal{V}(x) \twoheadrightarrow X(x))) \end{array}$$

commute. Here $c$ runs over all arrows $x \to y$ in $\mathcal{C}$ and $n$ over all non-negative integers.

**Remark 7.2.2.** Let us explain the commutativity condition:

It is by definition equivalent to the functoriality of the $\Lambda$-cellular filtration on the Čech diagram resulting from the refinement $\mathcal{V} \to \mathcal{U}$. Its foremost goal is to make Definition 7.1.5 functorial for $\Lambda$-cellular cover-filtrations by inducing morphisms

$$\mathrm{C}^i_{\mathcal{F}}(\check{\mathrm{C}}^j(\mathcal{U}^k \twoheadrightarrow X^k)) \to \mathrm{C}^i_{\mathcal{G}}(\check{\mathrm{C}}^j(\mathcal{V}^k \twoheadrightarrow X^k))$$

of triple complexes. Indeed, the other commutativities are already true by the definitions.

In light of Theorem 4.8.1, the condition means that $\mathcal{F}$ is coarse enough. Indeed, we have the eventual commutativity of $\Lambda$-cellular filtrations witnessed throughout Chapter 4.

Also note that if $\mathcal{U} = \mathcal{V}$, then the required commutativity is automatic by Lemma 4.6.3.

**Proposition 7.2.3.** *Let $k \subseteq \mathbb{C}$ be a field and let $\Lambda$ be a noetherian ring.*

*(a) Let $X^\bullet \in C^b(\mathrm{SmCor}(k, \Lambda)$ be a bounded complex of finite correspondences. Also let $u\colon \mathfrak{G} \to \mathfrak{F}$ be a functorial morphism of $\Lambda$-cellular cover-filtrations on $X^\bullet$.*

*Then it induces a quasi-isomorphism*

$$\varphi^\bullet_u\colon\ \mathrm{C}^\bullet_{\mathfrak{F}}(X^\bullet) \to \mathrm{C}^\bullet_{\mathfrak{G}}(X^\bullet)$$

*between Nori complexes in $C^+(\mathcal{MM}^{\mathrm{eff}}_{\mathrm{Nori}}(k, \Lambda))$. It is contravariantly functorial with respect to composition of functorial morphisms of $\Lambda$-cellular cover-filtrations.*

*(b) Assume that $\alpha^\bullet\colon X^\bullet \rightsquigarrow Y^\bullet$ is a finite correspondence of complexes over $\mathrm{SmCor}(k, \Lambda)$ and let $u\colon \mathfrak{G} \to \mathfrak{F}$ be a functorial morphism of $\Lambda$-cellular cover-filtrations on this diagram.*

*Then the diagram*

$$\begin{array}{ccc} \mathrm{C}^\bullet_{\mathfrak{F}}(Y^\bullet) & \xrightarrow{\mathrm{C}^\bullet_{\mathfrak{F}}(\alpha^\bullet)} & \mathrm{C}^\bullet_{\mathfrak{F}}(X^\bullet) \\ \varphi^\bullet_u \downarrow \text{q.-is.} & & \varphi^\bullet_u \downarrow \text{q.-is.} \\ \mathrm{C}^\bullet_{\mathfrak{G}}(Y^\bullet) & \xrightarrow{\mathrm{C}^\bullet_{\mathfrak{G}}(\alpha^\bullet)} & \mathrm{C}^\bullet_{\mathfrak{G}}(X^\bullet) \end{array}$$

*of Nori complexes in* $C^+(\mathcal{MM}^{\mathrm{eff}}_{\mathrm{Nori}}(k,\Lambda))$ *commutes.*

*Proof.* Let $\mathfrak{F} = (\mathcal{U}^\bullet \twoheadrightarrow X^\bullet, \mathcal{F})$ and $\mathfrak{G} = (\mathcal{V}^\bullet \twoheadrightarrow X^\bullet, \mathcal{G})$.

The functorial morphism $u\colon \mathfrak{G} \to \mathfrak{F}$ induces, as explained in Remark 7.2.2, a morphism

$$\psi_u^{i,-j,-k}\colon\ \mathrm{C}^i_{\mathcal{G}}(\check{\mathrm{C}}^j(\mathcal{V}^k \twoheadrightarrow X^k)) \to \mathrm{C}^i_{\mathcal{F}}(\check{\mathrm{C}}^j(\mathcal{U}^k \twoheadrightarrow X^k)$$

of triple complexes and hence a morphism $\varphi^\bullet_u\colon\ \mathrm{C}^\bullet_{\mathfrak{F}}(X) \to \mathrm{C}^\bullet_{\mathfrak{G}}(X)$ between their total complexes. Corollary 7.1.11 shows that it is a quasi-isomorphism if $X^\bullet = X \in \mathrm{Sm}_k$ is a single variety.

Now $\psi_u^{i,-j,-k}$ induces a morphism between the spectral sequences associated to the double complexes

$$\begin{aligned} F^{p,-q} &= \mathrm{Tot}^p_{i,-j}(\mathrm{C}^i_{\mathcal{F}}(\check{\mathrm{C}}^j(\mathcal{U}^q \twoheadrightarrow X^q))) \\ G^{p,-q} &= \mathrm{Tot}^p_{i,-j}(\mathrm{C}^i_{\mathcal{G}}(\check{\mathrm{C}}^j(\mathcal{V}^q \twoheadrightarrow X^q))) \end{aligned}$$

which shows that $\varphi_u$ is a quasi-isomorphism for all $X^\bullet \in C^b(\mathrm{SmCor}(k,\Lambda))$.

The functoriality follows by construction. The commutativity is a direct consequence of the commutativity condition for functorial morphisms of $\Lambda$-cellular cover-filtrations. □

**Theorem 7.2.4.** *Let* $\Lambda$ *be a noetherian ring. Let* $X^\bullet \in C^b(\mathrm{SmCor}(k,\Lambda))$ *be a complex of finite* $\Lambda$*-correspondences between smooth varieties over a field* $k \subseteq \mathbb{C}$*. Let* $\mathfrak{F}$ *and* $\mathfrak{G}$ *be two* $\Lambda$*-cellular cover-filtrations on* $X$*.*

*Then in* $D^+(\mathcal{MM}^{\mathrm{eff}}_{\mathrm{Nori}}(k,\Lambda))$ *there is an isomorphism*

$$\varphi_{\mathfrak{F},\mathfrak{G}}\colon\ \mathrm{C}^\bullet_{\mathfrak{G}}(X^\bullet) \xrightarrow{\cong} \mathrm{C}^\bullet_{\mathfrak{F}}(X^\bullet)$$

*of the corresponding Nori complexes.*

*Furthermore, this isomorphism is natural in both the complex and the* $\Lambda$*-cellular cover-filtrations, by which we mean:*

*(a)* $\varphi_{\mathfrak{F},\mathfrak{G}} \circ \varphi_{\mathfrak{G},\mathfrak{F}} = \mathrm{id}_{\mathrm{C}^\bullet_{\mathfrak{F}}(X)}$ *for any two* $\Lambda$*-cellular cover-filtrations* $\mathfrak{F}$ *and* $\mathfrak{G}$ *on* $X^\bullet$*,*

*(b)* $\varphi_{\mathfrak{F},\mathfrak{G}} \circ \varphi_{\mathfrak{G},\mathfrak{H}} = \varphi_{\mathfrak{F},\mathfrak{H}}$ *for any three* $\Lambda$*-cellular cover-filtrations* $\mathfrak{F}$, $\mathfrak{G}$ *and* $\mathfrak{H}$ *on* $X^\bullet$,

*(c) given any finite correspondence* $\alpha^\bullet\colon X^\bullet \rightsquigarrow Y^\bullet$ *in* $C^b(\mathrm{SmCor}(k,\Lambda))$ *and any two* $\Lambda$*-cellular cover-filtrations* $\mathfrak{F}$, $\mathfrak{G}$ *on this diagram, the diagram*

$$\begin{array}{ccc} \mathrm{C}^\bullet_{\mathfrak{G}}(Y^\bullet) & \xrightarrow{\mathrm{C}^\bullet_{\mathfrak{G}}(\alpha^\bullet)} & \mathrm{C}^\bullet_{\mathfrak{G}}(X^\bullet) \\ {\scriptstyle \varphi_{\mathfrak{F},\mathfrak{G}}}\downarrow\cong & & {\scriptstyle \varphi_{\mathfrak{F},\mathfrak{G}}}\downarrow\cong \\ \mathrm{C}^\bullet_{\mathfrak{F}}(Y^\bullet) & \xrightarrow{\mathrm{C}^\bullet_{\mathfrak{F}}(\alpha^\bullet)} & \mathrm{C}^\bullet_{\mathfrak{F}}(X^\bullet) \end{array}$$

*in* $D^+(\mathcal{MM}^{\mathrm{eff}}_{\mathrm{Nori}}(k,\Lambda))$ *commutes.*

*Proof.* If $\mathfrak{G}$ is a refinement of $\mathfrak{F}$ by a functorial morphism $u$, we set $\varphi_{\mathfrak{F},\mathfrak{G}} = \varphi^\bullet_u$ to be the morphism obtained in Proposition 7.2.3.

The general case is reduced to the above via appropriate zigzags, i.e. yoga of pro-ind-systems:

Let $\mathfrak{F} = (\mathcal{U}^\bullet \twoheadrightarrow X^\bullet, \mathcal{F})$ and $\mathfrak{G} = (\mathcal{V}^\bullet \twoheadrightarrow X^\bullet, \mathcal{G})$.

We interpret bounded complexes and finite acyclic diagrams thereof as finite acyclic diagrams. By Theorem 6.5.3 the two Nisnevich covers have a common refinement into an affine unifibrantly c-rigidified Nisnevich cover $\mathcal{S}^\bullet$ on $X^\bullet$. Then, using Proposition 6.6.3 and Theorem 4.8.4, we find a functorial $\Lambda$-cellular filtration $\mathcal{I}$ on the diagram

$$\begin{array}{ccccc} & & \check{\mathrm{C}}^\bullet(\mathcal{S}^\bullet \twoheadrightarrow X^\bullet) & & \\ & \swarrow & & \searrow & \\ \check{\mathrm{C}}^\bullet(\mathcal{U}^\bullet \twoheadrightarrow X^\bullet) & & & & \check{\mathrm{C}}^\bullet(\mathcal{V}^\bullet \twoheadrightarrow X^\bullet) \end{array}$$

coarser than the given $\mathcal{F}$ and $\mathcal{G}$ on the bottom row. Hence we got a diagram

$$\begin{array}{ccccccc} (\mathcal{U},\mathcal{F}) & & & (\mathcal{S},\mathcal{I}) & & & (\mathcal{V},\mathcal{G}) \\ & \searrow & \swarrow & & \searrow & \swarrow & \\ & (\mathcal{U},\mathcal{I}) & & & (\mathcal{V},\mathcal{I}) & & . \end{array}$$

of refinements of $\Lambda$-cellular cover-filtrations. Recall that, as mentioned in Remark 7.2.2, the outer arrows are indeed refinements due to Lemma 4.6.3.

This zigzag defines by Proposition 7.2.3 an isomorphism

$$\varphi_{\mathfrak{F},\mathfrak{G}}\colon\ \mathrm{C}^\bullet_{\mathfrak{G}}(X) \to \mathrm{C}^\bullet_{\mathfrak{F}}(X)$$

in $D^+(\mathcal{MM}^{\mathrm{eff}}_{\mathrm{Nori}}(k,\Lambda))$.

We need to check that this definition is independent of the choices we made, which we do individually for each part:

If $\mathcal{I}'$ is a second functorial $\Lambda$-cellular filtration akin to $\mathcal{I}$, we can pick a common coarsening $\widetilde{\mathcal{I}}$. This results in the diagram

$$\begin{tikzcd}
& (\mathcal{U},\mathcal{I}') \arrow[d] & (\mathcal{S},\mathcal{I}') \arrow[l] \arrow[r] \arrow[d] & (\mathcal{V},\mathcal{I}') \arrow[d] & \\
(\mathcal{U},\mathcal{F}) \arrow[ur] \arrow[r] \arrow[dr] & (\mathcal{U},\widetilde{\mathcal{I}}) & (\mathcal{S},\widetilde{\mathcal{I}}) \arrow[l] \arrow[r] & (\mathcal{V},\widetilde{\mathcal{I}}) & (\mathcal{V},\mathcal{G}) \arrow[ul] \arrow[l] \arrow[dl] \\
& (\mathcal{U},\mathcal{I}) \arrow[u] & (\mathcal{S},\mathcal{I}) \arrow[l] \arrow[r] \arrow[u] & (\mathcal{V},\mathcal{I}) \arrow[u] &
\end{tikzcd}$$

of refinements, which clearly commutes by the naturality of Proposition 7.2.3. Hence the isomorphism does not depend on the choice of $\mathcal{I}$.

Now let $\mathcal{S}'$ be a second choice of an affine unifibrantly c-rigidified Nisnevich cover. We pick by Theorem 6.5.3 a common refinement $\widetilde{\mathcal{S}}$ of the two covers. Proposition 6.6.3 and Theorem 4.8.4 give us a functorial $\Lambda$-cellular filtration $\mathcal{I}$ on the diagram

$$\begin{tikzcd}
& \check{\mathrm{C}}^\bullet(\mathcal{S}'^\bullet \twoheadrightarrow X^\bullet) \arrow[dl] \arrow[dr] & \\
\check{\mathrm{C}}^\bullet(\mathcal{U}^\bullet \twoheadrightarrow X^\bullet) & \check{\mathrm{C}}^\bullet(\widetilde{\mathcal{S}}^\bullet \twoheadrightarrow X^\bullet) \arrow[l] \arrow[r] \arrow[u] \arrow[d] & \check{\mathrm{C}}^\bullet(\mathcal{V}^\bullet \twoheadrightarrow X^\bullet) \\
& \check{\mathrm{C}}^\bullet(\mathcal{S}^\bullet \twoheadrightarrow X^\bullet) \arrow[ul] \arrow[ur] &
\end{tikzcd}$$

coarser than the given ones on $\check{\mathrm{C}}^\bullet(\mathcal{U}^\bullet \twoheadrightarrow X^\bullet)$ and $\check{\mathrm{C}}^\bullet(\mathcal{V}^\bullet \twoheadrightarrow X^\bullet)$. This yields the diagram

$$\begin{tikzcd}
& (\mathcal{U},\mathcal{I}) \arrow[dd, equal] & (\mathcal{S}',\mathcal{I}) \arrow[l] \arrow[r] & (\mathcal{V},\mathcal{I}) \arrow[dd, equal] & \\
(\mathcal{U},\mathcal{F}) \arrow[ur] \arrow[dr] & & (\widetilde{\mathcal{S}},\mathcal{I}) \arrow[ul] \arrow[u] \arrow[ur] \arrow[dl] \arrow[d] \arrow[dr] & & (\mathcal{V},\mathcal{G}) \arrow[ul] \arrow[dl] \\
& (\mathcal{U},\mathcal{I}) & (\mathcal{S},\mathcal{I}) \arrow[l] \arrow[r] & (\mathcal{V},\mathcal{I}) &
\end{tikzcd}$$

of refinements. Once more the independence follows from Proposition 7.2.3.

We are left to check the naturalities:

(a) This is immediate from the symmetry of the construction.

(b) Let $\mathfrak{H} = (\mathcal{W}, \mathcal{H})$ be a third $\Lambda$-cellular cover-filtration. Then we pick a common refinement $\mathcal{S}$ of all three Nisnevich covers and after that a functorial $\Lambda$-cellular filtration $\mathcal{I}$ on the diagram

$$\begin{array}{ccccc} & & \check{\mathrm{C}}^{\bullet}(\mathcal{W} \twoheadrightarrow X) & & \\ & & \uparrow & & \\ & & \check{\mathrm{C}}^{\bullet}(\mathcal{S} \twoheadrightarrow X) & & \\ & \swarrow & & \searrow & \\ \check{\mathrm{C}}^{\bullet}(\mathcal{V} \twoheadrightarrow X) & & & & \check{\mathrm{C}}^{\bullet}(\mathcal{U} \twoheadrightarrow X) \end{array}$$

which is coarser than the given ones on the outside. Hence we get a diagram

$$\begin{array}{ccccccc} & & & (\mathcal{W}, \mathcal{H}) & & & \\ & & & \downarrow & & & \\ & & & (\mathcal{W}, \mathcal{I}) & & & \\ & & & \uparrow & & & \\ & & & (\mathcal{S}, \mathcal{I}) & & & \\ & & \swarrow & & \searrow & & \\ & (\mathcal{V}, \mathcal{I}) & & & & (\mathcal{U}, \mathcal{I}) & \\ \nearrow & & & & & & \nwarrow \\ (\mathcal{V}, \mathcal{G}) & & & & & & (\mathcal{U}, \mathcal{F}) \end{array}$$

of refinements which shows the equality $\varphi_{\mathfrak{F},\mathfrak{G}} \circ \varphi_{\mathfrak{G},\mathfrak{H}} = \varphi_{\mathfrak{F},\mathfrak{H}}$.

(c) This follows, using Proposition 7.2.3 (b), by repeating the arguments in the construction of $\varphi_{\mathfrak{F},\mathfrak{G}}$, but this time using a finite correspondence $\alpha^{\bullet} \colon X^{\bullet} \rightsquigarrow Y^{\bullet}$, understood as a finite acyclic diagram, instead of the complex $X^{\bullet}$. □

## 7.3 Proof of the main theorem

We are now able to prove the effective version of our main theorem:

**Theorem 7.3.1** (Realization of Voevodsky's motives into and Nori motives, effective version)**.** *Assume that $k \subseteq \mathbb{C}$ is a field and let $\Lambda$ be a noetherian ring.*

*There exists a contravariant triangulated functor*

$$\mathrm{C}^{\mathrm{eff}} \colon \mathrm{DM}^{\mathrm{eff}}_{\mathrm{gm}}(k, \Lambda) \to D^b(\mathcal{MM}^{\mathrm{eff}}_{\mathrm{Nori}}(k, \Lambda))$$

*between Voevodsky's geometric motives and derived Nori motives which calculates singular cohomology $H_{\mathrm{sing}}$ in the following sense:*

*If* $\omega_{\text{sing}}\colon \mathcal{MM}^{\text{eff}}_{\text{Nori}}(k,\Lambda) \to \Lambda\text{-}\mathrm{Mod}$ *is the forgetful fibre functor of Nori motives, then there is a natural isomorphism*

$$\omega_{\text{sing}}\Big(H^n\big(\,\mathrm{C}^{\text{eff}}(X[0])\big)\Big) \cong H^n_{\text{sing}}(X^{\text{an}},\Lambda)$$

*for all smooth varieties* $X$ *over* $k$.

*If* $\Lambda$ *is a Dedekind domain or a field, then* $\mathrm{C}^{\text{eff}}$ *is a tensor functor.*

*Proof.* We first define a functor $\widetilde{\mathrm{C}}\colon C^b(\mathrm{SmCor}(k,\Lambda)) \to D^+(\mathcal{MM}^{\text{eff}}_{\text{Nori}}(k,\Lambda))$:

On objects it is given by

$$\widetilde{\mathrm{C}}(X^\bullet) := \mathrm{C}^\bullet_{\mathfrak{F}}(X^\bullet),$$

where $\mathfrak{F}$ is an arbitrary $\Lambda$-cellular cover-filtration on $X$, which exists by Theorem 7.1.3. By Theorem 7.2.4 this does, as an element of $D^+(\mathcal{MM}^{\text{eff}}_{\text{Nori}}(k,\Lambda))$, not depend on the choice of $\mathfrak{F}$.

On a morphism $\alpha^\bullet\colon X^\bullet \rightsquigarrow Y^\bullet$, the functor $\widetilde{\mathrm{C}}$ is defined by picking a $\Lambda$-cellular cover-filtration $\mathfrak{F}$ on this diagram and then taking the induced morphism

$$\mathrm{C}^\bullet_{\mathfrak{F}}(\alpha^\bullet)\colon\ \mathrm{C}^\bullet_{\mathfrak{F}}(Y^\bullet) \to \mathrm{C}^\bullet_{\mathfrak{F}}(X^\bullet).$$

This does, again by Theorem 7.2.4, not depend on the choice of $\mathfrak{F}$.

Compatibility with compositions follows immediately from choosing $\Lambda$-cellular cover-filtrations on the diagrams of the form $X^\bullet \rightsquigarrow Y^\bullet \rightsquigarrow Z^\bullet$. Hence $\widetilde{\mathrm{C}}$ is a functor.

Note that the diagram

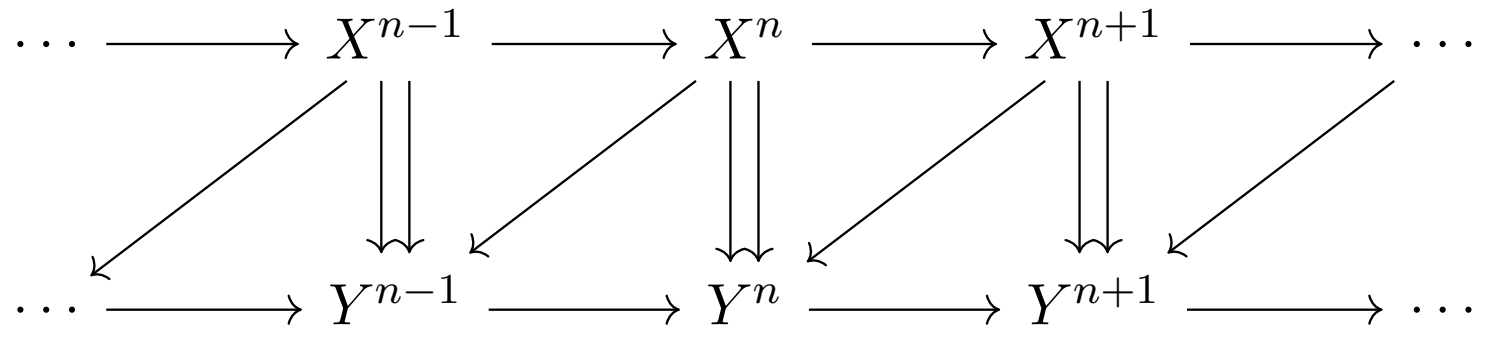

corresponding to a chain homotopy in $C^b(\mathrm{SmCor}(k,\Lambda))$ is finite and acyclic. Hence $\widetilde{\mathrm{C}}$ descends to a functor

$$K^b(\mathrm{SmCor}(k,\Lambda)) \to D^+(\mathcal{MM}^{\text{eff}}_{\text{Nori}}(k,\Lambda)).$$

This functor calculates singular cohomology due to Corollary 7.1.11. We have to show that it descends further to the localization at homotopy invariances HI and Mayer-Vietoris complexes MVN. This can by Lemma 7.1.8 be checked after application of the conservative functor $D(\omega_{\text{sing}})$, i.e. in singular cohomology. We have already seen in Corollary 7.1.11 that $\widetilde{\mathrm{C}}(X) = \mathrm{C}^\bullet_{\mathfrak{F}}(X)$ computes singular cohomology of $X \in \mathrm{Sm}_k$. Therefore the complexes of type $\widetilde{\mathrm{C}}(\mathrm{HI})$ vanish due to the homotopy invariance of singular cohomology, and those of type $\widetilde{\mathrm{C}}(\mathrm{MVN})$ vanish by the Čech spectral sequence, i.e. Theorem B.0.16.

Hence the functor $\widetilde{\mathrm{C}}$ descends further to a functor

$$\widetilde{\mathrm{C}}^{\mathrm{eff},\flat}\colon \underline{\mathrm{DM}}^{\mathrm{eff}}_{\mathrm{gm}}(k,\Lambda) \to D^+(\mathcal{MM}^{\mathrm{eff}}_{\mathrm{Nori}}(k,\Lambda)).$$

Note that $\mathcal{MM}^{\mathrm{eff}}_{\mathrm{Nori}}(k,\Lambda)$ is abelian, thus $D^+(\mathcal{MM}^{\mathrm{eff}}_{\mathrm{Nori}}(k,\Lambda))$ is pseudo-abelian by Lemma 2.4 of [BS01], therefore the functor finally descends to

$$\mathrm{C}^{\mathrm{eff}}\colon\ \mathrm{DM}^{\mathrm{eff}}_{\mathrm{gm}}(k,\Lambda) \to D^+(\mathcal{MM}^{\mathrm{eff}}_{\mathrm{Nori}}(k,\Lambda)).$$

The boundedness of singular cohomology and Lemma 7.1.8 show that this functor actually lands in $D^b(\mathcal{MM}^{\mathrm{eff}}_{\mathrm{Nori}}(k,\Lambda))$. Lastly, we postpone the tensor structure to the next Section 7.4, Theorem 7.4.10 in particular. □

**Remark 7.3.2.** It is clear from the construction that $\mathrm{C}^{\mathrm{eff}}$ extends the functor

$$C^b(\mathbb{Z}[\mathrm{Var}_k]) \to D^b(\mathcal{MM}^{\mathrm{eff}}_{\mathrm{Nori}}(k,\mathbb{Z}))$$

constructed in the proof of Theorem 9.2.21 in [HM16]. In particular their calculations remain valid. Let us give an example:

If $(X,Y,n)$ is a very good pair (cf. Definition 1.3.7), then

$$\mathrm{C}^{\mathrm{eff}}(\mathrm{Cone}(Y \hookrightarrow X)) \cong H^n_{\mathrm{Nori}}(X,Y)[-n].$$

Indeed, we can extend the given data to a $\Lambda$-cellular filtration $\mathcal{F}^\bullet$ with $\mathcal{F}^nX = X$ and $\mathcal{F}^iX = \mathcal{F}^iY$ for $i < n$. Then the claim follows from the construction of $\mathrm{C}^{\mathrm{eff}}$.

## 7.4 Tensor structure

**Definition 7.4.1.** Recall that an object $X$ of an abelian tensor category $\mathcal{A}$ is called *flat* if the endofunctor $-\otimes X\colon \mathcal{A}\to\mathcal{A}$ is exact.

We denote the full subcategory of $\mathcal{A}$ consisting of the flat objects by $\mathcal{A}^\flat$. It is automatically a tensor category.

**Remark 7.4.2.** Note that $K^b(\mathcal{A}^\flat)$, $K^+(\mathcal{A}^\flat)$, $K^-(\mathcal{A}^\flat)$ and $K(\mathcal{A}^\flat)$ are triangulated categories. By standard arguments, e.g. Proposition 10.4.1 of [Wei94], we find that quasi-isomorphisms form a multiplicative system. In particular, the derived categories $D^b(\mathcal{A}^\flat)$, $D^+(\mathcal{A}^\flat)$, $D^-(\mathcal{A}^\flat)$ and $D(\mathcal{A}^\flat)$ make sense.

**Lemma 7.4.3.** *Let $k \subseteq \mathbb{C}$ be a field and let $\Lambda$ be a noetherian ring.*

*Then $M \in \mathcal{MM}^{\mathrm{eff}}_{\mathrm{Nori}}(k,\Lambda)$ is flat if and only if $\omega_{\mathrm{sing}}(M)$ is a flat $\Lambda$-module.*

*Proof.* If $\omega_{\mathrm{sing}}(M)$ is flat, then so is $M$ because $\omega_{\mathrm{sing}}$ is a faithful exact tensor functor.

Let $\mathbb{1} = H^0_{\text{Nori}}(\text{Spec}(k), \emptyset)$. Then

$$\text{Hom}_{\mathcal{MM}^{\text{eff}}_{\text{Nori}}(k,\Lambda)}\left(\bigoplus_{i=1}^{m} \mathbb{1}, \bigoplus_{i=1}^{n} \mathbb{1}\right) \cong \text{Hom}_{\Lambda\text{-}\text{Mod}}(\Lambda^m, \Lambda^n).$$

This shows that the abelian subcategory of $\mathcal{MM}^{\text{eff}}_{\text{Nori}}(k, \Lambda)$ generated by $\mathbb{1}$ is equivalent to $\Lambda$-Mod. Hence flatness of $M$ implies that of $\omega_{\text{sing}}(M)$. □

**Corollary 7.4.4.** *Let $k \subseteq \mathbb{C}$ be a field and let $\Lambda$ be a Dedekind domain or a field.*

*Then every subobject of a flat object in $\mathcal{MM}^{\text{eff}}_{\text{Nori}}(k, \Lambda)$ is flat.*

*Proof.* This follows directly from Lemma 7.4.3 after observing that flat $\Lambda$-modules are exactly the torsion-free ones, a property clearly preserved by taking submodules. □

**Proposition 7.4.5.** *Let $k \subseteq \mathbb{C}$ be a field. Assume that the ring of coefficients $\Lambda$ is a Dedekind domain or a field.*

*Then the inclusion $i\colon \mathcal{MM}^{\text{eff}}_{\text{Nori}}(k, \Lambda)^\flat \to \mathcal{MM}^{\text{eff}}_{\text{Nori}}(k, \Lambda)$ induces an equivalence*

$$D^b(i)\colon D^b(\mathcal{MM}^{\text{eff}}_{\text{Nori}}(k, \Lambda)^\flat) \cong D^b(\mathcal{MM}^{\text{eff}}_{\text{Nori}}(k, \Lambda))$$

*of derived categories.*

*Proof.* Let $A \in \mathcal{MM}^{\text{eff}}_{\text{Nori}}(k, \Lambda)$. Proposition 1.3.10 says that $A$ is a quotient of a subobject $B$ of a finite direct sum $T$ of objects of the form $H^i_{\text{Nori}}(X, Y)$, with $(X, Y, i)$ a (very) good pair. Due to the definition of good pairs, $\omega_{\text{sing}}(T)$ is a finitely generated projective $\Lambda$-module.

Due to Corollary 7.4.4 have thus shown that every $A \in \mathcal{MM}^{\text{eff}}_{\text{Nori}}(k, \Lambda)$ admits an epimorphism $B \twoheadrightarrow A$ from a flat $B$. Induction then shows that every bounded complex $A^\bullet \in C^b(\mathcal{MM}^{\text{eff}}_{\text{Nori}}(k, \Lambda))$ has a flat left resolution, i.e. there is a quasi-isomorphism $B^\bullet \to A^\bullet$ with $B^\bullet \in C^-(\mathcal{MM}^{\text{eff}}_{\text{Nori}}(k, \Lambda)^\flat)$. By construction, $B^\bullet$ has only finitely many non-trivial cohomology groups. By Corollary 7.4.4 we can then truncate $B$ to the left to get a quasi-isomorphic complex in $C^b(\mathcal{MM}^{\text{eff}}_{\text{Nori}}(k, \Lambda)^\flat)$.

Hence $D(i)$ is essentially surjective. If $X^\bullet \to Y^\bullet$ is a quasi-isomorphism in $K^b(\mathcal{MM}^{\text{eff}}_{\text{Nori}}(k, \Lambda))$ with $Y^\bullet \in K^b(\mathcal{MM}^{\text{eff}}_{\text{Nori}}(k, \Lambda)^\flat)$, then we find a quasi-isomorphism $Z^\bullet \to X^\bullet$ such that $Z^\bullet \in K^b(\mathcal{MM}^{\text{eff}}_{\text{Nori}}(k, \Lambda)^\flat)$. By [Bor+87], Proposition 1.24, this shows that $D(i)$ is fully faithful, finishing the proof. □

**Lemma 7.4.6.** *Let $k \subseteq \mathbb{C}$ be a field and let $\Lambda$ be a noetherian ring of coefficients.*

*Then $D^b(\mathcal{MM}^{\text{eff}}_{\text{Nori}}(k, \Lambda)^\flat)$ carries a commutative tensor structure $- \overset{L}{\otimes} -$ such that the conservative forgetful functor*

$$D^b(\omega_{\text{sing}})\colon D^b(\mathcal{MM}^{\text{eff}}_{\text{Nori}}(k, \Lambda)^\flat) \to D^b(\Lambda\text{-}\text{Mod}^\flat)$$

*is a tensor functor.*

*Proof.* This is a standard result of homological algebra: given any $X^\bullet$ and $Y^\bullet$ in $K^b(\mathcal{MM}^{\mathrm{eff}}_{\mathrm{Nori}}(k,\Lambda)^\flat)$, we define their tensor product as $\mathrm{Tot}^\bullet(X^\bullet \otimes Y^\bullet)$. It descends to $D^b(\mathcal{MM}^{\mathrm{eff}}_{\mathrm{Nori}}(k,\Lambda)^\flat)$ by flatness.

Note that $D(\omega_{\mathrm{sing}})$ indeed lands in $D^b(\Lambda\text{-}\mathrm{Mod}^\flat)$ by Lemma 7.4.3. The tensor structure on the latter category is defined the same way, hence the final statement follows. □

**Proposition 7.4.7.** *Let $k \subseteq \mathbb{C}$ be a field. Assume that the ring of coefficients $\Lambda$ is a Dedekind domain or a field.*

*Then $D^b(\mathcal{MM}^{\mathrm{eff}}_{\mathrm{Nori}}(k,\Lambda))$ carries a commutative tensor structure $- \overset{L}{\otimes} -$ such that the conservative forgetful functor*

$$D^b(\omega_{\mathrm{sing}})\colon D^b(\mathcal{MM}^{\mathrm{eff}}_{\mathrm{Nori}}(k,\Lambda)) \to D^b(\Lambda\text{-}\mathrm{Mod})$$

*and the equivalence*

$$D^b(\mathcal{MM}^{\mathrm{eff}}_{\mathrm{Nori}}(k,\Lambda)^\flat) \cong D^b(\mathcal{MM}^{\mathrm{eff}}_{\mathrm{Nori}}(k,\Lambda))$$

*of Proposition 7.4.5 are tensor functors.*

*Proof.* This is immediate from Proposition 7.4.5 and Lemma 7.4.6. □

**Definition 7.4.8.** Recall that the functor $\mathrm{C}^{\mathrm{eff}}$ was defined using $\Lambda$-cellular filtrations, whose singular cohomologies are by definition projective $\Lambda$-modules. Hence the intermediate functor $\widetilde{\mathrm{C}}^{\mathrm{eff},\flat}$ seen in the proof of Theorem 7.3.1 defines actually a functor

$$\mathrm{C}^{\mathrm{eff},\flat}\colon \underline{\mathrm{DM}}^{\mathrm{eff}}_{\mathrm{gm}}(k,\Lambda) \to D^+(\mathcal{MM}^{\mathrm{eff}}_{\mathrm{Nori}}(k,\Lambda)^\flat),$$

even if $\Lambda$ is not a Dedekind domain.

**Proposition 7.4.9.** *Let $k \subseteq \mathbb{C}$ be a field and let $X, Y \in \mathrm{Sm}^{\mathrm{aff}}_k$. Let $\Lambda$ be a noetherian ring.*

*Then there is a natural isomorphism*

$$\mathrm{C}^{\mathrm{eff},\flat}(X \times Y) \cong \mathrm{C}^{\mathrm{eff},\flat}(X) \overset{L}{\otimes} \mathrm{C}^{\mathrm{eff},\flat}(Y)$$

*in $D^b(\mathcal{MM}^{\mathrm{eff}}_{\mathrm{Nori}}(k,\Lambda)^\flat)$.*

*Proof.* Choose $\Lambda$-cellular filtrations $\mathcal{F}$ and $\mathcal{G}$ on $X$ and $Y$, respectively. By Proposition 4.9.2 we get a chain

$$\begin{aligned}\mathrm{C}^{\mathrm{eff},\flat}(X \times Y) &\cong \mathrm{C}^\bullet_{\mathcal{F}\times\mathcal{G}}(X \times Y) \cong \\ &\cong \mathrm{Tot}^\bullet(\mathrm{C}^\bullet_{\mathcal{F}}(X) \otimes \mathrm{C}^\bullet_{\mathcal{G}}(Y)) \cong \mathrm{C}^{\mathrm{eff},\flat}(X) \overset{L}{\otimes} \mathrm{C}^{\mathrm{eff},\flat}(Y).\end{aligned}$$

□

**Theorem 7.4.10.** *Let $k \subseteq \mathbb{C}$ be a field and let $\Lambda$ be a noetherian ring. Then the functor*

$$\mathrm{C}^{\mathrm{eff},\flat}\colon \underline{\mathrm{DM}}^{\mathrm{eff}}_{\mathrm{gm}}(k,\Lambda) \to D^+(\mathcal{MM}^{\mathrm{eff}}_{\mathrm{Nori}}(k,\Lambda)^\flat)$$

*is a tensor functor.*

*Proof.* Let $X_1, X_2 \in \mathrm{Sm}_k$. Choose affine rigidified Nisnevich covers

$$f_1\colon \mathcal{U}_1 \twoheadrightarrow_{r_1} X,$$
$$f_2\colon \mathcal{U}_2 \twoheadrightarrow_{r_2} X_2.$$

We then use Theorem 4.9.7 to choose $\Lambda$-cellular filtrations $\mathcal{F}_1$, $\mathcal{F}_2$ and $\mathcal{F}$ on the Čech diagrams $\check{\mathrm{D}}(f)$, $\check{\mathrm{D}}(g)$ and $\check{\mathrm{D}}(f) \otimes \check{\mathrm{D}}(g)$. As before, this gives $\Lambda$-cellular filtrations on the (double) complexes $\check{\mathrm{C}}^\bullet(f_1)$, $\check{\mathrm{C}}^\bullet(f_2)$ and $\check{\mathrm{C}}^\bullet(f_1) \otimes \check{\mathrm{C}}^\bullet(f_2)$.

Also note that

$$f_1 \times f_2\colon \mathcal{U}_1 \times \mathcal{U}_2 \twoheadrightarrow_{r_1 \times r_2} X_1 \times X_2$$

is an affine pre-rigidified Nisnevich cover whose Čech complex $\check{\mathrm{C}}^\bullet(f_1 \times f_2)$ is the diagonal of $\check{\mathrm{C}}^\bullet(f_1) \otimes \check{\mathrm{C}}^\bullet(f_2)$. Let $f = (f_1 \times f_2)^{\mathrm{rig}}$ be its rigidification as in Definition 6.1.1.

Hence by construction we have for all integers $n_1, n_2$ quasi-isomorphisms

$$\mathrm{C}^\bullet_{\mathcal{F}}\left(\check{\mathrm{C}}^{n_1}(f_1) \times \check{\mathrm{C}}^{n_2}(f_2)\right) \to \mathrm{Tot}^\bullet\left(\mathrm{C}^\bullet_{\mathcal{F}_1}\left(\check{\mathrm{C}}^{n_1}(f_1)\right) \otimes \mathrm{C}^\bullet_{\mathcal{F}_2}\left(\check{\mathrm{C}}^{n_2}(f_2)\right)\right).$$

Looking at singular cohomology, i.e. using Proposition 7.1.9 and Lemma 7.1.8, we also get a quasi-isomorphism

$$\mathrm{Tot}^\bullet\left(\mathrm{C}^\bullet_{\mathcal{F}}\left(\check{\mathrm{C}}^\bullet(f)\right)\right) \to \mathrm{Tot}^\bullet\left(\mathrm{C}^\bullet_{\mathcal{F}}\left(\check{\mathrm{C}}^\bullet(f_1 \times f_2)\right)\right).$$

Recall that the Čech complexes are the associated complexes of simplicial objects $\Delta^{\mathrm{op}} \to \mathbb{Z}[\mathrm{Sm}_k]$, where $\Delta$ is the usual simplex category. Combined with the $\Lambda$-cellular filtrations we hence got for all integers $n$ co-simplicial objects

$$\mathrm{C}^n_{\mathcal{F}_1}, \mathrm{C}^n_{\mathcal{F}_2}, \mathrm{C}^n_{\mathcal{F}}\colon \Delta \to \mathcal{MM}^{\mathrm{eff}}_{\mathrm{Nori}}(k,\Lambda)^\flat.$$

Hence the Eilenberg-Zilber-Cartier Theorem (cf. Satz 2.9 of [DP61]) gives us for every integer $n$ a natural quasi-isomorphism

$$\mathrm{C}^n_{\mathcal{F}}\left(\check{\mathrm{C}}^\bullet(f_1 \times f_2)\right) \to \mathrm{Tot}^\bullet_{i,-j}\left(\mathrm{C}^n_{\mathcal{F}}\left(\check{\mathrm{C}}^i(f_1) \otimes \check{\mathrm{C}}^j(f_2)\right)\right).$$

Lastly, we note that by definition

$$\mathrm{C}(X_1) \cong \mathrm{Tot}^\bullet\left(\mathrm{C}^\bullet_{\mathcal{F}_1}\left(\check{\mathrm{C}}^\bullet(f_1)\right)\right),$$
$$\mathrm{C}(X_2) \cong \mathrm{Tot}^\bullet\left(\mathrm{C}^\bullet_{\mathcal{F}_2}\left(\check{\mathrm{C}}^\bullet(f_2)\right)\right),$$
$$\mathrm{C}(X_1 \times X_2) \cong \mathrm{Tot}^\bullet\left(\mathrm{C}^\bullet_{\mathcal{F}}\left(\check{\mathrm{C}}^\bullet(f)\right)\right).$$

We now put everything together and take total complexes. This then induces a natural zigzag that amounts to a natural isomorphism

$$\mathrm{C}(X) \otimes \mathrm{C}(Y) \cong \mathrm{C}(X \times Y)$$

in $D^+(\mathcal{MM}^{\mathrm{eff}}_{\mathrm{Nori}}(k,\Lambda)^\flat)$.

Taking rigidifications and tensor products are both compatible with rigid bridges as seen in Remark 6.4.2 and Theorem 6.4.5. Thus this argument extends naturally to complexes over $\mathrm{SmCor}(k,\Lambda)$, and to finite correspondences thereof. □

Combined with Proposition 7.4.5 and the obvious descent to pseudo-abelian envelopes, this shows the missing step in the proof of Theorem 7.3.1:

**Corollary 7.4.11.** *Let $k \subseteq \mathbb{C}$ be a field and assume that the ring of coefficients $\Lambda$ is a Dedekind domain or a field.*

*Then the functor*

$$\mathrm{C}^{\mathrm{eff}} \colon \mathrm{DM}^{\mathrm{eff}}_{\mathrm{gm}}(k,\Lambda) \to D^b(\mathcal{MM}^{\mathrm{eff}}_{\mathrm{Nori}}(k,\Lambda))$$

*of Theorem 7.3.1 is a tensor functor.*

Let us now switch to the tensor-localizations. All of the above works mutatis mutandis for $\mathcal{MM}_{\mathrm{Nori}}(k,\Lambda)$ instead of $\mathcal{MM}^{\mathrm{eff}}_{\mathrm{Nori}}(k,\Lambda)$. In particular we get:

**Proposition 7.4.12.** *Let $k \subseteq \mathbb{C}$ be a field. Assume that the ring of coefficients $\Lambda$ is a Dedekind domain or a field.*

*Then the inclusion $\mathcal{MM}_{\mathrm{Nori}}(k,\Lambda)^\flat \hookrightarrow \mathcal{MM}_{\mathrm{Nori}}(k,\Lambda)$ induces an equivalence*

$$D^b(\mathcal{MM}_{\mathrm{Nori}}(k,\Lambda)^\flat) \cong D^b(\mathcal{MM}_{\mathrm{Nori}}(k,\Lambda))$$

*of derived categories. Furthermore, $D^b(\mathcal{MM}_{\mathrm{Nori}}(k,\Lambda))$ carries a commutative tensor structure $- \overset{L}{\otimes} -$ such that the conservative forgetful functor*

$$D^b(\omega_{\mathrm{sing}}) \colon D^b(\mathcal{MM}_{\mathrm{Nori}}(k,\Lambda)) \to D^b(\Lambda\text{-}\mathrm{Mod})$$

*is a tensor functor.*

**Lemma 7.4.13.** *There is a natural equivalence*

$$D^b(\mathcal{MM}_{\mathrm{Nori}}(k,\Lambda)) \cong D^b(\mathcal{MM}^{\mathrm{eff}}_{\mathrm{Nori}}(k,\Lambda))[(\Lambda_{\mathrm{Nori}}(-1)[0])^{\otimes -1}]$$

*of categories.*

*Proof.* Note that $\omega_{\mathrm{sing}}(\Lambda_{\mathrm{Nori}}(-1)) \cong \Lambda$ is a flat $\Lambda$-module, hence $\Lambda_{\mathrm{Nori}}(-1)$ is flat by Lemma 7.4.3. Hence the result follows from Lemma 1.10.5. □

We also trivially have from their definition as 2-colimits:

**Lemma 7.4.14.** *Let $k \subseteq \mathbb{C}$ be a field and let $\Lambda$ be a noetherian ring.*

*Then there is a natural equivalence*

$$\mathcal{MM}_{\mathrm{Nori}}(k,\Lambda)^{\flat} \cong \mathcal{MM}^{\mathrm{eff}}_{\mathrm{Nori}}(k,\Lambda)^{\flat}[\Lambda_{\mathrm{Nori}}(-1)^{\otimes -1}]$$

*of triangulated tensor categories.*

**Lemma 7.4.15.** *Let $k \subseteq \mathbb{C}$ be a field and let $\Lambda$ be a noetherian ring.*

*The functor $\mathrm{C}^{\mathrm{eff}}$ sends the Tate motive $\Lambda_{\mathrm{DM}}(1) \in \mathrm{DM}^{\mathrm{eff}}_{\mathrm{gm}}(k,\Lambda)$ to the Lefschetz motive $\Lambda_{\mathrm{Nori}}(-1)[0]$.*

*Proof.* By Lemma 1.10.12, we can instead work with the complex $A^{\bullet} = (\mathbb{G}_m \to *)$, where $* := \mathrm{Spec}(k)$ and $\mathbb{G}_m$ is in degree 1. An explicit $\Lambda$-cellular filtration on this diagram is given by $A^{\bullet}_1 = A^{\bullet}$ and $A^{\bullet}_0 = (* \to *)$, where we used an arbitrary $k$-point $* \hookrightarrow \mathbb{G}_m$. Using this and Corollary 7.1.12 we find that

$$\mathrm{C}^{\mathrm{eff}}(\Lambda_{\mathrm{DM}}(1)) \cong (\Lambda \to \Lambda \to H^1_{\mathrm{Nori}}(\mathbb{G}_m, *))$$

with $H^1_{\mathrm{Nori}}(\mathbb{G}_m, *)$ in degree 0, where we simply wrote $\Lambda$ for $H^0_{\mathrm{Nori}}(*, \emptyset)$. This is clearly isomorphic to

$$H^1_{\mathrm{Nori}}(\mathbb{G}_m, *)[0] \cong H^1_{\mathrm{Nori}}(\mathbb{G}_m)[0] = \Lambda_{\mathrm{Nori}}(-1)[0].$$

□

**Theorem 7.4.16.** *Assume that $k \subseteq \mathbb{C}$ is a field and let $\Lambda$ be a noetherian ring.*

*Then there exists a contravariant triangulated tensor functor*

$$\mathrm{C}^{\flat} \colon \underline{\mathrm{DM}}_{\mathrm{gm}}(k,\Lambda) \to D^b(\mathcal{MM}_{\mathrm{Nori}}(k,\Lambda)^{\flat})$$

*between Voevodsky's triangulated premotives and derived flat Nori motives which calculates singular cohomology $H_{\mathrm{sing}}$ in the following sense:*

*If $\omega_{\mathrm{sing}} \colon \mathcal{MM}^{\mathrm{eff}}_{\mathrm{Nori}}(k,\Lambda) \to \Lambda\text{-}\mathrm{Mod}$ is the forgetful fibre functor of Nori motives, then there is a natural isomorphism*

$$\omega_{\mathrm{sing}}\Big(H^n\big(\mathrm{C}^{\flat}(X[0])\big)\Big) \cong H^n_{\mathrm{sing}}(X^{\mathrm{an}}, \Lambda)$$

*for all smooth varieties $X$ over $k$.*

*Furthermore, the diagram*

$$\begin{array}{ccccc}
\underline{\mathrm{DM}}^{\mathrm{eff}}_{\mathrm{gm}}(k,\Lambda) & \xrightarrow{\mathrm{C}^{\mathrm{eff},\flat}} & D^b(\mathcal{MM}^{\mathrm{eff}}_{\mathrm{Nori}}(k,\Lambda)^{\flat}) & \xrightarrow{\omega_{\mathrm{sing}}} & D^b(\Lambda\text{-}\mathrm{Mod}^{\flat}) \\
\downarrow{\scriptstyle -(0)} & & \downarrow{\scriptstyle -(0)} & & \| \\
\underline{\mathrm{DM}}_{\mathrm{gm}}(k,\Lambda) & \xrightarrow{\mathrm{C}^{\flat}} & D^b(\mathcal{MM}_{\mathrm{Nori}}(k,\Lambda)^{\flat}) & \xrightarrow{\omega_{\mathrm{sing}}} & D^b(\Lambda\text{-}\mathrm{Mod}^{\flat})
\end{array}$$

*commutes.*

*Proof.* By Lemma 7.4.13 it suffices to construct a functor

$$\underline{\mathrm{DM}}^{\mathrm{eff}}_{\mathrm{gm}}(k,\Lambda)[\Lambda_{\mathrm{DM}}(1)^{\otimes -1}] \to D^b(\mathcal{MM}^{\mathrm{eff}}_{\mathrm{Nori}}(k,\Lambda)^\flat)[\Lambda_{\mathrm{Nori}}(-1)^{\otimes -1}].$$

Due to their definition as 2-colimits it suffices to give functors

$$\mathrm{C}^\flat_n \colon \underline{\mathrm{DM}}^{\mathrm{eff}}_{\mathrm{gm}}(k,\Lambda) \to D^b(\mathcal{MM}^{\mathrm{eff}}_{\mathrm{Nori}}(k,\Lambda)^\flat)$$

such that the diagrams

$$\begin{array}{ccc} \underline{\mathrm{DM}}^{\mathrm{eff}}_{\mathrm{gm}}(k,\Lambda) & \xrightarrow{\mathrm{C}^\flat_n} & D^b(\mathcal{MM}^{\mathrm{eff}}_{\mathrm{Nori}}(k,\Lambda)^\flat) \\ \downarrow {\scriptstyle -\otimes\Lambda_{\mathrm{DM}}(1)} & & \downarrow {\scriptstyle -\otimes\Lambda_{\mathrm{Nori}}(-1)} \\ \underline{\mathrm{DM}}^{\mathrm{eff}}_{\mathrm{gm}}(k,\Lambda) & \xrightarrow{\mathrm{C}^\flat_{n+1}} & D^b(\mathcal{MM}^{\mathrm{eff}}_{\mathrm{Nori}}(k,\Lambda)^\flat) \end{array}$$

commute for all non-negative integers $n$. We pick $\mathrm{C}^\flat_n = \mathrm{C}^{\mathrm{eff},\flat}$, whereupon the compatibility follows immediately from Lemma 7.4.15 and Theorem 7.4.10. □

**Theorem 7.4.17** (Comparison between Voevodsky's and Nori's motives)**.** *Assume that $k \subseteq \mathbb{C}$ is a field and let $\Lambda$ be a field or a Dedekind domain.*

*Then there exists a contravariant triangulated tensor functor*

$$\mathrm{C}\colon \mathrm{DM}_{\mathrm{gm}}(k,\Lambda) \to D^b(\mathcal{MM}_{\mathrm{Nori}}(k,\Lambda))$$

*between Voevodsky's geometric motives and derived Nori motives which calculates singular cohomology $H_{\mathrm{sing}}$ in the following sense:*

*If $\omega_{\mathrm{sing}}\colon \mathcal{MM}^{\mathrm{eff}}_{\mathrm{Nori}}(k,\Lambda) \to \Lambda\text{-}\mathrm{Mod}$ is the forgetful fibre functor of Nori motives, then there is a natural isomorphism*

$$\omega_{\mathrm{sing}}\Big(H^n\big(\mathrm{C}(X[0])\big)\Big) \cong H^n_{\mathrm{sing}}(X^{\mathrm{an}},\Lambda)$$

*for all smooth varieties $X$ over $k$.*

*Furthermore, the diagram*

$$\begin{array}{ccccc} \mathrm{DM}^{\mathrm{eff}}_{\mathrm{gm}}(k,\Lambda) & \xrightarrow{\mathrm{C}^{\mathrm{eff}}} & D^b(\mathcal{MM}^{\mathrm{eff}}_{\mathrm{Nori}}(k,\Lambda)) & \xrightarrow{\omega_{\mathrm{sing}}} & D^b(\Lambda\text{-}\mathrm{Mod}) \\ \downarrow {\scriptstyle -(0)} & & \downarrow {\scriptstyle -(0)} & & \| \\ \mathrm{DM}_{\mathrm{gm}}(k,\Lambda) & \xrightarrow{\mathrm{C}} & D^b(\mathcal{MM}_{\mathrm{Nori}}(k,\Lambda)) & \xrightarrow{\omega_{\mathrm{sing}}} & D^b(\Lambda\text{-}\mathrm{Mod}) \end{array}$$

*commutes.*

*Proof.* This is immediate from Theorem 7.4.16 and Proposition 7.4.12, noting that $D^b(\mathcal{MM}_{\mathrm{Nori}}(k,\Lambda))$ is pseudo-abelian by Lemma 2.4 of [BS01]. □

## 7.5 An alternative proof

Ayoub (private communication) suggested a different approach to construct the functor C based on works by Choudhury and Gallauer Alves de Souza (cf. [CG15]). Their arguments are independent of ours and were developed in parallel. We want to sketch this approach, for which we need two more categories of triangulated motives.

The first one appeared in [Ayo14a]. We will, however, not repeat its definition here, as it is quite involved:

**Definition 7.5.1.** Let $\mathrm{DA}_{\text{ét}}(k, \Lambda)$ be the triangulated tensor category of *Ayoub motives* as defined in Section 3 of [Ayo14a].

The second one comes from [VSF00, chapter 5]:

**Definition 7.5.2.** Analogous to Definition 1.9.5 we have the category $\mathrm{Shv}_{\text{ét}}(\mathrm{SmCor}(S, \Lambda))$ of étale sheaves with $\Lambda$-transfers. Like its Nisnevich counterpart it comes with a Yoneda embedding

$$L_{\text{ét}}\colon \mathrm{SmCor}(S, \Lambda) \to \mathrm{Shv}_{\text{ét}}(\mathrm{SmCor}(S, \Lambda).$$

*Voevodsky's category of (unbounded) effective étale motivic complexes* is the Bousfield localization

$$\mathrm{DM}^{\text{eff}}_{\text{ét}}(S, \Lambda) := D(\mathrm{Shv}_{\text{ét}}(\mathrm{SmCor}(S, \Lambda)))/\langle L(\mathrm{HI})\rangle^{\oplus}$$

at the localizing subcategory generated by the images under the Yoneda embedding of the complexes $\mathbb{A}^1_X \to X$ of Definition 1.9.3. We let $\mathrm{DM}_{\text{ét}}(k, \Lambda)$ be its tensor-localization at the Tate motive.

*Sketch of an alternative proof of Theorem 7.4.17.* Adding transfers induces by Théorème B.1 of [Ayo14a] an equivalence

$$\mathrm{DA}_{\text{ét}}(k, \Lambda) \cong \mathrm{DM}_{\text{ét}}(k, \Lambda).$$

We furthermore have a change-of-topology functor

$$\mathrm{DM}(k, \Lambda) \to \mathrm{DM}^{\text{ét}}(k, \Lambda).$$

Proposition 7.1 of [CG15] gives a functor

$$\mathrm{DA}^{\text{ét}}(k, \Lambda) \to D^b(\text{ind-}\mathcal{HM}(k, \Lambda)),$$

where $\mathcal{HM}(k, \Lambda)$ is the abelian category of homological Nori motives, defined dually to $\mathcal{MM}_{\mathrm{Nori}}(k, \lambda)$ by using singular homology. Note that the argument of loc. cit. is similar to ours: it uses the yoga of filtrations and Zariski covers to reduce to affine objects. But it is simpler because it does not encounter finite correspondences.

An argument similar to that of loc. cit. also gives a functor

$$\mathrm{DA}^{\text{ét}}(k, \Lambda) \to D^b(\text{pro-}\mathcal{MM}_{\text{Nori}}(k, \lambda)).$$

We omit the details.

Composition thus gives a functor

$$\overline{\mathrm{C}}\colon \mathrm{DM}(k, \Lambda) \to D^b(\text{pro-}\mathcal{MM}_{\text{Nori}}(k, \lambda)).$$

In light of the full embedding $\mathrm{DM}^{\text{eff}}_{\text{gm}}(k, \Lambda) \hookrightarrow \mathrm{DM}^{\text{eff}}(k, \Lambda)$ of Theorem 1.9.10 it is sufficient to check that the Yoneda embedding of an $X \in \mathrm{Sm}_k$ lands in $D^b(\mathcal{MM}_{\text{Nori}}(k, \Lambda))$. Replacing $X$ with the Čech complex of an open cover allows us to assume that $X$ is affine. One now checks that by construction $\overline{\mathrm{C}}(L(X))$ is isomorphic to the Nori complex $\mathrm{C}^{\bullet}_{\mathcal{F}}(X)$ of any $\Lambda$-cellular filtration on $X$. □

**Remark 7.5.3.** A comparison of this construction with ours can be found as Proposition 7.8 of [CG15].

## 7.6 Consequences

We demonstrate a few consequences of Theorems 7.3.1 and 7.4.17. Some of them are not new results, but the functor C unifies them.

We start with the following well-known statement:

**Theorem 7.6.1.** *Let $k \subseteq \mathbb{C}$ be a field and let $\Lambda$ be a noetherian ring.*

*Then there exists a Betti realization*

$$H_{\text{sing}}(\Lambda)\colon \mathrm{DM}_{\text{gm}}(k, \Lambda) \to \Lambda\text{-Mod}$$

*with coefficients in $\Lambda$ such that the composition with*

$$-[0]\colon \mathrm{Var}_k \to \mathrm{DM}_{\text{gm}}(k, \Lambda)$$

*is singular cohomology with $\Lambda$-coefficients.*

*Proof.* This is immediate from setting $H_{\text{sing}}(\Lambda) = \omega_{\text{sing}} \circ \mathrm{C}$. □

**Remark 7.6.2.** Note that the existence of such a realization was not used in our proofs. Instead we only used the compatibility of singular cohomology with morphisms of varieties.

**Remark 7.6.3.** It is not difficult to check that for $\Lambda = \mathbb{Z}$ this agrees with Lecomte's realization in [Lec08]. By Théorème 1.1 of op. cit. it thus also agrees with Huber's realization in [Hub00] when $\Lambda = \mathbb{Q}$.

**Remark 7.6.4.** Ayoub recently outlined a promising approach to proving the conservativity conjecture, stating that the Betti realization $H_{\mathrm{sing}}(\mathbb{Q})$ is conservative, i.e. reflects isomorphisms. In particular it implies that the functor C is conservative as well.

**Theorem 7.6.5.** *Let $k \subseteq \mathbb{C}$ be a field.*

*Then the Betti realization $H_{\mathrm{sing}}(\mathbb{Q})\colon \mathrm{DM}_{\mathrm{gm}}(k,\mathbb{Q}) \to D^b(\mathbb{Q}\text{-}\mathrm{Mod})$ of Theorem 7.6.1 factors as*

$$\begin{aligned}\mathrm{DM}_{\mathrm{gm}}(k,\mathbb{Q}) &\to D^b(\mathcal{MM}_{\mathrm{Nori}}(k,\mathbb{Q})) \to D^b(\mathcal{MM}_{\mathrm{AH}}) \to \\ &\to \mathcal{D}_{\mathcal{MR}} \to D^b(\mathrm{MHS}(\mathbb{Q})) \to D^b(\mathbb{Q}\text{-}\mathrm{Mod}),\end{aligned}$$

*where*

- $\mathcal{MM}_{\mathrm{AH}}$ *are the absolute Hodge motives of Deligne and Jannsen,*
- $\mathcal{D}_{\mathcal{MR}}$ *are Huber's derived mixed realizations,*
- $\mathrm{MHS}(\mathbb{Q})$ *are the rational mixed Hodge structures.*

*Proof.* A conditional proof was given as Theorem 10.1.1 in [HM16], assuming the result of Theorem 7.4.17. □

**Theorem 7.6.6.** *Let $k \subseteq \mathbb{C}$ be a field.*

*Then every effective Nori motive $M \in \mathcal{MM}^{\mathrm{eff}}_{\mathrm{Nori}}(k,\mathbb{Q})$ carries a bounded increasing weight filtration $W_\bullet M$.*

*Similarly, every Nori motive in $\mathcal{MM}_{\mathrm{Nori}}(k,\mathbb{Q})$ carries a bounded increasing weight filtration. Both induce the usual weight filtrations on singular cohomology and are made unique by this property.*

*In view of Theorem 7.6.5, this weight filtration induces the weight filtration on Jannsen's mixed realizations $\mathcal{MR}$ and in particular the standard ones on singular, étale and $\ell$-adic cohomology. Furthermore, every morphism of Nori motives is strictly compatible with these filtrations.*

**Remark 7.6.7.** This result is originally due to Arapura (cf. Theorems 6.3.5 and 6.3.6 of [Ara13]). Unlike his more direct proof, ours takes use of the realization functor to push Bondarko's weight structures forward. This argument has the advantage of working without any essential change in positive characteristics, assuming that a reasonable theory of Nori motives together with an analogous realization functor has been constructed.

*Proof of Theorem 7.6.6.* We follow the conditional proof of [HM16], Proposition 10.2.5:

A faithful exact tensor functor $\mathcal{MM}^{\mathrm{eff}}_{\mathrm{Nori}}(k,\mathbb{Q}) \to \mathcal{MR}$ is constructed as part of Proposition 10.1.2 in [HM16]. In particular, there is a faithful exact mixed Hodge realization $h\colon \mathcal{MM}^{\mathrm{eff}}_{\mathrm{Nori}}(k,\mathbb{Q}) \to \mathrm{MHS}(\mathbb{Q})$ and by loc. cit. it calculates Deligne's mixed Hodge structures on smooth varieties.

By [Bon10], Proposition 6.5.3, there exists a weight structure $w$ in the sense of op. cit. Definition 1.1.1 on $\mathrm{DM}_{\mathrm{gm}}^{\mathrm{eff}}(k, \mathbb{Q})$. By Remark 2.4.3 of op. cit. and Remark 7.3.2 we hence have a weight filtration on every $H_{\mathrm{Nori}}^{n}(X, Y)$, where $(X, Y, n) \in \mathrm{VGood}^{\mathrm{eff}}$ is a very good pair. It induces the usual weight filtration on singular cohomology, and is made unique by this due to the faithful exactness of the forgetful functor $\omega_{\mathrm{sing}}$.

The weight filtration clearly extends to sums, subobjects and quotients. This is compatible with the same operations in $\mathrm{MHS}(\mathbb{Q})$ due to the faithful exactness of $h$. Therefore this extension to these three kinds of objects is again unique and well-defined. Proposition 1.3.10 now implies that we thus defined weight filtrations on all of $\mathcal{MM}_{\mathrm{Nori}}^{\mathrm{eff}}(k, \mathbb{Q})$.

The weight filtrations of Hodge theory on singular cohomology are furthermore strictly compatible with arbitrary morphisms (cf. [Del71b], Théorème 2.3.5 (iii)). Due to the faithful exactness of $\omega_{\mathrm{sing}}$ this translates into the strict compatibility of our weight filtration with all morphisms of $\mathcal{MM}_{\mathrm{Nori}}^{\mathrm{eff}}(k, \mathbb{Q})$.

The argument for non-effective Nori motives is completely analogous. The weight filtration induces that on $\mathcal{MR}$ because both are the unique ones inducing that on singular cohomology.

The interested reader should look at Proposition 10.2.5 of [HM16] for further details. □

**Remark 7.6.8.** We needed rational coefficients to have the weight filtration on $\mathrm{MHS}(\mathbb{Q})$ together with its strict compatibility at our disposal. Without that, it is not clear how to ensure the uniqueness, and hence well-definedness, of the weight filtrations induced on quotients and subobjects, as the same object might arise in several essentially different ways.

**Definition 7.6.9.** Let $\mathcal{C}$ be a symmetric tensor category with unit $\mathbb{1}$. An object $X \in \mathcal{C}$ is *rigid* if it has a *dual* $X^\vee \in \mathcal{C}$. This means that there are morphisms $\eta\colon \mathbb{1} \to X \otimes X^\vee$ and $\epsilon\colon X^\vee \otimes X \to \mathbb{1}$ such that

$$X \xrightarrow{\eta \otimes X} X \otimes X^\vee \otimes X \xrightarrow{X \otimes \epsilon} X$$

$$X^\vee \xrightarrow{X^\vee \otimes \eta} X^\vee \otimes X \otimes X^\vee \xrightarrow{\epsilon \otimes X^\vee} X^\vee$$

are the respective identities. The category $\mathcal{C}$ is *rigid* if all its objects are rigid.

An abelian $\Lambda$-linear symmetric tensor category is called *quasi-rigid* if its Tannaka dual, which automatically is an pro-algebraic monoid over $\Lambda$, is a pro-algebraic group over $\Lambda$.

The following is a result of Nori, whose proof can be found in [HM16]:

**Theorem 7.6.10.** *Assume that $k \subseteq \mathbb{C}$ is a field and let $\Lambda$ be a Dedekind domain or a field.*

*Then $\mathcal{MM}_{\mathrm{Nori}}(k, \Lambda)$ is quasi-rigid.*

*Proof.* Let $(X, Y, n)$ be a very good pair (cf. Definition 1.3.7). We consider

$$Z := (Y \hookrightarrow X)[-n] \in \mathrm{DM}_{\mathrm{gm}}(k, \Lambda).$$

Remark 7.3.2 computes $\mathrm{C}(Z)$ as $H^n_{\mathrm{Nori}}(X, Y)[0]$. Note that

$$\omega_{\mathrm{sing}}(H^n_{\mathrm{Nori}}(X, Y)) \cong H^n_{\mathrm{sing}}(X, Y)$$

is a finitely generated projective $\Lambda$-module and thus has a dual projective $\Lambda$-module $H^n_{\mathrm{sing}}(X, Y)^\vee$ in $\Lambda$- Mod.

Furthermore, $\mathrm{DM}_{\mathrm{gm}}(k, \Lambda)$ is rigid by Theorem 4.3.7 of [VSF00, chapter 5], hence $Z$ has a dual $Z^\vee$. The property of being a dual is preserved by tensor functors, hence $\mathrm{C}(Z^\vee)$ is a dual of $\mathrm{C}(Z)$. For the same reason, $D(\omega_{\mathrm{sing}})(\mathrm{C}(Z^\vee))$ is a dual of $D(\omega_{\mathrm{sing}})(\mathrm{C}(Z))$. But duals are unique up to unique isomorphism, hence

$$D(\omega_{\mathrm{sing}})(\mathrm{C}(Z^\vee)) \cong H^n_{\mathrm{sing}}(X, Y)^\vee[0].$$

By conservativity of $D(\omega_{\mathrm{sing}})$ (cf. Lemma 7.1.8) we conclude that $\mathrm{C}(Z^\vee) \cong H^0(\mathrm{C}(Z^\vee))[0]$.

From Proposition 7.4.12 we have the fully faithful tensor functor

$$\mathcal{MM}_{\mathrm{Nori}}(k, \Lambda)^\flat \xhookrightarrow{[0]} D^b(\mathcal{MM}_{\mathrm{Nori}}(k, \Lambda)^\flat) \cong D^b(\mathcal{MM}_{\mathrm{Nori}}(k, \Lambda)).$$

We have shown that $H^n_{\mathrm{Nori}}(X, Y)[0]$ and its dual $H^0(\mathrm{C}(Z^\vee))[0]$ lie in its essential image, therefore $H^n_{\mathrm{Nori}}(X, Y)$ has a dual in $\mathcal{MM}_{\mathrm{Nori}}(k, \Lambda)^\flat$ and in particular in $\mathcal{MM}_{\mathrm{Nori}}(k, \Lambda)$.

By definition we also have an inverse for the Lefschetz motive $\Lambda_{\mathrm{Nori}}(-1)$, thus every $H^n_{\mathrm{Nori}}(X, Y)(m)$, $m \in \mathbb{Z}$, is rigid. By Proposition 8.3.4 of [HM16] it therefore suffices to show that $\mathcal{MM}_{\mathrm{Nori}}(k, \Lambda)$ is generated by the $H^n_{\mathrm{Nori}}(X, Y)(m)$ as an abelian tensor category relative to $\omega_{\mathrm{sing}}$ in the sense of Definition 8.1.9 in [HM16]. As the $H^n_{\mathrm{Nori}}(X, Y)(m)$ are closed under tensor products, this by definition amounts to Lemma 1.10.17. □

**Remark 7.6.11.** The proof of quasi-rigidity of $\mathcal{MM}_{\mathrm{Nori}}(k, \Lambda)$ in [HM16] is significantly more technical: it explicitly constructs the tensor-duals using resolution of singularities. Our proof, on the other hand, avoids this. Both proofs, however, finish by using Nori's rigidity criterion, found as Proposition 8.3.4 in [HM16].

In Definition 1.10.13, we have furthermore defined $\mathcal{MM}_{\mathrm{Nori}}(k, \Lambda)$ as the tensor-localization of $\mathcal{MM}^{\mathrm{eff}}_{\mathrm{Nori}}(k, \Lambda)$ at the Lefschetz motive $\Lambda_{\mathrm{Nori}}(-1)$, avoiding the subtleties arising from the construction via a larger quiver seen in Proposition 8.2.5 of [HM16].

# Appendix A

# Finitistic Covers

By Lemma 5.4.3 the following definition makes sense:

**Definition A.0.1** (Finitistic covers)**.** Let $S$ be a noetherian scheme and let $X \in \mathrm{Sm}_S$.

An étale cover $\mathcal{U} \twoheadrightarrow X$ is *finitistic* if the complex

$$\cdots \longrightarrow c_S(-, \check{\mathrm{c}}^1(\mathcal{U} \twoheadrightarrow X)) \longrightarrow \underbrace{c_S(-, \check{\mathrm{c}}^0(\mathcal{U} \twoheadrightarrow X))}_{=c_S(-,\mathcal{U})} \longrightarrow c_S(-, X) \longrightarrow 0$$

of Nisnevich sheaves with transfers over $S$ is exact.

Recall that a Nisnevich sheaf with transfers is just a Nisnevich sheaf which, as a presheaf, has transfers. Therefore the Definition A.0.1 of finitistic covers is tantamount to being an exact complex of Nisnevich sheaves without transfers.

**Remark A.0.2.** Note that the given complex is the Yoneda embedding of the augmented reduced Čech complex.

The analogous variant using the non-reduced Čech complex is true for arbitrary Nisnevich covers by Theorem 10.3.3 of [CD12] (see Proposition 3.1.3 of [VSF00] for a special case). Theorem A.0.10 below shows that, at least over a regular base, exactly the unifibrant Nisnevich covers are finitistic.

**Remark A.0.3.** We consider a topological analogue in Appendix B. There we get similar results as in this section, Proposition B.0.11 in particular.

**Lemma A.0.4.** *Let $S$ be a regular noetherian scheme. Let morphisms $f\colon X \to Y$ and $g\colon Y \to Z$ in $\mathbb{Z}[\mathrm{Sm}_S]$ with $g \circ f = 0$ be given.*

*Let $\mathfrak{H}$ be the class of henselian local noetherian schemes which are finite and étale over a henselization $T_t^{\mathfrak{h}}$, where $t \in T \in \mathrm{Sm}_S$.*

*Then the following are equivalent:*

(a) *The complex*

$$\mathrm{SmCor}_S(-,X) \xrightarrow{-\circ f} \mathrm{SmCor}_S(-,Y) \xrightarrow{-\circ g} \mathrm{SmCor}_S(-,Z)$$

*is an exact sequence of Nisnevich sheaves with transfers.*

(b) *For every $\mathcal{H} \in \mathfrak{H}$ we get an exact complex*

$$c(\mathcal{H} \times_S X|\mathcal{H}) \xrightarrow{(\mathcal{H}\times_S f)_*} c(\mathcal{H} \times_S Y|\mathcal{H}) \xrightarrow{(\mathcal{H}\times_S g)_*} c(\mathcal{H} \times_S Z|\mathcal{H})$$

*of abelian groups.*

**Remark A.0.5.** The restriction on representable sheaves $\mathcal{F}$ in Lemma A.0.4 might sound arbitrary at first. We used it in the formulation to have an immediate extension of $\mathcal{F}$ to all regular noetherian henselian local schemes over $S$, which are not necessarily of finite type. This is, however, not the main reason.

The real problem lies in the formalism of points, or more specifically in Corollary 1.13 of [Nis89], which requires our presheaves to commute with certain projective limits as seen in the proof.

*Proof of Lemma A.0.4.* Note that $\mathrm{Sm}_S$ contains the Nisnevich site $X_{\mathrm{Nis}}$ of every $X \in \mathrm{Sm}_S$. Thus every Nisnevich sheaf $\mathcal{F} \in \mathrm{Shv}_{\mathrm{Nis}}(\mathrm{Sm}_S)$ restricts to a Nisnevich sheaf $\mathcal{F}_{|X_{\mathrm{Nis}}} \in \mathrm{Shv}_{\mathrm{Nis}}(X_{\mathrm{Nis}})$. It is then easily seen that a conservative family of points on $\mathrm{Sm}_S$ is given by a union $\bigcup_{X\in\mathrm{Sm}_S} P_X$ of conservative families of points $P_X$ on $X_{\mathrm{Nis}}$.

By Corollary 1.17 of [Nis89], the Nisnevich site $X_{\mathrm{Nis}}$ has enough points. They are the stalks $\mathcal{F} \mapsto \mathcal{F}_x$ induced by the field-valued points $\mathrm{Spec}(k) \to x \in X \in \mathrm{Sm}_S$ such that $k|\kappa(x)$ is finite and separable. Let us make this more precise:

Let $\overline{\mathrm{Sch}}_S$ be the category of noetherian schemes separated over $S$. If $T \in \mathrm{Sm}_S$, then the Nisnevich sheaf $\mathrm{SmCor}_S(-,T)$ on $\mathrm{Sm}_S$ extends to a Nisnevich sheaf $\widetilde{L}(T)\colon \overline{\mathrm{Sch}}_S \to \mathrm{Ab}$, $X \mapsto c(X \times_S T|X)$, where the presheaf structure comes from the pullback of relative cycles.

Proposition 9.3.9 of [CD12] tells us that $\widetilde{L}(T)$ commutes with projective limits over systems $D$ in $\mathrm{Sm}_S$ with flat, affine and dominant transition morphisms. Clearly, the directed system of affine and connected Nisnevich neighbourhoods of a point $x \in X \in \mathrm{Sm}_S$ satisfies this. Furthermore, they constitute a final subsystem of all pointed Nisnevich neighbourhoods of $(X, x)$. Hence

$$\varinjlim_{(\mathcal{U},u)\to(X,x)} \widetilde{L}(T)(\mathcal{U}) = c(X_x^{\mathfrak{h}} \times_S T|X_x^{\mathfrak{h}}).$$

Then by Corollary 1.13 (3) of [Nis89] and its proof, the aforementioned stalks, forming a conservative family, are given by $\mathcal{F}_y = c(\mathcal{H} \times_S T|\mathcal{H})$ where $\mathcal{H}$ is finite and étale over a henselization $X_x^{\mathfrak{h}}$. □

**Proposition A.0.6.** *Let $S$ be a regular noetherian scheme and let $X \in \mathrm{Sm}_S$. Then every unifibrant Nisnevich cover $(\mathcal{U}|I) \twoheadrightarrow X$ is finitistic.*

**Remark A.0.7.** We use an argument similar to the proof of Proposition 3.1.3 in [VSF00, chapter 5] and to the proof of its more general counterpart found as Theorem 10.3.3 in [CD12]. The main differences is the additional assumptions of unifibrancy, which is indeed necessary as seen in Theorem A.0.10 below. We have to assume regularity to separate Nisnevich points with respect to the point of $X$ their closed point lies over.

*Proof of Proposition A.0.6.* We fix an arbitrary noetherian regular henselian local scheme $\mathcal{H}$ over $S$. Due to Lemma A.0.4 it is sufficient to show the exactness of the complex

$$\to c(\mathcal{H} \times_S \check{c}^1(\mathcal{U} \twoheadrightarrow X)|\mathcal{H}) \to c(\mathcal{H} \times_S \check{c}^0(\mathcal{U} \twoheadrightarrow X)|\mathcal{H}) \to c(\mathcal{H} \times_S X|\mathcal{H}) \to 0.$$

Let $Y \in \mathrm{Sm}_S$. Because $\mathcal{H}$ is regular, $c(\mathcal{H} \times_S Y|\mathcal{H})$ is by Proposition 1.4.18 identified with the free abelian group generated by the closed integral subschemes of $\mathcal{H} \times_S Y$ which are finite and surjective over $\mathcal{H}$. We call them the *basic generators*. Because $\mathcal{H}$ is henselian local we get from Lemma 6.3.3 that each basic generator is also a henselian local scheme.

If furthermore $Y$ is a scheme over $X$ and $x \in X$ is a point, we let $\ell_x(Y) \subseteq c(\mathcal{H} \times_S Y|\mathcal{H})$ be the free abelian group generated by those basic generators whose closed point lies over $x$. Then the previous discussion shows that

$$c(\mathcal{H} \times_S Y|\mathcal{H}) \cong \bigoplus_{x \in X} \ell_x(Y)$$

and it is immediate that this decomposition is functorial with respect to morphisms of schemes over $X$. Thereby it is enough to show that for each $x \in X$ the complex

$$\ldots \to \ell_x(\check{c}^1(\mathcal{U} \twoheadrightarrow X)) \to \ell_x(\check{c}^0(\mathcal{U} \twoheadrightarrow X)) \to \ell_x(X) \to 0 \tag{A.1}$$

is exact.

By assumption there is an $i \in I$ such that $f_i^{-1}(x)$ consists of a single point $u \in \mathcal{U}_i$ where we have an isomorphism of residue fields $\kappa(u) \cong \kappa(x)$. By Lemma 5.4.3 we may assume that $i$ is the smallest element of $I$. We proceed to show that the complex (A.1) is homotopic to 0:

For a subset $J \subseteq I$ we set $\mathcal{U}_J := \mathcal{U}_{j_1} \times_X \ldots \times_X \mathcal{U}_{j_s}$, where $j_1 < j_2 < \ldots < j_s$ are the distinct elements of $J$. Assume for a moment that the projection $\mathrm{pr}_i \colon \mathcal{U}_{\{i\} \cup J} \to \mathcal{U}_J$ induces an isomorphism $(\mathrm{pr}_i)_* \colon \ell_x(\mathcal{U}_{\{i\} \cup J}) \cong \ell_x(\mathcal{U}_J)$. We now define a chain homotopy $h$ from the complex (A.1) to itself by

$$h_{|\ell_x(\mathcal{U}_J)} := \begin{cases} 0 \text{ if } i \in J \\ (\mathrm{pr}_i)_*^{-1} \text{ else.} \end{cases}$$

Let $d$ be the differential of the complex (A.1). We only have to check that $hd + dh$ is the identity, which can be done on the individual $\ell_x(\mathcal{U}_J)$:

If $i \notin J$ and $j \in J$, then $\mathrm{pr}_i \circ \mathrm{pr}_j = \mathrm{pr}_j \circ \mathrm{pr}_i$ gives $(\mathrm{pr}_i)_*^{-1} \circ (\mathrm{pr}_j)_* = (\mathrm{pr}_j)_* \circ (\mathrm{pr}_i)_*^{-1}$ and thus we get

$$\begin{aligned}
&(hd + dh)_{|\ell_x(\mathcal{U}_J)} = \\
=&(\mathrm{pr}_i)_*^{-1} \circ \sum_{k=0}^{s} (-1)^k (\mathrm{pr}_{j_k})_* + \\
&+ \left( (\mathrm{pr}_i)_* \circ (\mathrm{pr}_i)_*^{-1} - \sum_{k=0}^{s} (-1)^k (\mathrm{pr}_{j_k})_* \circ (\mathrm{pr}_i)_*^{-1} \right) = \\
=&\,\mathrm{id}_{\ell_x(\mathcal{U}_J)} + \sum_{k=0}^{s} \left( (-1)^k (\mathrm{pr}_i)_*^{-1} \circ (\mathrm{pr}_{j_k})_* - (-1)^k (\mathrm{pr}_{j_k})_* \circ (\mathrm{pr}_i)_*^{-1} \right) = \\
=&\,\mathrm{id}_{\ell_x(\mathcal{U}_J)} .
\end{aligned}$$

If instead $i \in J$, then $j_0 = i$ as $i$ is the smallest element of $I$. Hence $h \circ \mathrm{pr}_{j_k} = 0$ except if $k = 0$. Thus we find again that

$$(hd + dh)_{|\ell_x(\mathcal{U}_J)} = h \circ \sum_{k=0}^{s} (-1)^k (\mathrm{pr}_{j_k})_* = \sum_{k=0}^{s} (-1)^k h \circ (\mathrm{pr}_{j_k})_* = \mathrm{id}_{\ell_x(\mathcal{U}_J)} .$$

We are left to prove the postponed isomorphism. It is clearly enough to show that $\mathrm{pr}_i$ induces a bijection on the basic generators of $\ell_x(\mathcal{U}_{\{i\} \cup J})$ and $\ell_x(\mathcal{U}_J)$. We sort those generators by what point $u \in \mathcal{U}_J$ their closed point lies above, where $u$ runs through the points of $\mathcal{U}_J$ above $x$.

Due to our assumption on $i$ there is only a single point $u' \in \mathcal{U}_{\{i\} \cup J}$ above such a $u \in \mathcal{U}_J$ and $\mathrm{pr}_i$ induces an isomorphism $\kappa(u') \cong \kappa(u)$ of residue fields. The claim then follows immediately from the following technical Lemma A.0.8 by setting $(Y, y) = (\mathcal{U}_J, u)$ and $(\mathcal{V}, v) = (\mathcal{U}_{\{i\} \cup J}, u')$. □

**Lemma A.0.8.** *Let $\mathcal{H}$ be a noetherian henselian local scheme with closed point $\mathfrak{h}$. Let $Y$ be a scheme separated and of finite type over $\mathcal{H}$ and let $y \in Y$ lie over $\mathfrak{h}$. Let $f\colon (\mathcal{V}, v) \to (Y, y)$ be a pointed Nisnevich neighbourhood (cf. Definition 6.3.2).*

*Let $\mathfrak{b}(Y, y)$ be the set of closed integral subschemes of $Y$ which are henselian local with closed point $y$ as well as finite and surjective over $\mathcal{H}$. Similarly let $\mathfrak{b}(\mathcal{V}, v)$ be the set of closed integral subschemes of $\mathcal{V}$ which are henselian local with closed point $v$ as well as finite and surjective over $\mathcal{H}$.*

*Then taking images induces a bijection $f_*\colon \mathfrak{b}(\mathcal{V}, v) \to \mathfrak{b}(Y, y)$ and $f$ induces for each $R \in \mathfrak{b}(\mathcal{V}, v)$ an isomorphism $R \cong f_*(R)$.*

*Proof.* If $R \in \mathfrak{b}(\mathcal{V}, v)$, then the image $f(R)$ is a closed irreducible subset of $Y$ by [MVW06], Lemma 1.4, and using the reduced induced structure it is by loc. cit. finite and surjective over $\mathcal{H}$. By Lemma 6.3.3 it is also local and henselian. This defines the map $f_*\colon \mathfrak{b}(\mathcal{V}, v) \to \mathfrak{b}(Y, y)$.

Let now $T \in \mathfrak{b}(Y, y)$. As $y$ is the closed point of $T$ and because by assumption $\kappa(v) \cong \kappa(y)$ we get from [Stacks, Tag 08HQ] a lifting $l\colon T \to \mathcal{V} \times_Y T$ sending $y$ to $v$. As the projection $\mathcal{V} \times_Y T \to T$ is separated and étale we find by [Stacks, Tag 024T] that $l$ is an open and closed immersion, i.e. identifies $T$ with a connected component $\phi(T)$ of $\mathcal{V} \times_Y T$. It must be the unique connected component containing $v$.

By construction we got $f_*(\phi(T)) = T$ and $\phi(T) \cong T$. For $R \in \mathfrak{b}(\mathcal{V}, v)$ we also have $\phi(f_*(R)) \supseteq R$ because $R$ is connected. But $f$ induces an isomorphism of $\phi(f_*(R))$ onto the image $f(R) = f_*(R)$ of $R$, hence this inclusion is an equality, proving the lemma. □

**Corollary A.0.9.** *Let $S$ be a regular scheme of finite dimension and let $\Lambda$ be a ring. Let $X \in \mathrm{Sm}_S$ and let $f\colon \mathcal{U} \twoheadrightarrow X$ be a unifibrant Nisnevich cover.*

*Then the morphism $f[0]\colon \check{\mathrm{c}}^\bullet(\mathcal{U} \twoheadrightarrow X) \to X[0]$ becomes an isomorphism in $\mathrm{DM}^{\mathrm{eff}}_{\mathrm{gm}}(S, \Lambda)$.*

*Proof.* Let

$$L_{\mathbb{Z}}\colon \mathrm{SmCor}_{S,\mathbb{Z}} \to \mathrm{Shv}_{\mathrm{Nis}}(\mathrm{SmCor}(S, \mathbb{Z})),$$
$$L_{\Lambda}\colon \mathrm{SmCor}_{S,\Lambda} \to \mathrm{Shv}_{\mathrm{Nis}}(\mathrm{SmCor}(S, \Lambda))$$

be the Yoneda embeddings.

By Proposition A.0.6, the Nisnevich cover $f\colon \mathcal{U} \twoheadrightarrow X$ is finitistic. Hence $C^b(L_{\mathbb{Z}})(\underline{\check{\mathrm{c}}}^\bullet(\mathcal{U} \twoheadrightarrow X))$ is exact. For any $Y, Z \in \mathrm{Sm}_S$ we have by definition $\mathrm{SmCor}_{S,\Lambda}(Y, Z) = \mathrm{SmCor}_{S,\mathbb{Z}}(Y, Z) \otimes_{\mathbb{Z}} \Lambda$. Furthermore, the abelian groups $\mathrm{SmCor}_{S,\mathbb{Z}}(Y, Z)$ are free by Proposition 1.4.18, in particular flat. We thus have an isomorphism

$$C^b(L_{\Lambda})(\underline{\check{\mathrm{c}}}^\bullet(\mathcal{U} \twoheadrightarrow X)) \cong C^b(L_{\mathbb{Z}})(\underline{\check{\mathrm{c}}}^\bullet(\mathcal{U} \twoheadrightarrow X)) \otimes_{\mathbb{Z}} \Lambda$$

of exact complexes in $C^b(\mathrm{Shv}_{\mathrm{Nis}}(\mathrm{SmCor}(S, \Lambda)))$.

The result then follows from Theorem 1.9.10. □

**Theorem A.0.10.** *Let $S$ be a regular noetherian scheme and let $X \in \mathrm{Sm}_S$.*

*Then the following are equivalent for an ordered étale cover $f\colon (\mathcal{U}|I) \twoheadrightarrow X$:*

*(a) $f\colon (\mathcal{U}|I) \twoheadrightarrow X$ is a unifibrant Nisnevich cover.*

*(b) $f\colon (\mathcal{U}|I) \twoheadrightarrow X$ is finitistic.*

*Proof.* Proposition A.0.6 states that (a) implies (b), hence we only have to show the converse.

We fix an arbitrary $x \in X$ and a finite separable field extension $L|\kappa(x)$ such that $\mathrm{Spec}(L) \times_X \mathcal{U}$ splits into a disjoint union of copies of $\mathrm{Spec}(L)$. Such an $L$ exists because $f\colon \mathcal{U} \to X$ is étale and of finite type, thus in particular quasi-finite. Then by [EGA4-4], Proposition 18.5.15, there exists a (unique) finite étale pointed morphism $g\colon (\mathcal{H}, v) \to (X_x^{\mathfrak{h}}, x)$ of henselian local schemes such that $\kappa(v) \cong L$ and $v \to x$ corresponds to the chosen field extension. Note that $\mathcal{H}$ is as in Lemma A.0.4 (b).

As $f_i\colon \mathcal{U}_i \to X$ is quasi-finite, we get for every $i \in I$ a non-negative integer $a_i(x) := \dim_{\kappa(x)}(\mathcal{O}_{\mathcal{U}_i \times_X x}) = \dim_L(\mathcal{O}_{\mathcal{U}_i \times_X v})$ as the degree of the fibre of $f_i\colon \mathcal{U}_i \to X$ over $x$. Then the unifibrancy condition at $x$ states that $a_i = 1$ for at least one $i \in I$.

Assume now that $(\mathcal{U}|I) \twoheadrightarrow X$ is finitistic. The obvious embedding $\mathrm{Sm}_X \hookrightarrow \mathrm{Sm}_S$ allows us to restrict Nisnevich sheaves. This is clearly exact, hence we may now assume that $S = X$. Let

$$\ell_x(\check{c}^k(\mathcal{U} \twoheadrightarrow X)) \subseteq c(\mathcal{H} \times_X \check{c}^k(\mathcal{U} \twoheadrightarrow X)|\mathcal{H})$$

be as in the proof of Proposition A.0.6. It is freely generated by basic generators, which are identified with the $L$-points $v \to v \times_X \check{c}^k(\mathcal{U} \twoheadrightarrow X)$ by Lemma A.0.8. Due to their definition and the choice of $L$ there are exactly $a_{i_0} a_{i_1} \cdots a_{i_k}$ such $L$-points in

$$v \times_X \mathcal{U}_{i_0} \times_X \mathcal{U}_{i_1} \times_X \ldots \times_X \mathcal{U}_{i_k}.$$

The sum decomposition used in the proof of Proposition A.0.6 shows that the complex

$$\ldots \to \ell_x(\check{c}^1(\mathcal{U} \twoheadrightarrow X)) \to \ell_x(\check{c}^0(\mathcal{U} \twoheadrightarrow X)) \to \ell_x(X) \to 0 \qquad \text{(A.2)}$$

is exact. In particular, its Euler characteristic

$$1 - \sum_{i \in I} a_i + \sum_{\substack{i,j \in I \\ i<j}} a_i a_j - \sum_{\substack{i,j,k \in I \\ i<j<k}} a_i a_j a_k + \ldots = \prod_{i \in I}(1 - a_i)$$

vanishes. This can only happen if one of the $a_i$ is equal to 1, i.e. if the unifibrancy condition is satisfied at $x$. As $x$ was an arbitrary point we conclude that $(\mathcal{U}|I) \twoheadrightarrow X$ is unifibrant. □

# Appendix B

# Topologically Étale Covers

In this appendix we briefly explain the topological analogues of unifibrant Nisnevich covers. Most of this is standard.

**Definition B.0.1.** Let $f\colon U \to X$ be a continuous map between topological spaces with discrete fibres. We call $f$ *topologically étale* if it is a local homeomorphism: every $u \in U$ has an open neighbourhood $V$ which $f$ maps homeomorphically to an open neighbourhood $f(V)$ of $f(u)$.

**Remark B.0.2.** Note that this is not the same as being a covering space.

**Definition B.0.3** (Topological covers)**.** Let $X$ be a topological space and let $f\colon (\mathcal{U}|I) \twoheadrightarrow X$ be a jointly surjective pseudocover (cf. Definition 5.1.1).

- $f\colon (\mathcal{U}|I) \twoheadrightarrow X$ is called an *open cover* if each $\mathcal{U}_i \to X$ is the inclusion of an open subset.
- We call $f\colon (\mathcal{U}|I) \twoheadrightarrow X$ a *topologically étale cover* if the total map $\mathcal{U} \to X$ is topologically étale.
- Lastly, a topologically étale cover $f\colon (\mathcal{U}|I) \twoheadrightarrow X$ is called *unifibrant* if for each $x \in X$ there is an $i \in I$ such that $f^{-1}(x) \cap \mathcal{U}_i$ consists of a single point.

**Remark B.0.4.** As with the scheme-theoretic versions in Definition 5.3.1, one readily checks that these three covers form pretopologies on Top. We have the obvious and analogous inclusions between them.

A justification for the naming and an immediate consequence of the inverse function theorem is the following classical result of Grothendieck (cf. [SGA1], XII, Proposition 3.1):

**Proposition B.0.5.** *Let $k \subseteq \mathbb{C}$ be a field and let $X$ be a variety over $k$.*

*(a) Let $f\colon \mathcal{U} \to X$ be an étale morphism of varieties. Then the analytification $f^{\mathrm{an}}\colon \mathcal{U}^{\mathrm{an}} \to X^{\mathrm{an}}$ is topologically étale.*

*(b) Let $f\colon (\mathcal{U}|I) \twoheadrightarrow X$ be an étale cover of varieties. Then the analytification $f^{\mathrm{an}}\colon (\mathcal{U}^{\mathrm{an}}|I) \twoheadrightarrow X^{\mathrm{an}}$ is a topologically étale cover.*

*(c) Let $f\colon (\mathcal{U}|I) \twoheadrightarrow X$ be a unifibrant Nisnevich cover of varieties. Then the analytification $f^{\mathrm{an}}\colon (\mathcal{U}^{\mathrm{an}}|I) \twoheadrightarrow X^{\mathrm{an}}$ is a unifibrant topologically étale cover.*

**Definition B.0.6.** We use op and tét to denote the pretopologies of open and topologically étale covers, respectively.

The *(small) open site* is the category $X_{\mathrm{op}}$ of open subsets of $X$ with the open pretopology op. The *(small) étale site* over $X$ is the category $X_{\mathrm{tét}}$ of topological spaces which are topologically étale over $X$, provided with the pretopology tét.

Unlike the algebro-geometric setting, there is not much difference between open covers and topologically étale ones:

**Lemma B.0.7.** *Let $X$ be a topological space.*

*Then every open cover is topologically étale. Conversely, every topologically étale cover $g\colon \mathcal{V} \twoheadrightarrow X$ has an open cover $f\colon \mathcal{U} \twoheadrightarrow \mathcal{V}$ such that the composition $g \circ f\colon \mathcal{U} \twoheadrightarrow X$ is an open cover as well.*

*Proof.* This follows easily from the definitions. □

Note that this is slightly stronger than the pretopologies refining each other.

**Definition B.0.8.** Let $X$ be a topological space, let $t \in \{\mathrm{op}, \mathrm{tét}\}$ be a pretopology and let $U \in X_t$. We denote the $t$-sheafification of the presheaf $\mathrm{Hom}(-, U)\colon X_t \to \mathrm{Set}$ by $\rho_t(U)$.

Let now $\Lambda$ be a ring. Taking free $\Lambda$-modules with basis $\rho_t(U)(-)$ defines the *constant sheaf* $\Lambda^t_U = \Lambda\rho_t(U)\colon X_t \to \widetilde{\Lambda\text{-Mod}}$.

**Proposition B.0.9.** *Let $X$ be a locally contractible topological space, let $\Lambda$ be a ring and let $n$ be an integer.*

*Then we have natural isomorphisms*

$$H^n_{\mathrm{sing}}(X, \Lambda) \cong H^n_{\mathrm{op}}(X, \Lambda^{\mathrm{op}}_X) \cong H^n_{\mathrm{tét}}(X, \Lambda^{\mathrm{tét}}_X).$$

*Proof.* The fist isomorphism is standard, see for example Theorem 2.2.5 of [HM16] or Theorem 4.14 of [Ram05].

Lemma B.0.7 implies that on the category $X_{\mathrm{tét}}$ the sheaves for the open pretopology agree with those for the topologically étale one. Consider the full embedding $i\colon X_{\mathrm{op}} \hookrightarrow X_{\mathrm{tét}}$, by which the restriction

$$i^*\colon \mathrm{Shv}_{\mathrm{tét}}(X_{\mathrm{tét}}) = \mathrm{Shv}_{\mathrm{op}}(X_{\mathrm{tét}}) \to \mathrm{Shv}_{\mathrm{op}}(X_{\mathrm{op}})$$

is exact. As $\Lambda^{\mathrm{op}}_X = \Lambda^{\mathrm{tét}}_X \circ i = i^*(\Lambda^{\mathrm{tét}}_X)$, the result follows. □

Recall Definition 5.4.2. In accordance with Definition A.0.1 we define:

**Definition B.0.10.** A topologically étale cover $(\mathcal{U}|I) \twoheadrightarrow X$ is called *finitistic* if

$$\cdots \longrightarrow \mathbb{Z}\rho_{\text{ét}}(\check{c}^1(\mathcal{U} \twoheadrightarrow X)) \longrightarrow \underbrace{\mathbb{Z}\rho_{\text{ét}}(\check{c}^0(\mathcal{U} \twoheadrightarrow X))}_{=\mathbb{Z}\rho_{\text{ét}}(\mathcal{U})} \longrightarrow \mathbb{Z}\rho_t(X) \longrightarrow 0$$

is exact as a complex of topologically étale sheaves on $X_{\text{tét}}$.

Note that this does not depend on a choice of an ordering on $I$ by the obvious analogue of Lemma 5.4.3.

We then have the following topological analogue of Proposition A.0.6:

**Proposition B.0.11.** *Let $f\colon (\mathcal{U}|I) \twoheadrightarrow X$ be an ordered unifibrant topologically étale cover of topological spaces.*

*Then $f\colon (\mathcal{U}|I) \twoheadrightarrow X$ is finitistic.*

*Proof.* It suffices to check the exactness at the stalks, where we mimic the proof of Proposition A.0.6:

Fix an arbitrary $x \in X$. Let $i \in I$ be such that $f^{-1}(x) \cap U_i$ consists of a single point $u_i$ which has an open neighbourhood $V \hookrightarrow \mathcal{U}_i$ that is mapped homeomorphically to an open neighbourhood $W = f(V)$ of $x$. We may also assume that $i$ is the largest element of $I$ by Lemma 5.4.3. This gives us a section

$$\mathcal{U}_{j_1} \times_X \cdots \times_X \mathcal{U}_{j_r} \times_X W \to \mathcal{U}_{j_1} \times_X \cdots \times_X \mathcal{U}_{j_r} \times_X \mathcal{U}_i$$
$$(u_1, \ldots, u_r) \mapsto (u_1, \ldots, u_r, f_{|V}^{-1}(w))$$

for all $(u_1, \ldots, u_r)$ over a point $w \in W$. This induces a homotopy at the level of stalks as in the proof of Proposition A.0.6. □

**Remark B.0.12.** The counting argument in the proof of Theorem A.0.10 can be used to show that the converse holds when assuming that $f\colon \mathcal{U} \twoheadrightarrow X$ has finite fibres.

**Definition B.0.13.** Let $f\colon (\mathcal{U}|I) \twoheadrightarrow X$ be an ordered topologically étale cover of topological spaces and let $\mathcal{F}\colon X_{\text{tét}} \to \text{Ab}$ be a presheaf.

The *(full) Čech cohomology* of $\mathcal{F}$ with respect to $\mathcal{U} \twoheadrightarrow X$ is

$$\check{H}^n(\mathcal{U}, \mathcal{F}) := H^n(\mathcal{F}(\check{C}^\bullet(\mathcal{U} \twoheadrightarrow X))).$$

The *reduced Čech cohomology* of $\mathcal{F}$ with respect to $\mathcal{U} \twoheadrightarrow X$ is

$$\check{h}^n(\mathcal{U}, \mathcal{F}) := H^n(\mathcal{F}(\check{c}^\bullet(\mathcal{U} \twoheadrightarrow X))).$$

**Theorem B.0.14** (Čech spectral sequences)**.** *Let $f\colon (\mathcal{U}|I) \twoheadrightarrow X$ be an ordered topologically étale cover of topological spaces and let $\mathcal{F}$ be a sheaf on $X_{\text{tét}}$.*

1. *There exists a natural Čech spectral sequence*

$$E_2^{p,q} = \check{H}^q(\mathcal{U}, H^p_{\text{tét}}(-,\mathcal{F})) \implies H^{p+q}_{\text{tét}}(X,\mathcal{F}).$$

2. *Assume that $f\colon (\mathcal{U}|I) \twoheadrightarrow X$ is unifibrant. Then there exists a natural reduced Čech spectral sequence*

$$E_2^{p,q} = \check{h}^q(\mathcal{U}, H^p_{\text{tét}}(-,\mathcal{F})) \implies H^{p+q}_{\text{tét}}(X,\mathcal{F}).$$

**Remark B.0.15.** The first spectral sequence of Theorem B.0.14 is well known and holds in any site, see for example [Tam94], Theorem (3.4.4) i). The second one is usually only stated for open covers, which are clearly unifibrant.

*Proof of Theorem B.0.14.* By the preceding Remark B.0.15 we only have to show the existence of the second spectral sequence, but the following argument works, with small adaptations, also for the first case.

The cover $(\mathcal{U}|I) \twoheadrightarrow X$ is finitistic by Proposition B.0.11. Hence in the derived category $D(\mathrm{Shv}_{\text{tét}}(X,\mathbb{Z}))$ of topologically étale sheaves of abelian groups on $X_{\text{tét}}$ we have a natural isomorphism between $\mathbb{Z}\rho_{\text{tét}}(X)[0]$ and

$$\cdots \longrightarrow \mathbb{Z}\rho_{\text{tét}}(\check{c}^1(\mathcal{U} \twoheadrightarrow X)) \longrightarrow \mathbb{Z}\rho_{\text{tét}}(\check{c}^0(\mathcal{U} \twoheadrightarrow X)) \longrightarrow 0.$$

By construction we have $\mathrm{Hom}_{\mathrm{Shv}_{\text{tét}}}(\mathbb{Z}\rho_{\text{tét}}(V),\mathcal{F}) \cong \mathcal{F}(V)$ for all $V \in X_{\text{tét}}$, thus

$$\mathrm{Ext}^n_{\mathrm{Shv}_{\text{tét}}}(\mathbb{Z}\rho_{\text{tét}}(V),\mathcal{F}) \cong H^n_{\text{tét}}(V,\mathcal{F}).$$

This identifies the desired spectral sequence with the hypercohomology spectral sequence of the complex with respect to the functor $\mathrm{Hom}_{\mathrm{Shv}_{\text{tét}}}(-,\mathcal{F})$. □

Putting $\mathcal{F} = \mathbb{Z}\rho_{\text{tét}}(X)$ in Theorem B.0.14 and using Proposition B.0.9 we therefore get

**Theorem B.0.16** (Singular Čech spectral sequences)**.** *Let $f\colon (\mathcal{U}|I) \twoheadrightarrow X$ be an ordered topologically étale cover of a locally contractible topological space $X$.*

1. *There exists a natural Čech spectral sequence*

$$E_2^{p,q} = \check{H}^q(\mathcal{U}, H^p_{\text{sing}}(-)) \implies H^{p+q}_{\text{sing}}(X).$$

2. *Assume that $f\colon (\mathcal{U}|I) \twoheadrightarrow X$ is unifibrant. Then there exists a natural reduced Čech spectral sequence*

$$E_2^{p,q} = \check{h}^q(\mathcal{U}, H^p_{\mathrm{sing}}(-)) \implies H^{p+q}_{\mathrm{sing}}(X).$$

**Remark B.0.17.** A condition such as unifibrancy in the second version is easily seen to be necessary. For example, take the cover of any non-empty $X$ consisting of a single component $\mathcal{U}_1 = X \sqcup X$.

**Remark B.0.18.** On analytifications of varieties, both versions of Theorem B.0.16 can, using Proposition A.0.6, also be proven by invoking a Betti realization with $\Lambda$-coefficients of an adequate triangulated category of motives. For example, given one for $\mathrm{DM}^{\mathrm{eff}}_{\mathrm{gm}}(\mathbb{C}, \Lambda)$ or $\mathrm{DM}^{\mathrm{eff}}(\mathbb{C}, \Lambda)$, the first statement follows from Proposition 3.1.3 of [VSF00, chapter 5], and the second one from our Proposition A.0.6. In the case $\Lambda = \mathbb{Z}$ such a Betti realization was for example given in [Lec08].

# References

[AHP15] Giuseppe Ancona, Annette Huber, and Simon Pepin Lehalleur. *On the relative motive of a commutative group scheme.* 2015. arXiv: 1409.3401v2.

[AHW15] Aravind Asok, Marc Hoyois, and Matthias Wendt. *Affine representability results in $\mathbb{A}^1$-homotopy theory I: vector bundles.* 2015. arXiv: 1506.07093v2.

[Ara13] Donu Arapura. "An abelian category of motivic sheaves". In: *Adv. Math.* 233 (2013), pp. 135–195. DOI: 10.1016/j.aim.2012.10.004.

[Ayo07a] Joseph Ayoub. "Les six opérations de Grothendieck et le formalisme des cycles évanescents dans le monde motivique. I". In: *Astérisque* 314 (2007), x+466 pp. (2008).

[Ayo07b] Joseph Ayoub. "Les six opérations de Grothendieck et le formalisme des cycles évanescents dans le monde motivique. II". In: *Astérisque* 315 (2007), vi+364 pp. (2008).

[Ayo14a] Joseph Ayoub. "La réalisation étale et les opérations de Grothendieck". In: *Ann. Sci. Éc. Norm. Supér. (4)* 47.1 (2014), pp. 1–145.

[Ayo14b] Joseph Ayoub. "L'algèbre de Hopf et le groupe de Galois motiviques d'un corps de caractéristique nulle, II". In: *J. Reine Angew. Math.* 693 (2014), pp. 151–226. DOI: 10.1515/crelle-2012-0090.

[Beĭ87] A. A. Beĭlinson. "On the derived category of perverse sheaves". In: *K-theory, arithmetic and geometry (Moscow, 1984–1986).* Vol. 1289. Lecture Notes in Math. Springer, Berlin, 1987, pp. 27–41. DOI: 10.1007/BFb0078365.

[Bon10] M. V. Bondarko. "Weight structures vs. $t$-structures; weight filtrations, spectral sequences, and complexes (for motives and in general)". In: *J. K-Theory* 6.3 (2010), pp. 387–504. DOI: 10.1017/is010012005jkt083.

[Bor+87] A. Borel et al. *Algebraic D-modules.* Vol. 2. Perspectives in Mathematics. Academic Press, Inc., Boston, MA, 1987, pp. xii+355.

[BS01] Paul Balmer and Marco Schlichting. "Idempotent completion of triangulated categories". In: *J. Algebra* 236.2 (2001), pp. 819–834. DOI: `10.1006/jabr.2000.8529`.

[BS14] Manjul Bhargava and Matthew Satriano. "On a notion of "Galois closure" for extensions of rings". In: *J. Eur. Math. Soc. (JEMS)* 16.9 (2014), pp. 1881–1913. DOI: `10.4171/JEMS/478`.

[BV08] Alexander Beilinson and Vadim Vologodsky. "A DG guide to Voevodsky's motives". In: *Geom. Funct. Anal.* 17.6 (2008), pp. 1709–1787. DOI: `10.1007/s00039-007-0644-5`.

[CD12] Denis-Charles Cisinski and Frédéric Déglise. *Triangulated categories of mixed motives.* 2012. arXiv: `0912.2110v3`.

[CG15] Utsav Choudhury and Martin Gallauer Alves de Souza. *An isomorphism of motivic Galois groups.* 2015. arXiv: `1410.6104v2`.

[Del71a] Pierre Deligne. "Théorie de Hodge. I". In: *Actes du Congrès International des Mathématiciens (Nice, 1970), Tome 1.* Gauthier-Villars, Paris, 1971, pp. 425–430.

[Del71b] Pierre Deligne. "Théorie de Hodge. II". In: *Inst. Hautes Études Sci. Publ. Math.* 40 (1971), pp. 5–57.

[Del74] Pierre Deligne. "Théorie de Hodge. III". In: *Inst. Hautes Études Sci. Publ. Math.* 44 (1974), pp. 5–77.

[Del90] P. Deligne. "Catégories tannakiennes". In: *The Grothendieck Festschrift, Vol. II.* Vol. 87. Progr. Math. Birkhäuser Boston, Boston, MA, 1990, pp. 111–195.

[DP61] Albrecht Dold and Dieter Puppe. "Homologie nicht-additiver Funktoren. Anwendungen". In: *Ann. Inst. Fourier Grenoble* 11 (1961), pp. 201–312.

[EGA4-1] A. Grothendieck. "Éléments de géométrie algébrique. IV. Étude locale des schémas et des morphismes de schémas. I". In: *Inst. Hautes Études Sci. Publ. Math.* 20 (1964), p. 259.

[EGA4-2] A. Grothendieck. "Éléments de géométrie algébrique. IV. Étude locale des schémas et des morphismes de schémas. II". In: *Inst. Hautes Études Sci. Publ. Math.* 24 (1965), p. 231.

[EGA4-3] A. Grothendieck. "Éléments de géométrie algébrique. IV. Étude locale des schémas et des morphismes de schémas. III". In: *Inst. Hautes Études Sci. Publ. Math.* 28 (1966), p. 255.

[EGA4-4] A. Grothendieck. "Éléments de géométrie algébrique. IV. Étude locale des schémas et des morphismes de schémas IV". In: *Inst. Hautes Études Sci. Publ. Math.* 32 (1967), p. 361.

[Fer98] Daniel Ferrand. “Un foncteur norme”. In: *Bull. Soc. Math. France* 126.1 (1998), pp. 1–49. URL: http://www.numdam.org/item?id=BSMF_1998__126_1_1_0.

[Fri82] Eric M. Friedlander. *Étale homotopy of simplicial schemes.* Vol. 104. Annals of Mathematics Studies. Princeton University Press, Princeton, N.J.; University of Tokyo Press, Tokyo, 1982, pp. vii+190.

[Ful98] William Fulton. *Intersection theory.* Second. Vol. 2. Ergebnisse der Mathematik und ihrer Grenzgebiete. 3. Folge. A Series of Modern Surveys in Mathematics [Results in Mathematics and Related Areas. 3rd Series. A Series of Modern Surveys in Mathematics]. Springer-Verlag, Berlin, 1998, pp. xiv+470. DOI: 10.1007/978-1-4612-1700-8.

[Gov65] V. E. Govorov. “On flat modules”. In: *Sibirsk. Mat. Ž.* 6 (1965), pp. 300–304.

[GW10] Ulrich Görtz and Torsten Wedhorn. *Algebraic geometry I.* Advanced Lectures in Mathematics. Schemes with examples and exercises. Vieweg + Teubner, Wiesbaden, 2010, pp. viii+615. DOI: 10.1007/978-3-8348-9722-0.

[Han99] Masaki Hanamura. “Mixed motives and algebraic cycles. III”. In: *Math. Res. Lett.* 6.1 (1999), pp. 61–82. DOI: 10.4310/MRL.1999.v6.n1.a5.

[Hat02] Allen Hatcher. *Algebraic topology.* Cambridge University Press, Cambridge, 2002, pp. xii+544.

[HM14] Annette Huber and Stefan Müller-Stach. *On the relation between Nori Motives and Kontsevich Periods.* 2014. arXiv: 1105.0865v5.

[HM16] Annette Huber and Stefan Müller-Stach. *Periods and Nori Motives, version 3.* 2016, pp. xxiv+361. URL: http://home.mathematik.uni-freiburg.de/arithgeom/preprints/buch/buch-v3.pdf.

[Hub00] Annette Huber. “Realization of Voevodsky’s motives”. In: *J. Algebraic Geom.* 9.4 (2000), pp. 755–799.

[Ive70] Birger Iversen. “On linear determinants”. In: *Papers from the “Open House for Algebraists” (Aarhus Univ., Aarhus, 1970).* Mat. Inst., Aarhus Univ., Aarhus, 1970, pp. 12–22.

[Ivo05] Florian Ivorra. *Réalisation p-adique des motifs mixtes.* 2005. URL: https://perso.univ-rennes1.fr/florian.ivorra/PHD.pdf.

[Ivo07] Florian Ivorra. “Réalisation $l$-adique des motifs triangulés géométriques. I”. In: *Doc. Math.* 12 (2007), pp. 607–671.

[Ivo14] Florian Ivorra. *Perverse Nori motives.* 2014. URL: https://hal.archives-ouvertes.fr/hal-00978893.

[Jan90] Uwe Jannsen. *Mixed motives and algebraic K-theory.* Vol. 1400. Lecture Notes in Mathematics. With appendices by S. Bloch and C. Schoen. Springer-Verlag, Berlin, 1990, pp. xiv+246.

[Kel13] Shane Kelly. *Triangulated categories of motives in positive characteristic.* 2013. arXiv: 1305.5349v2.

[KSW99] Istvan Kovacs, Daniel S. Silver, and Susan G. Williams. "Determinants of commuting-block matrices". In: *Amer. Math. Monthly* 106.10 (1999), pp. 950–952. DOI: 10.2307/2589750.

[Lam99] T. Y. Lam. *Lectures on modules and rings.* Vol. 189. Graduate Texts in Mathematics. Springer-Verlag, New York, 1999, pp. xxiv+557. DOI: 10.1007/978-1-4612-0525-8.

[Laz69] Daniel Lazard. "Autour de la platitude". In: *Bull. Soc. Math. France* 97 (1969), pp. 81–128.

[Lec08] Florence Lecomte. "Réalisation de Betti des motifs de Voevodsky". In: *C. R. Math. Acad. Sci. Paris* 346.19-20 (2008), pp. 1083–1086. DOI: 10.1016/j.crma.2008.09.018.

[Lev98] Marc Levine. *Mixed motives.* Vol. 57. Mathematical Surveys and Monographs. American Mathematical Society, Providence, RI, 1998, pp. x+515. DOI: 10.1090/surv/057.

[LT07] Dan Laksov and Anders Thorup. "A determinantal formula for the exterior powers of the polynomial ring". In: *Indiana Univ. Math. J.* 56.2 (2007), pp. 825–845. DOI: 10.1512/iumj.2007.56.2937.

[Lun08] Christian Lundkvist. "Counterexamples regarding symmetric tensors and divided powers". In: *J. Pure Appl. Algebra* 212.10 (2008), pp. 2236–2249. DOI: 10.1016/j.jpaa.2008.03.024.

[Mil80] James S. Milne. *Étale cohomology.* Vol. 33. Princeton Mathematical Series. Princeton University Press, Princeton, N.J., 1980, pp. xiii+323.

[MVW06] Carlo Mazza, Vladimir Voevodsky, and Charles Weibel. *Lecture notes on motivic cohomology.* Vol. 2. Clay Mathematics Monographs. American Mathematical Society, Providence, RI; Clay Mathematics Institute, Cambridge, MA, 2006, pp. xiv+216.

[Nis89] Ye. A. Nisnevich. "The completely decomposed topology on schemes and associated descent spectral sequences in algebraic $K$-theory". In: *Algebraic $K$-theory: connections with geometry and topology (Lake Louise, AB, 1987)*. Vol. 279. NATO Adv. Sci. Inst. Ser. C Math. Phys. Sci. Kluwer Acad. Publ., Dordrecht, 1989, pp. 241–342.

[Nor02] Madhav V. Nori. "Constructible sheaves". In: *Algebra, arithmetic and geometry, Part I, II (Mumbai, 2000)*. Vol. 16. Tata Inst. Fund. Res. Stud. Math. Tata Inst. Fund. Res., Bombay, 2002, pp. 471–491.

[Ram05] S. Ramanan. *Global calculus*. Vol. 65. Graduate Studies in Mathematics. American Mathematical Society, Providence, RI, 2005, pp. xii+316.

[RG71] Michel Raynaud and Laurent Gruson. "Critères de platitude et de projectivité. Techniques de "platification" d'un module". In: *Invent. Math.* 13 (1971), pp. 1–89.

[Rob63] Norbert Roby. "Lois polynomes et lois formelles en théorie des modules". In: *Ann. Sci. École Norm. Sup. (3)* 80 (1963), pp. 213–348.

[Rob80] Norbert Roby. "Lois polynômes multiplicatives universelles". In: *C. R. Acad. Sci. Paris Sér. A-B* 290.19 (1980), A869–A871.

[Ryd07] David Rydh. "A minimal set of generators for the ring of multisymmetric functions". In: *Ann. Inst. Fourier (Grenoble)* 57.6 (2007), pp. 1741–1769. URL: `http://aif.cedram.org/item?id=AIF_2007__57_6_1741_0`.

[Ryd08a] David Rydh. *Families of zero-cycles and divided powers: I. Representability*. 2008. URL: `https://people.kth.se/~dary/famzerocyclesI-20080304.pdf`.

[Ryd08b] David Rydh. *Families of zero-cycles and divided powers: II. The universal family*. 2008. URL: `https://people.kth.se/~dary/famzerocyclesII-20080411.pdf`.

[Ryd08c] David Rydh. *Hilbert and Chow schemes of points, symmetric products and divided powers*. 2008. URL: `https://people.kth.se/~dary/hilbchowsymdiv-20080416.pdf`.

[Ryd08d] David Rydh. *Families of cycles*. 2008. URL: `https://people.kth.se/~dary/famofcycles20080518.pdf`.

[Ryd13] David Rydh. "Existence and properties of geometric quotients". In: *J. Algebraic Geom.* 22.4 (2013), pp. 629–669. DOI: `10.1090/S1056-3911-2013-00615-3`.

[SGA1] Alexander Grothendieck and Michele Raynaud. *Revêtements étales et groupe fondamental (SGA 1)*. 2004. arXiv: math/0206203v2.

[SS03] Michael Spieß and Tamás Szamuely. "On the Albanese map for smooth quasi-projective varieties". In: *Math. Ann.* 325.1 (2003), pp. 1–17. DOI: 10.1007/s00208-002-0359-8.

[Stacks] The Stacks Project Authors. *Stacks Project*. http://stacks.math.columbia.edu. 2015.

[SV96] Andrei Suslin and Vladimir Voevodsky. "Singular homology of abstract algebraic varieties". In: *Invent. Math.* 123.1 (1996), pp. 61–94. DOI: 10.1007/BF01232367.

[Tam94] Günter Tamme. *Introduction to étale cohomology*. Universitext. Translated from the German by Manfred Kolster. Springer-Verlag, Berlin, 1994, pp. x+186. DOI: 10.1007/978-3-642-78421-7.

[Tho80] R. W. Thomason. "Beware the phony multiplication on Quillen's $\mathscr{A}^{-1}\mathscr{A}$". In: *Proc. Amer. Math. Soc.* 80.4 (1980), pp. 569–573. DOI: 10.2307/2043425.

[Vac05] Francesco Vaccarino. "The ring of multisymmetric functions". In: *Ann. Inst. Fourier (Grenoble)* 55.3 (2005), pp. 717–731. URL: http://aif.cedram.org/item?id=AIF_2005__55_3_717_0.

[Vac06] Francesco Vaccarino. *Symmetric product as moduli space of linear representations*. 2006. arXiv: math/0608655v1.

[Voe10a] Vladimir Voevodsky. "Homotopy theory of simplicial sheaves in completely decomposable topologies". In: *J. Pure Appl. Algebra* 214.8 (2010), pp. 1384–1398. DOI: 10.1016/j.jpaa.2009.11.004.

[Voe10b] Vladimir Voevodsky. "Unstable motivic homotopy categories in Nisnevich and cdh-topologies". In: *J. Pure Appl. Algebra* 214.8 (2010), pp. 1399–1406. DOI: 10.1016/j.jpaa.2009.11.005.

[VSF00] Vladimir Voevodsky, Andrei Suslin, and Eric M. Friedlander. *Cycles, transfers, and motivic homology theories*. Vol. 143. Annals of Mathematics Studies. Princeton University Press, Princeton, NJ, 2000, pp. vi+254.

[Wei94] Charles A. Weibel. *An introduction to homological algebra*. Vol. 38. Cambridge Studies in Advanced Mathematics. Cambridge University Press, Cambridge, 1994, pp. xiv+450. DOI: 10.1017/CBO9781139644136.